\documentclass[11pt,twoside=semi,a4paper]{article}

\usepackage{layout_thesis}

\RequirePackage{bm}
\RequirePackage{bbm}
\RequirePackage{tabto}
\RequirePackage{authblk}
\RequirePackage{mathtools}
\RequirePackage{enumerate}
\RequirePackage{enumitem}
\RequirePackage{algorithm}
\RequirePackage{algpseudocode}

\newcommand{\N}{\mathbb{N}}
\newcommand{\Z}{\mathbb{Z}}
\newcommand{\Q}{\mathbb{Q}}
\newcommand{\R}{\mathbb{R}}
\newcommand{\C}{\mathbb{C}}
\newcommand{\E}{\mathbb{E}}

\newcommand{\bD}{\mathbb D}

\newcommand{\bP}{\mathbb P}

\newcommand{\bS}{\mathbb S}

\newcommand{\bV}{\mathbb V}

\newcommand{\ol}{\overline}
\newcommand{\ul}{\underline}
\newcommand{\tr}{\operatorname{tr}}

\renewcommand{\Re}{\operatorname{Re}}
\renewcommand{\Im}{\operatorname{Im}}

\newcommand{\cs}{\bm{s}}
\newcommand{\tz}{\tilde{z}}
\newcommand{\tJ}{\tilde{J}}

\newcommand{\diag}{\operatorname{diag}}
\newcommand{\dist}{\operatorname{dist}}
\newcommand{\supp}{\operatorname{supp}}

\newcommand{\Id}{\operatorname{Id}}
\newcommand{\Cov}{\operatorname{Cov}}

\makeatletter
\newcommand*\bigcdot{\mathpalette\bigcdot@{.5}}
\newcommand*\bigcdot@[2]{\mathbin{\vcenter{\hbox{\scalebox{#2}{$\m@th#1\bullet$}}}}}
\makeatother

\swapnumbers
\theoremstyle{plain}
\newtheorem{theorem}{Theorem}[section]

\newtheorem*{theorem*}{Theorem}
\newtheorem{lemma}[theorem]{Lemma}

\newtheorem{corollary}[theorem]{Corollary}
\theoremstyle{definition}
\newtheorem{definition}[theorem]{Definition}

\newtheorem{example}[theorem]{Example}
\newtheorem{assumption}[theorem]{Assumption}

\theoremstyle{remark}
\newtheorem{remark}[theorem]{Remark}

\swapnumbers
\theoremstyle{indented}

\title{Estimation of Population Linear Spectral Statistics by Marchenko--Pastur Inversion}
\date{}
\author{Ben Deitmar}
\affil{\small \textit{Department of Mathematical Stochastics, University of Freiburg \protect\\ Ernst-Zermelo-Str. 1, 79104 Freiburg, Germany \protect\\ E-mail: ben.deitmar@stochastik.uni-freiburg.de}}

\begin{document}
	\thispagestyle{empty}
	\maketitle
	
	\begin{abstract}
		\noindent
		A new method of estimating population linear spectral statistics from high-dimensional data is introduced. When the dimension $d$ grows with the sample size $n$ such that $\frac{d}{n} \to c>0$, the proposed method is the first with proven convergence rate of $\mathcal{O}(n^{\varepsilon - 1})$ for any $\varepsilon > 0$
		in a general nonparametric setting. For Gaussian data, a CLT for the estimation error with normalization factor $n$ is shown.
	\end{abstract}

	\section{Introduction}\label{Section_Introduction}
	The estimation of a high-dimensional covariance matrix $\Sigma_n \in \C^{d \times d}$ and its eigenvalues
	\begin{align*}
		& \lambda_1,\dots,\lambda_d \geq 0
	\end{align*}
	from i.i.d. samples $Y_1,\dots,Y_n \in \C^d$ with $\Cov[Y_1,Y_1] = \Sigma_n$ is a fundamental question in statistics. Often, quantities of interest are population linear spectral statistics (PLSS)
	\begin{align*}
		& L_{n}(g) \coloneq \frac{1}{d} \sum\limits_{j=1}^d g(\lambda_j) = \frac{1}{d} \tr\big( g(\Sigma_n) \big)
	\end{align*}
	for a function $g$ defined on an interval containing $\lambda_1,\dots,\lambda_d$. Influential examples include the log-determinant (see e.g. \cite{CaiZhou2015_MinimaxLogDet}), formulas for the limiting mean and covariance in spectral CLTs (see e.g. \cite{BaiCLT}), or deterministic equivalents to generalized linear spectral statistics of the form $\frac{1}{d} \tr\big( f(\bm{S}_n) g(\Sigma_n) \big)$ (see Remark~\ref{Remark_GLSS}). Contrary to the classical regime $d = o(n)$, the eigenvalues of the sample covariance matrix
	\begin{align*}
		& \bm{S}_n \coloneq \frac{1}{n} \sum\limits_{k=1}^n Y_k Y_k^* \in \C^{d \times d}
	\end{align*}
	are inconsistent estimators for $\lambda_1,\dots,\lambda_d$, when the dimension $d$ is comparable to the sample size $n$. In random matrix theory, high-dimensional data is modeled by the assumption that the dimension $d$ grows with $n$ such that $\frac{d}{n} \to c_\infty$ holds for some $c_\infty>0$. In this asymptotic regime, the celebrated Marchenko--Pastur law characterizes the first-order asymptotic behavior of the eigenvalues of $\bm{S}_n$. Due to their availability in practice, linear spectral statistics of the form $\frac{1}{d} \tr(f(\bm{S}_n))$ have been used extensively for hypothesis testing in high-dimensional settings, see e.g. \cite{BaiCLT, Testing3, Testing2, Testing1, Testing4}.
	\\[0.5em]
	This paper develops estimators for population linear spectral statistics $L_n(g)$ in the high-dimen\-sional regime with error rate of order $\mathcal{O}(n^{\varepsilon-1})$ for any $\varepsilon > 0$, when the function $g$ is holomorphic on a sufficiently large subset of $\C$. The proposed eigen-inference method thus improves upon the previously best error rate in this nonparametric setting of $\mathcal{O}(n^{-\frac{1}{2}})$, which was achieved in \cite{MPIKong} for estimation of $L_n(g)$ for the first two moments $g(\lambda)=\lambda$ and $g(\lambda) = \lambda^2$. In the parametric setting, where the population eigenvalues are assumed to take only $k$ distinct values for fixed $k \in \N$ and all $n \in \N$, a rate of $\mathcal{O}(n^{-1})$ was achieved in \cite{MPIBaiChenYao} and \cite{Speicher2008}.
	For the special case where $g(\lambda)=\log(\lambda)$, the population linear spectral statistic $L_n(g)$ coincides with the $\log$-determinant $\log(\det(\Sigma_n))$, for which an estimation method was introduced in \cite{CaiZhou2015_MinimaxLogDet}, which for Gaussian data achieves an error rate of order $\mathcal{O}(n^{-1})$.
	\\[0.5em]
	A key innovation of the proposed method is that it gives data-driven criteria for describing the domain on which a function $g$ must be holomorphic for the PLSS $L_n(g)$ to be estimateable with rate $\mathcal{O}(n^{\varepsilon-1})$, see Theorem~\ref{Thm_Consistency}, Remark~\ref{Remark_Replaceability} and Corollary~\ref{Cor_PLSS_Consistency}.

	\subsection{Model and notation}\label{Subsection_ModelAndNotation}
	This subsection gives a first description of the model and introduces the notation of the paper. A formal statement of assumptions will be given in Assumption \ref{Assumption_EigInf_Main}.
	\\[0.5em]
	Let $(d_n)_{n \in \N}$ be a sequence with values in $\N$ such that the quotient $\frac{d_n}{n}$ converges to a constant $c_\infty>0$. Suppressing the dependence of $d_n$ on $n$ in notation, write $c_n \coloneq \frac{d}{n} \to c_\infty > 0$.
	\\[0.5em]
	A fundamental model, commonly employed in random matrix theory, is that an observed data matrix $\bm{Y} = [Y_1,\dots,Y_n] \in \C^{d \times n}$ is  of the form
	\begin{align}\label{Eq_RMT_Fundamental_Assumption}
		& \bm{Y}_n \coloneq B_n \bm{X}_n \ ,
	\end{align}
	where $\bm{X}_n$ is a random $(d \times n)$ matrix with with independent centered variance-one entries, and $B_n$ is a deterministic $(d \times d)$ matrix satisfying $B_n B_n^* = \Sigma_n$.
	The sample covariance matrix is defined as
	\begin{align*}
		& \bm{S}_n \coloneq \frac{1}{n} \bm{Y}_n \bm{Y}_n^* = \frac{1}{n} B_n \bm{X}_n \bm{X}_n^* B_n^* \ .
	\end{align*}
	With the notation of $\lambda_1(A),\dots,\lambda_d(A)$ for the eigenvalues of a Hermitian matrix $A$, it is assumed that the \textit{population spectral distribution} (PSD)
	\begin{align}\label{Eq_Def_Hn}
		& H_n \coloneq \frac{1}{d} \sum\limits_{j=1}^d \delta_{\lambda_j(\Sigma_n)}
	\end{align}
	converges weakly to some limiting population spectral distribution $H_\infty \neq \delta_0$ with compact support on $[0,\infty)$.
	\\[0.5em]
	For any $n \in \N$, define the \textit{empirical spectral distribution} $\hat{\nu}_n \coloneq \frac{1}{d} \sum\limits_{j=1}^d \delta_{\lambda_j(\bm{S}_n)}$, and let $\nu_n$ denote the probability measure that arises from $H_n$ and $c_n$ by Lemma~\ref{Lemma_MP}. The measures $\nu_n$ will be referred to as \textit{deterministic equivalents} of $\hat{\nu}_n$.
	The symbols $\hat{\ul{\nu}}_n$ and $\ul{\nu}_n$ shall denote the probability measures
	\begin{align}\label{Eq_Def_ulNu}
		& \hat{\ul{\nu}}_n \coloneq \frac{1}{n} \sum\limits_{j=1}^n \delta_{\lambda_j(\frac{1}{n} \bm{X}_n^* \Sigma_n \bm{X}_n)} = (1-c_n) \delta_0 + c_n\hat{\nu}_n \ \ \text{ and } \ \ \ul{\nu}_n \coloneq (1-c_n) \delta_0 + c_n\nu_n \ ,
	\end{align}
	where $\ul{\nu}_\infty$ is defined analogously.
	Further, introduce the notation $\C^+ \coloneq \{z \in \C \mid \Im(z)>0\}$ for the open upper half plane and analogously $\C^- \coloneq \{z \in \C \mid \Im(z)<0\}$. The complex conjugate of a $z \in \C$ is denoted by $\ol{z}$, while the closure of a set $S \subset \C$ will be written as $\mathrm{closure}(S)$. The support of a measure $\mu$ on $\R$ is denoted by $\supp(\mu)$ and the distance between two sets $A,B \subset \C$ is defined as
	\begin{align*}
		& \dist(A,B) \coloneq \inf\big\{ |a-b| \ : \ a \in A, \, b \in B \big\} \ .
	\end{align*}
	Note that this distance function is not a metric and not to be confused with the Hausdorff distance.
	The open $\varepsilon$-neighborhood around a complex number $z \in \C^+$ will be denoted as $B^{\C}_\varepsilon(z) \coloneq \{v \in \C \, : \, |v-z| < \varepsilon\}$, and likewise $B^{\C^+}_\varepsilon(z) \coloneq \{v \in \C^+ \, : \, |v-z| < \varepsilon\}$. For a curve $\gamma : (0,1) \rightarrow \C$, let $\mathrm{image}(\gamma)$ denote the image $\{\gamma(t) \mid t \in (0,1)\}$. For the same curve $\gamma$, let $\operatorname{Hol}(\gamma)$ denote the space of functions $g : U_g \rightarrow \C$ holomorphic on some open set $U_g \subset \C$ that contains the convex hull of $\mathrm{image}(\gamma) \cup \ol{\mathrm{image}(\gamma)}$, equipped with the norm
	\begin{align*}
		& \|g\|_{\gamma} \coloneq \sup\limits_{\substack{z \in \mathrm{image}(\gamma) \\ \ \ \ \ \ \ \cup \ol{\mathrm{image}(\gamma)}}} |g(z)| \ .
	\end{align*}
	Functions $f : \R \rightarrow \R$ are applied to Hermitian matrices $M \in \C^{d \times d}$ canonically by $f(M) \coloneq U \diag(f(\lambda_1),\dots,f(\lambda_d)) U^*$, where $U \diag(\lambda_1,\dots,\lambda_d)U^*$ is the spectral decomposition of $M$. The spectral norm of a matrix $M$ is denoted by $\|M\|$.
	The maximum of two numbers $a,b \in \R$ is denoted as $a \lor b$, while $a \land b$ denotes their minimum.

	\subsection{The Marchenko--Pastur law and its inversion}\label{Subsection_MPlaw}
	The Marchenko--Pastur law is said to hold if the almost sure weak convergence of measures
	\begin{align}\label{Eq_MP_law}
		& 1 = \bP\Big( \hat{\nu}_n \coloneq \frac{1}{d} \sum\limits_{j=1}^d \delta_{\lambda_j(\bm{S}_n)} \xRightarrow{n \to \infty} \nu_\infty \Big)
	\end{align}
	is satisfied, where $\nu_\infty$ is uniquely defined by $H_\infty$ and $c_\infty$ through the Marchenko--Pastur equation (see Lemma~\ref{Lemma_MP}). For the model introduced in Subsection~\ref{Subsection_ModelAndNotation}, the Marchenko--Pastur law is known to hold under the additional condition that the entries of $\bm{X}_n$ are i.i.d. (see Theorem~2.14 of \cite{MPIYaoZhengBai}) or under the assumption that the fourth moments $\E[|(\bm{X}_n)_{j,k}|^4]$ are uniformly bounded (see \cite{MP_RowCorrelation_Bai} and Lemma \ref{Lemma_FourthMomentMP_law}).
	The Marchenko--Pastur equation is formulated in terms of Stieltjes transforms
	\begin{align}\label{Eq_DefStieltjes}
		& \cs_{\mu} : \C\setminus\supp(\mu) \rightarrow \C \ \ ; \ \ z \mapsto \int_\R \frac{1}{\lambda - z} \, d\mu(\lambda)
	\end{align}
	of measures $\mu$ on $\R$.
	The Stieltjes transform $\cs_{\mu}$ uniquely identifies the underlying probability measure $\mu$ on $\R$ and the asymptotic of $\cs_{\mu}(z)$ for $z$ approaching $\supp(\mu)$ is especially significant for reconstruction of $\mu$, as $\frac{1}{\pi} \Im(\cs_\mu(x+\bm{i}\eta))$ equates to the integral over the Cauchy kernel $\frac{\eta/\pi}{(\,\bigcdot\,-x)^2+\eta^2}$ with respect to $\mu$, for all $\eta>0$. A standard formulation of the Marchenko--Pastur equation 
	is as follows.
	
	\begin{lemma}[Marchenko--Pastur equation (Equation 1.4 in \cite{silverstein1995strong})]\label{Lemma_MP}\
		\\
		For any probability measure $H \neq \delta_0$ on $[0,\infty)$ with compact support and any constant $c>0$, there exists a probability measure $\nu \neq \delta_0$ on $[0,\infty)$ with compact support that is uniquely defined by the following property of its Stieltjes transform $\cs_\nu$.
		\\[0.5em]
		For all $\tilde{z} \in \C^+$, the Stieltjes transform $\cs_{\nu}(\tilde{z})$ is the unique solution to
		\begin{align}\label{Eq_MP_Equation}
			& s = \int_\R \frac{1}{\lambda(1-c\tilde{z}s-c) - \tilde{z}} \, dH(\lambda)
		\end{align}
		from the set
		\begin{align}\label{Eq_Def_Qzc}
			& \Big\{ s \in \C \ \Big| \ \Im\Big( cs + \frac{c-1}{\tilde{z}} \Big) > 0 \Big\} \ .
		\end{align}
	\end{lemma}
	
	The remainder of this subsection introduces an inversion formula for the Marchenko--Pastur equation to reconstruct $\cs_{H}(z)$ from $\nu$ and $c$ via a fixed-point equation, and defines the population Stieltjes transform estimator, which will be analyzed in this paper.
	\\[0.5em]
	First, the definition of the spectral domain on which Marchenko--Pastur inversion is possible. The most common choices for the parameter $\theta$ will be $1$ or $\infty$.
	
	\begin{definition}[Spectral domain]\
		\\
		Given a probability measure $H \neq \delta_0$ with compact support on $[0,\infty)$ and a constant $c>0$, for any $\theta \in (0,\infty]$ define the spectral domain
		\begin{align}\label{Eq_DefD}
			\bD_{H,c}(\theta) & \coloneq \Big\{ z \in \C^+ \ \Big| \ \Im\big( (1-cz\cs_H(z)-c)z \big) > 0 , \ \Big| \frac{cz \Im(z\cs_H(z))}{\Im((1-cz\cs_H(z)-c)z)} \Big| < \theta \Big\} \ .
		\end{align}
	\end{definition}
	
	Slightly weaker variations of the following lemma are known in the fields of random matrices and free probability. Chapter 6 in \cite{BaiSALDRM} employs identical maps for spectrum separation, while Section 2.3 of \cite{LedoitWolf1} describes how the same relationship may be used to characterize the density of $\nu$. Property (\ref{Eq_VarphiStieltjesProerty}) implies that $z \mapsto \frac{1}{\varphi_{\nu,c}(1/z)}$ is the analytic subordination function, originally introduced in \cite{Biane_Subordination}, for the free multiplicative convolution $H \boxtimes \mu_{\text{MP},c}$ with the standard Marchenko--Pastur distribution $\mu_{\text{MP},c}$. 
	
	\begin{lemma}[Transform between $H$-space and $\nu$-space]\label{Lemma_SpaceTransform}\
		\\
		For any probability measure $H\neq \delta_0$ with compact support on $[0,\infty)$ and any constant $c>0$, let $\nu = \nu(H,c)$ be the probability measure described in Lemma~\ref{Lemma_MP}.
		The holomorphic map
		\begin{align}\label{Eq_Def_Phi}
			& \Phi_{H,c} : \C^+ \rightarrow \C \ \ ; \ \ z \mapsto (1-cz\cs_H(z)-c)z
		\end{align}
		maps $\bD_{H,c}(\infty)$ to $\C^+$ bijectively. Its inverse function is
		\begin{align}\label{Eq_Def_varphi}
			& \varphi_{\nu,c} : \C^+ \rightarrow \bD_{H,c}(\infty) \ \ ; \ \ \tz \mapsto \frac{\tz}{1-c\tz\cs_{\nu}(\tz)-c} = \frac{-1}{\cs_{\ul{\nu}}(\tz)} \ ,
		\end{align}
		where $\ul{\nu} = c\nu + (1-c)\delta_0$ and thus $\cs_{\ul{\nu}}(\tz) = c\cs_{\nu}(\tz) - \frac{1-c}{\tz}$.
		It further holds that
		\begin{align}\label{Eq_PhiStieltjesProerty}
			& \forall z \in \bD_{H,c}(\infty) : \ \Phi_{H,c}(z) \cs_{\nu}(\Phi_{H,c}(z)) = z \cs_H(z) \ ,
		\end{align}
		or equivalently
		\begin{align}\label{Eq_VarphiStieltjesProerty}
			& \forall \tz \in \C^+ : \ \tz \cs_{\nu}(\tz) = \varphi_{\nu,c}(\tz) \cs_H(\varphi_{\nu,c}(\tz)) \ .
		\end{align}
	\end{lemma}
	
	\begin{figure}[H]
		\vspace{-0.3cm}
		\centering
		\includegraphics[width=\textwidth]{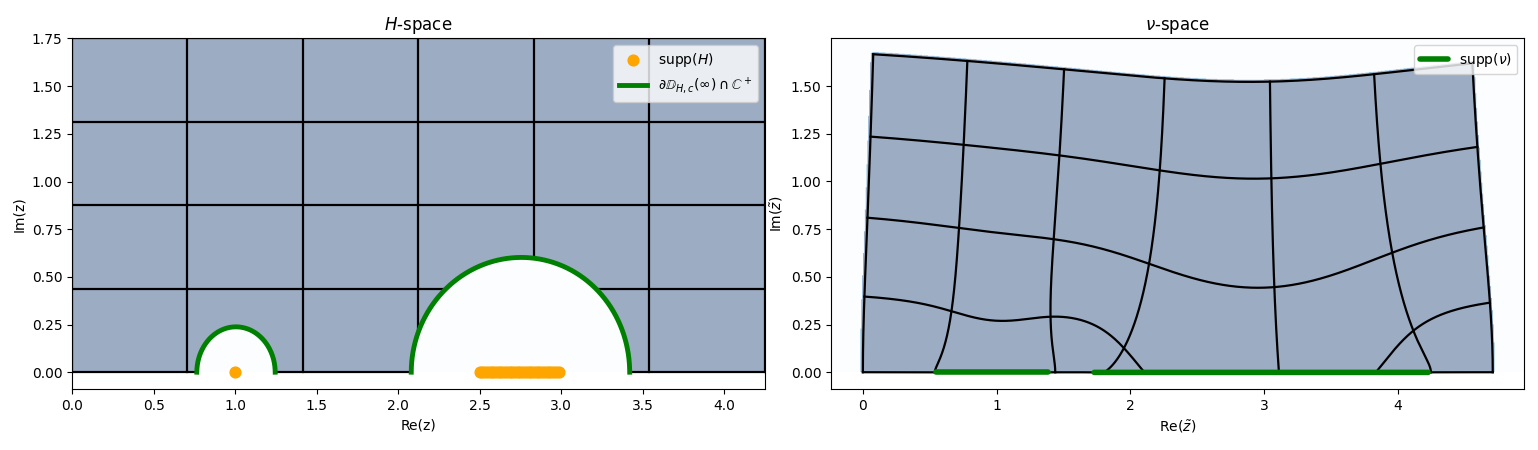}
		\vspace{-0.6cm}
		\caption{Left: Graphical representation of $\bD_{H,c}(\infty)$ (blue), where $H = \frac{1}{2} \delta_{1}+\frac{1}{2}\text{Uniform}([2.5,3])$ and $c=\frac{1}{10}$. Boundary points, with positive imaginary part, are marked green. The support of $H$ is marked orange.\\
		Right: The $\Phi_{H,c}$-image of the region of $\bD_{H,c}(\infty)$ shown on the left. The support of $\nu$ is marked green.}\label{Fig_SpaceTransform}
	\end{figure}
	
	The fact that $\varphi_{\nu,c}$ may be calculated without knowledge of $H$ ensures that the map $\cs_H$ on $\bD_{H,c}(\infty)$ may by (\ref{Eq_VarphiStieltjesProerty}) be reconstructed from only $\nu$ and $c$. The following lemma turns the above relationship into a method of reconstructing $\cs_{H}$ on the spectral domain $\bD_{H,c}(\infty)$ via a fixed-point equation.
	
	\begin{lemma}[Marchenko--Pastur inversion in the infinite-sample case]\label{Lemma_MPI}\
		\\
		For any probability measure $H\neq \delta_0$ with compact support on $[0,\infty)$ and any constant $c>0$, let $\nu$ be the probability measure described in Lemma~\ref{Lemma_MP}. The following two statements hold.
		\begin{itemize}
			\item[a)]
			For every $z \in \C^+$ satisfying
			\begin{align}\label{Eq_LemmaMPI_Assumptions_x}
				& \Im\big( (1-cz\cs_H(z)-c)z \big) > 0 \ ,
			\end{align}
			the Stieltjes transform $\cs_H(z)$ is the only solution $s \in \C^+$ to the equation
			\begin{align}\label{Eq_ReverseMPEquation_re}
				& zs + 1 = \int_\R \frac{\lambda}{\lambda - (1 - czs - c)z} \, d\nu(\lambda) \ .
			\end{align}
			
			\item[b)]
			For any $z \in \C^+$ satisfying
			\begin{align}\label{Eq_LemmaMPI_Assumptions1}
				& \Im\big( (1-cz\cs_H(z)-c)z \big) \leq 0 \ ,
			\end{align}
			there is no solution $s \in \C^+$ to (\ref{Eq_ReverseMPEquation_re}) that also satisfies $\Im\big( (1-czs-c)z \big) > 0$.
		\end{itemize}
	\end{lemma}
	
	The question naturally arises whether the connection explored in Lemmas~\ref{Lemma_SpaceTransform} and~\ref{Lemma_MPI} may be used for eigen-inference in the finite sample setting, i.e. to infer $\cs_{H_n}$ when only $\hat{\nu}_n$ and $c_n$ are known. A first hurdle to overcome is the fact that $\varphi_{\hat{\nu}_n,c_n} : \C^+ \rightarrow \C^+$ may no longer be injective. While in the infinite-sample setting, $\cs_{H}(z)$ was for $z \in \bD_{H,c}(\infty)$ uniquely defined by either $\cs_{H}(z) = \frac{\varphi_{\nu,c}^{-1}(z) \cs_{\nu}(\varphi_{\nu,c}^{-1}(z))}{z}$ or as the unique solution to (\ref{Eq_ReverseMPEquation_re}), in the finite sample case, there may be multiple pre-images $\varphi_{\hat{\nu}_n,c_n}^{-1}(z)$ or equivalently, multiple solutions $\hat{s} \in \C^+$ with $\Im((1-cz\hat{s}-c)z) > 0$ to the equation
	\begin{align}\label{Eq_ReverseMPEquation_Empirical}
		& z\hat{s} + 1 = \int_\R \frac{\lambda}{\lambda - (1 - cz\hat{s} - c)z} \, d\hat{\nu}_n(\lambda) \ .
	\end{align}
	The additional assumption $\big| \frac{c z \Im(z\hat{s})}{\Im((1-cz\hat{s}-c)z)} \big| < 1$ will ensure that $\hat{s}$ is unique.
	
	\begin{figure}[H]
		\vspace{-0.2cm}
		\centering
		\includegraphics[width=0.95\textwidth]{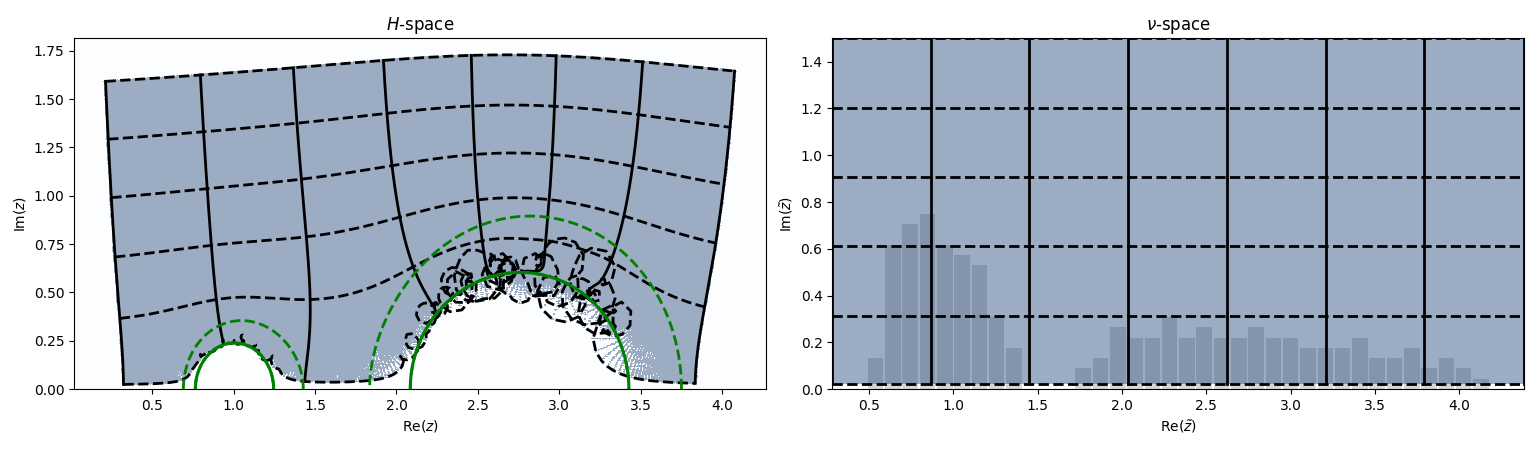}
		
		\vspace{-0.2cm}
		\caption{Right: Graphical representation of $\C^+$ with a blue grid covering the space $\{\Im(\tz) \geq 0.02\}$. Histograms of sample eigenvalues of a $(d=200,n=2000)$ data matrix with population spectral distribution $H_n \approx \frac{1}{2} \delta_{1}+\frac{1}{2}\text{Uniform}([2.5,3])$ and $c_n=\frac{1}{20}$ are added to the background.\\
		Left: Images of the right-hand grids under the maps $\varphi_{\hat{\nu}_n,c_n}$. The boundaries of the sets $\bD_{H_n,c_n}(\infty)$ (solid) and $\bD_{H_n,c_n}(1)$ (dashed) are added as green lines.\\
		Note that the map $\varphi_{\hat{\nu}_n,c_n}$ is not injective, as is demonstrated by the overlapping grid lines in the left-hand image.
		}\label{Fig_SpaceTransform_FiniteSample}
	\end{figure}

	\begin{definition}[Population Stieltjes transform estimator]\label{Def_StieltjesEstimator}\
		\\
		For all $z \in \C^+$, if a solution $\hat{s}_n(z) \in \C^+$ of
		\begin{align}\label{Eq_ReverseMPEquation_Empirical_again}
			& z\hat{s} + 1 = \int_\R \frac{\lambda}{\lambda - (1 - c_nz\hat{s} - c_n)z} \, d\hat{\nu}_n(\lambda)
		\end{align}
		exists that also satisfies
		\begin{align}\label{Eq_StilEstimator_DefiningConditions}
			& \Im\big( (1-c_nz\hat{s}-c_n)z \big) > 0 \ \ \text{ and } \ \ 
			\Big| \frac{c_n z \Im(z\hat{s})}{\Im((1-c_nz\hat{s}-c_n)z)} \Big| < 1 \ ,
		\end{align}
		then call $\hat{s}_n(z)$ a \textit{population Stieltjes transform estimator} to $\cs_{H_n}(z)$.
	\end{definition}

	It will be seen in Lemma~\ref{Lemma_PopStilEst_Uniqueness} that the population Stieltjes transform estimator is unique.\\
	The above definition may be extended canonically to $z \in \C^-$ by $\hat{s}_n(z) = \ol{\hat{s}_n(\ol{z})}$.
	\\[0.5em]
	The Stieltjes transform $\cs_{H_n}$ determines population linear spectral statistics by the formula
	\begin{align*}
		& \frac{-1}{2\pi \bm{i}} \oint_\gamma f(z) \cs_{H_n}(z) \, dz = \int_{\R \cap \operatorname{int}(\gamma)} f(\lambda) \, dH_n(\lambda)
	\end{align*}
	for any closed curve $\gamma$ enclosing (while not intersecting) a portion of $\supp(H_n)$ with a counter-clockwise orientation and $f \in \operatorname{Hol}(\gamma)$. This naturally suggests estimators of the form
	\begin{align*}
		& \frac{-1}{2\pi \bm{i}} \oint_\gamma f(z) \hat{s}_n(z) \, dz
	\end{align*}
	for population linear spectral statistics, provided that $\hat{s}_n(z)$ exists for all $z \in \C^+ \cap \operatorname{image}(\gamma)$ and $z \in \ol{\C^- \cap \operatorname{image}(\gamma)}$. This will be made precise in Subsection \ref{Subsection_MainResults_PLSS}.
	
	\subsection{Overview of the literature}
	\subsubsection{Marchenko--Pastur laws and extensions}
	The Marchenko--Pastur law was established for $B_n = \Id_d$ and $(\bm{X}_n)_{j,k} \sim \mathcal{N}(0,1)$ in \cite{MP_original}. The generalization to arbitrary centered i.i.d. entries of $\bm{X}_n$ with variance one was achieved in \cite{MP_Yin} under mild moment-conditions on $H_\infty$. The limiting spectral distribution of $A + \frac{1}{n} \bm{X}_n^* T \bm{X}_n$ for a deterministic matrix $A$ and possibly non-positive-semidefinite $T$ was first characterized by \cite{MP_Bai}. The assumption of independence between rows of $\bm{X}_n$ was weakened in \cite{MP_RowCorrelation_Bai} and, in the isotropic case $B_n = \Id_d$, and \cite{MP_Heiny} allows for correlations between rows and columns of $\bm{X}_n$ provided they go to zero sufficiently quickly as $n \to \infty$. A series of papers \cite{MP_NessSuff_Yaskov, MP_NessSuff_Nina, MP_NessSuff_Yao} deal with necessary and sufficient conditions for the Marchenko--Pastur law to hold in the isotropic case. The paper \cite{StrongMP} relaxes the assumption (\ref{Eq_RMT_Fundamental_Assumption}) and the data matrix $\bm{Y}_n$ is allowed to have more general independent columns, while still assuming the covariance matrices of said columns to be simultaneously diagonalizable. Marchenko--Pastur laws arising in the context of high-dimensional time series are studied in the papers \cite{TimeMP4, TimeMPLast, TimeMP1, Aue, TimeMP2}.
	\\[0.5em]
	Local laws present a quantitative generalization of Marchenko--Pastur laws. They describe the behavior of the Stieltjes transforms $\cs_{\hat{\nu}_n}(z)$ with precision depending on the position $z$ relative to the support of the LSD $\nu_\infty$,  allowing for more detailed analysis of eigenvalues at the edge of the spectrum, such as the largest or smallest eigenvalues. Among the most influential and comprehensive works on local laws for sample covariance matrices are \cite{BloemendalIsotropicLocalLaws} in the isotropic case and \cite{KnowlesAnisotropicLocalLaws} in the anisotropic case.
	
	\subsubsection{Spectral CLTs}
	A well-known phenomenon of high-dimensional random matrix theory is the prevalence of fast convergence rates of order $\frac{1}{n}$ instead of $\frac{1}{\sqrt{n}}$. Similarly to the standard central limit theorem (CLT), one can for a measurable function $g:\R \rightarrow \R$ observe the difference between the empirical integral $\int_\R g \, d\hat{\nu}_n$ and the limiting integral $\int_\R g \, d\nu_\infty$. Spectral central limit theorems (spectral CLTs) for certain functions $g$ describe the weak convergence of the rescaled difference
	\begin{align*}
		& n \Big( \int_\R g \, d\hat{\nu}_n - \int_\R g \, d\nu_n \Big)
	\end{align*}
	to Gaussian distributions. 
	The earliest high-dimensional spectral CLT for the setting $T_n = \Id_d$ and $X_{j,k} \sim \mathcal{N}(0,1)$ goes back to \cite{JonssonCLT} by Jonsson. In the celebrated paper \cite{BaiCLT}, Bai and Silverstein first formulated a spectral CLT for general $\Sigma_n$ and $B_n = \Sigma_n^{\frac{1}{2}}$ in the case where the entries of $\bm{X}_n$ are i.i.d. and have the same fourth moment as a standard normal distribution. The latter condition was removed, and the class of functions $g$ for which the CLT holds was expanded by Najim and Yao in \cite{NajimYao}. In the paper \cite{CLTniceCov}, improved formulas for the limiting mean and covariance are given. Generalizations of the spectral CLT to the case $\frac{d}{n} \to \infty$ or to columns taken from a high-dimensional time series were developed in \cite{DetteSequential, YaoUltraCLT}. Novel spectral CLTs for generalized linear spectral statistics of the form $\frac{1}{d} \tr\big( f(\bm{S}_n) g(\Sigma_n) \big)$ were obtained in \cite{CLT_generalizedLSS}.
	
	\subsubsection{Eigen-inference}
	Eigen-inference refers to the inference of properties of population eigenvalues $\big(\lambda_j(\Sigma_n)\big)_{j \leq d}$ from the observable sample eigenvalues $\big(\lambda_j(\bm{S}_n)\big)_{j \leq d}$ or data matrix $\bm{Y}_n$. Eigen-inference methods may: a) estimate the population eigenvalues directly, b) construct measures $\hat{H}_n$, which attempt to approximate $H_n$, c) estimate the Stieltjes transforms $\cs_{H_n}(z)$, or d) derive estimators for population linear spectral statistics $L_{n}(g)$ for various functions $g: \R \rightarrow \C$.
	\\[0.5em]
	An early work on the estimation of $H_n$ by solving a convex optimization problem is \cite{MPIKaroui}. El Karoui proves consistency of the resulting estimator $\hat{H}_n$ in the sense that the weak convergence $\hat{H}_n \xRightarrow{n \to \infty} H_\infty$ holds with probability one, but does not provide bounds for the rate of convergence. Building on \cite{Speicher2008}, Bai, Chen and Yao in \cite{MPIBaiChenYao} construct a moment-based estimator under the assumption $H_\infty = t_1\delta_{\theta_1}+\dots+t_k\delta_{\theta_k}$ for model parameters $(t_1,\dots,t_k,\theta_1,\dots,\theta_k)$. They were also able to show asymptotic normality of the estimation error with rate $\frac{1}{n}$. Further work on parametric models of this type was done in \cite{MPILiBaiYao} and the textbook \cite{MPIYaoZhengBai}. To the best of the author's knowledge, there are no previous results achieving an error rate of $\mathcal{O}(\frac{1}{n})$ in a nonparametric setting.
	\\[0.5em]
	The papers \cite{LedoitWolf2} and \cite{LedoitWolfNumerical} by Ledoit and Wolf present a minimization algorithm which solves the discretized Marchenko--Pastur equation numerically. The corresponding argmin-estimators $(\hat{\lambda}_j)_{j \leq d}$ for $(\lambda_j(\Sigma_n))_{j \leq d}$ are shown to satisfy the consistency property
	\begin{align*}
		& \frac{1}{d} \sum\limits_{j=1}^d \big( \hat{\lambda}_j - \lambda_j(\Sigma_n) \big)^2 \xrightarrow{n \to \infty}_{\text{a.s.}} 0 \ .
	\end{align*}
	For estimation accuracy, Ledoit and Wolf's method and its extensions are widely regarded as state-of-the-art in high-dimensional population eigenvalue estimation. To the best of the author's knowledge, to date no theoretically guaranteed error rates are known for any method of estimating population eigenvalues directly.
	\\[0.5em]
	In \cite{LedoitWolf_OptShrink}, Ledoit and Wolf forgo the estimation of $H_n$ entirely in favor of estimating the deterministic equivalent $\nu_n$, which, by the same connection explored in Lemma~\ref{Lemma_SpaceTransform}, has applications to optimal shrinkage of the sample covariance matrix.
	\\[0.5em]
	The paper \cite{DobribanSpectrode} by Dobriban specializes in fast estimation of the population eigenvalues, while \cite{MPIKong} by Kong and Valiant develops a new ansatz for estimating the moments of the population distribution $H_n$ from the data matrix $\bm{Y}_n$. An achievement of \cite{MPIKong} was the derivation of explicit convergence rates for the estimation error in the nonparametric setting. These rates depend on the moment to be estimated, but are bounded from below by $\mathcal{O}(\frac{1}{\sqrt{n}})$.
	\\[0.5em]
	With formulas from free probability, Arizmendi, Tarrago and Vargas in \cite{MPIDeconvolution} develop methods for inverting free convolutions, which may be applied to the calculation of the limiting population distribution $H_\infty$ from the limiting distribution $\nu_\infty$. While an asymptotic theory for the resulting finite-sample estimators is currently missing, the high generality of their setting makes this direction a promising area of research.
	
	\subsection{Inherent limitations to the estimation of $H_n$}
	This subsection highlights two limitations that are intrinsic to the use of the
	Marchenko--Pastur equation for population eigenvalue inference. First, the
	equation only gives direct access to the population Stieltjes transform
	$\cs_H(z)$ on the domain $\bD_{H,c}(\infty)$. Values of $\cs_H$ outside this
	domain do not enter the equation in a directly invertible way. Second, the
	numerical example below suggests that, for methods which aim to estimate the
	entire population spectral distribution $H_n$, substantially faster rates than
	$\mathcal{O}(1/\log n)$ may in general be difficult to obtain.
	
	
	\begin{remark}[Fundamental limits of the Marchenko--Pastur equation for inference]\
		\\
		The Marchenko--Pastur equation from Lemma \ref{Lemma_MP} may be equivalently re-arranged to show that the Stieltjes transform $\cs_{\nu}$ is uniquely defined by
		\begin{align*}
			& \forall z \in \C^+ : \ \cs_{\nu}(\tz) (1-c\tz\cs_{\nu}(\tz)-c) = \cs_H\Big( \frac{\tz}{1-c\tz\cs_{\nu}(\tz)-c} \Big)
		\end{align*}
		and $\Im\big( c \cs_{\nu}(\tz) + \frac{c-1}{\tz} \big)>0$, which demonstrates that $\cs_{\nu}$ is directly determined through the values $\cs_H(z)$ for
		\begin{align*}
			& z \in \Big\{ \frac{\tz}{1-c\tz\cs_{\nu}(\tz)-c} \ \Big| \ \tz \in \C^+ \Big\} \overset{\text{(\ref{Eq_Def_varphi})}}{=} \bD_{H,c}(\infty) \ .
		\end{align*}
		Conversely, knowledge of $\nu$, will through the Marchenko--Pastur equation itself only directly yield $\cs_H(z)$ for points $z \in \bD_{H,c}(\infty)$. The characterization
		\begin{align*}
			& \bD_{H,c}(\infty) = \Big\{ z \in \C^+ \ \Big| \ 1 > c\int_\R \frac{\lambda^2}{|\lambda-z|^2} \, dH(\lambda) \Big\}
		\end{align*}
		(cf. proof of Lemma \ref{Lemma_SpaceTransform}) shows that
		$\bD_{H,c}(\infty)$ will (except for some edge-cases) be bounded away from $\supp(H)$.
		\\[0.5em]
		Existing methods for the estimation of $H$ (cf. \cite{MPIKaroui, LedoitWolf2}) work around this restriction by employing minimization algorithms and are ultimately consistent, since knowledge of $\cs_H(z)_{z \in U}$ for any open $U \subset \C^+$ is already sufficient to uniquely determine $H$. To-date, no error rates for such methods have been shown and the following numerical example suggests that error rates below $\mathcal{O}\big( \frac{1}{\log(n)} \big)$ might not be achievable when estimating the cumulative distribution functions (CDF) of $H$.
	\end{remark}
	
\begin{figure}[H]
		\begin{center}
			\includegraphics[width=\textwidth]{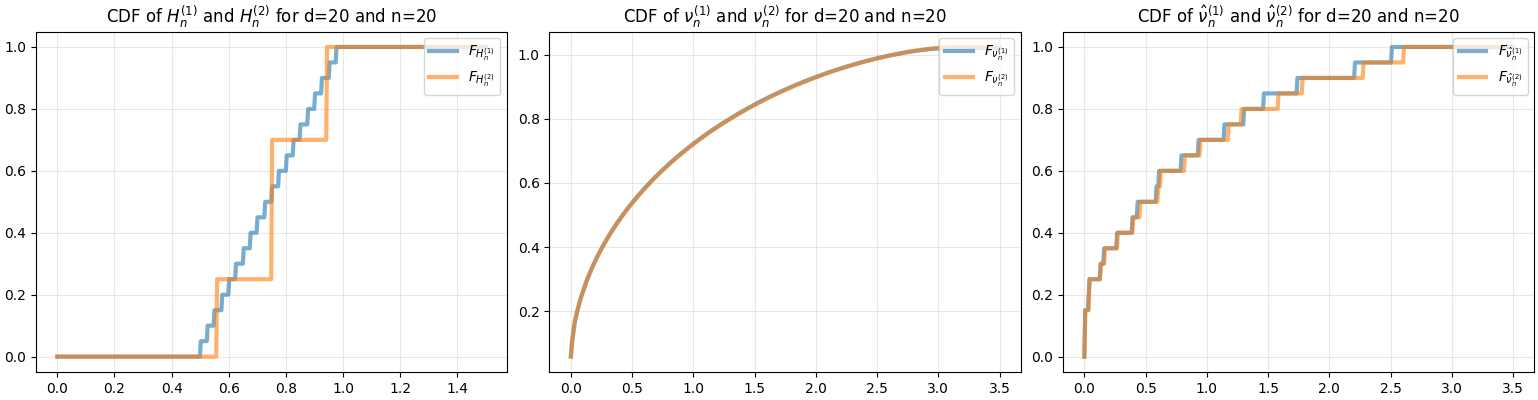}\\
			\includegraphics[width=\textwidth]{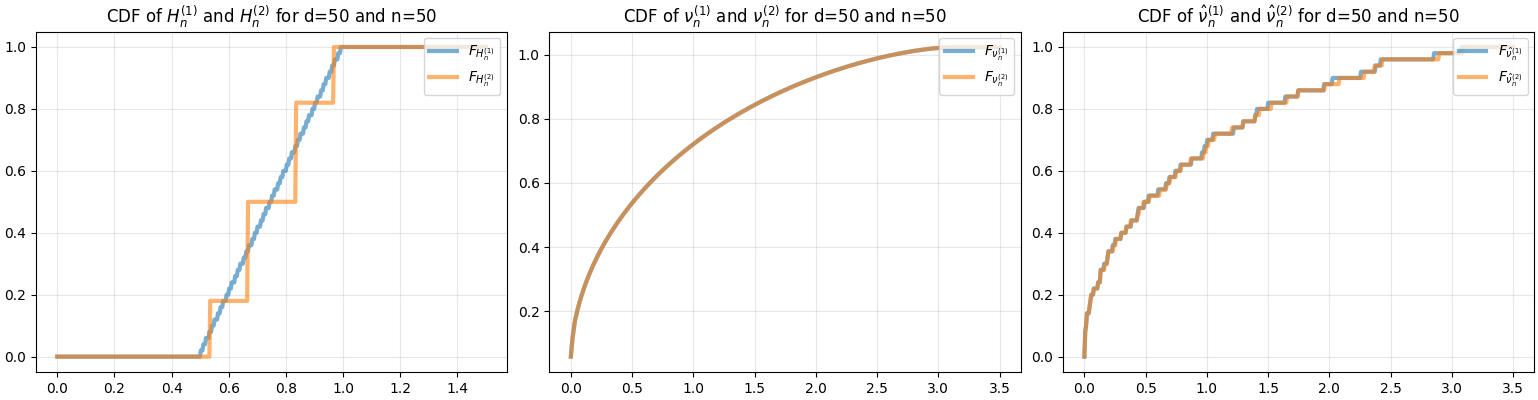}
		\end{center}
		\vspace{-0.4cm}
		\caption{\footnotesize Cumulative distributions functions of $H_n^{(i)}$ (left), $\nu_n^{(i)}$ (middle) and $\hat{\nu}_n^{(i)}$ (right) for $i=1$ (blue) and $i=2$ (orange). The values of $d=n$ are $20$ (top) and $50$ (bottom).
		}\label{Fig_CDFs}
		\vspace{-0.2cm}
	\end{figure}

\begin{figure}[H]
		\begin{center}
			\includegraphics[width=\textwidth]{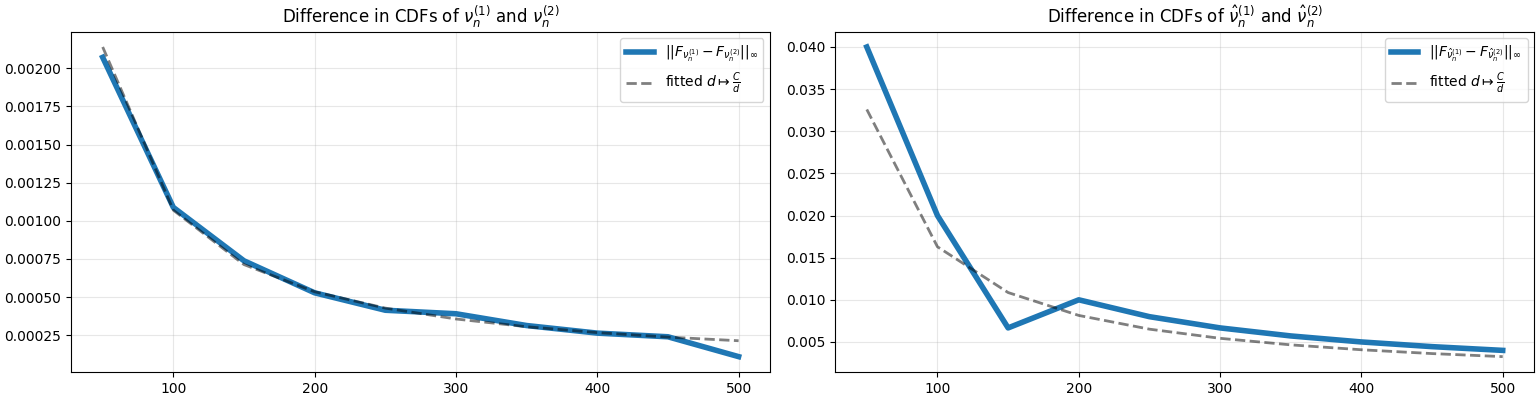}
		\end{center}
		\vspace{-0.4cm}
		\caption{\footnotesize Differences $\sup\limits_{x \in \R} |F_{\nu_n^{(1)}}(x)-F_{\nu_n^{(2)}}(x)|$ (left, blue) and $\sup\limits_{x \in \R} |F_{\hat{\nu}_n^{(1)}}(x)-F_{\hat{\nu}_n^{(2)}}(x)|$ (right, blue) for $d=n$ ranging from $50$ to $500$. Fitted curves $d \mapsto \frac{C}{d}$ are added as dashed gray lines.}\label{Fig_CDF_rates}
	\end{figure}

	\begin{example}[Very similar sample spectral distributions]\
		\\
		Let $\Sigma_n^{(1)} = \diag\big(\frac{1+1/d}{2},\ldots,\frac{1+d/d}{2}\big)$, such that the associated population spectral measure $H_n^{(1)}$ closely approximates the uniform distribution on $[\frac{1}{2},1]$. For $N=\lceil \log n\rceil$, let $x_1,\ldots,x_N$ denote the Gauss--Legendre nodes scaled to the interval $[\frac{1}{2},1]$ and let $w_1,\ldots,w_N>0$ be the corresponding weights. Define
		\begin{align*}
			& \Sigma_n^{(2)}=\diag\big(\underbrace{x_1,\ldots,x_1}_{\times [ w_1*d ]} \ ; \ \underbrace{x_2,\ldots,x_2}_{\times [ w_1*d ]} \ ; \ \ldots \ ; \ \hspace{-0.4cm} \underbrace{x_N,\ldots,x_N}_{\times (d - [ w_1*d ] - \ldots - [ w_{N-1}*d ])} \hspace{-0.4cm} \big)
		\end{align*}
		such that the associated
		spectral measure $H_n^{(2)}$ approximates
		$2\sum\limits_{k=1}^N w_k \delta_{x_k}$.
		\\[0.5em]
		The two population spectral measures $H_n^{(1)}$ and $H_n^{(2)}$ constructed in this way
		satisfy
		\begin{align*}
			& \sup_{x\in\R}\big|F_{H_n^{(1)}}(x)-F_{H_n^{(2)}}(x)\big|
			\asymp \frac{1}{\log n} \ .
		\end{align*}
		Figure~\ref{Fig_CDFs} illustrates that the corresponding empirical spectral distributions $\hat{\nu}_n^{(1)}$ and $\hat{\nu}_n^{(2)}$, as well as their deterministic equivalents
		$\nu_n^{(1)}$ and $\nu_n^{(2)}$, quickly become indistinguishable for
		$n\to\infty$. Additionally, Figure~\ref{Fig_CDF_rates} suggests that the empirical CDF distance $\sup\limits_{x \in \R} |F_{\hat{\nu}_n^{(1)}}(x)-F_{\hat{\nu}_n^{(2)}}(x)|$ appears to fall with rate $\mathcal{O}(1/n)$.
	\end{example}
	
	\subsection{Organization of the paper}
	The remainder of the paper is organized as follows. Section~\ref{Section_MainResults} presents the main results and is divided into four subsections. Subsection~\ref{Subsection_MainResults_ExistenceConsistency} gives first results on existence and basic properties of the population Stieltjes transform estimator, before Subsections~\ref{Subsection_MainResults_StatisticalGuaranteesForInference} and~\ref{Subsection_MainResults_PLSS} develop a method for the inference of population linear spectral statistics (PLSS) with error rate $\mathcal{O}(n^{\varepsilon-1})$. Subsection~\ref{Subsection_MainResults_SpectralCLTs} presents a method for inverting spectral CLTs into CLTs for the error of the population Stieltjes transform estimator, which is also applied to derive a CLT for the PLSS estimation error in the Gaussian case. Section \ref{Section_GLSS} introduces generalized linear spectral statistics and presents a method of estimation with error rate $\mathcal{O}(n^{\varepsilon-\frac{1}{2}})$.
	\\[0.5em]
	Section~\ref{Section_ProofsOfLemmas} proves the introductory results of Subsection~\ref{Subsection_MainResults_ExistenceConsistency}, while Section~\ref{Section_ConsistencyProof} proves Theorem~\ref{Thm_Consistency}, the main contribution of this paper to the field of eigen-inference. Section~\ref{Section_Numerical} presents a numerical implementation of the method described in Subsections~\ref{Subsection_MainResults_StatisticalGuaranteesForInference} and~\ref{Subsection_MainResults_PLSS}, while comparing results with four commonly used eigen-inference methods.
	\\[0.5em]
	The proofs of Lemmas~\ref{Lemma_SpaceTransform} and~\ref{Lemma_MPI}, the proofs of Theorems~\ref{Thm_CLT_Inversion} and \ref{Thm_GLSS_consistency}, as well as the proofs of Corollaries~\ref{Cor_PLSS_Consistency} and~\ref{Cor_GaussPLSSCLT} are deferred to the appendix, which also contains a list of symbols.

	\section{Main results}\label{Section_MainResults}
	
	\begin{assumption}[Main assumptions]\label{Assumption_EigInf_Main}\
		\begin{enumerate}[label=(A\arabic*)]
			\item\label{EI_ItemAssumption_Asymptotics}
			Suppose $d$ and $n$ tend to infinity simultaneously such that
			\begin{align}\label{Eq_Assumption0_AsymptoticRegime}
				&c_n \coloneq \frac{d}{n} \to c_{\infty} > 0 \ .
			\end{align}
			
			\item\label{EI_ItemAssumption_X_Structure}
			Suppose the sample covariance matrix is of the form
			\begin{align}\label{Eq_Assumption-2_ShapeS}
				& \bm{S}_n \coloneq \frac{1}{n} B_n \bm{X}_n \bm{X}_n^* B_n^*
			\end{align}
			for a $(d\times d)$ matrix $B_n$ with $B_nB_n^* = \Sigma_n$ and a random $(d \times n)$ matrix $\bm{X}_n$ with independent entries that satisfy
			\begin{align}\label{Eq_Assumption1_EntriesBasic}
				& \E[(\bm{X}_n)_{j,k}] = 0 \ \ \text{ and } \ \ \E[|(\bm{X}_n)_{j,k}|^2] = 1 \ .
			\end{align}
			
			\item\label{EI_ItemAssumption_PopConv}
			Suppose the weak convergence
			\begin{align}\label{Eq_Assumption2_Hconvergence}
				& H_n \coloneq \frac{1}{d} \sum\limits_{j=1}^d \delta_{\lambda_j(\Sigma_n)} \xRightarrow{n \to \infty} H_\infty
			\end{align}
			holds for a probability measure $H_\infty \neq \delta_0$ with compact support on $[0,\infty)$.
			
			\item\label{EI_ItemAssumption_sigmaBound}
			Suppose there exists a constant $\sigma^2 > 0$ such that
			\begin{align}\label{Eq_Assumption4_sigmaBound}
				& \forall n \in \N: \ \|\Sigma_n\| \leq \sigma^2 \ .
			\end{align}
			
			\item\label{EI_ItemAssumption_MomentBound}
			Suppose that for every $p \in \N$ there exists a constant $C_p>0$ such that
			\begin{align}\label{Eq_Assumption5_BoundedMoments}
				& \forall n \in \N, \, j\leq d, \, k \leq n : \ \E\big[ |(\bm{X}_n)_{j,k}|^p \big] \leq C_p \ .
			\end{align}
			
		\end{enumerate}
	\end{assumption}
	
	\begin{remark}[Discussion of assumptions]\
		\\
		Assumptions~\ref{EI_ItemAssumption_Asymptotics}--\ref{EI_ItemAssumption_PopConv} are standard in the field of random matrices (see e.g. \cite{MPIBaiChenYao, BaiCLT, KnowlesAnisotropicLocalLaws}). Assumption~\ref{EI_ItemAssumption_sigmaBound} also appears in most works on eigen-inference, including \cite{DingFan_SpikedShrinkage, MPIKaroui, MPIKong, LedoitWolf2}. 
		Lastly, assumption~\ref{EI_ItemAssumption_MomentBound} is stronger than moment assumptions needed for some other eigen-inference methods (for example \cite{MPIKong} only requires finite fourth moments), however, it is needed only for the statistical guarantees on inference from Subsections~\ref{Subsection_MainResults_StatisticalGuaranteesForInference} and~\ref{Subsection_MainResults_PLSS}. The results of Lemma~\ref{Lemma_ConsistencyBasic} and Theorem~\ref{Thm_CLT_Inversion} hold under minimal moment assumptions.
	\end{remark}
	
	\subsection{Existence and consistency of the population Stieltjes transform estimator}\label{Subsection_MainResults_ExistenceConsistency}
	
	The first lemma is only concerned with properties of the map $(z,c_n,\hat{\nu}_n) \mapsto \hat{s}_n(z)$ and its statements hold for any $c_n>0$ and probability measure $\hat{\nu}_n$ with compact support on $[0,\infty)$. None of the items from Assumption \ref{Assumption_EigInf_Main} are required at this point.
	
	\begin{lemma}[Existence and uniqueness of $\hat{s}_n(z)$]\label{Lemma_PopStilEst_Uniqueness}\
		\\
		The following statements on the population Stieltjes transform estimator $\hat{s}_n(z)$ from Definition~\ref{Def_StieltjesEstimator} hold deterministically for any $c_n>0$ and realization $\hat{\nu}_n$, which is by construction a probability measure with compact support on $[0,\infty)$.
		\begin{itemize}
			\item[a)] The population Stieltjes transform estimator $\hat{s}_n(z)$ is uniquely defined.
			
			\item[b)] The set $U = U(\hat{\nu}_n,c_n) = \{z \in \C^+ \mid \hat{s}_n(z) \text{ exists}\}$ is open and contains $\C^+\setminus B_{\kappa}^{\C^+}(0)$ for some $\kappa = \kappa(\hat{\nu}_n,c_n) > 0$.
			
			\item[c)] For $U$ as in (b), the map $\hat{s}_n : U \rightarrow \C^+$ is holomorphic.
		\end{itemize}
	\end{lemma}
	
	Next, the consistency of the estimator $\hat{s}_n(z)$ on the domain $\bD_{H_\infty,c_\infty}(1)$ is shown to be true, whenever the Marchenko--Pastur law holds.
	
	\begin{lemma}[Consistency of $\hat{s}_n(z)$]\label{Lemma_ConsistencyBasic}\
		\\
		Suppose~\ref{EI_ItemAssumption_Asymptotics}--\ref{EI_ItemAssumption_sigmaBound} of Assumption~\ref{Assumption_EigInf_Main} hold. If the Marchenko--Pastur law:
		\begin{align}\label{Eq_MP_Basic}
			& 1 = \bP\big( \hat{\nu}_n \xRightarrow{n \to \infty} \nu_\infty \big) \ ,
		\end{align}
		holds, the Stieltjes transform estimator is consistent in the sense
		\begin{align}\label{Eq_ConsistencyBasic}
			& 1 = \bP\Big( \forall J \subset \bD_{H_\infty,c_\infty}(1) \text{ compact, } \, \exists N_J>0 \, \forall n \geq N_J : \nonumber\\
			& \hspace{2cm} \hat{s}_{n} \text{ exists on $J$ and } \ \ \sup\limits_{z \in J} \big|\hat{s}_n(z) - \cs_{H_\infty}(z)\big| \xrightarrow{n \to \infty} 0 \Big) \ .
		\end{align}
	\end{lemma}
	
\begin{figure}[H]
		\centering
		\includegraphics[width=0.325\textwidth]{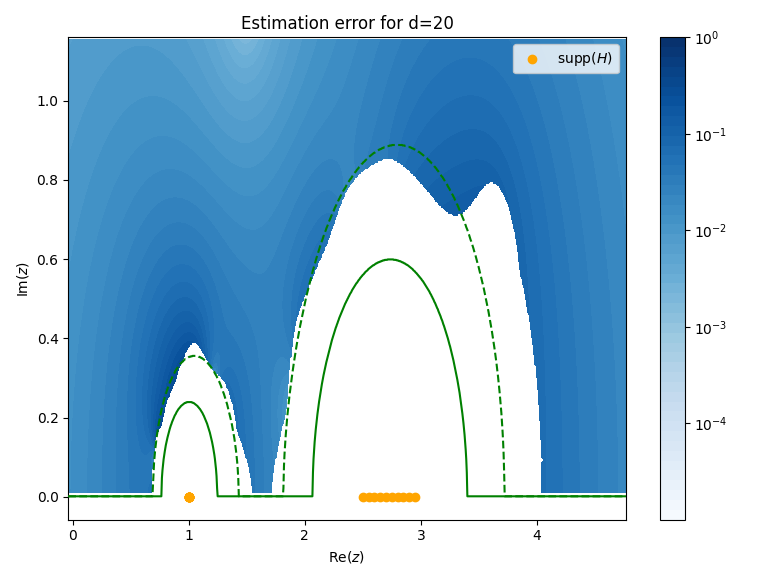} \includegraphics[width=0.325\textwidth]{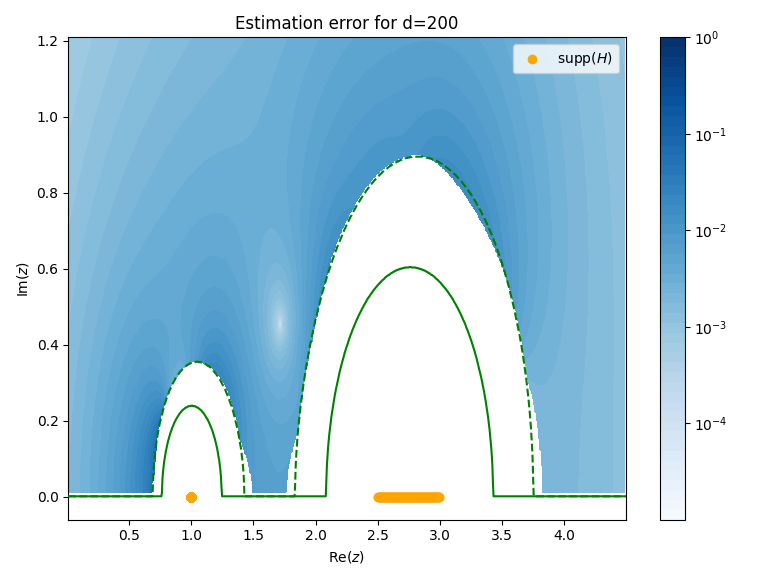} \includegraphics[width=0.325\textwidth]{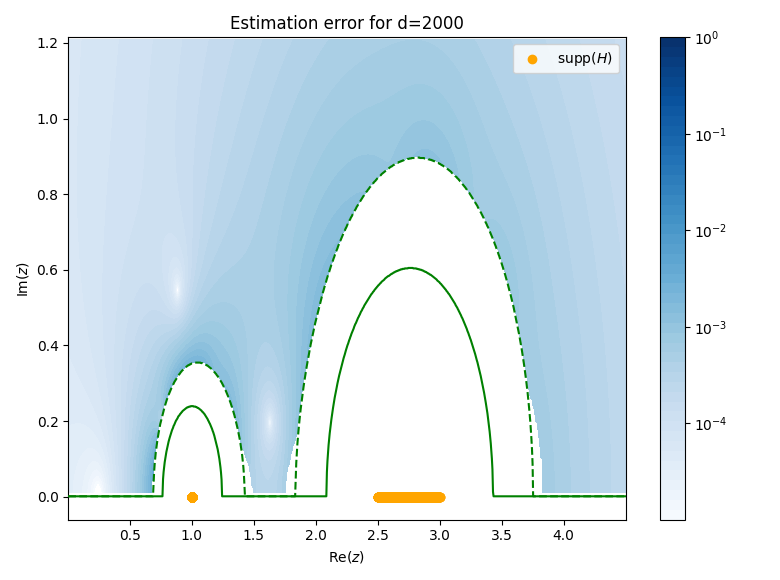}
		
		\vspace{-0.2cm}
		\caption{Three contour plots of the estimation error $|\hat{s}_n(z)-\cs_{H_n}(z)|$, where $\hat{s}_n(z)$ exists (blue), based on the example data (Ex.1) given in Definition~\ref{Def_Examples}. The dimension $d$ takes values $20$ (left), $200$ (middle) and $2000$ (right), while $n = 10d$ such that $c_n = \frac{1}{10}$. The color scale is logarithmic and the boundaries of $\bD_{H_n,c_n}(\infty)$ (solid) and $\bD_{H_n,c_n}(1)$ (dashed) are added as green lines.}\label{Fig_StieltjesEstimationError}
	\end{figure}

	Figure \ref{Fig_StieltjesEstimationError} illustrates the validity of Lemma \ref{Lemma_ConsistencyBasic}. The region in which $\hat{s}_n$ may be found numerically is seen to cover $\bD_{H_n,c_n}(1)$ almost entirely for large $n$, while the estimation error is seen to fall quickly with growing $d$ and $n$.
	
	\subsection{Inference of population Stieltjes transforms}\label{Subsection_MainResults_StatisticalGuaranteesForInference}
	
	The following theorem establishes the consistency of $\hat{s}_n$ with a rate of $\mathcal{O}(n^{\varepsilon-1})$ on a large, mostly data-driven, spectral domain.
	
	\begin{theorem}[Consistency of $\hat{s}_n$ with rate $\mathcal{O}(n^{\varepsilon-1})$]\label{Thm_Consistency}\
		\\
		Suppose Assumption~\ref{Assumption_EigInf_Main} holds.
		Dependent on constants $\tau,\kappa > 0$,
		define the empirical spectral domain $\hat{\bD}(\tau,\kappa,n)$ as the set of all $z \in \C^+$ for which the conditions
		\begin{align}
			& \tau < |z| < \kappa \label{Eq_Consistency_EmpDomain1_1} \ , \\
			& \hat{s}_n(z) \text{ as in Def.~\ref{Def_StieltjesEstimator} exists} \label{Eq_Consistency_EmpDomain1_2}\\
			& \tau < \big| \hat{\Phi}_n(z)  \big| < \kappa \label{Eq_Consistency_EmpDomain2}\\
			& \tau < \dist\big( \hat{\Phi}_n(z), \{\lambda_j(\bm{S}_n) \mid j \leq d\} \big) \label{Eq_Consistency_EmpDomain3}\\
			& \tau < \dist\big( \hat{\Phi}_n(z), \supp(\nu_n) \big) \label{Eq_Consistency_EmpDomain4}
		\end{align}
		hold with the notation $\hat{\Phi}_n(z) \coloneq (1-c_nz\hat{s}_n(z)-c_n)z \in \C^+$.
		Then for every $\varepsilon \in (0,1)$ and $D>0$ there exists a constant $C=C(\tau,\kappa,\varepsilon,D)>0$ such that
		\begin{align}\label{Eq_Consistency_Result}
			& \bP\Big( \forall z \in \hat{\bD}(\tau,\kappa,n) : \ \big| \hat{s}_n(z) - \cs_{H_n}(z) \big| \leq n^{\varepsilon-1} \Big) \geq 1 - \frac{C}{n^D} \ .
		\end{align}
	\end{theorem}
	
\begin{figure}[H]
		\centering
		\includegraphics[width=0.325\textwidth]{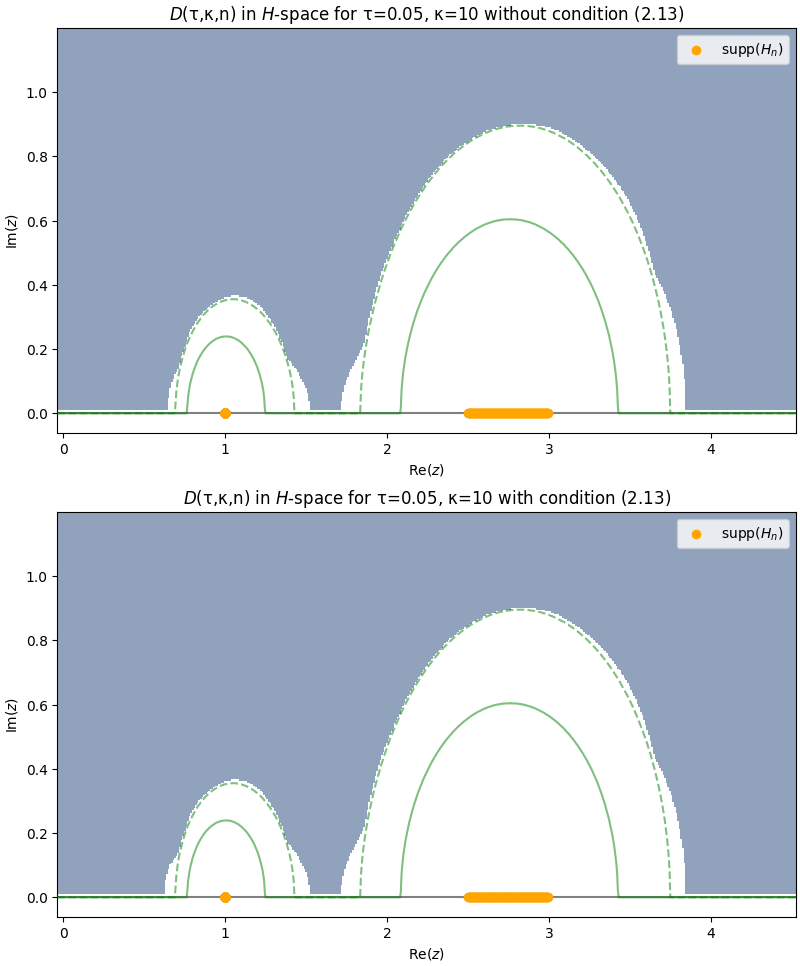} \includegraphics[width=0.325\textwidth]{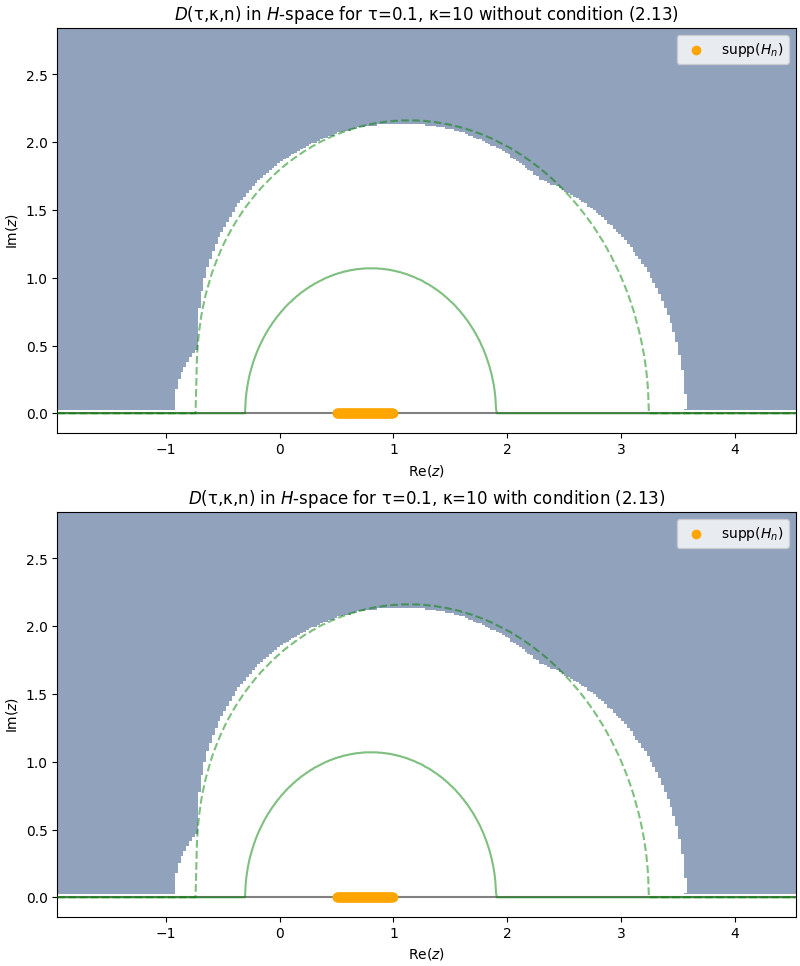} \includegraphics[width=0.325\textwidth]{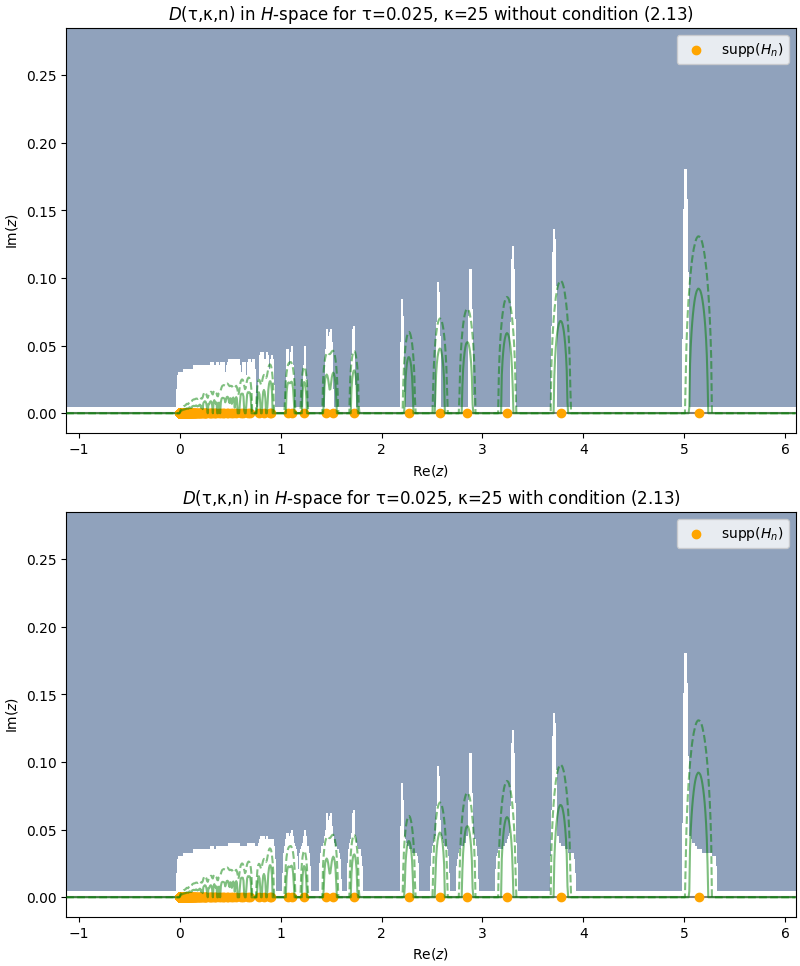}
		
		\vspace{-0.2cm}
		\caption{Plots of the set $\hat{\bD}(\tau,\kappa,n)$ (blue), \underline{without} enforcing condition (\ref{Eq_Consistency_EmpDomain4}) (top) and \underline{with} enforcing condition (\ref{Eq_Consistency_EmpDomain4}) (bottom) for specified values of $\tau,\kappa$ and example data matrices introduced in Definition~\ref{Def_Examples}.\\ Left: Example (Ex.1) with $d=200$. Middle: Example (Ex.2) with $d=100$. Right: Example (Ex.3) with $d=784$.\\
			The boundaries of $\bD_{H_n,c_n}(\infty)$ and $\bD_{H_n,c_n}(1)$ are marked green as in Figure~\ref{Fig_StieltjesEstimationError}.
		}\label{Fig_CondComparison}
		\vspace{-0.3cm}
	\end{figure}

	\begin{remark}[Replaceability of condition (\ref{Eq_Consistency_EmpDomain4})]\label{Remark_Replaceability}\
		\\
		A caveat of the above theorem is that condition (\ref{Eq_Consistency_EmpDomain4}) is not directly observable from data, but requires approximate knowledge of the support of the deterministic equivalent $\nu_n$.
		Under quite general assumptions on the structure of $H_n$, which amount to restricting the number of outlier eigenvalues, existing eigenvalue rigidity results such as Theorem 3.12 of \cite{KnowlesAnisotropicLocalLaws} or Theorem S.3.11 of \cite{DingSpikedRigidity} may be employed to remove condition (\ref{Eq_Consistency_EmpDomain4}) from the above theorem.
		\\
		If no assumptions on $\supp(\nu_n)$ or $H_n$ other than $\supp(H_n) \in [0,\sigma^2]$ are permissible, it follows from (b) of Lemma~\ref{Lemma_BasicStieltjesConvergence} that condition (\ref{Eq_Consistency_EmpDomain4}) may be replaced with
		\begin{align}
			& \dist\big(\hat{\Phi}_n(z) , \, [0,\sigma^2(1+\sqrt{c_n})^2]\big) \geq \tau \ . \label{Eq_NoFalsePositives_Assumption3_ReplacementSafe}
		\end{align}
		Figure \ref{Fig_CondComparison} suggests that condition (\ref{Eq_Consistency_EmpDomain4}) is only noticeable around outlier eigenvalues.
	\end{remark}

	\subsection{Inference of population linear spectral statistics}\label{Subsection_MainResults_PLSS}
\begin{figure}[H]
		\centering
		\includegraphics[width=0.325\textwidth]{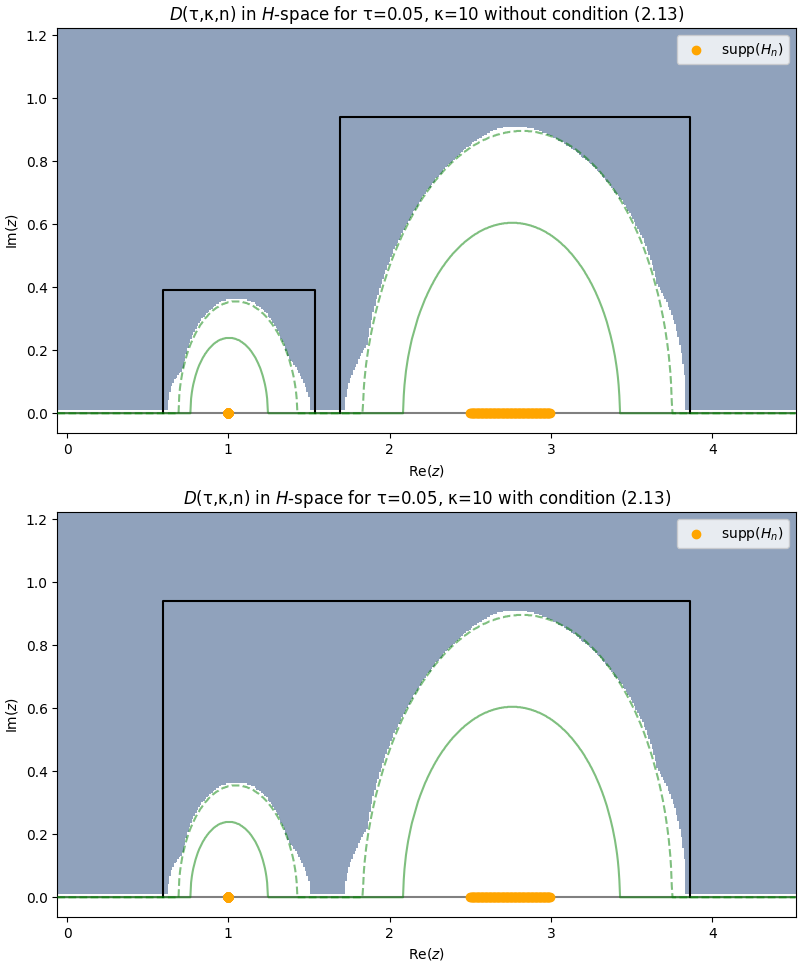} \includegraphics[width=0.325\textwidth]{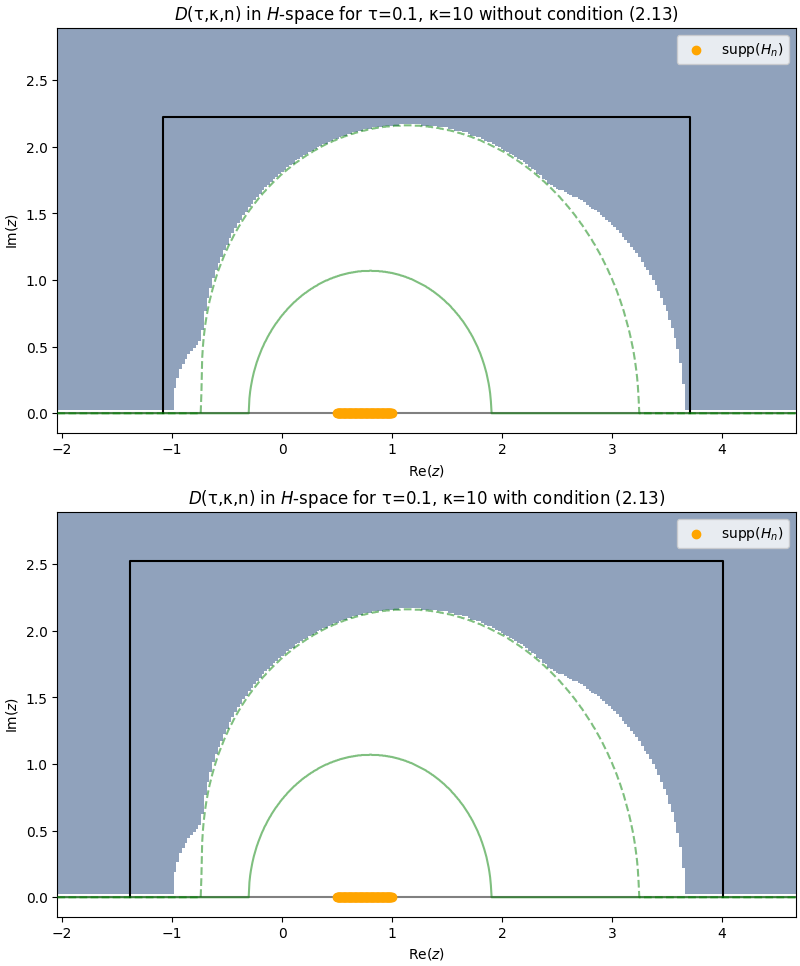} \includegraphics[width=0.325\textwidth]{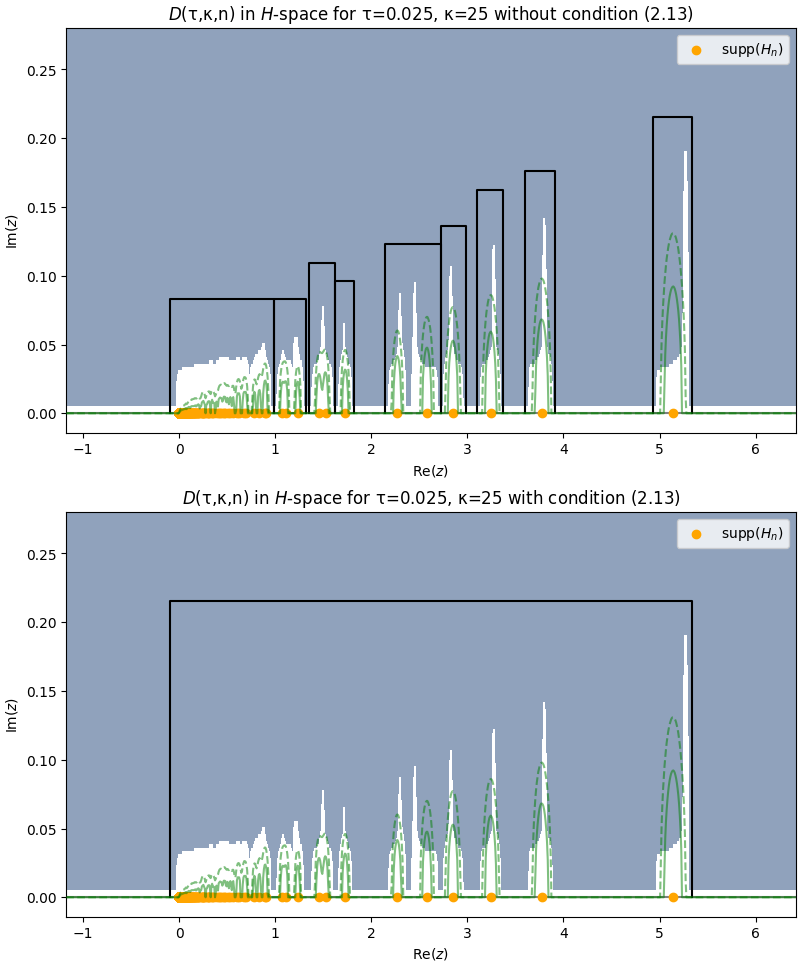}
		
		\vspace{-0.2cm}
		\caption{Plots of the set $\hat{\bD}(\tau,\kappa,n)$ (blue) \underline{with} enforcing condition (\ref{Eq_Consistency_EmpDomain4}) for specified values of $\tau,\kappa$ and example data matrices introduced in Definition~\ref{Def_Examples}. Left: Example (Ex.1) with $d=200$. Middle: Example (Ex.2) with $d=100$. Right: Example (Ex.3) with $d=784$.
		The boundaries of $\bD_{H_n,c_n}(\infty)$ and $\bD_{H_n,c_n}(1)$ are marked green as in Figure~\ref{Fig_StieltjesEstimationError}. Admissible curve candidates are added in black.}\label{Fig_CurveDiscovery}
	\end{figure}

	Figure~\ref{Fig_CurveDiscovery} illustrates the possibility of plotting curves $\gamma_n$ which surround a portion of $\supp(H_n)$ while remaining in the domain $\hat{\bD}(\tau,\kappa,n)$. For sufficiently analytic functions $g$, Cauchy's integral formula then allows for the estimation of population linear spectral statistics $\int g \, dH_n$, where the integral goes over the region of $\R$ encapsulated by the curve.
	This notion will be made precise in Corollary~\ref{Cor_PLSS_Consistency}.
	
	\begin{definition}[Admissible curves]\label{Def_AdmissibleCurve}\
		\\
		For given $\tau,\kappa>0$, 
		call a piecewise $C^1$-curve $\gamma_n : (0,1) \rightarrow \C^+$ \textit{admissible}, if
		\begin{align}\label{Eq_AdmissibleCurve_Cond1}
			& \mathrm{image}(\gamma_n) \subset \hat{\bD}(\tau,\kappa,n)
		\end{align}
		and
		\begin{align}\label{Eq_AdmissibleCurve_Cond2}
			& \gamma_n(t) \xrightarrow{t \nearrow 1} a_{\gamma_n} \in \R \ \ ; \ \ \gamma_n(t) \xrightarrow{t \searrow 0} b_{\gamma_n} \in \R
		\end{align}
		hold for some $a_{\gamma_n} < b_{\gamma_n}$ from $\R$.
	\end{definition}
	
	\begin{definition}[Population linear spectral statistic estimator]\label{Def_PLSS_estimator}\
		\\
		For given $\tau,\kappa>0$, let $\gamma_n$ be an admissible curve. 
		For any $g \in \operatorname{Hol}(\gamma_n)$, define the \textit{population linear spectral statistic estimator} (PLSS estimator) as
		\begin{align}\label{Eq_DefPLSS_estimator_old}
			\hat{L}_{n,\gamma_n}(g) \coloneq & \frac{-1}{2\pi \bm{i}} \int_{\gamma_n} g(z) \hat{s}_n(z) - g(\ol{z}) \ol{\hat{s}_n(z)} \, dz \ .
		\end{align}
	\end{definition}
	
	\begin{remark}[Homotopy invariance]\label{Remark_Homotopy}\
		\\
		Since the population Stieltjes transform estimator $\hat{s}_n$ is by Lemma \ref{Lemma_PopStilEst_Uniqueness} holomorphic, standard results on the homotopy invariance of curve integrals are applicable to the previously defined PLSS estimator. For example, if $g: \C \rightarrow \C$ is an entire function and $\gamma_{n,1},\dots,\gamma_{n,m}$ is a sequence of admissible curves with the property that every component of $\supp(H_n)$ is contained in precisely one interval $(a_{\gamma_{n,k}},b_{\gamma_{n,k}})$ (compare Figure \ref{Fig_CurveDiscovery}), then the sum $\hat{L}_{n,\gamma_{n,1}}(g) + \dots + \hat{L}_{n,\gamma_{n,m}}(g)$ does not depend on the precise choice of the sequence $\gamma_{n,1},\dots,\gamma_{n,m}$.
	\end{remark}
	
	Figure \ref{Fig_CurveDiscovery} shows some (algorithmically discovered) admissible curves. An algorithm for discovering admissible curves is included in the author's Github repository\footnote{\href{https://github.com/BenDeitmar/PLSS_estimation}{{https://github.com/BenDeitmar/PLSS\_estimation}}}.
	
	\begin{corollary}[Consistency of the PLSS estimator with rate $\mathcal{O}(n^{\varepsilon-1})$]\label{Cor_PLSS_Consistency}\
		\\
		Suppose Assumption~\ref{Assumption_EigInf_Main} holds and fix constants $\tau,\kappa, \varepsilon > 0$. On the high-probability event
		\begin{align*}
			& \forall z \in \hat{\bD}(\tau,\kappa,n) : \ \big| \hat{s}_n(z) - \cs_{H_n}(z) \big| \leq n^{\varepsilon-1}
		\end{align*}
		from Theorem~\ref{Thm_Consistency}, it holds that
		\begin{align}\label{Eq_PLSS_consistency_Result}
			& \hspace{0cm}  \forall \gamma_n \text{ admissible curve}, \, \forall g \in \operatorname{Hol}(\gamma_n) :\nonumber\\[-0.3cm]
			& \hspace{1cm} \bigg| \hat{L}_{n,\gamma_n}(g) - \int_{a_{\gamma_n}}^{b_{\gamma_n}} g \, dH_n \bigg| \leq n^{\varepsilon-1} \frac{\|g\|_{\gamma_n}}{\pi} \overbrace{\int_{\gamma_n} |dz|}^{= \text{length of }\gamma_n}
		\end{align}
		holds for all $n \in \N$.
	\end{corollary}
	
	The fact that the points $a_{\gamma_n}, b_{\gamma_n} \in \R$ may be chosen with $\gamma_n$ makes the previous result stronger than if only standard population linear spectral statistics $L_n(g) = \int_\R g \, dH_n$ were estimated. The necessary assumption that $g$ be holomorphic on the region enclosed by $\gamma_n$ is considerably weaker, if $\gamma_n$ may be chosen to enclose only certain portions of $\supp(H_n)$.
	

	\begin{remark}[Generalized linear spectral statistics]\label{Remark_GLSS}\
		\\
		Generalized linear spectral statistics (GLSS) are quantities of the form $\frac{1}{d} \tr\big( f(\bm{S}_n) g(\Sigma_n) \big)$ for some functions $f,g : \R \rightarrow \C$. As they encode information on the relationship between sample eigenvectors and population eigenvectors, GLSS are of particular interest for applications in random matrix theory, leading to their examination in \cite{LedoitPeche} or \cite{CLT_generalizedLSS}. Similarly to how the sample Stieltjes transforms $\cs_{\hat{\nu}_n}(\tz)$ are known to approximate their deterministic equivalents $\cs_{\nu_n}(\tz)$, it is under certain conditions known that
		\begin{align}\label{Eq_GLSS_DeterministicEquivalent}
			& \frac{1}{d} \tr\Big( (\bm{S}_n-\tz\Id_d)^{-1} g(\Sigma_n) \Big) \approx \frac{1}{d} \tr\Big( -\frac{1}{\tz} (\Id_d+\cs_{\ul{\nu}_n}(\tz)\Sigma_n)^{-1} g(\Sigma_n) \Big) \ ,
		\end{align}
		where the error of this approximation may vary between $\mathcal{O}(\frac{1}{n})$ and $\mathcal{O}(\frac{1}{\sqrt{n}})$ depending on the behavior of $g(\Sigma_n)$ (see \cite{CLT_generalizedLSS}). Since the right-hand side of (\ref{Eq_GLSS_DeterministicEquivalent}) is a population linear spectral statistic, it is estimateable by the methods introduced here. Section~\ref{Section_GLSS} goes into further detail and presents an analogue to Corollary~\ref{Cor_PLSS_Consistency} for generalized linear spectral statistics with error rate $\mathcal{O}(n^{\varepsilon-\frac{1}{2}})$.
	\end{remark}

	\subsection{Inversion of spectral CLTs into CLTs for the PLSS estimation error}\label{Subsection_MainResults_SpectralCLTs}
	
	\begin{theorem}[Inversion of spectral CLTs]\label{Thm_CLT_Inversion}\
		\\
		Suppose~\ref{EI_ItemAssumption_Asymptotics}--\ref{EI_ItemAssumption_sigmaBound} of Assumption~\ref{Assumption_EigInf_Main} hold, and that there exists a constant $C_4$ bounding the fourth moments of $\bm{X}_n$, i.e. $\E[|(\bm{X}_n)_{j,k}|^4] \leq C_4$ for all $n,j,k$. Further assume that a spectral CLT holds uniformly on $\tilde{U} = \Phi_{H_\infty,c_\infty}(U)$ for some open connected set $U \subset \bD_{H_\infty,c_\infty}(1)$, i.e. the convergence in distribution
		\begin{align}\label{Eq_CLTInversion_Cond2}
			& n\big( \cs_{\hat{\ul{\nu}}_n}(\tz) - \cs_{\ul{\nu}_n}(\tz) \big)_{\tz \in \tilde{U}} \xrightarrow{n \to \infty}_{\mathcal{D}} \tilde{G}(\tz)_{\tz \in \tilde{U}} \ \ \text{ in } \ell^{\infty}(\tilde{U})
		\end{align}
		holds for a Gaussian process $\tilde{G}(\tz)_{\tz \in \tilde{U}}$. Then for any compact subset $J \subset U$, the uniform convergence in distribution
		\begin{align}\label{Eq_CLTInversion_Result}
			& n\big( \hat{s}_n^{(0)}(z) - \cs_{H_n}(z) \big)_{z \in J} \xrightarrow{n \to \infty}_{\mathcal{D}} \Big(-\frac{\tilde{G}(\Phi_{H_\infty,c_\infty}(z))}{c_\infty z^2\cs_{\ul{\nu}_\infty}'(\Phi_{H_\infty,c_\infty}(z))}\Big)_{z \in J} \ \ \text{ in } \ell^{\infty}(J)
		\end{align}
		holds, where $\hat{s}_n^{(0)}(z) \coloneq \hat{s}_n(z)$ whenever the latter exists and $\hat{s}_n^{(0)}(z) \coloneq 0$ otherwise.
	\end{theorem}
	
	The assumption (\ref{Eq_CLTInversion_Cond2}) has a slightly different form from spectral CLTs in the literature, which are commonly of the form
	\begin{align}\label{Eq_CLTLiterature_Example1}
		& n\big( \cs_{\hat{\ul{\nu}}_n}(\tz) - \cs_{\ul{\nu}_n}(\tz) \big)_{\tz \in \mathcal{C}} \xrightarrow{n \to \infty}_{\mathcal{D}} \tilde{G}(\tz)_{\tz \in \mathcal{C}} \ \ \text{ in } \ell^\infty(\mathcal{C})
	\end{align}
	for some contour $\mathcal{C}$ surrounding $\supp(\ul{\nu}_\infty)$ (compare Lemma 1.1 of \cite{BaiCLT} and Theorem 1 of \cite{NajimYao}). The assumption (\ref{Eq_CLTInversion_Cond2}) of Theorem~\ref{Thm_CLT_Inversion} follows for any $\tilde{U}$ lying outside the contour $\mathcal{C}$ by the continuous mapping theorem.
	\\[0.5em]
	As an example illustrating the applicability of Theorem~\ref{Thm_CLT_Inversion}, the following corollary gives a CLT for the rescaled estimation error of the PLSS estimator, in a Gaussian setting.
	
	\begin{corollary}[CLT for the PLSS estimation error]\label{Cor_GaussPLSSCLT}\
		\\
		Suppose Assumption~\ref{Assumption_EigInf_Main} holds and that additionally the entries of $\bm{X}_n$ are i.i.d. real ($\beta=1$) or complex ($\beta=2$) standard normal.
		Let $\gamma : (0,1) \rightarrow \bD_{H_\infty,c_\infty}(1)$ be a piecewise $C^1$-curve with
		\begin{align}\label{Eq_GaussPLSSCLT_Cond1}
			& \lim\limits_{t \nearrow 1} \gamma(t) \in \big(-\infty, \hspace{-0.3cm} \underbrace{\varphi_{\nu_\infty,c_\infty}(-\delta)}_{\substack{\text{exists by (b) of Lemma~\ref{Lemma_StandardBounds}} \\ \in \R \text{ by (\ref{Eq_Def_varphi})}}} \hspace{-0.3cm} \big) \ \text{ and } \ \lim\limits_{t \searrow 0} \gamma(t) \in \big(\underbrace{\varphi_{\nu_\infty,c_\infty}(\sigma^2(1+\sqrt{c_\infty})^2+\delta)}_{\substack{\text{exists by (b) of Lemmas~\ref{Lemma_BasicStieltjesConvergence} and~\ref{Lemma_StandardBounds}} \\ \in \R \text{ by (\ref{Eq_Def_varphi})}}},\infty\big)
		\end{align}
		for some $\delta>0$.
		For any $g \in \operatorname{Hol}(\gamma)$, the \textit{population linear spectral statistic (PLSS) estimator}
		\begin{align}\label{Eq_DefPLSS_estimator}
			\hat{L}_{n,\gamma}^{(0)}(g) \coloneq & \frac{-1}{2\pi \bm{i}} \int_{\gamma} g(z) \hat{s}_n^{(0)}(z) - g(\ol{z}) \ol{\hat{s}_n^{(0)}(z)} \, dz
		\end{align}
		satisfies the convergence in distribution
		\begin{align}\label{Eq_GaussPLSSCLT_Result2}
			& n\Big( \hat{L}_{n,\gamma}^{(0)}(g) - \int_\R g \, dH_n \Big) \xrightarrow{n \to \infty}_{\mathcal{D}} Z(g)
		\end{align}
		for a (complex) Gaussian $Z(g)$. The mean and covariance structure of $Z(g)$ is given by
		\begin{align*}
			& \E[Z(g)] = \frac{-1}{2\pi \bm{i}} \int_\gamma g(z) \bm{e}(z) - g(\ol{z}) \bm{e}(\ol{z}) \, dz
		\end{align*}
		as well as
		\begin{align*}
			& \Cov[Z(g),\ol{Z(g)}]\\
			& = \frac{1}{4\pi^2} \int_{\gamma} \int_{\gamma} g(z_1) \ol{g(z_2)} \bm{c}(z_1,\ol{z}_2) \, dz_1 \, dz_2 - \frac{1}{4\pi^2} \int_{\gamma} \int_{\gamma} g(z_1) \ol{g(\ol{z}_2)} \bm{c}(z_1,z_2) \, dz_1 \, dz_2\\
			& \hspace{0.5cm} - \frac{1}{4\pi^2} \int_{\gamma} \int_{\gamma} g(\ol{z}_1) \ol{g(z_2)} \bm{c}(\ol{z}_1,\ol{z}_2) \, dz_1 \, dz_2 + \frac{1}{4\pi^2} \int_{\gamma} \int_{\gamma} g(\ol{z}_1) \ol{g(\ol{z}_2)} \bm{c}(\ol{z}_1,z_2) \, dz_1 \, dz_2\\
			& \Cov[Z(g),Z(g)]\\
			& = \frac{-1}{4\pi^2} \int_{\gamma} \int_{\gamma} g(z_1) g(z_2) \bm{c}(z_1,z_2) \, dz_1 \, dz_2 + \frac{1}{4\pi^2} \int_{\gamma} \int_{\gamma} g(z_1) g(\ol{z}_2) \bm{c}(z_1,\ol{z}_2) \, dz_1 \, dz_2\\
			& \hspace{0.5cm} + \frac{1}{4\pi^2} \int_{\gamma} \int_{\gamma} g(\ol{z}_1) g(z_2) \bm{c}(\ol{z}_1,z_2) \, dz_1 \, dz_2 - \frac{1}{4\pi^2} \int_{\gamma} \int_{\gamma} g(\ol{z}_1) g(\ol{z}_2) \bm{c}(\ol{z}_1,\ol{z}_2) \, dz_1 \, dz_2 \ ,
		\end{align*}
		where $\bm{e}(z) \coloneq - \mathbbm{1}_{\beta=1} \, \frac{\Phi_{H_\infty,c_\infty}''(z)}{2c_\infty \Phi_{H_\infty,c_\infty}'(z)}$ and
		\begin{align*}
			& \bm{c}(z_1,z_2) \coloneq \frac{1}{\beta c_\infty^2} \Big(\frac{1}{(z_1-z_2)^2} - \frac{\Phi_{H_\infty,c_\infty}'(z_1) \Phi_{H_\infty,c_\infty}'(z_2)}{(\Phi_{H_\infty,c_\infty}(z_1)-\Phi_{H_\infty,c_\infty}(z_2))^2}\Big) \ .
		\end{align*}
	\end{corollary}

	\section{Inference of generalized linear spectral statistics}\label{Section_GLSS}
	A generalized linear spectral statistic for functions $f,g : \R \rightarrow \C$ has the form
	\begin{align}\label{Eq_DefGLSS_theoretical}
		& L_n(f,g) \coloneq \frac{1}{d} \tr\big( f(\bm{S}_n) g(\Sigma_n) \big) \ .
	\end{align}
	By simple spectral decomposition of $\bm{S}_n$ and $\Sigma_n$,
	it may be written as
	\begin{align}\label{Eq_GLSS_interpretation0}
		L_n(f,g) 
		& = \frac{1}{d} \sum\limits_{j=1}^d f(\lambda_j(\bm{S}_n)) \sum\limits_{k=1}^d |u_j^*v_k|^2 \, g(\lambda_k(\Sigma_n)) \ ,
	\end{align}
	where $u_j$ is the (sample) eigenvector of $\bm{S}_n$ corresponding to $\lambda_j(\bm{S}_n)$ and $v_k$ is the (population) eigenvector of $\Sigma_n$ corresponding to $\lambda_k(\Sigma_n)$. Estimation of such generalized linear spectral statistics has clear utility, for example with regards to high-dimensional PCA algorithms. Two restrictions new to this generalized case are that the function $g$ must be deterministic, except for multiplication with indicator functions of certain random intervals (see Definition~\ref{Def_GLSS_estimator} and Theorem~\ref{Thm_GLSS_consistency}), and that the rate of consistency proven in Theorem~\ref{Thm_GLSS_consistency} is $\mathcal{O}(n^{\varepsilon-\frac{1}{2}})$ instead of $\mathcal{O}(n^{\varepsilon-1})$ as in Corollary~\ref{Cor_PLSS_Consistency}.
	
	\begin{definition}[Admissible pair of curves]\label{Def_AdmissiblePair}\
		\\
		For given $\tau,\kappa>0$ and $K \in \N$,
		call a pair of curves $(\gamma_n^{(g)},\gamma_n^{(f)})$ an \textit{admissible pair}, if the following properties hold.
		\begin{itemize}
			\item[a)]
			The piecewise $C^1$-curve $\gamma_n^{(g)}$ is admissible as in Definition~\ref{Def_AdmissibleCurve}. Additionally, it is assumed that the points $a_{\gamma_n^{(g)}}$ and $b_{\gamma_n^{(g)}}$ lie on the deterministic grid
			\begin{align}\label{Eq_GLSS_consistency_DefJn_new}
				& J_{n,K} \coloneq \Big\{\frac{k}{n^K} \ \Big| \ k \in \Z , \, \frac{k}{n^K} \in [-\kappa,\kappa] \Big\} \ .
			\end{align}
			
			\item[b)]
			The piecewise $C^1$-curve $\gamma_n^{(f)} : (0,1) \rightarrow B_{\kappa}^{\C^+}(0)$ satisfies
			\begin{align}\label{Eq_GLSS_gammafCondition1}
				& a_{\gamma_n^{(f)}} \coloneq \lim\limits_{t \nearrow 1} \gamma_n^{(f)}(t) \in \R \ \ \text{ and } \ \ b_{\gamma_n^{(f)}} \coloneq \lim\limits_{t \searrow 0} \gamma_n^{(f)}(t) \in \R
			\end{align}
			as well as
			\begin{align}\label{Eq_GLSS_gammafCondition2}
				& \dist\big(\mathrm{image}(\gamma_n^{(f)}),\supp(\hat{\ul{\nu}}_n)\cup\supp(\ul{\nu}_n)\cup\{0\}\big) \geq \tau
			\end{align}
			and
			\begin{align}\label{Eq_GLSS_gammafCondition2_2}
				& \forall \tz \in \mathrm{image}(\gamma_n^{(f)}) : \ \tau \leq |\varphi_{\hat{\nu}_n,c_n}(\tz)| \leq \kappa \ .
			\end{align}
			
			\item[c)] The curve images must avoid each other in the sense
			\begin{align}\label{Eq_GLSS_gammafCondition3}
				& \dist\big( \mathrm{image}(\gamma_n^{(g)}), \varphi_{\hat{\nu}_n,c_n}(\mathrm{image}(\gamma_n^{(f)})) \big) \geq \tau \ ,
			\end{align}
			where $\varphi_{\hat{\nu}_n,c_n}(\tz) = \frac{-1}{\cs_{\hat{\ul{\nu}}_n}(\tz)}$, similarly to (\ref{Eq_Def_varphi}).
		\end{itemize}
	\end{definition}
	
	\begin{definition}[GLSS estimator]\label{Def_GLSS_estimator}\
		\\
		To an admissible pair $(\gamma_n^{(g)},\gamma_n^{(f)})$ and functions $f \in \operatorname{Hol}(\gamma_n^{(f)})$, $g \in \operatorname{Hol}(\gamma_n^{(g)})$, define the generalized linear spectral statistic estimator as the sum of double curve integrals
		\begin{align}\label{Eq_DefGLSS_estimator}
			& \hat{L}_{n,\gamma_n^{(f)},\gamma_n^{(g)}}(f,g) \coloneq \frac{1}{4\pi^2} \int_{\gamma_n^{(f)}} \int_{\gamma_n^{(g)}} g(z_1) f(\tz_2) \hat{k}_n(z_1,\tz_2) \, dz_1 \, d\tz_2 \nonumber\\
			& \hspace{3.5cm} - \frac{1}{4\pi^2} \int_{\gamma_n^{(f)}} \int_{\gamma_n^{(g)}} g(\ol{z_1}) f(\tz_2) \hat{k}_n(\ol{z_1},\tz_2) \, dz_1 \, d\tz_2 \nonumber\\
			& \hspace{3.5cm} - \frac{1}{4\pi^2} \int_{\gamma_n^{(f)}} \int_{\gamma_n^{(g)}} g(z_1) f(\ol{\tz_2}) \hat{k}_n(z_1,\ol{\tz_2}) \, dz_1 \, d\tz_2 \nonumber\\
			& \hspace{3.5cm} + \frac{1}{4\pi^2} \int_{\gamma_n^{(f)}} \int_{\gamma_n^{(g)}} g(\ol{z_1}) f(\ol{\tz_2}) \hat{k}_n(\ol{z_1},\ol{\tz_2}) \, dz_1 \, d\tz_2 \ ,
		\end{align}
		where
		\begin{align}\label{Eq_DefGLSS_kernel}
			& \hat{k}_n(z_1,\tz_2) \coloneq \frac{\tz_2 \cs_{\hat{\ul{\nu}}_n}(\tz_2)^2 + (1-c_n)\cs_{\hat{\ul{\nu}}_n}(\tz_2) + c_n\hat{s}_n(z_1)}{c_n(z_1\cs_{\hat{\ul{\nu}}_n}(\tz_2)+1) \tz_2} \ .
		\end{align}
		Evaluation of $\hat{k}_n(\ol{z_1},\tz_2)$, $\hat{k}_n(z_1,\ol{\tz_2})$ or $\hat{k}_n(\ol{z_1},\ol{\tz_2})$ may be done with the canonical extensions $\cs_{\hat{\ul{\nu}}_n}(\ol{\tz_2}) = \ol{\cs_{\hat{\ul{\nu}}_n}(\tz_2)}$ and $\hat{s}_n(\ol{z_1}) = \ol{\hat{s}_n(z_1)}$.
	\end{definition}

	\begin{theorem}[Consistency of the GLSS estimator with error rate $\mathcal{O}(n^{\varepsilon-\frac{1}{2}})$]\label{Thm_GLSS_consistency}\
		\\
		Suppose Assumption~\ref{Assumption_EigInf_Main} holds and fix constants $\tau,\kappa>0$. Let $g : U_0 \rightarrow \C$ be a fixed holomorphic function on some symmetric (with respect to conjugation), open and simply connected set $U_0 \subset \C$ containing $[0,\sigma^2]$ such that $g(\R) \subset \R$.
		Then for every $\varepsilon \in (0,1)$ and $D>0$, there exists a constant $C=C(\tau,\kappa,\varepsilon,D)>0$ such that
		\begin{align}
			& \hspace{0cm} \bP_*\bigg( \forall (\gamma_n^{(f)}, \gamma_n^{(g)}) \text{ as in Def.~\ref{Def_AdmissiblePair}, where } \mathrm{image}(\gamma_n^{(g)}) \subset U_0, \, \forall f \in \operatorname{Hol}(\gamma_n^{(f)}) : \nonumber\\
			& \hspace{1cm} \Big| \hat{L}_{n,\gamma_n^{(f)},\gamma_n^{(g)}}(f,g) - L_n\big( f \mathbbm{1}_{[a_{\gamma_n^{(f)}},b_{\gamma_n^{(f)}}]} \, , \, g \mathbbm{1}_{[a_{\gamma_n^{(g)}},b_{\gamma_n^{(g)}}]} \big) \Big| \leq K_{\gamma_n^{(f)},\gamma_n^{(g)}}(f,g) \, n^{\varepsilon-\frac{1}{2}} \bigg) \nonumber\\
			& \hspace{10cm} \geq 1 - \frac{C}{n^D}
		\end{align}
		holds for all $n \in \N$, where
		\begin{align*}
			& K_{\gamma_n^{(f)},\gamma_n^{(g)}}(f,g) \coloneq \|f\|_{\gamma_n^{(f)}} \ell(\gamma_n^{(f)}) \, \|g\|_{\gamma_n^{(g)}} \big(\ell(\gamma_n^{(g)})+1\big)
		\end{align*}
		and $\ell(\gamma)$ denotes the length of the curve $\gamma$.
	\end{theorem}

	\section{Proofs of Lemmas~\ref{Lemma_PopStilEst_Uniqueness} and~\ref{Lemma_ConsistencyBasic}}\label{Section_ProofsOfLemmas}
	
	Two auxiliary lemmas are introduced prior to proving Lemmas~\ref{Lemma_PopStilEst_Uniqueness} and~\ref{Lemma_ConsistencyBasic}. The first analyzes the behavior of the deterministic equivalents $\nu_n$, while the second contains simple bounds, variations of which will be required throughout the rest of the paper. Both lemmas are proved in Section~\ref{Section_LemmaProofs}.
	
	\begin{lemma}[Convergence of deterministic Stieltjes transforms]\label{Lemma_BasicStieltjesConvergence}\
		\\
		Conditions~\ref{EI_ItemAssumption_Asymptotics},~\ref{EI_ItemAssumption_PopConv} and~\ref{EI_ItemAssumption_sigmaBound} of Assumption~\ref{Assumption_EigInf_Main} jointly imply the following properties of the deterministic measures $H_n$ and $\nu_n$.
		\begin{itemize}
			\item[a)] The convergence $\cs_{H_n} \xrightarrow{n \to \infty} \cs_{H_\infty}$ holds uniformly on compact subsets of $\C\setminus[0,\sigma^2]$.
			
			\item[b)] The supports of the measures $\nu_\infty$ and $\nu_n$ are contained in the intervals $[0,\sigma^2(1+\sqrt{c_\infty})^2]$ and $[0,\sigma^2(1+\sqrt{c_n})^2]$ respectively.
			
			\item[c)] The convergence $\cs_{\nu_n} \xrightarrow{n \to \infty} \cs_{\nu_\infty}$ holds uniformly on compact subsets of the domain $\C\setminus[0,\sigma^2(1+\sqrt{c_\infty})^2]$.
			
			\item[d)] The convergence $\nu_n \xRightarrow{n \to \infty} \nu_\infty$ holds.
			
			\item[e)] The convergence $\varphi_{\nu_n,c_n} \xrightarrow{n \to \infty} \varphi_{\nu_\infty,c_\infty}$ holds uniformly on compact subsets of the domain $\C\setminus[0,\sigma^2(1+\sqrt{c_\infty})^2]$.
		\end{itemize}
	\end{lemma}
	
	\begin{lemma}[Standard bounds]\label{Lemma_StandardBounds}\
		The following statements hold and (a)--(c) are also true if every occurrence of $\nu_n$ and $\ul{\nu}_n$ is replaced with $\hat{\nu}_n$ and $\hat{\ul{\nu}}_n$ respectively. Furthermore, statement (b) is always applicable with $[0,\infty)$ instead of $[0,K]$.
		\begin{itemize}
			\item[a)]
			The support $\supp(\ul{\nu}_n)$ is contained in $\supp(\nu_n) \cup \{0\}$.
			
			\item[b)]
			Suppose a compact interval $[0,K]$ for some $K>0$ satisfies $\supp(\nu_n) \subset [0,K]$ for all $n \in \N$, then the implication
			\begin{align}\label{Eq_StandardBounds_b1}
				& \big( \nu_n \xRightarrow{n \to \infty} \nu_\infty \big) \Rightarrow \big( \supp(\nu_\infty), \supp(\ul{\nu}_\infty) \subset [0,K] \big)
			\end{align}
			holds. Additionally, for all $\tz \in \C$ with
			\begin{align}\label{Eq_StandardBounds_b2}
				& \dist(\tz,[0,K]) \geq \tau
			\end{align}
			for some $\tau > 0$, the bound
			\begin{align}\label{Eq_StandardBounds_b3}
				& |\cs_{\nu}(\tz)| \geq \int_\R \frac{\tau/2}{|\lambda-\tz|^2} \, d\nu(\lambda)
			\end{align}
			holds for both $\nu \in \{\nu_n,\ul{\nu}_n\}$ and also for $\nu \in \{\nu_\infty,\ul{\nu}_\infty\}$, if the right-hand side of (\ref{Eq_StandardBounds_b1}) is true.
			The map $\varphi_{\nu_n,c_n}(\tz) = \frac{-1}{\cs_{\ul{\nu}_n}(\tz)}$ is holomorphic on $\C\setminus[0,K]$, and if the right-hand side of (\ref{Eq_StandardBounds_b1}) holds, then so is $\varphi_{\nu_\infty,c_\infty}(\tz) = \frac{-1}{\cs_{\ul{\nu}_\infty}(\tz)}$.
			
			\item[c)] For all $\tz \in \C$ with
			\begin{align}\label{Eq_StandardBounds_c1}
				& \dist(\tz,\supp(\nu_n)\cup\{0\}) \geq \tau
			\end{align}
			the bounds
			\begin{align}\label{Eq_StandardBounds_c2}
				& |\cs_{\nu_n}(\tz)| \leq \frac{1}{\tau} \ \ \text{ and } \ \ |\cs_{\ul{\nu}_n}(\tz)| \leq \frac{1}{\tau}
			\end{align}
			hold. Additionally, for any two $\tz_1,\tz_2 \in \C$ both satisfying (\ref{Eq_StandardBounds_c1}) it holds that
			\begin{align}\label{Eq_StandardBounds_c3}
				& |\cs_{\nu_n}(\tz_1) - \cs_{\nu_n}(\tz_2)| \leq \frac{|\tz_1-\tz_2|}{\tau^2} \ \ \text{ and } \ \ |\cs_{\ul{\nu}_n}(\tz_1) - \cs_{\ul{\nu}_n}(\tz_2)| \leq \frac{|\tz_1-\tz_2|}{\tau^2} \ .
			\end{align}
			
			\item[d)]
			The deterministic maps $\varphi_{\nu_n,c_n}$ and $\varphi_{\nu_\infty,c_\infty}$ satisfy
			\begin{align}\label{Eq_StandardBounds_d1}
				& \frac{|\varphi_{\nu_n,c_n}(\tz_1)-\varphi_{\nu_n,c_n}(\tz_2)|}{|\tz_1-\tz_2|} \geq \frac{1}{2} \ \ \text{ and } \ \ \frac{|\varphi_{\nu_\infty,c_\infty}(\tz_1)-\varphi_{\nu_\infty,c_\infty}(\tz_2)|}{|\tz_1-\tz_2|} \geq \frac{1}{2}
			\end{align}
			for all $\tz_1 \neq \tz_2$ from $\C\setminus\R$. In the setting of (b), continuity of the involved maps extends (\ref{Eq_StandardBounds_d1}) to hold for all $\tz_1 \neq \tz_2$ from $\C\setminus[0,K]$.
		\end{itemize}
	\end{lemma}

	\subsection{Proof of Lemma~\ref{Lemma_PopStilEst_Uniqueness}}\label{Proof_Lemma_PopStilEst_Uniqueness}
	\begin{itemize}
		\item[i)] \textit{Proof of property (a) for $\hat{\nu}_n \neq \delta_0$}:\\
		As a preparatory step, observe that any solution $s \in \C$ with $\Im((1 - czs - c)z)>0$ to an equation of the form
		\begin{align}\label{Eq_InversionFormulaClean}
			& zs + 1 = \int_\R \frac{\lambda}{\lambda - (1 - czs - c)z} \, d\nu(\lambda)
		\end{align}
		must satisfy
		\begin{align}\label{Eq_Imag_Calc0}
			& \Im(zs) = \Im(zs + 1) = \int_\R \Im\Big(\frac{\lambda}{\lambda - (1 - czs - c)z}\Big) \, d\nu(\lambda) \nonumber\\
			& = -\int_\R \frac{\lambda \Im(\lambda - (1 - czs - c)z)}{|\lambda - (1 - czs - c)z|^2} \, d\nu(\lambda) \nonumber\\
			& = \underbrace{\Im((1 - czs - c)z)}_{>0} \int_\R \frac{\lambda}{|\lambda - (1 - czs - c)z|^2} \, d\nu(\lambda) > 0 \ ,
		\end{align}
		which may be rearranged as
		\begin{align}\label{Eq_Imag_Calc}
			& \Big| \frac{c z \Im(zs)}{\Im((1-czs-c)z)} \Big| = \frac{c |z| \Im(zs)}{\Im((1-czs-c)z)} = c|z| \int_\R \frac{\lambda}{|\lambda - (1 - czs - c)z|^2} \, d\nu(\lambda) \ .
		\end{align}
		Suppose there are two such solutions $s,\ul{s} \in \C$ both satisfying
		\begin{align*}
			& \Big| \frac{c z \Im(zs)}{\Im((1-czs-c)z)} \Big| < 1 \ \ \text{ and } \ \ \Big| \frac{c z \Im(z\ul{s})}{\Im((1-cz\ul{s}-c)z)} \Big| < 1 \ ,
		\end{align*}
		then
		\begin{align}\label{Eq_UniquenessCalc_new}
			& s-\ul{s} = \frac{1}{z} \int_\R \frac{\lambda}{\lambda - (1 - czs - c)z} - \frac{\lambda}{\lambda - (1 - cz\ul{s} - c)z} \, d\nu(\lambda) \nonumber\\
			& = \frac{1}{z} \int_\R \lambda \frac{(1 - czs - c)z - (1 - cz\ul{s} - c)z}{(\lambda - (1 - czs - c)z)(\lambda - (1 - cz\ul{s} - c)z)} \, d\nu(\lambda) \nonumber\\
			& = \int_\R \lambda \frac{cz(\ul{s}-s)}{(\lambda - (1 - czs - c)z)(\lambda - (1 - cz\ul{s} - c)z)} \, d\nu(\lambda) \nonumber\\
			& = (s-\ul{s}) \int_\R \frac{- cz \lambda}{(\lambda - (1 - czs - c)z)(\lambda - (1 - cz\ul{s} - c)z)} \, d\nu(\lambda) \ .
		\end{align}
		Using Cauchy--Schwarz and (\ref{Eq_Imag_Calc}), one may bound the right-hand factor by
		\begin{align*}
			& \bigg| \int_\R \frac{- c z \lambda}{(\lambda - (1 - czs - c)z)(\lambda - (1 - cz\ul{s} - c)z)} \, d\nu(\lambda) \bigg|\\
			& \leq \bigg( c |z| \int_\R \frac{\lambda}{|\lambda - (1 - czs - c)z|^2} \, d\nu(\lambda) \bigg)^{\frac{1}{2}} \bigg( c |z| \int_\R \frac{\lambda}{|\lambda - (1 - cz\ul{s} - c)z|^2} \, d\nu(\lambda) \bigg)^{\frac{1}{2}}\\
			& \hspace{-0.12cm} \overset{\text{(\ref{Eq_Imag_Calc})}}{=} \bigg( \Big|\frac{cz \Im(zs)}{\Im((1-czs-c)z)} \Big| \bigg)^{\frac{1}{2}} \bigg( \Big| \frac{cz \Im(z\ul{s})}{\Im((1-cz\ul{s}-c)z)} \Big| \bigg)^{\frac{1}{2}} \ ,
		\end{align*}
		which is strictly less than one by assumption, thus proving $s=\ul{s}$.
		
		\item[(ii)] \textit{Proof of property (a) for $\hat{\nu}_n=\delta_0$}:\\
		For $\nu=\delta_0$, equation (\ref{Eq_InversionFormulaClean}) under the assumption $\Im((1-czs-c)z)>0$ reduces to
		\begin{align*}
			& zs + 1= 0 \ ,
		\end{align*}
		so the only possible candidate for a solution to (\ref{Eq_InversionFormulaClean}) is $s=\frac{-1}{z}$. This implies uniqueness of the solution, when it exists. In this case, properties (b) and (c) also follow immediately. Hence, one may for the rest of this proof assume $\hat{\nu}_n\neq\delta_0$.
		
		\item[iii)] \textit{Proof of property (b)}:\\
		As a preparatory step, observe that for any $z, s \in \C$ with $(1 - czs - c)z \notin \supp(\nu)$ and $0 \neq \cs_{\ul{\nu}}((1 - czs - c)z)$, the calculation
		\begin{align*}
			& zs + 1 = \int_\R \frac{\lambda}{\lambda - (1 - czs - c)z} \, d\nu(\lambda) = \int_\R \frac{(1 - czs - c)z}{\lambda - (1 - czs - c)z} \, d\nu(\lambda) + 1 \nonumber\\
			\Leftrightarrow \ & \frac{zs}{(1 - czs - c)z} = \int_\R \frac{1}{\lambda - (1 - czs - c)z} \, d\nu(\lambda) \overset{\text{(\ref{Eq_DefStieltjes})}}{=} \cs_\nu\big( (1 - czs - c)z \big) \nonumber\\
			\Leftrightarrow \ & \frac{1-czs-c}{(1 - czs - c)z} = -c\cs_\nu\big( (1 - czs - c)z \big) + \frac{1-c}{(1 - czs - c)z} \overset{\text{(\ref{Eq_Def_ulNu})}}{=} -\cs_{\ul{\nu}}\big( (1 - czs - c)z \big)\\
			\Leftrightarrow \ & z = \frac{-1}{\cs_{\ul{\nu}}((1 - czs - c)z)} \overset{\text{(\ref{Eq_Def_varphi})}}{=} \varphi_{\nu,c}\big( (1 - czs - c)z \big)
		\end{align*}
		shows that
		\begin{align}\label{Eq_RevMPE_Equivalence1}
			& \forall z, s \in \C \text{ with } \ (1 - czs - c)z \notin \supp(\nu) \ \text{ and } \ 0 \neq \cs_{\ul{\nu}}((1 - czs - c)z) : \nonumber\\
			& zs + 1 = \int_\R \frac{\lambda}{\lambda - (1 - czs - c)z} \, d\nu(\lambda) \ \ \Leftrightarrow \ \ z = \varphi_{\nu,c}\big( (1 - czs - c)z \big) \ .
		\end{align}
		Define
		\begin{align*}
			& U \coloneq \big\{ z \in \C^+ \ \big| \ \hat{s}_n(z) \text{ exists} \big\} \ .
		\end{align*}
		By (\ref{Eq_RevMPE_Equivalence1}) and (a), $\hat{s}_n(z)$ is also the unique solution to
		\begin{align}\label{Eq_RevMPE_Equivalence1_copy}
			& z = \varphi_{\hat{\nu}_n,c_n}\big( \underbrace{(1 - c_nz\hat{s}_n(z) - c_n)z}_{\eqcolon  \tz} \big) \ ,
		\end{align}
		which also satisfies (\ref{Eq_StilEstimator_DefiningConditions}). The first property in (\ref{Eq_StilEstimator_DefiningConditions}) is equivalent to the assertion $\tz \in \C^+$. The second property is by (\ref{Eq_Imag_Calc0}) equivalent to
		\begin{align*}
			& c_n |\overbrace{\varphi_{\hat{\nu}_n,c_n}(\tz)}^{\overset{\text{(\ref{Eq_RevMPE_Equivalence1_copy})}}{=}z}| \int_\R \frac{\lambda}{|\lambda - \tz|^2} \, d\hat{\nu}_n(\lambda) < 1 .
		\end{align*}
		The set
		\begin{align*}
			& \tilde{U} \coloneq \big\{ (1 - c_nz\hat{s}_n(z) - c_n)z \ \big| \ z \in U \big\} \overset{\text{(\ref{Eq_RevMPE_Equivalence1_copy})}}{=} \varphi_{\hat{\nu}_n,c_n}^{-1}(U)
		\end{align*}
		may thus be characterized as
		\begin{align}\label{Eq_UniquenessLemma_Step2_1}
			& \tilde{U} = \Big\{ \tz \in \C^+ \ \Big| \ c_n |\varphi_{\hat{\nu}_n,c_n}(\tz)| \int_\R \frac{\lambda}{|\lambda - \tz|^2} \, d\hat{\nu}_n(\lambda) < 1 \Big\} \ ,
		\end{align}
		which is an open subset of $\C^+$.
		As $\varphi_{\hat{\nu}_n,c_n}(\tilde{U})$ is holomorphic on all of $\C^+$, the open mapping theorem guarantees that $U = \varphi_{\hat{\nu}_n,c_n}(\tilde{U})$ is also open.
		\\[0.5em]
		It remains to show the existence of a $\kappa(\hat{\nu}_n,c_n)>0$ such that $\C^+ \setminus B^{\C^+}_{\kappa(\hat{\nu}_n,c_n)}(0)$ lies in $U$. As this step is somewhat technical without introducing additional conceptual ideas, it has been moved to the appendix (see Subsection \ref{Proof_Lemma_PopStilEst_Uniqueness_Continuation}).

		\item[iv)] \textit{Proof of property (c)}:\\
		As already seen in (\ref{Eq_RevMPE_Equivalence1_copy}), the map
		\begin{align*}
			& \hat{\Phi}_n : U \rightarrow \tilde{U} \ \ ; \ \ z \mapsto (1 - c_nz\hat{s}_n(z) - c_n)z
		\end{align*}
		is the local inverse of $\varphi_{\hat{\nu}_n,c_n}$. Since $\varphi_{\hat{\nu}_n,c_n}$ is by construction holomorphic on all of $\C \setminus \supp(\hat{\ul{\nu}}_n)$, the holomorphic inverse function theorem applied locally at every point in $U$, to guarantee that $\hat{\Phi}_n$ and thus also $\hat{s}_n$ are holomorphic on $U$. \qed
	\end{itemize}

	Next, the proof of Lemma~\ref{Lemma_ConsistencyBasic}. A slightly stronger version of the lemma will be proved, which allows the set $J$ to approach the real line, provided it stays sufficiently far away from $\supp(H_\infty)$. This is realized by letting $J$ be a compact subset of $\bD_{H_\infty,c_\infty}(1) \cup (\R\setminus I)$, for a certain interval $I$ defined in (\ref{Eq_StrongerVersion_Def_I}).
	\\[0.5em]
	The core part of the proof will be an application of Rouch\'{e}'s theorem to show that $\hat{\Phi}_n(z) = (1-c_nz\hat{s}_n(z)-c_n)z$ exists and lies in a small neighborhood of $\Phi_{H_n,c_n}(z)$. Statement (d) of Lemma~\ref{Lemma_StandardBounds} plays a key role for the applicability of Rouch\'{e}'s theorem.
	
	\subsection{Proof of Lemma~\ref{Lemma_ConsistencyBasic}}\label{Proof_Lemma_ConsistencyBasic}
	Under the additional assumption that there exists a constant $K>0$ such that
	\begin{align}\label{Eq_StrongerVersion_AdditionalConditionK}
		& \forall n \in \N : \supp(\hat{\nu}_n) \cup \supp(\nu_n) \subset [0,K] \ ,
	\end{align}
	it will be shown here that
	\begin{align}\label{Eq_ConsistencyBasic_Stronger}
		& 1 = \bP\Big( \forall J \subset \bD_{H_\infty,c_\infty}(1) \cup (\R\setminus I) \text{ compact } \, \exists N_J>0 \, \forall n \geq N_J : \nonumber\\
		& \hspace{2cm} \hat{s}_{n} \text{ exists on $J\setminus\R$ and } \ \ \sup\limits_{z \in J\setminus\R} \big|\hat{s}_n(z) - \cs_{H_\infty}(z)\big| \xrightarrow{n \to \infty} 0 \Big) \ ,
	\end{align}
	where
	\begin{align}\label{Eq_StrongerVersion_Def_I}
		& I \coloneq \big[\varphi_{\nu_\infty,c_\infty}^{-1}(-\delta), \varphi_{\nu_\infty,c_\infty}^{-1}(K+\delta) \big] \subset \R
	\end{align}
	for some $\delta>0$. It follows from (b) of Lemma~\ref{Lemma_StandardBounds} and the definition $\varphi_{\nu_\infty,c_\infty}(\tz) = \frac{-1}{\cs_{\ul{\nu}_\infty}(\tz)}$ that the interval $I$ is well-defined.
	To prove the original version of Lemma~\ref{Lemma_ConsistencyBasic}, whose assumptions do not yield (\ref{Eq_StrongerVersion_AdditionalConditionK}), one may replace both $[0,K]$ and $I$ in the following proof with $\R$.
	\\[0.5em]
	The proof of measurability of the event in (\ref{Eq_ConsistencyBasic_Stronger}) is deferred to Section~\ref{Section_Measurability}.
	\begin{itemize}
		\item[i)] \textit{Preliminaries}:\\
		The entirety of the following proof will take place on the event
		\begin{align}\label{Eq_BasicConsistency_Event}
			& \mathcal{E}_{\text{MP}} \coloneq \{\hat{\nu}_n \xRightarrow{n \to \infty} \nu_\infty\} \ .
		\end{align}
		Assumption (\ref{Eq_MP_Basic}) then guarantees that all following statements hold almost surely.
		\\[0.5em]
		Let $J \subset \bD_{H_\infty,c_\infty}(1) \cup (\R\setminus I)$ denote an arbitrary compact subset. The map $\varphi_{\nu_\infty,c_\infty}$ is by (b) of Lemma~\ref{Lemma_StandardBounds} holomorphic on $\C\setminus[0,K]$ and it follows from Lemma~\ref{Lemma_SpaceTransform} that it has the holomorphic inverse
		\begin{align*}
			& \Phi_{H_\infty,c_\infty} : \varphi_{\nu_\infty,c_\infty}\big( \C\setminus[0,K] \big) \rightarrow \C\setminus[0,K] \ \ ; \ \ z \mapsto (1-c_\infty z \cs_{H_\infty} - c_\infty)z \ .
		\end{align*}
		By construction of $I$ and Lemma~\ref{Lemma_SpaceTransform} it also holds that
		\begin{align*}
			& \bD_{H_\infty,c_\infty}(1) \cup (\R\setminus I) \subset \varphi_{\nu_\infty,c_\infty}\big( \C\setminus[0,K] \big) \ ,
		\end{align*}
		making $\Phi_{H_\infty,c_\infty}$ continuous on $J \subset \bD_{H_\infty,c_\infty}(1) \cup (\R\setminus I)$. Hence, the image $\Phi_{H_\infty,c_\infty}(J) \subset \C\setminus[0,K]$ must also be compact.
		Let $V_\varepsilon$ denote the closed $\varepsilon$-neighborhood of $\Phi_{H_\infty,c_\infty}(J)$, i.e.
		\begin{align}\label{Eq_BasicConsistency_UDef}
			& V_\varepsilon \coloneq \big\{ \tz \in \C \ \big| \ \dist(\tz, \Phi_{H_\infty,c_\infty}(J)) \leq \varepsilon \big\} \ .
		\end{align}
		Choose $\varepsilon>0$ small enough to satisfy
		\begin{align}\label{Eq_BasicConsistency_UProp1}
			& \dist(V_\varepsilon,[0,K]) > \tau
		\end{align}
		for some $\tau>0$, which is possible by compactness of $\Phi_{H_\infty,c_\infty}(J) \subset \C\setminus[0,K]$. By (\ref{Eq_StrongerVersion_AdditionalConditionK}) and (a) of Lemma~\ref{Lemma_StandardBounds}, this implies
		\begin{align}\label{Eq_BasicConsistency_UProp2}
			& \dist\big(V_\varepsilon,\supp(\hat{\ul{\nu}}_n) \cup \supp(\ul{\nu}_n)\big) > \tau \ .
		\end{align}
		
		\item[ii)] \textit{Showing that $\cs_{\hat{\nu}_n} \xrightarrow{n \to \infty} \cs_{\nu_\infty}$ uniformly on $V_\varepsilon$}:\\
		The pointwise convergence $\cs_{\hat{\nu}_n}(\tz) \xrightarrow{n \to \infty} \cs_{\nu_\infty}(\tz)$ for all $\tz \in V_\varepsilon$ follows directly from (\ref{Eq_StrongerVersion_AdditionalConditionK}), (\ref{Eq_BasicConsistency_Event}) and the fact that the function $f_{\tz} : [0,K] \rightarrow \C \ ; \ \lambda \mapsto \frac{1}{\lambda-\tz}$ is for every $\tz \in V_\varepsilon$ continuous and bounded by $\frac{1}{\tau}$ due to (\ref{Eq_BasicConsistency_UProp1}). Moreover, one by (\ref{Eq_BasicConsistency_UProp2}) and (c) of Lemma~\ref{Lemma_StandardBounds} also has
		\begin{align}
			& \hspace{1cm} \forall \tz \in V_\varepsilon : \ |\cs_{\hat{\nu}_n}(\tz)| \overset{\text{(\ref{Eq_StandardBounds_c2})}}{\leq} \frac{1}{\tau} \label{Eq_BasicConsistency_Step2_02}\\
			& \hspace{-0.5cm} \forall \tz_1,\tz_2 \in V_\varepsilon : \ |\cs_{\hat{\nu}_n}(\tz_1) - \cs_{\hat{\nu}_n}(\tz_2)| \overset{\text{(\ref{Eq_StandardBounds_c3})}}{\leq} \frac{|\tz_1-\tz_2|}{\tau^2} \label{Eq_BasicConsistency_Step2_03} \ ,
		\end{align}
		thus showing that the family of functions $(\cs_{\hat{\nu}_n})_{n \in \N}$ is uniformly bounded and equicontinuous on $V_\varepsilon \subset \C$, which is compact by construction. Arzel\`a--Ascoli's theorem gives the existence of a sub-sequence $(\cs_{\hat{\nu}_{n_k}})_{k \in \N}$ uniformly convergent on $V_\varepsilon$. The fact that the limit can only be the pointwise limit $\cs_{\nu_\infty}$, by standard topological arguments implies that the original sequence must have already converged uniformly to $\cs_{\nu_\infty}$ on $V_\varepsilon$.
		
		\item[iii)] \textit{Showing that $\varphi_{\hat{\nu}_n,c_n} \xrightarrow{n \to \infty} \varphi_{\nu_\infty,c_\infty}$ uniformly on $V_\varepsilon$}:\\
		By construction, $V_\varepsilon$ is bounded, so there exists a $C>0$ such that
		\begin{align}\label{Eq_BasicConsistency_Step3_01}
			& \forall \tz \in V_\varepsilon : \ |\tz| < C \ .
		\end{align}
		Using (\ref{Eq_BasicConsistency_UProp1}) and (b) of Lemma~\ref{Lemma_StandardBounds}, one may calculate
		\begin{align*}
			& |\cs_{\hat{\ul{\nu}}_n}(\tz)| \overset{\text{(\ref{Eq_StandardBounds_b3})}}{\geq} \int_\R \frac{\tau/2}{|\lambda-\tz|^2} \, d\hat{\ul{\nu}}_n(\lambda) \overset{\text{(\ref{Eq_BasicConsistency_Step3_01})}}{\geq} \int_\R \frac{\tau/2}{(\lambda+C)^2} \, d\hat{\ul{\nu}}_n(\lambda) \ .
		\end{align*}
		As the right-hand side does not depend on $\tz$ and by (\ref{Eq_BasicConsistency_Event}) converges to $\int_\R \frac{\tau/2}{(\lambda+C)^2} \, d\ul{\nu}_\infty(\lambda)>0$, the convergence shown in (ii) directly implies
		\begin{align}\label{Eq_BasicConsistency_Step3_02}
			& \sup\limits_{\tz \in V_\varepsilon} \big| \varphi_{\hat{\nu}_n,c_n}(\tz) - \varphi_{\nu_\infty,c_\infty}(\tz) \big| \overset{\text{(\ref{Eq_Def_varphi})}}{=} \sup\limits_{\tz \in V_\varepsilon} \Big| \frac{-1}{\cs_{\hat{\ul{\nu}}_n}(\tz)} - \frac{-1}{\cs_{\ul{\nu}_\infty}(\tz)} \Big| \xrightarrow{n \to \infty} 0 \ .
		\end{align}
		
		\item[iv)] \textit{Uniform application of Rouch\'{e}'s theorem}:\\
		It is shown here that for each $\varepsilon > 0$ small enough to satisfy (\ref{Eq_BasicConsistency_UProp1}) and all $n$ greater than some $N_\varepsilon'>0$ there exist maps $\hat{s}_{n,\varepsilon} : J \rightarrow \C$, which satisfy
		\begin{align}\label{Eq_BasicConsistency_Step4_Prop1}
			& \big| (1-c_nz\hat{s}_{n,\varepsilon}(z)-c_n)z - \Phi_{H_\infty,c_\infty}(z) \big| \leq \varepsilon
		\end{align}
		and
		\begin{align}\label{Eq_BasicConsistency_Step4_Prop2}
			& z = \varphi_{\hat{\nu}_n,c_n}\big( (1-c_nz\hat{s}_{n,\varepsilon}(z)-c_n)z \big)
		\end{align}
		for all $z \in J$.
		\\[0.5em]
		In (i), the set $V_\varepsilon$ was defined as the closed $\varepsilon$-neighborhood of the compact set $\Phi_{H_\infty,c_\infty}(J)$, i.e.
		\begin{align*}
			& V_\varepsilon \coloneq \big\{ \tz \in \C \ \big| \ \dist(\tz, \Phi_{H_\infty,c_\infty}(J)) \leq \varepsilon \big\} \subset \C\setminus[0,K] \ .
		\end{align*}
		By (d) of Lemma~\ref{Lemma_StandardBounds}, it holds that
		\begin{align}\label{Eq_BasicConsistency_Step4_4}
			& \forall z_* \in J \, \forall \tz \in \underbrace{\partial B_\varepsilon^{\C}(\Phi_{H_\infty,c_\infty}(z_*))}_{\subset V_\varepsilon} : \ \big| \varphi_{\nu_\infty,c_\infty}(\tz) - \underbrace{\varphi_{\nu_\infty,c_\infty}(\Phi_{H_\infty,c_\infty}(z_*))}_{= z_* \text{ by Lemma~\ref{Lemma_SpaceTransform}}} \big| \geq \frac{\varepsilon}{2} \ .
		\end{align}
		By (iii), there exists an $N_{\varepsilon}' > 0$ such that
		\begin{align}\label{Eq_BasicConsistency_Step4_5}
			& \sup\limits_{\tz \in V_\varepsilon} \big| \overbrace{\varphi_{\nu_\infty,c_\infty}(\tz) - \varphi_{\hat{\nu}_n,c_n}(\tz)}^{\eqcolon  G_n(\tz)} \big| < \frac{\varepsilon}{4}
		\end{align}
		holds for all $n \geq N_{\varepsilon}'$. Consequently, one may for all $z_* \in J$ and $\tz \in \partial B_\varepsilon^{\C}(\Phi_{H_\infty,c_\infty}(z_*))$ bound
		\begin{align}\label{Eq_BasicConsistency_Step4_6}
			& \big| \overbrace{\varphi_{\hat{\nu}_n,c_n}(\tz) - z_*}^{\eqcolon  F_{n,z_*}} \big| \geq \big| \varphi_{\nu_\infty,c_\infty}(\tz) - z_* \big| - \big| \varphi_{\hat{\nu}_n,c_n}(\tz) - \varphi_{\nu_\infty,c_\infty}(\tz) \big| > \frac{\varepsilon}{4} \ ,
		\end{align}
		which with (\ref{Eq_BasicConsistency_Step4_5}) proves that the map $G_n$ is dominated by $F_{n,z_*}$ on the boundary. As both are holomorphic on $\C\setminus[0,K] \supset V_\varepsilon$ by (b) of Lemma~\ref{Lemma_StandardBounds}, Rouch\'{e}'s theorem then yields that $F_{n,z_*}$ and $F_{n,z_*}+G_n$ must have the same number of zeros/roots in $B_\varepsilon^{\C}(\Phi_{H_\infty,c_\infty}(z_*))$ when counting multiplicity. Lemma~\ref{Lemma_SpaceTransform} yields
		\begin{align*}
			& (F_{n,z_*}+G_n)(\Phi_{H_\infty,c_\infty}(z_*)) = \varphi_{\nu_\infty,c_\infty}(\Phi_{H_\infty,c_\infty}(z_*)) - z_* = 0
		\end{align*}
		(technically only for $z_* \in \bD_{H_\infty,c_\infty}(\infty)$, but the construction of $I$ and the expansion of $\Phi_{H_\infty,c_\infty}$ done in (i) gives the same identity for all $z_* \in J \subset \bD_{H_\infty,c_\infty}(\infty) \cup (\R\setminus I)$), so by Rouch\'{e}'s theorem there must exist at least one $\tz_0 \in B_\varepsilon^{\C}(\Phi_{H_\infty,c_\infty}(z_*))$ such that
		\begin{align}\label{Eq_BasicConsistency_Step4_07}
			& 0 = F_{n,z_*}(\tz_0) = \varphi_{\hat{\nu}_n,c_n}(\tz_0) - z_* \ .
		\end{align}
		The value $\hat{s}_{n,\varepsilon}(z_*) \coloneq \frac{1-\frac{\tz_0}{z_*}-c_n}{c_nz_*}$ by construction satisfies
		\begin{align}\label{Eq_BasicConsistency_Step4_7}
			& \big| \underbrace{(1-c_n z_*\hat{s}_{n,\varepsilon}(z_*)-c_n)z_*}_{= \tz_0} - \underbrace{(1-c_\infty z_*\cs_{H_\infty}(z_*)-c_\infty)z_*}_{= \Phi_{H_\infty,c_\infty}(z_*)} \big| < \varepsilon
		\end{align}
		as well as
		\begin{align*}
			& \varphi_{\hat{\nu}_n,c_n}\big( \underbrace{(1-c_n z_*\hat{s}_{n,\varepsilon}(z_*)-c_n)z_*}_{= \tz_0} \big) \overset{\text{(\ref{Eq_BasicConsistency_Step4_07})}}{=} z_* \ ,
		\end{align*}
		which immediately yields properties (\ref{Eq_BasicConsistency_Step4_Prop1}) and (\ref{Eq_BasicConsistency_Step4_Prop2}).
		
		\item[v)] \textit{Constructing $\hat{s}_n$ from $\hat{s}_{n,\varepsilon}$}:\\
		Let $\varepsilon_0 > 0$ be small enough to satisfy (\ref{Eq_BasicConsistency_UProp1}).
		For each $\varepsilon \in (0,\varepsilon_0)$, let $N_\varepsilon$ denote the smallest integer such that the statement of (iv) still holds.
		Define the sequence $(\varepsilon_n)_{n \in \N}$ by
		\begin{align}
			& \varepsilon_n = 2^{-i} \land \varepsilon_0 \text{ for all } N_{2^{-i}} \leq n < N_{2^{-(i+1)}}
		\end{align}
		for every $i \in \N$ and $\varepsilon_n = \varepsilon_0$ for all $n<N_{\frac{1}{2}}$. The sequence converges to zero and by construction satisfies $n \geq N_{\varepsilon_n}$ whenever $\varepsilon_n < \varepsilon_0$. By (\ref{Eq_BasicConsistency_Step4_Prop1}), the maps $z \mapsto (1-c_nz\hat{s}_{n,\varepsilon_n}(z)-c_n)z$ converge to $\Phi_{H_\infty,c_\infty}$ uniformly in $z \in J$. As the compact set $J\subset\C^+$ is bounded away from zero, the maps $\hat{s}_{n,\varepsilon_n}$ thus by (\ref{Eq_Def_Phi}) converge to $\cs_{H_\infty}$ uniformly on $J$. The compactness of $J \subset \bD_{H_\infty,c_\infty}(1)\cup(\R\setminus I)$ by definition of $\bD_{H_\infty,c_\infty}(1)$ also yields
		\begin{align*}
			& \forall z \in J\setminus\R : \ \Big| \frac{c_\infty z \Im(z\cs_{H_\infty}(z))}{\Im((1-c_\infty z\cs_{H_\infty}(z)-c_\infty)z)} \Big| < \theta
		\end{align*}
		for some $\theta \in (0,1)$ and the previous two uniform convergences will guarantee
		\begin{align}\label{Eq_BasicConsistency_Step5_01}
			& \forall z \in J\setminus\R : \ \Big| \frac{c_n z \Im(z\hat{s}_{n,\varepsilon_n}(z))}{\Im((1-c_nz\hat{s}_{n,\varepsilon_n}(z)-c_n)z)} \Big| < 1
		\end{align}
		for sufficiently large $n$. Furthermore, it follows from the definition of $\varphi_{\hat{\nu}_n,c_n}$ that the sign of $\Im(\varphi_{\hat{\nu}_n,c_n}(\tz))$ will coincide with the sign of $\Im(\tz)\overset{\text{(\ref{Eq_BasicConsistency_UProp1})}}{>}0$, which by (\ref{Eq_BasicConsistency_Step4_Prop2}) yields
		\begin{align}\label{Eq_BasicConsistency_Step5_02}
			& \forall z \in J\setminus\R : \ \Im((1-c_nz\hat{s}_{n,\varepsilon_n}(z)-c_n)z)>0 \ .
		\end{align}
		Finally, (\ref{Eq_BasicConsistency_Step5_02}) and (\ref{Eq_BasicConsistency_Step4_Prop2}) together with (\ref{Eq_RevMPE_Equivalence1}) prove
		\begin{align}\label{Eq_BasicConsistency_Step5_03}
			& \forall z \in J\setminus\R : \ z\hat{s}_{n,\varepsilon_n}(z) + 1 = \int_\R \frac{\lambda}{\lambda - (1 - c_nz\hat{s}_{n,\varepsilon_n}(z) - c_n)z} \, d\hat{\nu}_n(\lambda) \ .
		\end{align}
		By (\ref{Eq_BasicConsistency_Step5_01}), (\ref{Eq_BasicConsistency_Step5_02}) and (\ref{Eq_BasicConsistency_Step5_03}), the value $\hat{s}_{n,\varepsilon_n}(z)$ will for large $n$ satisfy the defining properties of $\hat{s}_n(z)$ from Definition~\ref{Def_StieltjesEstimator} for all $z \in J\setminus\R$. By (a) of Lemma~\ref{Lemma_PopStilEst_Uniqueness} it follows that $\hat{s}_n = \hat{s}_{n,\varepsilon_n}$ on $J\setminus\R$ for sufficiently large $n$, i.e. there exists an $N_J>0$ such that $\hat{s}_n(z) = \hat{s}_{n,\varepsilon_n}(z)$ for all $z \in J$ and $n \geq N_J$. The previously observed convergence of $\hat{s}_{n,\varepsilon_n}$ to $\cs_{H_\infty}$ becomes
		\begin{align}\label{Eq_BasicConsistency_Step5_02_old}
			& \sup\limits_{z \in J\setminus\R} \big| \hat{s}_n(z) - \cs_{H_\infty}(z) \big| \xrightarrow{n \to \infty} 0 \ ,
		\end{align}
		thus proving (\ref{Eq_ConsistencyBasic_Stronger}). \qed
	\end{itemize}

	\section{Consistency for inference of population Stieltjes transforms}\label{Section_ConsistencyProof}
	This section is dedicated to proving Theorem~\ref{Thm_Consistency}, the main contribution of this paper to the field of eigen-inference. A key ingredient is the following special case of an anisotropic local law shown in \cite{KnowlesAnisotropicLocalLaws}. Since the spectral domain $\hat{\bS}(\tau,n)$ is not allowed to approach the support of $\nu_n$, as opposed to classical local laws, the convergence rate $n^{\tilde{\varepsilon}-1}$ holds uniformly on the entire domain. A full derivation of this theorem from results in \cite{KnowlesAnisotropicLocalLaws} is included in Section~\ref{Section_OuterLaws}.
	
	\begin{theorem}[Knowles-Yin: Outer local law]\label{Corollary_OuterLaw}\
		\\
		Suppose Assumption~\ref{Assumption_EigInf_Main} holds and define the spectral domain
		\begin{align}\label{Eq_Def_SInfty}
			& \hat{\bS}(\tau,n) \coloneq \big\{ \tz \in \C^+ \ \big| \ \dist(\tz,\supp(\ul{\nu}_n)\cup\supp(\hat{\ul{\nu}}_n)\cup\{0\}) \geq 2\tau \big\} \ .
		\end{align}
		For any $\tilde{\varepsilon},D > 0$ there exists a constant $C = C(\tilde{\varepsilon},\tau,D) > 0$ such that
		\begin{align}\label{Eq_OuterLaw_StieltjesInfty}
			& \bP\Big( \forall \tz \in \hat{\bS}(\tau,n) : \ \big| \cs_{\hat{\ul{\nu}}_n}(\tz) - \cs_{\ul{\nu}_n}(\tz) \big| \leq n^{\tilde{\varepsilon}-1} \Big) \geq 1 - \frac{C}{n^D}
		\end{align}
		holds for all $n\in \N$.
	\end{theorem}
	
	With this, the proof of Theorem~\ref{Thm_Consistency} commences. The critical step is part (iii), where again Rouch\'{e}'s theorem is applied, but in the opposite direction to that in the proof of Lemma~\ref{Lemma_ConsistencyBasic}, to show that the true $\Phi_{H_n,c_n}(z)$ must lie in an $\varepsilon_n$-neighborhood of $\hat{\Phi}_n(z) = (1-c_nz\hat{s}_n(z)-c_n)z$.

	\subsection{Proof of Theorem~\ref{Thm_Consistency}}\label{Proof_Thm_Consistency}
	The proof of measurability of the event from (\ref{Eq_Consistency_Result}) is deferred to Section~\ref{Section_Measurability}.
	\begin{itemize}
		\item[i)] \textit{Basic definitions and assumptions}:\\
		Define
		\begin{align*}
			& \tilde{\varepsilon} \coloneq \frac{\varepsilon}{4} \ \ \text{ and } \ \ \varepsilon_n \coloneq n^{3\tilde{\varepsilon}-1} \ .
		\end{align*}
		Without loss of generality, assume $\kappa$ and $n$ to be large enough that the following technical conditions hold:\\
		\begin{minipage}{0.43\textwidth}
			\begin{align}
				& \varepsilon_n < \frac{\tau}{2} < \kappa \label{Eq_Consistency_LargeNAssumption1}\\
				& \frac{\varepsilon_n}{\tau^2/4} \leq \frac{1}{2\kappa} \label{Eq_Consistency_LargeNAssumption2}\\
				& n^{\tilde{\varepsilon}-1} \leq \frac{1}{4\kappa} \label{Eq_Consistency_LargeNAssumption3}
			\end{align}
		\end{minipage}
		\hspace{0.3cm}\vrule\hspace{0.3cm}
		\begin{minipage}{0.43\textwidth}
			\begin{align}
				& 8\kappa^2 \leq n^{\tilde{\varepsilon}} \label{Eq_Consistency_LargeNAssumption4}\\
				& 2n^{2\tilde{\varepsilon}-1} \leq \frac{\varepsilon_n}{4} \label{Eq_Consistency_LargeNAssumption5}\\
				& \frac{\varepsilon_n}{c_n \tau^2} \leq n^{4\tilde{\varepsilon}-1} \label{Eq_Consistency_LargeNAssumption6}
			\end{align}
		\end{minipage}
		\\
		Let $\hat{U}(\tau,\kappa,n,\varepsilon_n)$ denote the closed $\varepsilon_n$-neighborhood of $\hat{\Phi}_n(\hat{\bD}(\tau,\kappa,n))$ in the sense
		\begin{align*}
			& \hat{U}(\tau,\kappa,n,\varepsilon_n) = \big\{ \tz \in \C \ \big| \ \dist(\tz, \hat{\Phi}_n(\hat{\bD}(\tau,\kappa,n))) \leq \varepsilon_n \big\} \ .
		\end{align*}
		While $\hat{\Phi}_n(\hat{\bD}(\tau,\kappa,n))$ by construction lies in $\C^+$, the set $\hat{U}(\tau,\kappa,n,\varepsilon_n)$ need not necessarily lie in $\C^+$, which will require an extra argument at the end of step (ii).
		Assumptions (\ref{Eq_Consistency_EmpDomain2})--(\ref{Eq_Consistency_EmpDomain4}) with (\ref{Eq_Consistency_LargeNAssumption1}) nonetheless for all $\tz \in \hat{U}(\tau,\kappa,n,\varepsilon_n)$ guarantee:
		\begin{align}
			& \frac{\tau}{2} < |\tz| < \kappa+\frac{\tau}{2} \label{Eq_Consistency_UDomain2}\\
			& \frac{\tau}{2} < \dist\big( \tz, \{\lambda_j(\bm{S}_n) \mid j \leq d\} \big) \label{Eq_Consistency_UDomain3}\\
			& \frac{\tau}{2} < \dist\big( \tz, \supp(\nu_n) \big) \label{Eq_Consistency_UDomain4}
		\end{align}

		\item[ii)] \textit{Bounding the difference $|\varphi_{\hat{\nu}_n,c_n}(\tz) - \varphi_{\nu_n,c_n}(\tz)|$ on $\hat{U}(\tau,\kappa,n,\varepsilon_n)$}:\\
		Theorem~\ref{Corollary_OuterLaw} yields the existence of a constant $C'=C'(\tilde{\varepsilon},\tau/4,D)>0$ such that
		\begin{align}\label{Eq_Consistency_Event1}
			& \bP\Big( \forall \tz \in \hat{\bS}(\tau/4,n) : \ \big| \cs_{\hat{\ul{\nu}}_n}(\tz) - \cs_{\ul{\nu}_n}(\tz) \big| \leq n^{\tilde{\varepsilon}-1} \Big) \geq 1 - \frac{C'}{n^{D}} \ .
		\end{align}
		As properties (\ref{Eq_Consistency_UDomain2})--(\ref{Eq_Consistency_UDomain4}) imply the inclusion $\hat{U}(\tau,\kappa,n,\varepsilon_n) \subset \hat{\bS}(\tau/4,n)$, it on the same event holds that
		\begin{align}\label{Eq_Consistency_Event1_2}
			& \forall \tz \in \hat{U}(\tau,\kappa,n,\varepsilon_n) : \ \big| \cs_{\hat{\ul{\nu}}_n}(\tz) - \cs_{\ul{\nu}_n}(\tz) \big| \leq n^{\tilde{\varepsilon}-1} \ .
		\end{align}
		For every $\tz \in \hat{U}(\tau,\kappa,n,\varepsilon_n)$, there by construction exists a $\tz_0 = \hat{\Phi}_n(z_0) \in \hat{\Phi}_n(\hat{\bD}(\tau,\kappa,n))$ such that
		\begin{align}\label{Eq_Consistency_Step2_01}
			& |\tz - \tz_0| \leq \varepsilon_n \ .
		\end{align}
		By (\ref{Eq_RevMPE_Equivalence1}) and Definition~\ref{Def_StieltjesEstimator}, it holds that $z_0 = \varphi_{\hat{\nu}_n,c_n}(\tz_0) \overset{\text{(\ref{Eq_Def_varphi})}}{=} \frac{-1}{\cs_{\hat{\ul{\nu}}_n}(\tz_0)}$, which implies
		\begin{align}\label{Eq_Consistency_Step2_02}
			& |\cs_{\hat{\ul{\nu}}_n}(\tz_0)| = \frac{1}{|z_0|} \overset{\text{(\ref{Eq_Consistency_EmpDomain1_1})}}{\geq} \frac{1}{\kappa} \ .
		\end{align}
		This bound is extended to $\tz \in \hat{U}(\tau,\kappa,n,\varepsilon_n)$ with (b) of Lemma~\ref{Lemma_StandardBounds} applied to $\hat{\ul{\nu}}_n$ and the calculation
		\begin{align}\label{Eq_Consistency_Step2_03}
			& |\cs_{\hat{\ul{\nu}}_n}(\tz)| \geq |\cs_{\hat{\ul{\nu}}_n}(\tz_0)| - |\cs_{\hat{\ul{\nu}}_n}(\tz) - \cs_{\hat{\ul{\nu}}_n}(\tz_0)| \overset{\text{(\ref{Eq_Consistency_Step2_02})}}{\geq} \frac{1}{\kappa} - |\cs_{\hat{\ul{\nu}}_n}(\tz) - \cs_{\hat{\ul{\nu}}_n}(\tz_0)| \nonumber\\
			& \overset{\text{(\ref{Eq_StandardBounds_c3})}}{\geq} \frac{1}{\kappa} - \frac{|\tz-\tz_0|}{\tau^2/4} \overset{\text{(\ref{Eq_Consistency_Step2_01})}}{\geq} \frac{1}{\kappa} - \frac{\varepsilon_n}{\tau^2/4} \overset{\text{(\ref{Eq_Consistency_LargeNAssumption2})}}{\geq} \frac{1}{2\kappa}
		\end{align}
		and one may, on the event $\mathcal{E}_{n,\text{(\ref{Eq_Consistency_Event1})}}$ from (\ref{Eq_Consistency_Event1}), transfer the bound to $|\cs_{\ul{\nu}_n}(\tz)|$ by
		\begin{align}\label{Eq_Consistency_Step2_04}
			& |\cs_{\ul{\nu}_n}(\tz)| \geq |\cs_{\hat{\ul{\nu}}_n}(\tz)| - |\cs_{\hat{\ul{\nu}}_n}(\tz) - \cs_{\ul{\nu}_n}(\tz)| \nonumber\\
			& \overset{\text{(\ref{Eq_Consistency_Step2_03})}}{\geq} \frac{1}{2\kappa} - |\cs_{\hat{\ul{\nu}}_n}(\tz) - \cs_{\ul{\nu}_n}(\tz)| \overset{\text{(\ref{Eq_Consistency_Event1_2})}}{\geq} \frac{1}{2\kappa} - n^{\tilde{\varepsilon}-1} \overset{\text{(\ref{Eq_Consistency_LargeNAssumption3})}}{\geq} \frac{1}{4\kappa} \ .
		\end{align}
		On the event $\mathcal{E}_{n,\text{(\ref{Eq_Consistency_Event1})}}$ it then holds that
		\begin{align}\label{Eq_Consistency_Step2_05}
			& \big|\varphi_{\hat{\nu}_n,c_n}(\tz) - \varphi_{\nu_n,c_n}(\tz)\big| = \Big| \frac{-1}{\cs_{\hat{\ul{\nu}}_n}(\tz)} - \frac{-1}{\cs_{\ul{\nu}_n}(\tz)} \Big| = \Big| \frac{\cs_{\hat{\ul{\nu}}_n}(\tz)-\cs_{\ul{\nu}_n}(\tz)}{\cs_{\hat{\ul{\nu}}_n}(\tz)\cs_{\ul{\nu}_n}(\tz)} \Big| \nonumber\\
			& \overset{\text{(\ref{Eq_Consistency_Step2_03}), (\ref{Eq_Consistency_Step2_04})}}{\leq} 8\kappa^2 |\cs_{\hat{\ul{\nu}}_n}(\tz)-\cs_{\ul{\nu}_n}(\tz)| \overset{\text{(\ref{Eq_Consistency_Event1_2})}}{\leq} 8\kappa^2 n^{\tilde{\varepsilon}-1} \overset{\text{(\ref{Eq_Consistency_LargeNAssumption4})}}{\leq} n^{2\tilde{\varepsilon}-1}
		\end{align}
		for all $\tz \in \hat{U}(\tau,\kappa,n,\varepsilon_n)$. Note also, that (\ref{Eq_Consistency_Step2_03}) and (\ref{Eq_Consistency_Step2_04}) together with (\ref{Eq_Consistency_UDomain2})--(\ref{Eq_Consistency_UDomain4}) guarantee that the maps $\varphi_{\hat{\nu}_n,c_n}(\tz) = \frac{-1}{\cs_{\hat{\ul{\nu}}_n}(\tz)}$ and $\varphi_{\nu_n,c_n}(\tz) = \frac{-1}{\cs_{\ul{\nu}_n}(\tz)}$ are holomorphic on $\hat{U}(\tau,\kappa,n,\varepsilon_n)$. This required separate verification, since $\hat{U}(\tau,\kappa,n,\varepsilon_n)$ may contain segments of $\R$ between elements of $\supp(\hat{\nu}_n) \cup \supp(\nu_n)$, which makes (b) of Lemma~\ref{Lemma_StandardBounds} no longer applicable.
		
		\item[iii)] \textit{Application of Rouch\'{e}'s theorem}:\\
		By (d) of Lemma~\ref{Lemma_StandardBounds} it holds that
		\begin{align*}
			& \forall \tz_1\neq\tz_2 \in \C\setminus\R : \ |\varphi_{\nu_n,c_n}(\tz_1) - \varphi_{\nu_n,c_n}(\tz_2)| \geq \frac{|\tz_1-\tz_2|}{2} \ .
		\end{align*}
		At the end of (ii) it was shown that $\varphi_{\nu_n,c_n}$ is holomorphic on $\hat{U}(\tau,\kappa,n,\varepsilon_n)$ and the implied continuity extends the above bound to
		\begin{align}\label{Eq_Consistency_Step3_3}
			& \forall \tz_1\neq\tz_2 \in \hat{U}(\tau,\kappa,n,\varepsilon_n) : \ |\varphi_{\nu_n,c_n}(\tz_1) - \varphi_{\nu_n,c_n}(\tz_2)| \geq \frac{|\tz_1-\tz_2|}{2} \ .
		\end{align}
		For every $\tz_* = \hat{\Phi}_n(z_*) \in \hat{\Phi}_n(\hat{\bD}(\tau,\kappa,n))$, the boundary $\partial B_{\varepsilon_n}^{\C}(\tz_*)$ is by construction contained in $\hat{U}(\tau,\kappa,n,\varepsilon_n)$ and one for each $\tz \in \partial B_{\varepsilon_n}^{\C}(\tz_*)$ on the event $\mathcal{E}_{n,\text{(\ref{Eq_Consistency_Event1})}}$ has
		\begin{align}\label{Eq_Consistency_Step4_1}
			& |\overbrace{\varphi_{\hat{\nu}_n,c_n}(\tz) - \varphi_{\hat{\nu}_n,c_n}(\tz_*)}^{\eqcolon  F_{n,\tz_*}(\tz)}| \overset{\text{(\ref{Eq_Consistency_Step2_05})}}{\geq} |\varphi_{\nu_n,c_n}(\tz) - \varphi_{\nu_n,c_n}(\tz_*)| - 2n^{2\tilde{\varepsilon}-1} \nonumber\\
			& \overset{\text{(\ref{Eq_Consistency_Step3_3})}}{\geq} \frac{|\tz-\tz_*|}{2} - 2n^{2\tilde{\varepsilon}-1} = \frac{\varepsilon_n}{2} - 2n^{2\tilde{\varepsilon}-1} \overset{\text{(\ref{Eq_Consistency_LargeNAssumption5})}}{\geq} \frac{\varepsilon_n}{4}
		\end{align}
		as well as
		\begin{align}\label{Eq_Consistency_Step4_2}
			& |\overbrace{\varphi_{\nu_n,c_n}(\tz) - \varphi_{\hat{\nu}_n,c_n}(\tz)}^{\eqcolon  G_{n}(\tz)}| \overset{\text{(\ref{Eq_Consistency_Step2_05})}}{\leq} n^{2\tilde{\varepsilon}-1} \overset{\text{(\ref{Eq_Consistency_LargeNAssumption5})}}{\leq} \frac{\varepsilon_n}{8} \overset{\text{(\ref{Eq_Consistency_Step4_1})}}{<} |F_{n,\tz_*}(\tz)| \ .
		\end{align}
		The functions $F_{n,\tz_*}$ and $G_{n}$ are by (ii) holomorphic on $\hat{U}(\tau,\kappa,n,\varepsilon_n)$ and thus Rouch\'{e}'s theorem with (\ref{Eq_Consistency_Step4_2}) implies that $F_{n,\tz_*}$ and $F_{n,\tz_*}+G_{n}$ have the same number of zeros/roots in $B_{\varepsilon_n}^{\C}(\tz_*)$, when counting multiplicity. Since $F_{n,\tz_*}(\tz_*) = 0$ holds by construction, there must be at least one $\tz_0 \in B_{\varepsilon_n}^{\C}(\tz_*)$ with
		\begin{align*}
			& 0 = F_{n,\tz_*}(\tz_0)+G_{n}(\tz_0) = \varphi_{\nu_n,c_n}(\tz_0) - \varphi_{\hat{\nu}_n,c_n}(\overbrace{\tz_*}^{= \hat{\Phi}_n(z_*)}) \overset{\text{(\ref{Eq_RevMPE_Equivalence1})}}{=} \varphi_{\nu_n,c_n}(\tz_0) - z_*
		\end{align*}
		and Lemma~\ref{Lemma_SpaceTransform} guarantees that $\tz_0 = \Phi_{H_n,c_n}(z_*)$. Since $\tz_0 \in B_{\varepsilon_n}^{\C}(\tz_*)$, it was thus shown that
		\begin{align}\label{Eq_Consistency_Step4_3}
			& |\hat{\Phi}_n(z_*) - \Phi_{H_n,c_n}(z_*)| \leq \varepsilon_n
		\end{align}
		holds on the event $\mathcal{E}_{n,\text{(\ref{Eq_Consistency_Event1})}}$ for all $z_* \in \hat{\bD}(\tau,\kappa,n)$.
		
		\item[iv)] \textit{Choice of $C$}:\\
		On the event $\mathcal{E}_{n,\text{(\ref{Eq_Consistency_Event1})}}$ it follows from (\ref{Eq_Consistency_Step4_3}) that every $z_* \in \hat{\bD}(\tau,\kappa,n)$ satisfies
		\begin{align*}
			& |\hat{s}_n(z) - \cs_{H_n}(z)| = \frac{1}{c_n |z|^2} \big| \overbrace{(1-c_nz\hat{s}_n(z)-c_n)z}^{= \hat{\Phi}_n(z)} - \overbrace{(1-c_nz\cs_{H_n}(z)-c_n)z}^{= \Phi_{H_n,c_n}(z)} \big|\\
			& \overset{\text{(\ref{Eq_Consistency_Step4_3})}}{\leq} \frac{\varepsilon_n}{c_n |z|^2} \overset{\text{(\ref{Eq_Consistency_EmpDomain1_1})}}{\leq} \frac{\varepsilon_n}{c_n \tau^2} \overset{\text{(\ref{Eq_Consistency_LargeNAssumption6})}}{\leq} n^{4\tilde{\varepsilon}-1} = n^{\varepsilon-1} \ .
		\end{align*}
		By choosing $C>0$ larger than $C'$ from (\ref{Eq_Consistency_Event1}) and large enough that (\ref{Eq_Consistency_LargeNAssumption1})--(\ref{Eq_Consistency_LargeNAssumption6}) hold for all $n \geq C^{\frac{1}{D}}$ the assertion of Theorem~\ref{Thm_Consistency} follows. \qed
	\end{itemize}

	\section{Numerical applications and comparison to existing methods}\label{Section_Numerical}
	The following section showcases numerical applications of the proposed method and their performance when compared with existing methods. All algorithms in this section will forgo enforcing condition (\ref{Eq_Consistency_EmpDomain4}) to ensure that the algorithms do not rely on prior knowledge of $H_n$. By Remark~\ref{Remark_Replaceability} and the homotopy invariance of Remark~\ref{Remark_Homotopy}, the consistency results are still valid for all consequent estimation tasks involving examples (Ex.1) and (Ex.2).
	\\[0.5em]
	Example (Ex.3) is not intended as an example covered by Theorem \ref{Thm_Consistency}. Since it is constructed from real data, assumptions~\ref{EI_ItemAssumption_PopConv} and~\ref{EI_ItemAssumption_sigmaBound} cannot be verified. Instead, estimation tasks involving example (Ex.3) serve as a stress test for the proposed method. Its large number of outlier eigenvalues makes it a natural setting in which condition (\ref{Eq_Consistency_EmpDomain4}) is expected to become relevant.
	\\[0.5em]
	A full \texttt{Python} implementation of the proposed method is available from the author's GitHub repository \footnote{\href{https://github.com/BenDeitmar/PLSS_estimation}{{https://github.com/BenDeitmar/PLSS\_estimation}}}. The code to produce Figures \ref{Fig_SpaceTransform}--\ref{Fig_Intervals} may be found in the same repository.
	
	\newpage
	\begin{definition}[Example data] \label{Def_Examples} The following examples are formulated in the notation of~\ref{EI_ItemAssumption_Asymptotics}--\ref{EI_ItemAssumption_PopConv} from Assumption~\ref{Assumption_EigInf_Main}.
		\begin{itemize}
			\item[Ex.1)] Let
			\begin{align*}
				& B_n \coloneq \diag\Big(\underbrace{1,\dots,1}_{\times \lfloor d/2 \rfloor},\sqrt{2.5+\frac{1}{d}},\sqrt{2.5+\frac{2}{d}},\dots,\sqrt{2.5+\frac{d - \lfloor d/2 \rfloor}{d}}\Big) \ ,
			\end{align*}
			such that $H_n = \frac{\lfloor d/2 \rfloor}{d} \delta_{1} + \frac{1}{d} \sum\limits_{j=1}^{d-\lfloor d/2 \rfloor} \delta_{2.5+\frac{j}{d}}$ and $H_\infty = \frac{1}{2} \delta_{1} + \frac{1}{2} \operatorname{Uniform}([2.5,3])$. Suppose $\bm{X}_n$ has i.i.d. standard normal entries. Unless stated otherwise, let $n \coloneq 10d$, such that $c_n = \frac{1}{10}$.
			
			\item[Ex.2)] Suppose the i.i.d. entries of $\bm{X}_n$ satisfy $\bP((\bm{X}_n)_{j,k}=1)=\frac{1}{2}=\bP((\bm{X}_n)_{j,k}=-1)$ and for a $\operatorname{Haar}(O(d))$-distributed orthogonal matrix $V$ let
			\begin{align*}
				& B_n \coloneq V\diag\Big(\sqrt{\frac{1}{2}+\frac{1}{2d}},\sqrt{\frac{1}{2}+\frac{2}{2d}},\dots,1\Big)
			\end{align*}
			such that $H_n = \frac{1}{d} \sum\limits_{j=1}^{d} \delta_{\frac{1}{2}+\frac{j}{2d}}$ and $H_\infty = \operatorname{Uniform}\big(\big[\frac{1}{2},1\big]\big)$. Let $n \coloneq \lceil \frac{d}{2} \rceil$, such that $c_n \approx 2$.
			
			\item[Ex.3)] Suppose $d \leq 784$ and $n \leq 70{,}000$. Let $\bm{Y}_n$ be a $(d \times n)$ matrix constructed by randomly selecting $d$ rows and $n$ columns (without replacement) from the $(784 \times 70{,}000)$ matrix whose columns are vectorized, normalized and centralized images of the MNIST dataset \cite{MNIST}. The sample covariance matrix under consideration is $\bm{S}_n=\frac{1}{n}\bm{Y}_n\bm{Y}_n^\top$, while the full sample of size $70{,}000$ leads to a surrogate population covariance matrix $\Sigma_n \in \R^{d \times d}$ and also a surrogate population spectral distribution $H_n$. Let $n \coloneq 4d$, such that $c_n = \frac{1}{4}$.
		\end{itemize}
	\end{definition}

	\subsection{Numerical estimation of population Stieltjes transforms}
\begin{figure}[H]
		\centering
		\includegraphics[width=0.3\textwidth]{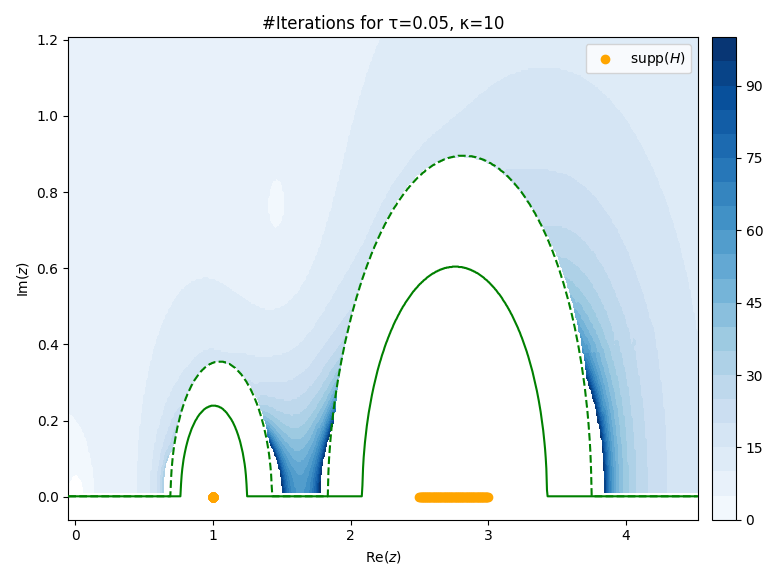}
		\includegraphics[width=0.3\textwidth]{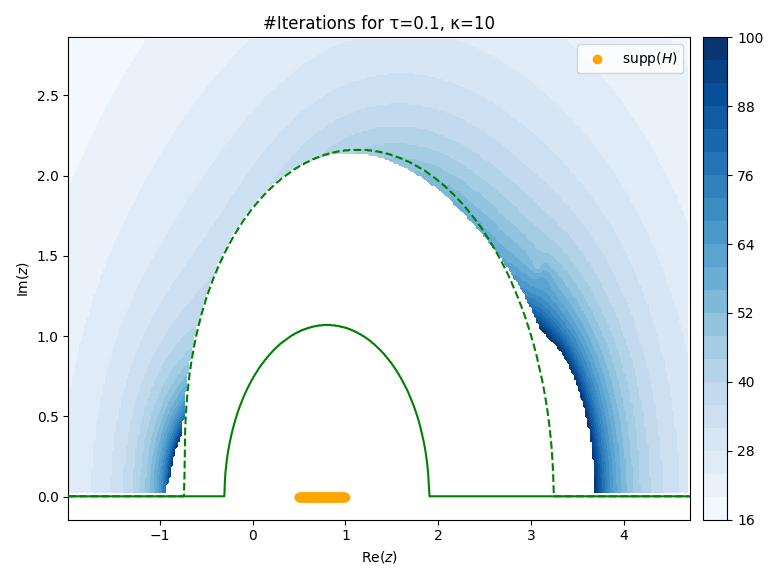}
		\includegraphics[width=0.3\textwidth]{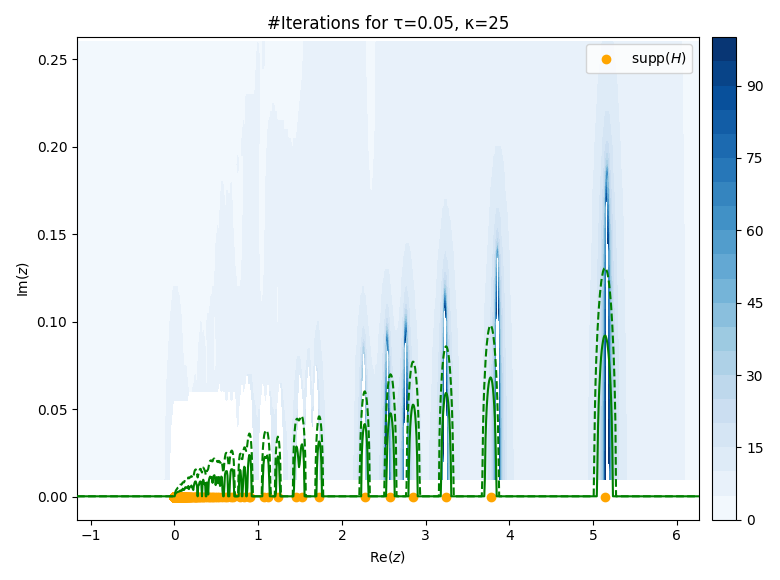}
		\vspace{-0.25cm}
		\caption{Contour plots showing the number of iterations of (\ref{Eq_StilEstimator_iteration}) required to reach $|v_{i}-v_{i-1}| \leq 10^{-9}$ for Examples (Ex.1) with $d=200$ (left), (Ex.2) with $d=100$ (middle) and (Ex.3) with $d=784$ (right).}\label{Fig_Iterations}
		\vspace{-0.3cm}
	\end{figure}

	The population Stieltjes transform estimator $\hat{s}_n(z)$ of Definition~\ref{Def_StieltjesEstimator} may be computed efficiently using fixed-point iterations. With the notation $v_* = z\hat{s}_n(z)+1$, equation (\ref{Eq_ReverseMPEquation_Empirical_again}) translates to
	\begin{align*}
		& v_* = \int_\R \frac{\lambda}{\lambda-(1-c_nv_*)z} \, d\hat{\nu}_n(\lambda) \ .
	\end{align*}
	By iterative application of the operator
	\begin{align}\label{Eq_StilEstimator_iteration}
		& v_{i+1} = \hat{T}_{n,z}(v_i) \coloneq \int_\R \frac{\lambda}{\lambda-(1-c_nv_i)z} \, d\hat{\nu}_n(\lambda)
	\end{align}
	to any starting point $v_0 \in \C^+$, the solution $v_*=z\hat{s}_n(z)+1$ is typically found with fewer than twenty iterations (compare Figure~\ref{Fig_Iterations}).

	\subsection{Comparison of PLSS estimation methods}
	Several estimators for $H_n$, and by extension for population linear spectral statistics $L_n(g) = \int_\R g(\lambda) \, dH_n(\lambda)$ have been proposed in the literature. Two of the most prominent nonparametric approaches are:
	\begin{itemize}
		\item[a)] El Karoui's method from \cite{MPIKaroui} for estimating $H_n$ by solving a convex optimization problem induced by the Marchenko--Pastur equation. The implementation was done in \texttt{Python} according to Section 3.2 and the appendix of \cite{MPIKaroui}. Population linear spectral statistic estimators are then found by $\int_\R g(\lambda) \, d\hat{H}_n$, where $\hat{H}_n$ is the estimator for $H_n$.
		
		\item[b)] Ledoit--Wolf's method from \cite{LedoitWolf2} and \cite{LedoitWolfNumerical} for the estimation of population eigenvalues, i.e. $\lambda_1(\Sigma_n),\dots,\lambda_d(\Sigma_n)$, 
		by solving a non-convex optimization problem arising from the connection between $H_n$ and $\nu_n$ explored in Lemma~\ref{Lemma_SpaceTransform}. The following comparison utilizes the implementation from the \texttt{R} package \texttt{nlshrink}. The estimators $\hat{\lambda}_1,\dots,\hat{\lambda}_d$ induce an estimator for population linear spectral statistics by $\frac{1}{d} \sum\limits_{j=1}^d g(\hat{\lambda}_j)$.
	\end{itemize}
	\vspace{-0.3cm}
	Several methods of consistently estimating population spectral moments $\int_\R \lambda^k \, dH_n(\lambda)$ have also been developed in the literature, with the aim of reconstructing $H_n$ from its moments. 
	Consistency of the resulting estimators for $H_n$ has only been shown for parametric models, where $H_n = \theta_1 \delta_{t_1} + \dots + \theta_m \delta_{t_m}$ for fixed $m \in \N$ and parameters $\theta_1,...,\theta_m,t_1,...,t_m \geq 0$. The estimators for the population spectral moments $\int_\R \lambda^k \, dH_n(\lambda)$ employed in the following two papers, are however equally valid in non-parametric models.
	\begin{itemize}
		\item[c)] Bai, Chen and Yao's method from \cite{MPIBaiChenYao} for estimation of population spectral moments by inverting the combinatorial formulas linking the moments of $\nu_n$ to those of $H_n$. The implementation in \texttt{Python} follows Lemma 1 in \cite{MPIBaiChenYao}.
		
		\item[d)] Kong and Valiant's method from \cite{MPIKong} for estimation of population spectral moments by calculating traces of modified sample covariance matrices. The implementation in \texttt{Python} follows Algorithm 1 in \cite{MPIKong}.
	\end{itemize}

	\subsubsection{Estimation of Stieltjes transforms}
\begin{figure}[H]
		\centering
		\includegraphics[width=0.325\textwidth]{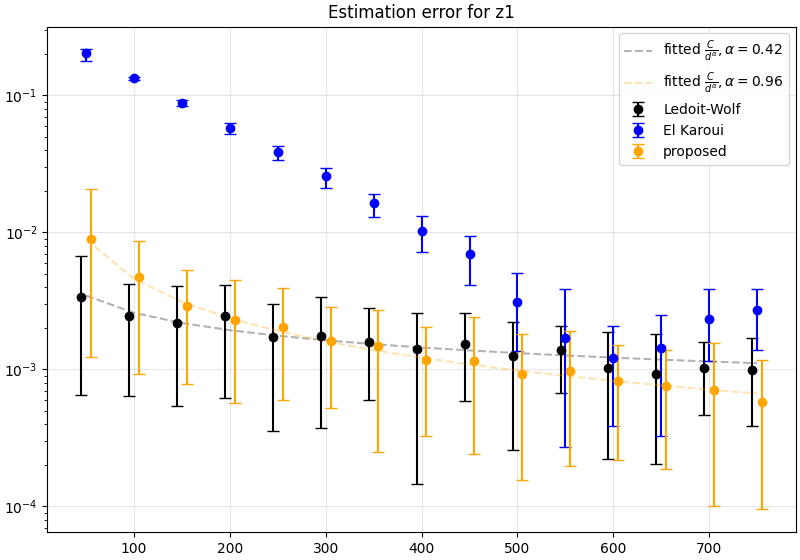}
		\includegraphics[width=0.325\textwidth]{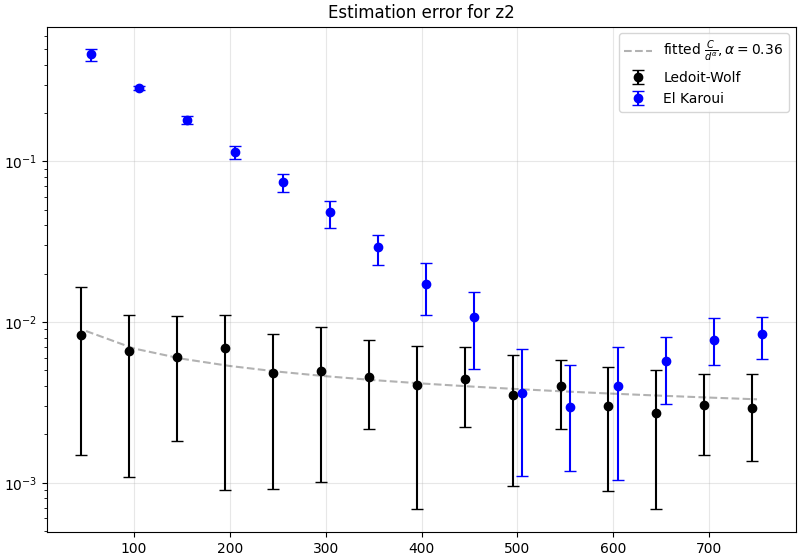}
		\includegraphics[width=0.325\textwidth]{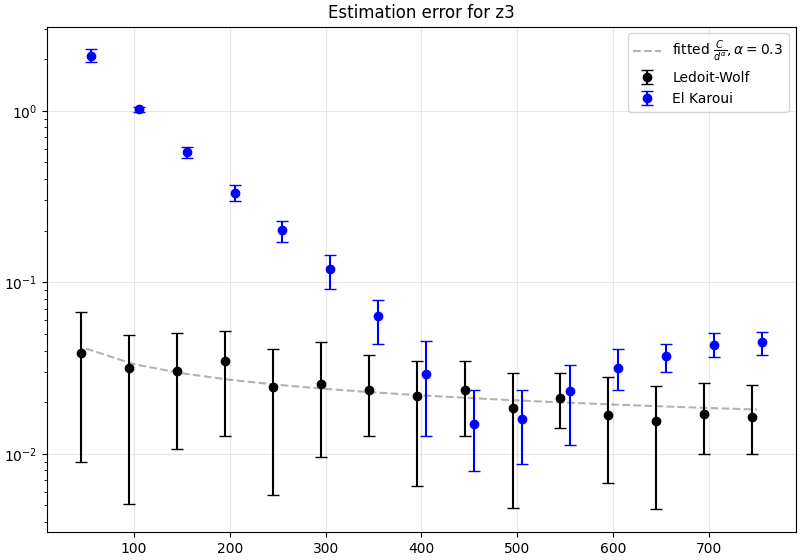}
		\begin{minipage}{0.7\textwidth}
			\caption{Average absolute errors for estimators of $\cs_{H_n}(z)$ for $z_1=-1+0.1\bm{i}$ (left), $z_2=-\frac{1}{2}+0.1\bm{i}$ (middle) and $z_3=0+0.1\bm{i}$ (right). The data is from 50 realizations of example (Ex.2) where the dimension ranges from 50 to 750. The points are chosen such that $z_1 \in \bD_{H_n,c_n}(1)$, $z_2 \in \bD_{H_n,c_n}(\infty)\setminus\bD_{H_n,c_n}(1)$ and $z_3 \notin \bD_{H_n,c_n}(\infty)$ (see bottom right). The bars in the plots represent the range spanned by the empirical 10$\%$ and $90\%$ quantiles and the green lines are as in Figure~\ref{Fig_StieltjesEstimationError}. The scale for the error axis is logarithmic.}\label{Fig_Stieltjes_comparison}
		\end{minipage}
		\begin{minipage}{0.29\textwidth}
			\includegraphics[width=\textwidth]{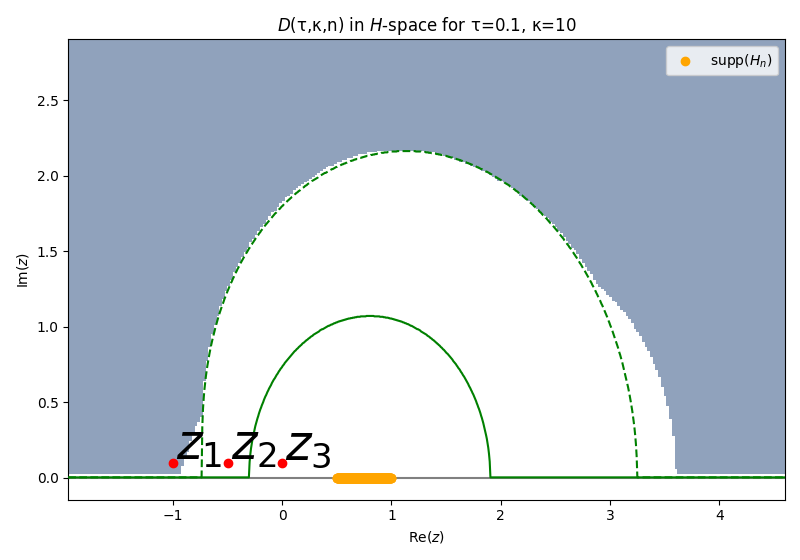}
		\end{minipage}
		\vspace{-0.5cm}
	\end{figure}

	Figure~\ref{Fig_Stieltjes_comparison} illustrates that the estimation methods of El Karoui and Ledoit--Wolf do not necessarily achieve an error rate of order $\mathcal{O}(n^{-1})$ for the estimation of population Stieltjes transforms $\cs_{H_n}(z)$ even when $z$ is well separated from the support of $H_n$. In exchange, they are less reliant on $z$ lying in the domain $\bD_{H_n,c_n}(1)$ or $\bD_{H_n,c_n}(\infty)$. The proposed method is observed to achieve an error rate near $\mathcal{O}(n^{-1})$, when applicable.
	
	\subsubsection{Estimation of monomials}
\begin{figure}[H]
		\centering
		\includegraphics[width=0.99\textwidth]{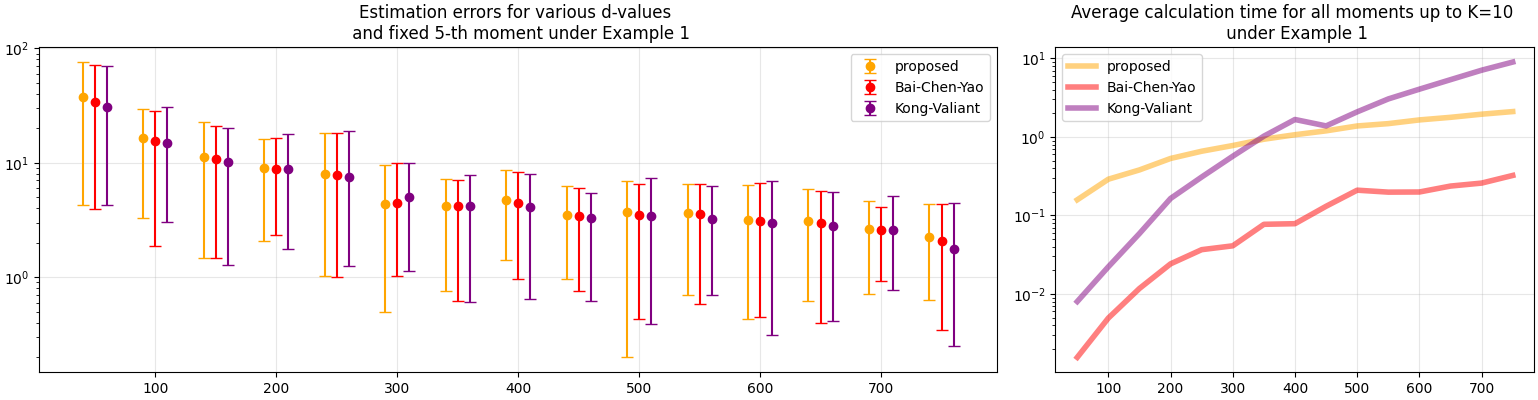}\\
		\includegraphics[width=0.99\textwidth]{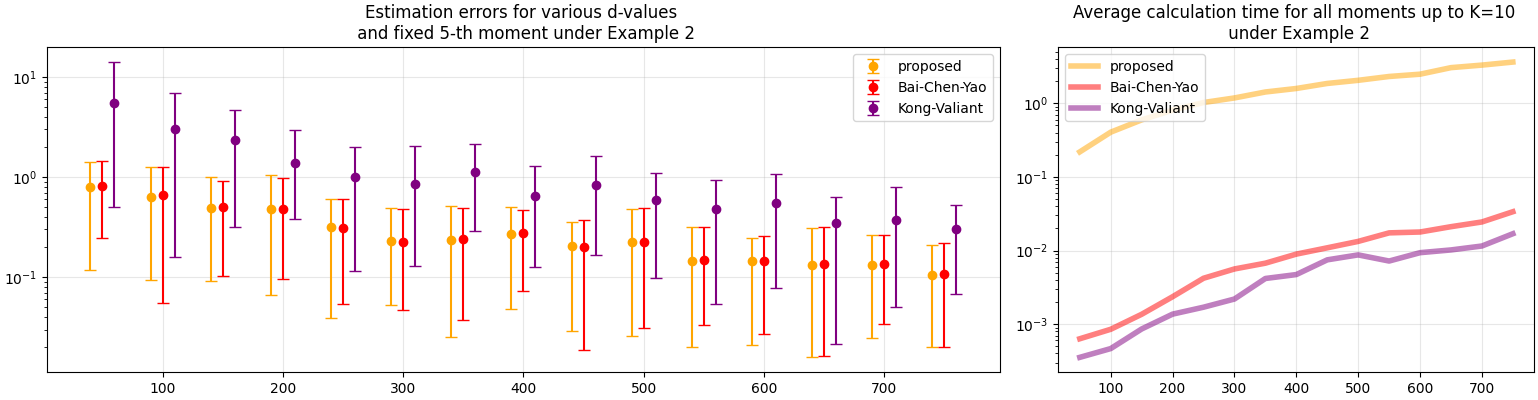}\\
		\includegraphics[width=0.99\textwidth]{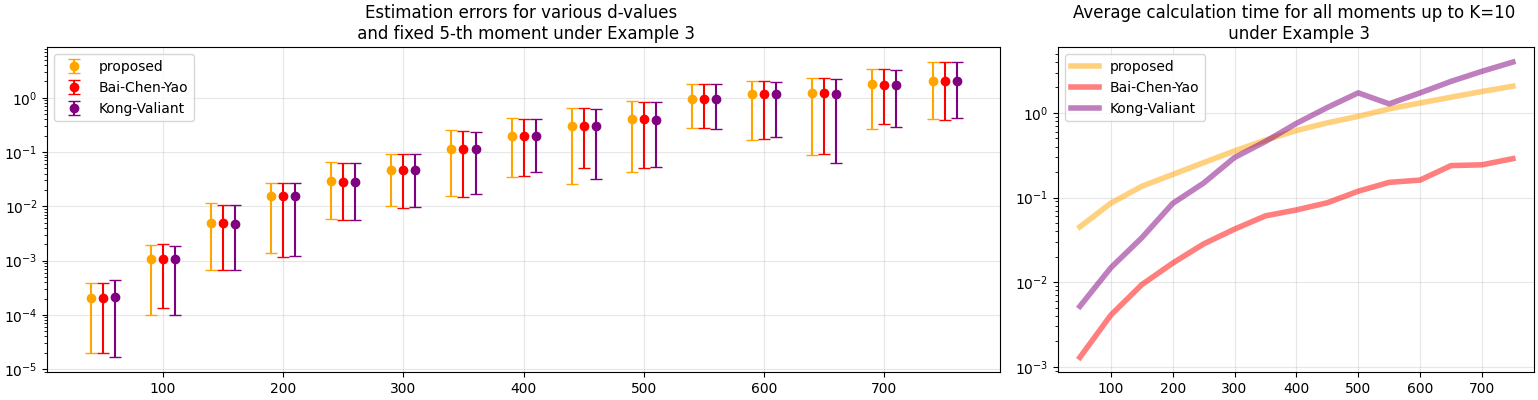}
		\vspace{-0.25cm}
		\caption{Average errors (left) and average computation times (right) for estimators of $\int_\R \lambda^5 \, dH_n(\lambda)$. The data is from 50 realizations of examples (Ex.1) with $c_n = \frac{1}{5}$: (top), (Ex.2) with $c_n = 2$: (middle) and (Ex.3) with $c_n = \frac{1}{4}$: (bottom), where the dimension ranges from 50 to 750.
			All plots use a logarithmic scale for the $y$-axis.}\label{Fig_Monomial_comparison}
		\vspace{-0.3cm}
	\end{figure}

\begin{figure}[H]
		\centering
		\includegraphics[width=0.325\textwidth]{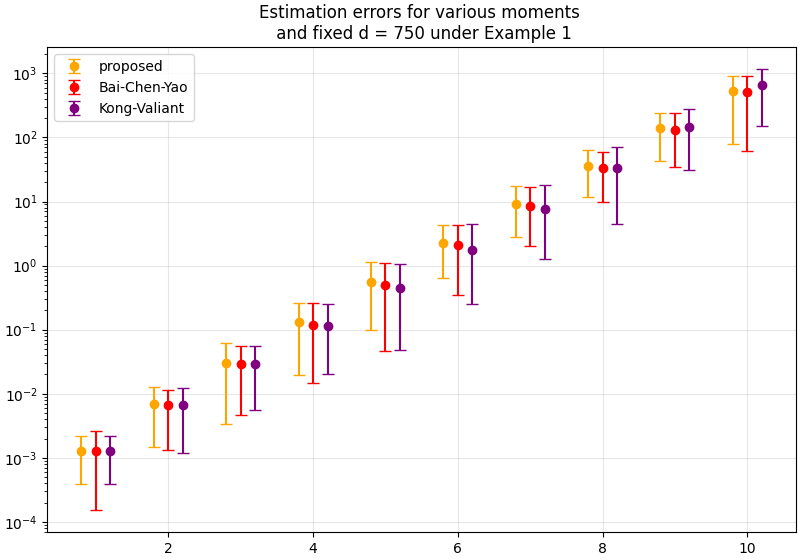}
		\includegraphics[width=0.325\textwidth]{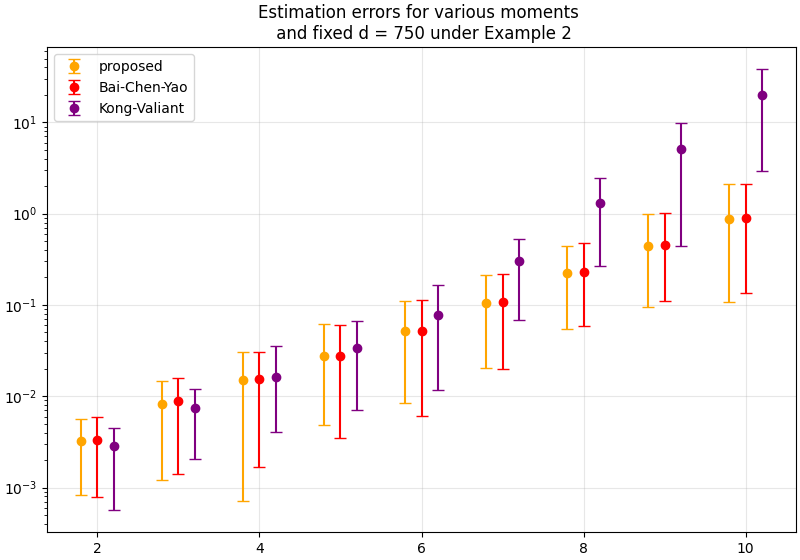}
		\includegraphics[width=0.325\textwidth]{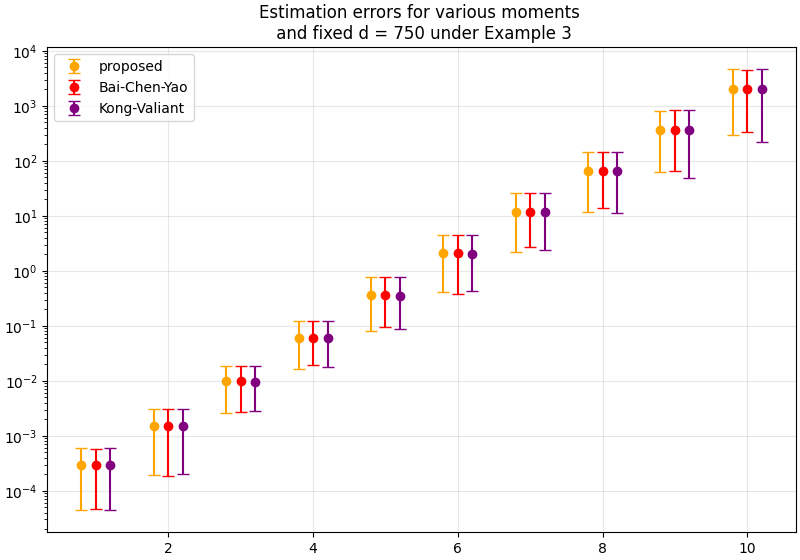}
		\vspace{-0.25cm}
		\caption{Average errors estimators of $\int_\R \lambda^k \, dH_n(\lambda)$ with $k$ ranging from $1$ to $10$. The data is from 50 realizations of examples (Ex.1) with $c_n = \frac{1}{5}$ (left), (Ex.2) with $c_n = 2$ (middle) and (Ex.3) with $c_n = \frac{1}{4}$ (right), where the dimension is $d=750$.
			All plots use a logarithmic scale for the $y$-axis.}\label{Fig_Monomial_powers}
		\vspace{-0.4cm}
	\end{figure}

	Specialized estimators for population spectral moments $\int_\R \lambda^k \, dH_n(\lambda) = \frac{1}{d} \tr(\Sigma_n^k)$ were introduced in \cite{MPIBaiChenYao} and \cite{MPIKong} (see (c) and (d) above).
	Figures~\ref{Fig_Monomial_comparison} and~\ref{Fig_Monomial_powers} compare the proposed method to the specialized estimators of Bai--Chen--Yao and Kong--Valiant in terms of:
	\begin{itemize}
		\item[1)] estimation error when estimating $\int_\R \lambda^5 \, dH_n(\lambda)$ for increasing $d$ and $n$,
		
		\item[2)] computation time when estimating $\big(\int_\R \lambda^k \, dH_n(\lambda)\big)_{k \in \{1,...,10\}}$ for increasing $d$ and $n$.
		
		\item[3)] estimation error for increasing power $k$ for fixed dimension $d=750$,
	\end{itemize}
	The proposed method is observed to be very close to the specialized method of Bai--Chen--Yao in terms of estimation error, which may be explained by the fact that both methods invert the same relationship between $H_n$ and the deterministic equivalent $\nu_n$ to construct the estimator.
	\\[0.5em]
	In the bottom plot of Figure \ref{Fig_Monomial_comparison} (Example (Ex.3)), the estimation errors are not observed to go to zero for growing $d$ and $n$. This is not an error of the proposed method, but rather an artifact from the way example (Ex.3) is constructed from real data by sub-sampling both columns (dimension) and rows (samples) from an underlying data-matrix. Conditions \ref{EI_ItemAssumption_PopConv} and \ref{EI_ItemAssumption_sigmaBound} can not be verified for this example and the lack of convergence is likely due to the range of the eigenvalues $[\lambda_{\min}(\Sigma_n),\lambda_{\max}(\Sigma_n)]$ growing with increasing dimension (number of sub-sampled columns).
	\\[0.5em]
	Figure \ref{Fig_Monomial_powers} demonstrates how the estimation errors of all methods grow exponentially with the power $k$ of the monomial.

	\begin{remark}[Moment relationships and the Marchenko--Pastur inversion formula]\
		\\
		The relationship between the population moments $H_n(\bullet^k) = \int_\R \lambda^k \, dH_n(\lambda)$ and the deterministic equivalents $\nu_n(\bullet^k) = \int_\R \lambda^k \, d\nu_n(\lambda)$ may by combinatorial methods be shown to satisfy
		\begin{align*}
			\nu_n(\bullet^K) = \frac{K!}{c_n} \sum_{b=1}^K \frac{c_n^b}{b!(K+1-b)!} \sum_{\substack{k_1,\ldots,k_b \geq 1 \\ k_1+\cdots+k_b = K}} H_n(\bullet^{k_1}) \cdots H_n(\bullet^{k_b}) \ ,
		\end{align*}
		or (by Theorem~2.2 of~\cite{MottelsonINCP}) equivalently
		\begin{align*}
			H_n(\bullet^K) = -\frac{1}{c_n} \sum_{b=1}^K \frac{(-c_n)^b (K+b-2)!}{b!(K-1)!} \sum_{\substack{k_1,\ldots,k_b \geq 1 \\ k_1+\cdots+k_b = K}} \nu_n(\bullet^{k_1}) \cdots \nu_n(\bullet^{k_b}) \ .
		\end{align*}
		The above equation may, together with Lagrange's inversion formula, be used to show that $\cs_{H_n}(z)$ for sufficiently large $|z|$ satisfies (\ref{Eq_ReverseMPEquation_re}), thus providing an alternate way to arrive at the Marchenko--Pastur inversion formula from Lemma \ref{Lemma_MPI}.
	\end{remark}

	\subsubsection{Estimation of generalized linear spectral statistics}
	\begin{minipage}{0.5\textwidth}
		In Theorem \ref{Thm_GLSS_consistency}, the presented estimator for generalized linear spectral statistics (GLSS) of the form $\frac{1}{d} \tr\big( f(\bm{S}_n) g(\Sigma_n) \big)$ was shown to be consistent with error rate $\mathcal{O}(n^{\varepsilon-\frac{1}{2}})$. While GLSS were already consistently estimateable by combining the approximation (\ref{Eq_GLSS_DeterministicEquivalent}) with the Ledoit--Wolf estimator for $H_n$, the proposed method intrinsically provides criteria for deciding when eigenvalues are sufficiently well separated for their eigenvectors to be discernable.
		\\[0.5em]
		Figure~\ref{Fig_LargestEigVec} compares the proposed method with that of Ledoit--Wolf when estimating $|u_1^\top v_1|^2$, where $u_1$ is the sample eigenvector to the largest eigenvalue of $\bm{S}_n$ and $v_1$ is the population eigenvector to the largest eigenvalue of $\Sigma_n$.
		\\[0.5em]
		If $f = \mathbbm{1}_{\bullet \geq x}$ for $x \in (\lambda_2(\bm{S}_n),\lambda_1(\bm{S}_n))$ and $g = \mathbbm{1}_{\bullet \geq y}$ for $y \in (\lambda_2(\Sigma_n),\lambda_1(\Sigma_n))$, then the quantity $|u_1^\top v_1|^2$ is equal to $\tr\big( f(\bm{S}_n) g(\Sigma_n) \big)$ and may thus be estimated as a GLSS.
	\end{minipage}
	\hspace{0.3cm}
	\begin{minipage}{0.49\textwidth}
		\vspace{-0.5cm}
		\begin{figure}[H]
			\includegraphics[width=\textwidth]{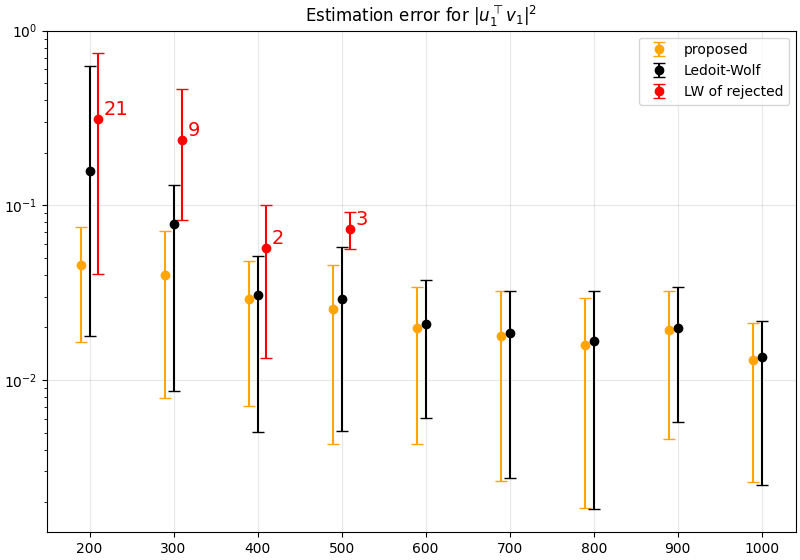}
			\vspace{-0.6cm}
			\caption{ Average absolute estimation errors of $|u_1^\top v_1|^2$, computed for $50$ realizations of example (Ex.3) with $d=784$ and $n$ ranging from $200$ to $1000$. The bars represent the range spanned by the empirical $10\%$ and $90\%$ quantiles. Separate error-bars for the Ledoit--Wolf errors to only those realizations rejected by the proposed method, together with the number of rejections, are included in red.
			}\label{Fig_LargestEigVec}
		\end{figure}
	\end{minipage}
	\\[0.5em]
	The proposed method rejects the estimation-task, if the largest eigenvalue $\lambda_1(\Sigma_n)$ may not be separated from the remainder of $\supp(H_n)$ by an admissible curve. In Figure~\ref{Fig_LargestEigVec}, it is observed that this rejection does seem to filter out cases in which $|u_1^\top v_1|^2$ may not be well estimated. Note that the scale of the $y$-axis in Figure~\ref{Fig_LargestEigVec} is logarithmic.

	\subsection{CLT and confidence intervals}\label{Subsection_InferenceCLT}
\begin{figure}[bp]
		\centering
		\includegraphics[width=0.32\textwidth]{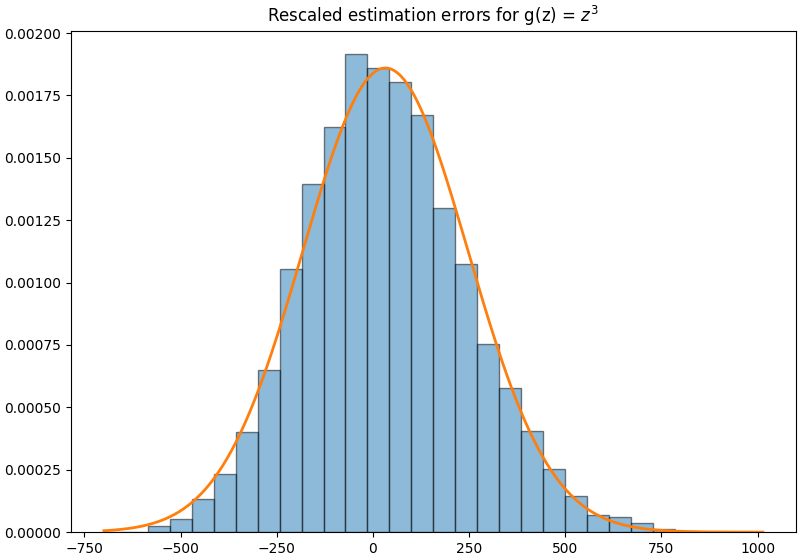}
		\includegraphics[width=0.32\textwidth]{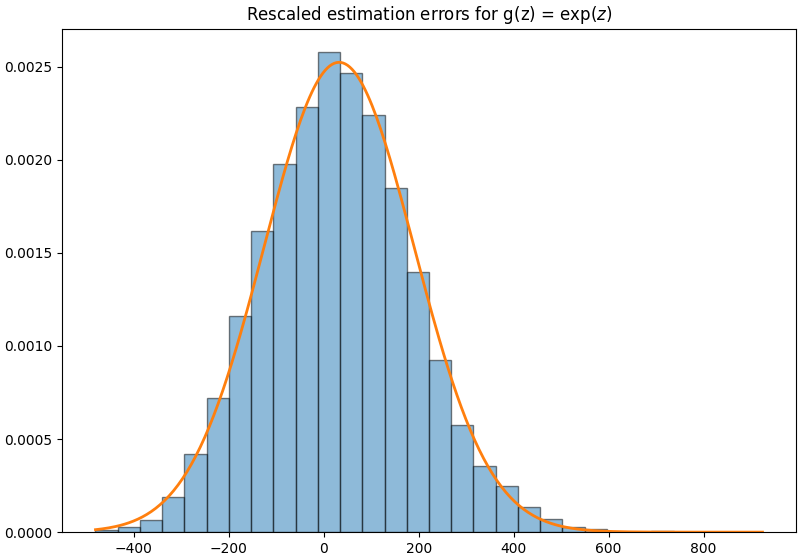}
		\includegraphics[width=0.32\textwidth]{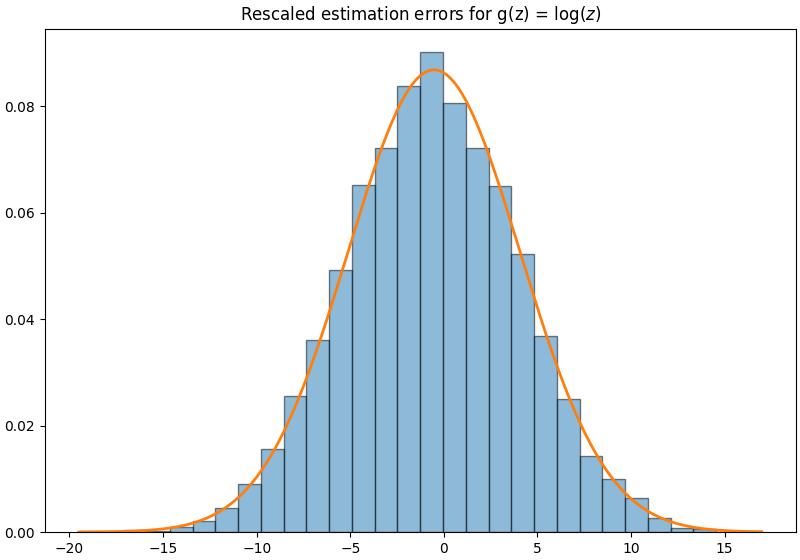}
		
		\vspace{-0.2cm}
		\caption{Histograms of the rescaled estimation errors $n(\hat{L}_{n,\gamma}(g)-L_n(g))$ with $g(z)=z^3$ (left), $g(z) = \exp(z)$ (middle) and $g(z)=\log(z)$ (right) for $10{,}000$ realizations of example (Ex.1) and dimensions $d=20$ (left), $d=50$ (middle) and $d=100$ (right). The density of the Gaussian distribution with mean and covariance given in Corollary~\ref{Cor_GaussPLSSCLT} is overlaid in orange.}\label{Fig_CLT}
		\vspace{-0.3cm}
	\end{figure}

\begin{figure}[H]
		\centering
		\includegraphics[width=0.32\textwidth]{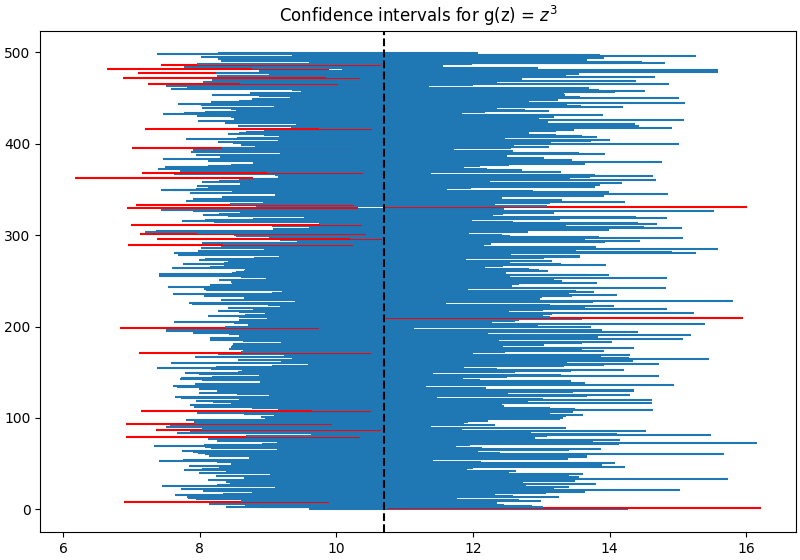}
		\includegraphics[width=0.32\textwidth]{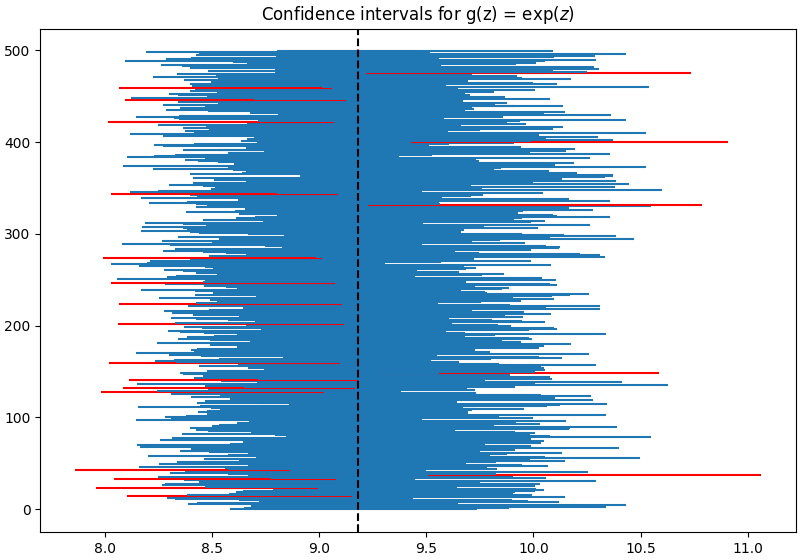}
		\includegraphics[width=0.32\textwidth]{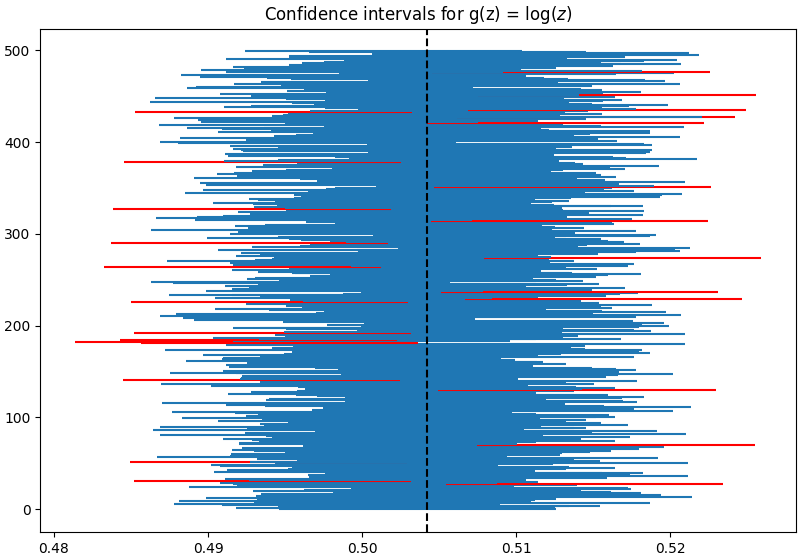}
		\vspace{-0.2cm}
		\caption{Confidence intervals $\hat{I}_{n,\alpha}(g)$ with $\alpha=0.05$ for 500 realizations of example (Ex.1) with dimension $d$ and function $g$ as in Figure~\ref{Fig_CLT}. Accurate intervals, where $L_n(g) \in \hat{I}_{n,\alpha}(g)$, are marked blue, while inaccurate intervals are marked red. The percentages of accurate intervals are $95.0\%$ (left),  $95.8\%$ (middle) and  $94.8\%$ (right).}\label{Fig_Intervals}
		\vspace{-0.2cm}
	\end{figure}

	The central limit theorem for the rescaled estimation error $n(\hat{L}_{n,\gamma}(g)-L_n(g))$, shown in Corollary~\ref{Cor_GaussPLSSCLT} under the additional assumption that $\bm{X}_n$ has i.i.d. Gaussian entries, is observed in Figure~\ref{Fig_CLT} to generalize well to the finite-sample setting.
	\\[0.5em]
	If one further assumes that $g$ maps $\R$ into itself, the same CLT may also be used to construct confidence intervals. To a given confidence level $1-\alpha \in (0,1)$ define
	\begin{align*}
		& \hat{I}_{n,\alpha}(g) \coloneq \Big[ Q_{\mathcal{N}(\mu_{n}(g),\sigma^2_n(g))}\Big(\frac{\alpha}{2}\Big), Q_{\mathcal{N}(\mu_{n}(g),\sigma^2_n(g))}\Big(1-\frac{\alpha}{2}\Big) \Big] \ ,
	\end{align*}
	where $Q_{\mathcal{N}(\mu_{n}(g),\sigma^2_n(g))}(x)$ is the $x$-quantile of the Gaussian distribution with mean
	\begin{align*}
		& \mu_n(g) \coloneq \hat{L}_{n,\gamma}(g) + \frac{1}{2\pi \bm{i} n} \int_{\gamma} g(z) \hat{\bm{e}}_n(z) - g(\ol{z}) \hat{\bm{e}}_n(\ol{z}) \, dz
	\end{align*}
	and variance
	\begin{align*}
		& \sigma^2_n(g) \coloneq \frac{1}{\pi^2n^2} \Re\bigg( \int_{\gamma} \int_{\gamma}g(z_1)g(\ol{z}_2)\hat{\bm{c}}_n(z_1,\ol{z}_2) - g(z_1)g(z_2)\hat{\bm{c}}_n(z_1,z_2) \, dz_1 \, dz_2 \bigg) \ .
	\end{align*}
	Here, $\hat{\bm{e}}_n$ and $\hat{\bm{c}}_n$ are defined through $\hat{\Phi}_n(z) = (1-c_nz\hat{s}_n(z)-c_n)z$ as the estimated counterparts to $\bm{e}$ and $\bm{c}$ from Corollary~\ref{Cor_GaussPLSSCLT}, i.e.
	\begin{align*}
		& \hat{\bm{e}}_n(z) = -\mathbbm{1}_{\beta=1} \frac{\hat{\Phi}_n''(z)}{2c_n\hat{\Phi}_n'(z)} \ \text{ and } \ \hat{\bm{c}}_n(z_1,z_2) = \frac{1}{\beta c_n^2} \Big(\frac{1}{(z_1-z_2)^2} - \frac{\hat{\Phi}_n'(z_1) \hat{\Phi}_n'(z_2)}{(\hat{\Phi}_n(z_1)-\hat{\Phi}_n(z_2))^2}\Big) \ .
	\end{align*}
	In the Gaussian case and for sufficiently smooth $g$, Corollary~\ref{Cor_GaussPLSSCLT} together with Lemma~\ref{Lemma_ConsistencyBasic} may be used to show that the coverage probability of $\hat{I}_{n,\alpha}(g)$ will converge to $1-\alpha$ for $d,n \to \infty$. Figure~\ref{Fig_Intervals} illustrates the efficacy of these intervals for a finite sample setting.
	
	\subsection{Conclusion of the numerical applications}
	The preceding section illustrates the main theoretical and practical features of the proposed
	eigen-inference method. Its primary theoretical advantage is the error rate
	$\mathcal{O}(n^{\varepsilon-1})$ for the estimation of population Stieltjes transforms, and
	therefore of population linear spectral statistics, in a general non-parametric setting. The
	numerical experiment in Figure~\ref{Fig_Stieltjes_comparison} supports this rate prediction.
	\\[0.5em]
	For the estimation of population spectral moments
	\[
	\int_\R \lambda^k \, dH_n(\lambda)
	=
	\frac{1}{d}\tr(\Sigma_n^k) \ ,
	\]
	Figure~\ref{Fig_Monomial_comparison} shows that the proposed method performs almost
	identically to the specialized estimator of Bai--Chen--Yao. Thus, for this class of examples,
	the general contour-based approach does not appear to lose accuracy compared with a method
	designed specifically for moment estimation.
	\\[0.5em]
	The generalized linear spectral statistic estimator introduced in
	Section~\ref{Section_GLSS} was also seen to be competitive with the Ledoit--Wolf method over
	a wide range of aspect ratios $\frac{d}{n}$ (see Figure~\ref{Fig_LargestEigVec}). In this setting,
	the requirement that admissible curves exist plays an additional diagnostic role. In particular,
	the numerical results suggest that this requirement provides an effective criterion for detecting
	whether the largest sample eigenvalue is sufficiently separated from the bulk for the quantity
	$|u_1^\top v_1|^2$ to be estimable as a generalized linear spectral statistic.
	\\[0.5em]
	For Gaussian data, Figure~\ref{Fig_CLT} illustrates the central limit theorem for the estimation
	error from Corollary~\ref{Cor_GaussPLSSCLT}. Figure~\ref{Fig_Intervals} then demonstrates how the
	CLT may be used to construct asymptotically valid confidence intervals for population linear
	spectral statistics.
	\\[0.5em]
	In summary, the proposed method combines generality with competitive finite-sample
	performance. For the estimation of population spectral statistics
	\[
	L_n(g)
	=
	\int_\R g(\lambda)\,dH_n(\lambda)
	\]
	under general non-parametric population spectral distributions, it is the first method with a
	proven error rate of $\mathcal{O}(n^{\varepsilon-1})$. At the same time, the simulations indicate
	that it remains competitive with specialized eigen-inference methods in settings where those
	methods apply. Moreover, the method naturally provides diagnostics for deciding when separated
	spectral components can be analyzed independently of the bulk, and, in the Gaussian case, it
	also yields a practical route to confidence intervals for population linear spectral statistics.

	\newpage
	\appendix
	\section{Proofs of introductory and auxiliary lemmas}\label{Section_LemmaProofs}
	
	\subsection{Proof of Lemma~\ref{Lemma_SpaceTransform}}\label{Proof_Lemma_SpaceTransform}
	For ease of notation, write $\bD = \bD_{H,c}(\infty)$, $\Phi = \Phi_{H,c}$ as well as $\varphi = \varphi_{\nu,c}$. Note that $\varphi(\C^+) \subset \C^+$ is a consequence of the construction of $\varphi$.
	\begin{itemize}
		\item[i)] \textit{Proving bijectivity of $\Phi|_{\varphi(\C^+)}$ and $\Phi \circ \varphi = \Id$:}\\
		Observe that the Marchenko--Pastur equation yields
		\begin{align}\label{Eq_SubOrd_Calc1}
			& \cs_{\nu}(\tz) \overset{\text{(\ref{Eq_MP_Equation})}}{=} \int_\R \frac{1}{\lambda(1-c\tz\cs_{\nu}(\tz)-c) - \tz} \, dH(\lambda) = \frac{-1}{\tz} \int_\R \frac{1}{1 + \lambda \cs_{\ul{\nu}}(\tz)} \, dH(\lambda) \ ,
		\end{align}
		where the second equality follows from $\cs_{\ul{\nu}}(\tz) = c\cs_{\nu}(\tz) - \frac{1-c}{\tz}$.
		It follows that
		\begin{align*}
			& \Phi(\varphi(\tz)) = \Big(1-c\frac{-1}{\cs_{\ul{\nu}}(\tz)}\cs_H\Big(\frac{-1}{\cs_{\ul{\nu}}(\tz)}\Big)-c\Big)\frac{-1}{\cs_{\ul{\nu}}(\tz)}\\
			& \overset{\text{(\ref{Eq_DefStieltjes})}}{=} \bigg( 1 + c\int_\R \frac{1}{\lambda\cs_{\ul{\nu}}(\tz)+1} \, dH(\lambda) - c \bigg) \frac{-1}{\cs_{\ul{\nu}}(\tz)}\\
			& \overset{\text{(\ref{Eq_SubOrd_Calc1})}}{=} \big( 1 - c\tz\cs_{\nu}(\tz) - c \big) \frac{-1}{\cs_{\ul{\nu}}(\tz)} = -\tz\cs_{\ul{\nu}}(\tz) \frac{-1}{\cs_{\ul{\nu}}(\tz)} = \tz \ ,
		\end{align*}
		which proves $\Phi \circ \varphi = \Id$, so both $\varphi : \C^+ \rightarrow \varphi(\C^+)$ and $\Phi : \varphi(\C^+) \rightarrow \C^+$ are bi-holomorphic and $\varphi(\C^+) \subset \bD$.
		
		\item[ii)] \textit{Proving that $\Phi$ is injective on $\bD$ and $\bD = \varphi(\C^+)$:}\\
		For every $z \in \bD$, the calculation
		\begin{align*}
			& \overbrace{\Im(\Phi(z))}^{>0} = \Im\big( (1 - cz\cs_H(z) - c) z \big) \overset{\text{(\ref{Eq_DefStieltjes})}}{=} \Im\bigg( \bigg(1 - c\int_\R \frac{\lambda}{\lambda-z} dH(\lambda)\bigg) z \bigg) \nonumber\\
			& = \Im(z) - c \int_\R \lambda \frac{\Im(z (\lambda - \ol{z}))}{|\lambda-z|^2} \, dH(\lambda) = \Im(z) - c \Im(z) \int_\R \frac{\lambda^2}{|\lambda-z|^2} \, dH(\lambda) \nonumber\\
			& 
			= \underbrace{\Im(z)}_{>0} \bigg( 1 - c \int_\R \frac{\lambda^2}{|\lambda-z|^2} \, dH(\lambda) \bigg)
		\end{align*}
		proves
		\begin{align}\label{Eq_bDChar}
			& 1 > c \int_\R \frac{\lambda^2}{|\lambda-z|^2} \, dH(\lambda) \ .
		\end{align}
		For any $z_1\neq z_2 \in \bD$, the Cauchy-Schwarz bound
		\begin{align}\label{Eq_CauchySchwarz}
			& \Big| c\int_\R \frac{\lambda^2}{(\lambda-z_1)(\lambda-z_2)} \, dH(\lambda) \Big| \nonumber\\
			& \leq \Big( \underbrace{c\int_\R \frac{\lambda^2}{|\lambda-z_1|^2} \, dH(\lambda)}_{\overset{\text{(\ref{Eq_bDChar})}}{<} 1} \Big)^{\frac{1}{2}} \Big( \underbrace{c\int_\R \frac{\lambda^2}{|\lambda-z_2|^2} \, dH(\lambda)}_{\overset{\text{(\ref{Eq_bDChar})}}{<} 1} \Big)^{\frac{1}{2}} < 1
		\end{align}
		may be used to see
		\begin{align*}
			& \Phi(z_1) - \Phi(z_2) = (z_1-z_2) - c \int_\R \frac{\lambda z_1}{\lambda-z_1} - \frac{\lambda z_2}{\lambda-z_2} \, dH(\lambda)\\
			& = (z_1-z_2) - c\int_\R \frac{\lambda^2 (z_1-z_2)}{(\lambda-z_1)(\lambda-z_2)} \, dH(\lambda)\\
			& = (z_1-z_2) \Big( \underbrace{1 - \overbrace{c\int_\R \frac{\lambda^2}{(\lambda-z_1)(\lambda-z_2)} \, dH(\lambda)}^{|\cdot| \overset{\text{(\ref{Eq_CauchySchwarz})}}{<} 1}}_{\neq 0} \Big) \ , 
		\end{align*}
		thus proving that $\Phi$ is injective on $\bD$.
		In the previous part, we had already seen $\varphi(\C^+) \subset \bD$, so it remains to show $\bD \subset \varphi(\C^+)$. If a $z_0 \in \bD \setminus \varphi(\C^+)$ were to exist, the known injectivity of $\varphi$ guarantees that $\Phi(z_0) = \Phi(z_1)$ for $z_1 = \varphi(\Phi(z_0)) \in \varphi(\C^+)$. Injectivity of $\Phi$ would then yield $\varphi(\C^+)\not\ni z_0 = z_1 \in \varphi(\C^+)$, which is not possible. This concludes the proof of $\bD = \varphi(\C^+)$.
		
		\item[iii)] \textit{Checking (\ref{Eq_PhiStieltjesProerty}) and (\ref{Eq_VarphiStieltjesProerty}):}\\
		First, (\ref{Eq_PhiStieltjesProerty}) will be shown by verifying the equivalent formula
		\begin{align}\label{Eq_PhiStieltjesProerty2}
			& \forall z \in \bD : \ \cs_{\nu}\big( (1-cz\cs_H(z)-c)z \big) = \frac{\cs_H(z)}{1-cz\cs_H(z)-c}
		\end{align}
		which is straightforward to check with Lemma~\ref{Lemma_MP}. For $\tilde{z} \coloneq (1 - cz\cs_H(z) - c)z \in \C^+$, the calculation
		\begin{align*}
			& \Im\Big( c \frac{\cs_H(z)}{1 - cz\cs_H(z) - c} + \frac{c-1}{\tz} \Big) = \Im\Big( -\frac{1-cz \cs_H(z) - c}{(1 - cz\cs_H(z) - c)z} \Big) = \Im\Big( \frac{-1}{z} \Big) > 0
		\end{align*}
		shows $\frac{\cs_H(z)}{1 - cz\cs_H(z) - c} \in \text{(\ref{Eq_Def_Qzc})}$.  Further observe
		\begin{align*}
			& \int_\R \frac{1}{\lambda(1-c\tilde{z}\frac{\cs_H(z)}{1 - cz\cs_H(z) - c}-c) - \tilde{z}} \, dH(\lambda)\\
			& = \int_\R \frac{1}{\lambda(1-cz\cs_H(z)-c) - (1 - cz\cs_H(z) - c)z} \, dH(\lambda)\\
			& = \frac{1}{1 - cz\cs_H(z) - c} \int_\R \frac{1}{\lambda -z} \, dH(\lambda) = \frac{\cs_H(z)}{1 - cz\cs_H(z) - c} \ ,
		\end{align*}
		which shows that $\frac{\cs_H(z)}{1 - cz\cs_H(z) - c}$ also satisfies the defining property of $\cs_{\nu}(\tz)$ from Lemma~\ref{Lemma_MP}. The equality (\ref{Eq_PhiStieltjesProerty2}) and thus (\ref{Eq_PhiStieltjesProerty}) follows. 
		Property (\ref{Eq_VarphiStieltjesProerty}) is then an immediate consequence of the fact that $\varphi$ is the inverse function to $\Phi : \bD \rightarrow \C^+$. \qed
	\end{itemize}

	\subsection{Proof of Lemma~\ref{Lemma_MPI}}\label{Proof_Lemma_MPI}
	For ease of notation, write $\bD = \bD_{H,c}(\infty)$, $\Phi = \Phi_{H,c}$ as well as $\varphi = \varphi_{\nu,c}$.
	\\[0.5em]
	In part (iii) of the proof of Lemma~\ref{Lemma_PopStilEst_Uniqueness} it is without any additional assumptions shown that
	\begin{align}\label{Eq_RevMPE_Equivalence2}
		& \forall z, s \in \C \text{ with } \ (1 - czs - c)z \notin \supp(\nu) \ \text{ and } \ 0 \neq \cs_{\ul{\nu}}((1 - czs - c)z) : \nonumber\\
		& zs + 1 = \int_\R \frac{\lambda}{\lambda - (1 - czs - c)z} \, d\nu(\lambda) \ \ \Leftrightarrow \ \ z = \varphi_{\nu,c}\big( (1 - czs - c)z \big) \ .
	\end{align}
	For any $z \in \bD_{H,c}(\infty)$, one may set $s=\cs_{H}(z)$, since $z \in \bD_{H,c}(\infty)$ implies $(1 - cz\cs_{H}(z) - c)z \in \C^+$ and $\cs_{\ul{\nu}}((1 - czs - c)z) \in \C^+$ holds by positivity of Stieltjes transforms.
	\begin{itemize}
		\item[a)]
		By the injectivity of $\varphi : \C^+ \rightarrow \bD$ shown in Lemma~\ref{Lemma_SpaceTransform}, there can be precisely one pre-image $\varphi^{-1}(z) \in \C^+$ to $z \in \bD$, which by $\varphi\circ\Phi=\Id$ must be $\Phi(z) = (1-cz\cs_H(z)-c)z$. Combining this with (\ref{Eq_RevMPE_Equivalence2}) proves (a).
		
		\item[b)]
		By (\ref{Eq_RevMPE_Equivalence2}), any solution $s \in \C^+$ to (\ref{Eq_ReverseMPEquation_re}) would have to satisfy $z = \varphi\big( (1 - czs - c)z \big)$. The additional assumption $\Im\big( (1-czs-c)z \big) > 0$ thus guarantees $z \in \varphi(\C^+)$. It was however shown in Lemma ~\ref{Lemma_SpaceTransform} that $\bD = \varphi(\C^+)$, and (\ref{Eq_LemmaMPI_Assumptions1}) implies $z \notin \bD$. Consequently, there can be no such solutions, proving (b). \qed
	\end{itemize}

	\subsection{Proof of Lemma~\ref{Lemma_BasicStieltjesConvergence}}\label{Proof_Lemma_BasicStieltjesConvergence}
	\begin{itemize}
		\item[i)] \textit{Proof of (a)}:\\
		For any compact set $J \subset \C\setminus[0,\sigma^2]$ it must hold that
		\begin{align}\label{Eq_BasicStieltjesConvergence_Step1_1}
			& \tau \coloneq \min\limits_{z \in J} \dist(z,[0,\sigma^2]) > 0 \ ,
		\end{align}
		which means the functions
		\begin{align*}
			& \Big\{ f_{z} : [0,\sigma^2] \rightarrow \C \ ; \ \lambda \mapsto \frac{1}{\lambda - z} \ \Big| z \in J \Big\}
		\end{align*}
		are all continuous and uniformly bounded by $\frac{1}{\tau}$. The assumptions~\ref{EI_ItemAssumption_PopConv} and~\ref{EI_ItemAssumption_sigmaBound} then give the pointwise convergence
		\begin{align}\label{Eq_BasicStieltjesConvergence_Step1_2}
			& \forall z \in J : \ \cs_{H_n}(z) = \int_{[0,\sigma^2]} f_{z} \, dH_n \xrightarrow{n \to \infty} \int_{[0,\sigma^2]} f_{z} \, dH_\infty = \cs_{H_\infty}(z) \ .
		\end{align}
		Analogously to (\ref{Eq_StandardBounds_c2}) and (\ref{Eq_StandardBounds_c3}) of Lemma~\ref{Lemma_StandardBounds}, one may with (\ref{Eq_BasicStieltjesConvergence_Step1_1}) show
		\begin{align*}
			& \forall z \in J : \ |\cs_{H_n}(z)| \leq \frac{1}{\tau}\\
			& \forall z_1,z_2 \in J : \ \big| \cs_{H_n}(z_1) - \cs_{H_n}(z_2) \big| \leq \frac{|z_1-z_2|}{\tau^2} \ ,
		\end{align*}
		so the sequence $(\cs_{H_n})_{n \in \N}$ is uniformly bounded and equicontinuous on the compact set $J$. Arzel\`a--Ascoli gives the existence of a sub-sequence $(\cs_{H_{n_k}})_{k \in \N}$ uniformly convergent on $S$. The fact that the limit can only be the pointwise limit $\cs_{H_\infty}$, by standard topological arguments implies that the original sequence must have already converged uniformly to $\cs_{H_\infty}$ on $S$.
		
		\item[ii)] \textit{Proof of (b)}:\\
		It will first be shown that $\supp(\nu_\infty) \subset [0,\sigma^2(1+\sqrt{c_\infty})^2]$. To this end, one may also assume that~\ref{EI_ItemAssumption_X_Structure} of Assumption~\ref{Assumption_EigInf_Main} holds, where $\bm{X}_n$ is assumed to have i.i.d. entries, and also that~\ref{EI_ItemAssumption_MomentBound} of Assumption~\ref{Assumption_EigInf_Main} holds for $p=4$. This is permissible, since such a meta-model is sure to exist and the additional assumptions do not influence $\nu_n$ or $\nu_\infty$. Let $B_n = U_n \diag(\sigma_{1,n},\dots,\sigma_{d,n}) V_n$ be the singular value decomposition of $B_n$. By assumption (\ref{Eq_Assumption4_sigmaBound}) one has $\sigma_{1,n}^2 \leq \sigma^2$. Since the difference
		\begin{align*}
			& \frac{1}{n} \bm{X}_n^* V_n^* \diag(\sigma^2,\dots,\sigma^2) V_n \bm{X}_n - \frac{1}{n} \bm{Y}_n^* \bm{Y}_n\\
			& = \frac{1}{n} \bm{X}_n^* V_n^* \diag(\sigma^2-\sigma_{1,n}^2,\dots,\sigma^2-\sigma_{d,n}^2) V_n \bm{X}_n
		\end{align*}
		is positive semi-definite, it must hold that
		\begin{align*}
			& \lambda_{\max}(\bm{S}_n) = \lambda_{\max}\Big( \frac{1}{n} \bm{Y}_n^* \bm{Y}_n \Big)\\
			& \leq \lambda_{\max}\Big( \frac{1}{n} \bm{X}_n^* V_n^* \diag(\sigma^2,\dots,\sigma^2) V_n \bm{X}_n \Big) = \sigma^2 \lambda_{\max}\Big( \frac{1}{n} \bm{X}_n \bm{X}_n^* \Big) \ .
		\end{align*}
		It was in \cite{YinBaiKrishnaiah} for the above meta-model shown that
		\begin{align*}
			& 1 = \bP\Big( \lambda_{\max}\Big( \frac{1}{n} \bm{X}_n \bm{X}_n^* \Big) \xrightarrow{n \to \infty} (1+\sqrt{c_\infty})^2 \Big) \ .
		\end{align*}
		Since the standard Marchenko--Pastur law
		\begin{align*}
			& 1 = \bP\big( \hat{\nu}_n \xRightarrow{n \to \infty} \nu_\infty \big)
		\end{align*}
		also holds under these conditions (see \cite{MP_Bai}), it for the deterministic limiting distribution $\nu_\infty$ must hold that
		\begin{align*}
			& \supp(\nu_\infty) \subset [0,\sigma^2(1+\sqrt{c_\infty})^2] \ .
		\end{align*}
		This is extended to hold for $(\nu_n,c_n)$ instead of $\nu_\infty,c_\infty$ by a simple meta-model argument, since $(H_n,c_n,\nu_n)$ are themselves valid values for $(H_\infty,c_\infty,\nu_\infty)$.
		
		\item[iii)] \textit{Proof of (c)}:\\
		Property (\ref{Eq_PhiStieltjesProerty}) for every $n \in \N \cup \{\infty\}$ and all $z \in \bD_{H_n,c_n}(\infty)$ implies
		\begin{align}\label{Eq_InversionUse0_new}
			& \cs_{\nu_n}\big( (1 - c_nz\cs_{H_n}(z) - c_n)z \big) = \frac{\cs_{H_n}(z)}{1 - c_nz\cs_{H_n}(z) - c_n}
		\end{align}
		and the map
		\begin{align*}
			& \Phi_n = \Phi_{H_n,c_n} : \bD_{H_n,c_n}(\infty) \rightarrow \C^+ \ \ ; \ \ z \mapsto (1 - c_nz\cs_{H_n}(z) - c_n)z
		\end{align*}
		is a bi-holomorphism by Lemma~\ref{Lemma_SpaceTransform}. Thus, for every $\tz \in \C\setminus[0,\sigma^2(1+\sqrt{c_\infty})^2]$ and $n \in \N \cup \{\infty\}$ exists a $z_n \in D^+_{H_n,c_n}(0_+,\infty)$ such that
		\begin{align}\label{Eq_tz_Injectivity_new}
			& \tz = \Phi_n(z_n) \ .
		\end{align}
		Observe
		\begin{align*}
			& |\cs_{\nu_n}(\tz) - \cs_{\nu_\infty}(\tz)| = |\cs_{\nu_n}(\Phi_\infty(z_\infty)) - \cs_{\nu_\infty}(\Phi_\infty(z_\infty))|\\
			& \leq |\cs_{\nu_n}(\Phi_\infty(z_\infty)) - \cs_{\nu_n}(\Phi_n(z_\infty))| + \underbrace{|\cs_{\nu_n}(\Phi_n(z_\infty)) - \cs_{\nu_\infty}(\Phi_\infty(z_\infty))|}_{\to 0 \text{ , by (\ref{Eq_InversionUse0_new}) and (a)}}\\
			& \leq \int_\R \Big| \frac{1}{\lambda - \Phi_\infty(z_\infty)} - \frac{1}{\lambda - \Phi_n(z_\infty)} \Big| \, d\nu_n(\lambda) + o(1)\\
			& = \int_\R \frac{|\Phi_\infty(z_\infty) - \Phi_n(z_\infty)|}{|\lambda - \Phi_\infty(z_\infty)| \, |\lambda - \Phi_n(z_\infty)|} \, d\nu_n(\lambda) + o(1)\\
			& \leq \frac{|\Phi_\infty(z_\infty) - \Phi_n(z_\infty)|}{\Im(\Phi_\infty(z_\infty)) \, \Im(\Phi_n(z_\infty))} + o(1) \ .
		\end{align*}
		Since (a) implies $\Phi_n(z_\infty) \to \Phi_\infty(z_\infty) = \tz \in \C\setminus[0,\sigma^2(1+\sqrt{c_\infty})^2]$, it follows that $\cs_{\nu_n} \xrightarrow{n \to \infty} \cs_{\nu_\infty}$ pointwise on $\C\setminus[0,\sigma^2(1+\sqrt{c_\infty})^2]$. With Arzel\`a--Ascoli, one can analogously to the proof of (a) get uniform convergence on compact sets.
		
		\item[iv)] \textit{Proof of (d)}:\\
		It is well known, and shown for example in Theorem 5.8 of \cite{ProofMethodsRMT}, that pointwise convergence of Stieltjes transforms on $\C^+$ implies weak convergence of the underlying probability measures. Thus, (d) follows directly from result (c).
		
		\item[v)] \textit{Proof of (e)}:\\
		For every compact subset $\tJ \subset \C\setminus[0,\sigma^2(1+\sqrt{c_\infty})^2]$, it must hold that
		\begin{align}\label{Eq_BasicStieltjesConvergence_Step5_1}
			& \tau \coloneq \min\limits_{\tz \in \tJ} \dist(\tz,[0,\sigma^2(1+\sqrt{c_\infty})^2]) > 0 \ ,
		\end{align}
		Use (\ref{Eq_BasicStieltjesConvergence_Step5_1}), (b) and (d) for the applicability of (b) of Lemma~\ref{Lemma_StandardBounds} to see
		\begin{align}\label{Eq_BasicStieltjesConvergence_Step5_3}
			& \forall \tz \in \tJ : \ |\cs_{\ul{\nu}_\infty}(\tz)| \geq \underbrace{\int_\R \frac{\tau/2}{|\lambda-\tz|^2} \, d\ul{\nu}_\infty(\lambda)}_{>0} \eqcolon  2\delta \ .
		\end{align}
		By (c) and definition (\ref{Eq_Def_ulNu}), it also holds that
		\begin{align}\label{Eq_BasicStieltjesConvergence_Step5_4}
			& \sup\limits_{\tz \in \tJ} \big| \cs_{\ul{\nu}_n}(\tz) - \cs_{\ul{\nu}_\infty}(\tz) \big| \xrightarrow{n \to \infty} 0 \ ,
		\end{align}
		which with (\ref{Eq_BasicStieltjesConvergence_Step5_3}) implies
		\begin{align}\label{Eq_BasicStieltjesConvergence_Step5_5}
			& \exists N > 0 \, \forall n \geq N \, \forall \tz \in \tJ : \ |\cs_{\ul{\nu}_n}(\tz)| \geq \delta \ .
		\end{align}
		Combining (\ref{Eq_BasicStieltjesConvergence_Step5_3})--(\ref{Eq_BasicStieltjesConvergence_Step5_5}) with definition (\ref{Eq_Def_varphi}) yields
		\begin{align*}
			& \sup\limits_{\tz \in \tJ} \big| \varphi_{\nu_n,c_n}(\tz) - \varphi_{\nu_\infty,c_\infty}(\tz) \big| \xrightarrow{n \to \infty} 0 \ .
		\end{align*}
		\qed
	\end{itemize}

	\subsection{Proof of Lemma~\ref{Lemma_StandardBounds}}\label{Proof_Lemma_StandardBounds}
	\begin{itemize}
		\item[a)]
		This holds immediately by definition (\ref{Eq_Def_ulNu}) and likewise for $\hat{\nu}_n$ and $\hat{\ul{\nu}}_n$ instead of $\nu_n$ and $\ul{\nu}_n$.
		
		\item[b)] 
		The implication
		\begin{align}\label{Eq_StandardBounds_b1_copy}
			& \big( \nu_n \xRightarrow{n \to \infty} \nu_\infty \big) \Rightarrow \big( \supp(\nu_\infty) \subset [0,K] \big)
		\end{align}
		is a direct consequence of the Portmanteau theorem and (a) turns this into (\ref{Eq_StandardBounds_b1}). By (\ref{Eq_StandardBounds_b2}), each $\tz$ must fall into one of the three cases
		\begin{align*}
			& \Re(\tz) < -\tau/2 \ \ ; \ \ |\Im(\tz)| > \tau/2 \ \ ; \ \ \Re(\tz) > K+\tau/2 \ .
		\end{align*}
		In the first and third case, one has
		\begin{align*}
			& \big| \Re\big( \cs_{\nu}(\tz) \big) \big| \overset{\text{(\ref{Eq_DefStieltjes})}}{=} \Big| \int_\R \frac{\lambda-\Re(\tz)}{|\lambda-\tz|^2} \, d\nu(\lambda) \Big| \geq \int_\R \frac{\tau/2}{|\lambda-\tz|^2} \, d\nu(\lambda) \ ,
		\end{align*}
		while in the second case one has
		\begin{align*}
			& \big| \Im\big( \cs_{\nu}(\tz) \big) \big| \overset{\text{(\ref{Eq_DefStieltjes})}}{=} \Big| \int_\R \frac{\Im(\tz)}{|\lambda-\tz|^2} \, d\nu(\lambda) \Big| \geq \int_\R \frac{\tau/2}{|\lambda-\tz|^2} \, d\nu(\lambda) \ ,
		\end{align*}
		thus proving (\ref{Eq_StandardBounds_b3}).
		\\[0.5em]
		It remains to show that the map $\varphi_{\nu_n,c_n}(\tz)$ is holomorphic in $\C \setminus [0,K]$. Without loss of generality, assume $\nu_n \neq \delta_0$, since in this case the map has the trivial form $\varphi_{\nu_n,c_n}(\tz) = \tz$. It is clear that $\cs_{\ul{\nu}_n}$ (by (a)) will be holomorphic on $\C \setminus [0,K]$ and by the already shown property (\ref{Eq_StandardBounds_b3}) for $\tau = \dist(\tz,[0,K])>0$, it must hold that $0 \neq\cs_{\ul{\nu}_n}(\tz)$ for all $\tz \in \C\setminus[0,K]$. Thus, $\varphi_{\nu_n,c_n}(\tz)$ is holomorphic in $\C \setminus [0,K]$. The argument is analogous for $\varphi_{\nu_\infty,c_\infty}$, if the right-hand side of (\ref{Eq_StandardBounds_b1}) holds.
		\\[0.5em]
		For $\hat{\nu}_n$ and $\hat{\ul{\nu}}_n$ instead of $\nu_n$ and $\ul{\nu}_n$, all arguments are analogous.
		
		\item[c)]
		These bounds are shown with the notation $\nu$ for either $\nu_n$ or $\ul{\nu}_n$ by the simple calculation
		\begin{align*}
			& |\cs_{\nu}(\tz)| \overset{\text{(\ref{Eq_DefStieltjes})}}{\leq} \int_\R \frac{1}{|\lambda-\tz|} \, d\nu(\lambda) \overset{\text{(\ref{Eq_StandardBounds_c1})}}{\leq} \frac{1}{\tau} \ ,
		\end{align*}
		where (a) goes into the application of (\ref{Eq_StandardBounds_c1}), if $\nu = \ul{\nu}_n$. The bound (\ref{Eq_StandardBounds_c3}) is shown analogously by \begin{align*}
			& |\cs_{\nu}(\tz_1) - \cs_{\nu}(\tz_2)| \overset{\text{(\ref{Eq_DefStieltjes})}}{\leq} \int_\R \frac{|\tz_1-\tz_2|}{|\lambda-\tz_1| \, |\lambda-\tz_2|} \, d\nu(\lambda) \overset{\text{(\ref{Eq_StandardBounds_c1})}}{\leq} \frac{|\tz_1-\tz_2|}{\tau^2} \ .
		\end{align*}
		For $\hat{\nu}_n$ and $\hat{\ul{\nu}}_n$ instead of $\nu_n$ and $\ul{\nu}_n$, all arguments are again analogous.
		
		\item[d)]
		The map $\Phi_{H_n,c_n} : \bD_{H_n,c_n}(\infty) \rightarrow \C^+$ is by Lemma~\ref{Lemma_SpaceTransform} the inverse of $\varphi_{\nu_n,c_n} : \C^+ \rightarrow \bD_{H_n,c_n}(\infty)$. Since $\Phi_{H_n,c_n}(z) = (1-c_nz\cs_{H_n}(z)-c_n)z$ is defined (and holomorphic) on all of $\C\setminus\supp(H_n)$, basic symmetry allows for
		\begin{align}\label{Eq_BoundLemma_dProof_1}
			& \Phi_{H_n,c_n}|_{\bD_{H_n,c_n}(\infty) \cup \ol{\bD_{H_n,c_n}(\infty)}} = \varphi|_{\C\setminus\R}^{-1} \ .
		\end{align}
		The definition of $\bD_{H_n,c_n}(\infty)$ yields
		\begin{align}\label{Eq_BoundLemma_dProof_2}
			& \forall z \in \bD_{H_n,c_n}(\infty) \cup \ol{\bD_{H_n,c_n}(\infty)} : \ \operatorname{sign}(\Im(z)) = \operatorname{sign}(\Im(\Phi_{H_n,c_n}(z))) \ ,
		\end{align}
		so the calculation
		\begin{align}\label{Eq_imPhiCalc_copy}
			& \Im(\Phi_{H_n,c_n}(z)) = \Im\big( (1 - c_nz\cs_{H_n}(z) - c_n) z \big) \overset{\text{(\ref{Eq_DefStieltjes})}}{=} \Im\bigg( \bigg(1 - c_n\int_\R \frac{\lambda}{\lambda-z} dH_n(\lambda)\bigg) z \bigg) \nonumber\\
			& = \Im(z) - c_n \int_\R \lambda \frac{\Im(z (\lambda - \ol{z}))}{|\lambda-z|^2} \, dH_n(\lambda) = \Im(z) - c_n \Im(z) \int_\R \frac{\lambda^2}{|\lambda-z|^2} \, dH_n(\lambda) \nonumber\\
			& = \Im(z) \bigg( 1 - c_n \int_\R \frac{\lambda^2}{|\lambda-z|^2} \, dH_n(\lambda) \bigg)
		\end{align}
		shows that
		\begin{align}\label{Eq_BoundLemma_dProof_3}
			& 0 < 1 - c_n \int_\R \frac{\lambda^2}{|\lambda-z|^2} \, dH_n(\lambda)
		\end{align}
		holds for all $z \in \bD_{H_n,c_n}(\infty) \cup \ol{\bD_{H_n,c_n}(\infty)}$. The equality
		\begin{align}\label{Eq_PhiDiffCalc_copy}
			& \Phi_{H_n,c_n}(z_1) - \Phi_{H_n,c_n}(z_2) \overset{\text{(\ref{Eq_DefStieltjes})}}{=} (z_1-z_2) - c_n \int_\R \frac{\lambda z_1}{\lambda - z_1} - \frac{\lambda z_2}{\lambda - z_2} \, dH_n(\lambda) \nonumber\\
			& = (z_1-z_2) - c_n \int_\R \frac{\lambda z_1(\lambda - z_2) - \lambda z_2 (\lambda - z_1)}{(\lambda - z_1)(\lambda - z_2)} \, dH_n(\lambda) \nonumber\\
			& = (z_1-z_2) \bigg( 1 - c_n \int_\R \frac{\lambda^2}{(\lambda - z_1)(\lambda - z_2)} \, dH_n(\lambda) \bigg) \ .
		\end{align}
		thus together with the bound
		\begin{align*}
			& \Big| c_n \int_\R \frac{\lambda^2}{(\lambda - z_1)(\lambda - z_2)} \, dH_n(\lambda) \Big| \nonumber\\
			& \leq \Big( c_n\int_\R \frac{\lambda^2}{|\lambda-z_1|^2} \, dH_n(\lambda) \Big)^{\frac{1}{2}} \Big( c_n\int_\R \frac{\lambda^2}{|\lambda-z_2|^2} \, dH_n(\lambda) \Big)^{\frac{1}{2}} \overset{\text{(\ref{Eq_BoundLemma_dProof_3})}}{<} 1
		\end{align*}
		thus shows that $\Phi_{H_n,c_n}$ is Lipschitz-continuous on $\bD_{H_n,c_n}(\infty) \cup \ol{\bD_{H_n,c_n}(\infty)}$ with Lipschitz-constant $2$. By (\ref{Eq_BoundLemma_dProof_1}) it follows that (\ref{Eq_StandardBounds_d1}) holds for all $\tz_1 \neq \tz_2 \in \C\setminus\R$. 
		The proof for $\nu_\infty$ is completely analogous, while the proof does not hold for $\hat{\nu}_n$, since in this case Lemma~\ref{Lemma_SpaceTransform} is not applicable. \qed
	\end{itemize}
	
	\subsection{Proof of (b) from Lemma \ref{Lemma_PopStilEst_Uniqueness} (continuation)}\label{Proof_Lemma_PopStilEst_Uniqueness_Continuation}
	The following proof continuation ties directly into the notation of the proof of part (b) from Lemma \ref{Lemma_PopStilEst_Uniqueness} in Subsection \ref{Proof_Lemma_PopStilEst_Uniqueness}.
	\\[0.5em]
	It remains to show the existence of a $\kappa(\hat{\nu}_n,c_n)>0$ such that $\C^+ \setminus B^{\C^+}_{\kappa(\hat{\nu}_n,c_n)}(0)$ lies in $U$. This will be done by proving the existence of a $\tilde{\kappa}>0$ such that $\tilde{U}$ contains $\C^+\setminus B_{\tilde{\kappa}}^{\C^+}(0)$ and then showing that $\varphi_{\hat{\nu}_n,c_n}$ is for large $|\tz|$ sufficiently close to the identity that $U=\varphi_{\hat{\nu}_n,c_n}(\tilde{U})$ must contain $\C^+\setminus B_{\kappa}^{\C^+}(0)$ for some $\kappa>0$.
	\\[0.5em]
	By plugging in definitions (\ref{Eq_Def_varphi}) and (\ref{Eq_DefStieltjes}) into (\ref{Eq_UniquenessLemma_Step2_1}), one sees that a $\tz \in \C^+$ lies in $\tilde{U}$ if and only if
	\begin{align}\label{Eq_UniquenessLemma_Step2_2}
		& c_n \int_\R \frac{\lambda}{|\lambda - \tz|^2} \, d\hat{\nu}_n(\lambda) < \Big| \int_\R \frac{1}{\lambda-\tz} \, d\hat{\ul{\nu}}_n(\lambda) \Big| \ .
	\end{align}
	As $\hat{\nu}_n$ is assumed to have compact support, there exists a $K=K(\hat{\nu}_n)>0$ such that $\supp(\hat{\nu}_n) \subset [0,K]$. Result (b) of Lemma~\ref{Lemma_StandardBounds} with $\tau_{\text{Lemma~\ref{Lemma_StandardBounds}}} = \dist(\tz,[0,K]) \eqcolon  d_{\tz}$ then shows
	\begin{align}\label{Eq_UniquenessLemma_Step2_3}
		& \Big|\int_\R \frac{1}{\lambda-\tz} \, d\hat{\ul{\nu}}_n(\lambda)\Big| > \int_\R \frac{d_{\tz}/2}{|\lambda-\tz|^2} \, d\hat{\ul{\nu}}_n(\lambda)
	\end{align}
	for all $\tz \in \C^+$.
	The bound
	\begin{align}\label{Eq_UniquenessLemma_Step2_4}
		& \int_\R \frac{1}{|\lambda-\tz|^2} \, d\hat{\ul{\nu}}_n(\lambda) \geq \frac{1}{K} \int_\R \frac{\lambda}{|\lambda-\tz|^2} \, d\hat{\ul{\nu}}_n(\lambda) \overset{\text{(\ref{Eq_Def_ulNu})}}{=} \frac{c_n}{K} \int_\R \frac{\lambda}{|\lambda-\tz|^2} \, d\hat{\ul{\nu}}_n(\lambda)
	\end{align}
	may be combined with (\ref{Eq_UniquenessLemma_Step2_3}) to see that $d_{\tz} \geq 2K$ is sufficient for (\ref{Eq_UniquenessLemma_Step2_2}). Consequently, it by $\supp(\hat{\nu}_n) \subset [0,K]$ and basic triangle inequality must holds that
	\begin{align}\label{Eq_UniquenessLemma_Step2_5}
		& \C^+\setminus B_{\tilde{\kappa}}^{\C^+}(0) \subset \tilde{U}
	\end{align}
	for any $\tilde{\kappa} \geq 3K$.
	\\[0.5em]
	The convergence
	\begin{align*}
		& \frac{\Im(\varphi_{\hat{\nu}_n,c_n}(\tz))}{\Im(\tz)} \overset{\text{(\ref{Eq_Def_varphi})}}{=} -\frac{\Im(\frac{1}{\cs_{\hat{\ul{\nu}}_n}(\tz)})}{\Im(\tz)} = \frac{\Im(\cs_{\hat{\ul{\nu}}_n}(\tz))}{\Im(\tz) |\cs_{\hat{\ul{\nu}}_n}(\tz)|^2} \overset{\text{(\ref{Eq_DefStieltjes})}}{=} \frac{\int_\R \frac{1}{|\lambda-\tz|^2} \, d\hat{\ul{\nu}}_n(\lambda)}{|\int_\R \frac{1}{\lambda-\tz} \, d\hat{\ul{\nu}}_n(\lambda)|^2}\\
		& \leq \frac{\frac{1}{d_{\tz}^2}}{(\frac{1}{K+|\tz|})^2} \leq \frac{(K+|\tz|)^2}{d_{\tz}^2} \xrightarrow{|\tz| \to \infty} 1
	\end{align*}
	proves that
	\begin{align}\label{Eq_UniquenessLemma_Step2_6}
		& \frac{\Im(\varphi_{\hat{\nu}_n,c_n}(\tz))}{\Im(\tz)} \text{ remains bounded on $\C^+\setminus B_{\tilde{\kappa}}^{\C^+}(0)$ for large enough $\tilde{\kappa}$,}
	\end{align}
	while the convergence
	\begin{align*}
		& \Big|1 - \frac{\tz}{\varphi_{\hat{\nu}_n,c_n}(\tz)}\Big| \overset{\text{(\ref{Eq_Def_varphi})}}{=} |1 + \tz\cs_{\hat{\ul{\nu}}_n}(\tz)| \overset{\text{(\ref{Eq_DefStieltjes})}}{=} \Big|\int_\R \frac{\lambda}{\tz-\lambda} \, d\hat{\ul{\nu}}_n(\lambda)\Big| \leq \frac{K}{d_{\tz}} \xrightarrow{|\tz| \to \infty} 0
	\end{align*}
	in turn shows that $\varphi_{\hat{\nu}_n,c_n}(\tz)$ approaches the identity for large $\kappa$, which together with (\ref{Eq_UniquenessLemma_Step2_5}), (\ref{Eq_UniquenessLemma_Step2_6}) and the continuity of $\varphi_{\hat{\nu}_n,c_n}$ ensures the existence of some $\kappa>0$ such that $\C^+\setminus B_{\kappa}^{\C^+}(0) \subset U$. \qed

	\section{Measurability}\label{Section_Measurability}
	This section proves measurability of the events described in Lemma~\ref{Lemma_ConsistencyBasic} and Theorem~\ref{Thm_Consistency}. The first step will be to show that the map, which maps the sample covariance matrix $\bm{S}_n$ to the population Stieltjes transform estimator $\hat{s}_n(z)$ is Borel-measurable for fixed $z \in \C^+$.
	
	\begin{lemma}[Borel-measurability of $\hat{s}_n$]\label{Lemma_BorelMeasurability}\
		\\
		Let $\textrm{PD}_d$ denote the cone of positive semi-definite matrices in $\C^{d \times d}$ equipped with the topology inherited from $\C^{d \times d}$.
		For fixed $c>0$ and $z \in \C^+$, define the map $\xi_{c,z} : \textrm{PD}_d \rightarrow \C$ such that $s = s(A) = \xi_{c,z}(A)$ is the (unique by parts (i) and (ii) in the proof of Lemma~\ref{Lemma_PopStilEst_Uniqueness}) solution to 
		\begin{align}\label{Eq_ReverseMPEquation_Empirical_copy}
			& zs + 1 = \int_\R \frac{\lambda}{\lambda - (1 - czs - c)z} \, d\underbrace{\Big( \frac{1}{d} \sum\limits_{j=1}^d \delta_{\lambda_j(A)} \Big)}_{\eqcolon  \nu(A)}(\lambda) \ ,
		\end{align}
		which also satisfies
		\begin{align}\label{Eq_StilEstimator_DefiningConditions_copy}
			& \Im\big( (1-czs-c)z \big) > 0 \ \ \text{ and } \ \ 
			\Big| \frac{c z \Im(zs)}{\Im((1-czs-c)z)} \Big| < 1 \ .
		\end{align}
		If such a solution $s(A)$ does not exist, define $\xi_{c,z}(A)=0$. Then:
		\begin{itemize}
			\item[a)] For any $A \in\textrm{PD}_d$, a solution $s(A) \in\C$ to (\ref{Eq_ReverseMPEquation_Empirical_copy}), which satisfies (\ref{Eq_StilEstimator_DefiningConditions_copy}), is never zero.
			
			\item[b)] The map $\xi_{c,z} : \textrm{PD}_d \rightarrow \C$ is Borel-measurable.
		\end{itemize}
	\end{lemma}
	As a consequence of (a), the set $\xi_{c,z}^{-1}(\{0\})$ is precisely the set of all positive semi-definite $(d \times d)$-matrices $A$, for which a solution $s(A)$ does not exist.
	\begin{proof}\
		Proof of (a):
		Suppose for some choice of $A \in \textrm{PD}_d$, $c>0$ and $z \in \C^+$ it holds that $s=0$ is a solution to (\ref{Eq_ReverseMPEquation_Empirical_copy}), which satisfies (\ref{Eq_StilEstimator_DefiningConditions_copy}). The fixed-point equation (\ref{Eq_ReverseMPEquation_Empirical_copy}) becomes
		\begin{align*}
			& 1 = \int_\R \frac{\lambda}{\lambda-(1-c)z} \, d\nu(A)(\lambda) \ .
		\end{align*}
		Since $z \in \C^+$, this is only possible for $c=1$. However, for $(s,c)=(0,1)$ the left hand side of (\ref{Eq_StilEstimator_DefiningConditions_copy}) can not hold, which brings the assumption of $s=0$ to contradiction.
		\\[0.5em]
		The proof of (b) is carried out in the following three steps:
		\begin{itemize}
			\item[i)] \textit{The graph of $\xi_{c,z}$ is a Borel set}:\\
			For fixed $c>0$ and $z \in \C^+$, let $U_{c,z}$ denote the open set of all $s \in \C$, which satisfy (\ref{Eq_StilEstimator_DefiningConditions_copy}). By (\ref{Eq_RevMPE_Equivalence1}) the equivalence
			\begin{align*}
				& \Big(s \in U_{c,z} \ \text{ and $s$ solves (\ref{Eq_ReverseMPEquation_Empirical_copy})}\Big) \ \Leftrightarrow \ \Big(s \in U_{c,z} \ \text{ and } \ z = \varphi_{\nu(A),c}\big( \underbrace{(1-cz\hat{s}-c)z}_{\in \C^+ \text{ by } s \in U_{c,z}} \big)\Big)
			\end{align*}
			holds, which yields that the set
			\begin{align*}
				& G_{c,z} \coloneq \big\{(A,s(A)) \ \big| \ A \in \textrm{PD}_d , \, s(A) \in \C^+ \text{ solves (\ref{Eq_ReverseMPEquation_Empirical_copy}) and satisfies (\ref{Eq_StilEstimator_DefiningConditions_copy})}\big\}
			\end{align*}
			is equal to the pre-image to $\Psi_{c,z}^{-1}(\{z\})$ for the map
			\begin{align*}
				& \Psi_{c,z} : \textrm{PD}_d \times U_{c,z} \rightarrow \C \ \ ; \ \ (A,s) \mapsto \varphi_{\nu(A),c}\big( (1-czs-c)z \big) \ ,
			\end{align*}
			with $\varphi_{\nu(A),c}$ as defined in (\ref{Eq_Def_varphi}).
			This map $\Psi_{c,z}$ is continuous in $s$ and continuous in $A$, since the eigenvalues of Hermitian matrices depend continuously on the entries. Consequently, the set $G_{c,z} = \Psi_{c,z}^{-1}(\{z\})$ is a Borel set.
			
			\item[ii)] \textit{The set $\xi_{c,z}^{-1}(\{0\})$ is a Borel set}:\\
			In (i) it was shown that $G_{c,z} \subset \textrm{PD}_d \times U_{c,z}$ is a Borel set, so the projection onto the first component:
			\begin{align*}
				& \{A \in \textrm{PD}_d \mid \exists s(A) \in \C^+ \text{ solution to (\ref{Eq_ReverseMPEquation_Empirical_copy}) and satisfies (\ref{Eq_StilEstimator_DefiningConditions_copy})} \} \overset{\text{(a)}}{=} \xi_{c,z}^{-1}(\C\setminus\{0\})
			\end{align*}
			must also be a Borel set. Since $\xi_{c,z}^{-1}(\C\setminus\{0\})$ is a Borel set, so must be $\xi_{c,z}^{-1}(\{0\}) = \textrm{PD}_d \setminus \xi_{c,z}^{-1}(\C\setminus\{0\})$.
			
			\item[iii)] \textit{The map $\xi_{c,z}$ is Borel-measurable}:\\
			It remains to show that $\xi_{c,z}^{-1}(B) \subset \textrm{PD}_d$ is a Borel set for every Borel set $B \subset \C \setminus\{0\}$. In part (i) is was shown that the graph $G_{c,z}$ is a Borel set, so the intersection
			\begin{align*}
				& G_{c,z} \cap (\textrm{PD}_d \times B) = \big\{(A,s(A)) \ \big| \ A \in \textrm{PD}_d , \, s(A) \in B \text{ solves (\ref{Eq_ReverseMPEquation_Empirical_copy}) and satisfies (\ref{Eq_StilEstimator_DefiningConditions_copy})}\big\}
			\end{align*}
			is also a Borel set. The projection onto the first component:
			\begin{align*}
				& \{A \in \textrm{PD}_d \mid \exists s(A) \in B \text{ solution to (\ref{Eq_ReverseMPEquation_Empirical_copy}) and satisfies (\ref{Eq_StilEstimator_DefiningConditions_copy})} \} = \xi_{c,z}^{-1}(B)
			\end{align*}
			is then also a Borel set. \qedhere
		\end{itemize}
	\end{proof}
	
	Since it holds by construction that $\hat{s}_n(z) = \xi_{c_n,z}(\bm{S}_n)$, when the left-hand side exists, the above Lemma proves that $\hat{s}_n^{(0)}(z)$ as defined in Theorem~\ref{Thm_CLT_Inversion} is measurable with respect to $\sigma(\bm{X}_n)$.

	\subsection{Proof of measurability for Lemma~\ref{Lemma_ConsistencyBasic}}\label{Proof_ConsistencyBasic_Measurability}
	With the notation $(\Omega,\mathcal{A},\bP)$ for the measure space on which the random matrices $\bm{X}_n$ are defined, it is shown here that the expression in (\ref{Eq_ConsistencyBasic_Stronger}) indeed describes an event from $\mathcal{A}$.
	Let $(J_k)_{k \in \N}$ be a deterministic increasing sequence of compact sets $J_k \subset \C$ such that $\bigcup\limits_{k=1}^\infty J_k = \bD_{H_\infty,c_\infty}(1) \cup (\R\setminus I)$.
	The event from (\ref{Eq_ConsistencyBasic_Stronger}) is equivalent to
	\begin{align}\label{Eq_BasicConsistency_Measurability_1}
		& \bigcap\limits_{k=1}^\infty \bigcap\limits_{m=1}^\infty \bigcup\limits_{N=1}^\infty \bigcap\limits_{n=N}^\infty \Big(\Big\{ \omega \in \Omega \ \Big| \ \hat{s}_{n} \text{ exists on $J_k\setminus\R$} \Big\} \nonumber\\
		& \hspace{3cm} \cap \Big\{ \omega \in \Omega \ \Big| \ \sup\limits_{z \in J_k\setminus\R} \big|\hat{s}_n(z) - \cs_{H_\infty}(z)\big| \leq \frac{1}{m} \Big\} \Big) \ .
	\end{align}
	For each $k \in \N$, let $Q_k \subset J_k\setminus\R$ be countable and dense.
	In terms of $\xi_{c_n,z}$ from Lemma~\ref{Lemma_BorelMeasurability}, one may write
	\begin{align*}
		& \Big\{ \omega \in \Omega \ \Big| \ \hat{s}_{n} \text{ exists on $J_k\setminus\R$} \Big\} = \Big\{ \omega \in \Omega \ \Big| \ J_k\setminus\R \subset \overbrace{\{z \in \C^+ \mid \hat{s}_n(z) \text{ exists}\}}^{\text{open by (b) of Lemma~\ref{Lemma_PopStilEst_Uniqueness}}} \Big\}\\
		& = \bigcap\limits_{q \in Q_k} \underbrace{\Big\{ \omega \in \Omega \ \Big| \ \xi_{c_n,q}(\bm{S}_n(\omega)) \neq 0 \Big\}}_{\in \mathcal{A} \text{ by Lemma~\ref{Lemma_BorelMeasurability}}} \ ,
	\end{align*}
	so the upper set in (\ref{Eq_BasicConsistency_Measurability_1}) is also in $\mathcal{A}$. Similarly, one may use the continuity from (c) of Lemma~\ref{Lemma_PopStilEst_Uniqueness} to see
	\begin{align*}
		& \Big\{ \omega \in \Omega \ \Big| \ \hat{s}_{n} \text{ exists on $J_k\setminus\R$} \Big\} \cap \Big\{ \omega \in \Omega \ \Big| \ \sup\limits_{z \in J_k\setminus\R} \big|\hat{s}_n(z) - \cs_{H_\infty}(z)\big| \leq \frac{1}{m} \Big\}\\
		& = \bigcap\limits_{q \in Q_k} \underbrace{\Big\{ \omega \in \Omega \ \Big| \ \xi_{c_n,q}(\bm{S}_n(\omega)) \neq 0 \Big\}}_{\in \mathcal{A} \text{ by Lemma~\ref{Lemma_BorelMeasurability}}} \cap \bigcap\limits_{q \in Q_k} \underbrace{\Big\{ \omega \in \Omega \ \Big| \ |\xi_{c_n,q}(\bm{S}_n(\omega)) - \cs_{H_\infty}(z)| \leq \frac{1}{m} \Big\}}_{\in \mathcal{A} \text{ by Lemma~\ref{Lemma_BorelMeasurability}}} \ ,
	\end{align*}
	which proves measurability of (\ref{Eq_BasicConsistency_Measurability_1}) thus of the event in (\ref{Eq_ConsistencyBasic_Stronger}). \qed

	\subsection{Proof of measurability for Theorem~\ref{Thm_Consistency}}\label{Proof_Consistency_Measurability}
	With the notation $(\Omega,\mathcal{A},\bP)$ for the measure space on which the random matrices $\bm{X}_n$ are defined, it is shown here that the expression in (\ref{Eq_Consistency_Result}) indeed describes an event from $\mathcal{A}$.
	\\[0.5em]
	For each $\omega \in \Omega$, the set $\hat{D}(\tau,\kappa,n;\omega) \subset \C^+$ is open by (b) and (c) of Lemma~\ref{Lemma_PopStilEst_Uniqueness}. Since both $\hat{s}_n(\bullet;\omega)$ and $\cs_{H_\infty}$ are continuous on $\hat{D}(\tau,\kappa,n;\omega)$, the expression in (\ref{Eq_Consistency_Result}) is equivalent to
	\begin{align}\label{Eq_ConistencyMeasurability_Equiv1}
		& \Big\{ \omega \in \Omega \ \Big| \ \forall z \in \hat{D}(\tau,\kappa,n;\omega) \cap (\Q \times i\Q_{>0}) : \ \big| \hat{s}_n(z) - \cs_{H_\infty}(z) \big| \leq n^{\varepsilon-1} \Big\} \ .
	\end{align}
	Further, the random (open) set $\hat{D}(\tau,\kappa,n) \subset \C^+$ may by Theorem~\ref{Thm_Consistency} and Lemma~\ref{Lemma_BorelMeasurability} be characterized by
	\begin{align*}
		& \mathbbm{1}_{z \in \hat{D}(\tau,\kappa,n)}\\
		& = \overbrace{\mathbbm{1}_{\tau < |z| < \kappa}}^{\text{deterministic}} \times \overbrace{\mathbbm{1}_{\xi_{c_n,z}(\bm{S}_n) \neq 0}}^{\mathcal{A}\text{-measurable by Lemma~\ref{Lemma_BorelMeasurability}}} \times \overbrace{\mathbbm{1}_{\tau < |(1-c_nz\xi_{c_n,z}(\bm{S}_n)-c_n)z| < \kappa}}^{\mathcal{A}\text{-measurable by Lemma~\ref{Lemma_BorelMeasurability}}}\\
		& \hspace{1cm} \times \mathbbm{1}_{\forall j \leq d : \, \tau < |(1-c_nz\xi_{c_n,z}(\bm{S}_n)-c_n)z - \lambda_j(\bm{S}_n)|} \times \underbrace{\mathbbm{1}_{\tau < \dist((1-c_nz\xi_{c_n,z}(\bm{S}_n)-c_n)z, \supp(\nu_n))}}_{\mathcal{A}\text{-measurable by Lemma~\ref{Lemma_BorelMeasurability}}} \ .
	\end{align*}
	The remaining indicator function is also $\mathcal{A}$-measurable by the fact that
	\begin{align*}
		& \textrm{PD}_d \times \textrm{PD}_d \rightarrow \C \ \ ; \ \ (A,B) \mapsto \min\limits_{j \leq d} |(1-c_nz\xi_{c_n,z}(A)-c_n)z - \lambda_j(B)|
	\end{align*}
	is continuous in $B$ and measurable in $A$ by Lemma~\ref{Lemma_BorelMeasurability}. 
	For each $z \in \C^+$, the random variable
	\begin{align*}
		& \hat{s}_n^{(H_\infty)}(z) \coloneq \mathbbm{1}_{z \in \hat{\bD}(\tau,\kappa,n)} \xi_{c_n,z}(\bm{S}_n) + \mathbbm{1}_{z \notin \hat{\bD}(\tau,\kappa,n)} \cs_{H_\infty}(z) = \begin{cases}
			\hat{s}_n(z) & \text{ if } z \in \hat{\bD}(\tau,\kappa,n)\\
			\cs_{H_\infty}(z) & \text{ else}
		\end{cases} \ ,
	\end{align*}
	is thus $\mathcal{A}$-measurable and one observes that the expression in (\ref{Eq_Consistency_Result}), which was already seen to be equivalent to (\ref{Eq_ConistencyMeasurability_Equiv1}), is equivalent to
	\begin{align*}
		& \bigcap\limits_{z \in \Q \times i\Q_{>0}} \Big\{ \omega \in \Omega \ \Big| \ \big| \hat{s}_n^{(H_\infty)}(z) - \cs_{H_\infty}(z) \big| \leq n^{\varepsilon-1} \Big\} \ ,
	\end{align*}
	which must lie in $\mathcal{A}$ due to the measurability of $\hat{s}_n^{(H_\infty)}(z)$. \qed

	\section{Outer local laws}\label{Section_OuterLaws}
	
	\begin{theorem}[Knowles-Yin: General outer local law]\label{Thm_OuterLaw}\
		\\
		Suppose Assumption~\ref{Assumption_EigInf_Main} holds.
		For a fixed $\tau > 0$ define
		\begin{align*}
			& \bm{D}(\tau) \coloneq \big\{ \tz \in \C^+ \ \big| \ 0 < \Im(\tz) \leq \tau^{-1} , \, |\Re(\tz)| \leq \tau^{-1}, \, \tau \leq |\tz| \big\}\\
			& \bS(\tau,n) \coloneq \big\{ \tz \in \bm{D}(\tau) \ \big| \ \dist(\tz,\supp(\ul{\nu}_n)) \geq \tau \big\} \ .
		\end{align*}
		For any $\tilde{\varepsilon},D,\tau>0$ there exists a constant $C = C(\tilde{\varepsilon},D,\tau)>0$, which additionally depends on $\inf_{n \in \N} c_n$, $\sup_{n \in \N} c_n$, $\sigma^2$ and the constants $(C_p)_{p \in \N}$ (but not on the explicit distributions of the entries of $\bm{X}_n$ or the covariances $\Sigma_n$), such that
		\begin{align}
			& \bP\Big( \forall \tz \in \bS(\tau,n) : \ \big| \cs_{\hat{\ul{\nu}}_n}(\tz) - \cs_{\ul{\nu}_n}(\tz) \big| \leq \frac{n^{\tilde{\varepsilon}}}{n \Im(\tz)} \Big) \geq 1 - \frac{C}{n^D} \label{Eq_OuterLaw_Stieltjes}
		\end{align}
		and
		\begin{align}\label{Eq_OuterLaw_Matrix}
			& \forall g : [0,\infty) \rightarrow \R \ \text{ bounded and measurable} : \nonumber\\
			& \bP\bigg( \forall \tz \in \bS(\tau,n) : \ \Big| \frac{1}{d} \tr\big( g(\Sigma_n) \bm{R}_n(\tz) \big) - \frac{1}{d} \tr\Big(- \frac{1}{\tz}g(\Sigma_n)(\Id_d + \cs_{\ul{\nu}_n}(\tz)\Sigma_n)^{-1} \Big) \Big| \nonumber\\
			& \hspace{6cm} \leq  2\frac{\sigma^2 \|g\|_\infty}{\tau^2} n^{\tilde{\varepsilon}-\frac{1}{2}} \bigg) \geq 1 - d\frac{C}{n^D} \ .
		\end{align}
		hold for all $n \in \N$, where $\bm{R}_n(\tz) = (\bm{S}_n - \tz\Id_d)^{-1}$.
	\end{theorem}
	
	Importantly, the above theorem does not require $B_n = B_n^* = \Sigma^{\frac{1}{2}} > 0$ as assumed in (2.9) of \cite{KnowlesAnisotropicLocalLaws}, since this assumption is only used to simplify the stability analysis near the spectrum, and is removed in Section 11 of \cite{KnowlesAnisotropicLocalLaws}.. The theorem also does not require the regularity assumptions on the eigenvalues of $\Sigma_n$ from Definition 2.7 of \cite{KnowlesAnisotropicLocalLaws}, since in the present theorem the spectral domain $\bS(\tau,n)$ is uniformly separated from $\supp(\ul{\nu}_n)$, so no edge or cusp regularity assumptions are needed.
	
	\subsection{Proof of Theorem~\ref{Thm_OuterLaw}}\label{Proof_Thm_OuterLaw}
	
	\subsubsection{Checking assumptions}
	Before proving the theorem, it is briefly shown that Assumption 2.1 of \cite{KnowlesAnisotropicLocalLaws} follows from Assumption~\ref{Assumption_EigInf_Main}. Note the notational difference that \cite{KnowlesAnisotropicLocalLaws} examines a sample covariance matrix $TXX^*T$, where the normalization factor $\frac{1}{n}$ is already part of $X$, i.e. it is assumed that $X$ has centered independent entries with variance $\frac{1}{n}$. Since this paper denotes the sample covariance matrix $\bm{S}_n$ as $\frac{1}{n} B_n \bm{X}_n \bm{X}_n^* B_n^*$ for a matrix $\bm{X}_n$ with centered independent entries with variance $1$, one may apply results of \cite{KnowlesAnisotropicLocalLaws} with the translation $X = \frac{1}{\sqrt{n}} \bm{X}_n$ and $T = B_n$.
	\begin{itemize}
		\item Checking (2.1) of \cite{KnowlesAnisotropicLocalLaws}:\\
		By~\ref{EI_ItemAssumption_X_Structure} of Assumption~\ref{Assumption_EigInf_Main}, equation (2.1) of \cite{KnowlesAnisotropicLocalLaws} is in our notation equivalent to $d \asymp n$, which holds by~\ref{EI_ItemAssumption_Asymptotics} of Assumption~\ref{Assumption_EigInf_Main}.
		
		\item Checking (2.4) of \cite{KnowlesAnisotropicLocalLaws}:\\
		The equalities $\E[X_{i \, \mu}] = 0$ and $\E[|X_{i \, \mu}|^2]=\frac{1}{N}$ are $\frac{1}{\sqrt{n}}\E[(\bm{X}_n)_{j,k}]=0$ and $\frac{1}{n}\E[|(\bm{X}_n)_{j,k}|^2]=\frac{1}{n}$, when translated into the notation of this paper. They are thus equivalent to (\ref{Eq_Assumption1_EntriesBasic}) from Assumption~\ref{Assumption_EigInf_Main}.
		
		\item Checking (2.5) of \cite{KnowlesAnisotropicLocalLaws}:\\
		This is exactly~\ref{EI_ItemAssumption_MomentBound} of Assumption~\ref{Assumption_EigInf_Main} after translation into the notation of this paper.
		
		\item Checking (2.7) of \cite{KnowlesAnisotropicLocalLaws}:\\
		Equation (2.7) of \cite{KnowlesAnisotropicLocalLaws} assumes the existence of a sufficiently small $\tau'>0$ that $\lambda_{\max}(\Sigma_n) \leq (\tau')^{-1}$ for all $n \in \N$. By choosing $\tau' < \frac{1}{\sigma^2}$, this follows from~\ref{EI_ItemAssumption_sigmaBound} of Assumption~\ref{Assumption_EigInf_Main}.
		
		\item Checking (2.8) of \cite{KnowlesAnisotropicLocalLaws}:\\
		Equation (2.8) of \cite{KnowlesAnisotropicLocalLaws} assumes the existence of a sufficiently small $\tau'>0$ that $H_n([0,\tau']) \leq 1-\tau'$ holds for every $n \in \N$. As all results of \cite{KnowlesAnisotropicLocalLaws} are asymptotic in nature, it suffices to show the existence of an $N_0>0$ such that $H_n([0,\tau']) \leq 1-\tau'$ holds for all $n \geq N_0$, which follows from~\ref{EI_ItemAssumption_PopConv} of Assumption~\ref{Assumption_EigInf_Main} by the following argument. Since $H_\infty\neq 0$ is a probability measure on $[0,\infty)$, there exists a $\tau'>0$ such that $H_\infty([0,\tau']) \leq 1-2\tau'$. The Portmanteau Theorem then guarantees $\limsup\limits_{n \to \infty} H_n([0,\tau']) \leq H_\infty([0,\tau']) \leq 1-2\tau'$ and there thus must exist an $N_0>0$ such that $H_n([0,\tau']) \leq \limsup\limits_{n \to \infty} H_n([0,\tau']) + \tau' \leq 1-\tau'$ holds for all $n \geq N_0$.
	\end{itemize}
	
	\subsubsection{Proof structure}
	The main part of the proof will work under the additional assumption
	\begin{align}\label{Eq_Assumption6_TSig}
		& B_n = B_n^* = \Sigma_n^{\frac{1}{2}} > 0 \ ,
	\end{align}
	which will be removed at the end using arguments described in Section 11 of \cite{KnowlesAnisotropicLocalLaws}. It will be shown in Subsection~\ref{Subsubsection_EntrywiseLocalLaw_to_GeneralizedStil} that (\ref{Eq_OuterLaw_Matrix}) follows from the anisotropic outer local law:
	\begin{align}\label{Eq_AnisotropicOuterLaw}
		& \forall v,w \in \C^{d} , \, \|v\|_2,\|w\|_2 \leq 1 : \nonumber\\
		& \bP\bigg( \exists \tz \in \bS(\tau,n) : \ \Big| \Big\langle v, \Sigma_n^{-1} \big( \tz\Sigma_n^{\frac{1}{2}}\bm{R}_n(\tz)\Sigma_n^{\frac{1}{2}} + \Sigma_n (\Id_d + \cs_{\ul{\nu}_n}(\tz)\Sigma_n)^{-1} \big) \Sigma_n^{-1} w \Big\rangle \Big| \nonumber\\
		& \hspace{4cm} \geq 2n^{\tilde{\varepsilon}} \sqrt{\frac{\Im(\cs_{\ul{\nu}_n}(\tz))}{n \Im(\tz)}} \bigg) \leq \frac{C}{n^D} \ , \hspace{0.3cm} 
	\end{align}
	which will be proved here in parallel with (\ref{Eq_OuterLaw_Stieltjes}).
	As the spectral domain $\bS(\tau,n)$ is a subset of $\bS_1(\tau,n) \cup \bS_2(\tau,n)$ for
	\begin{align*}
		& \bS_1(\tau,n) \coloneq \{\tz \in \bm{D}(\tau) \mid \dist(\Re(\tz),\supp(\ul{\nu}_n)) \geq \tau/2 \}\\
		& \bS_2(\tau,n) \coloneq \{\tz \in \bm{D}(\tau) \mid \Im(\tz) \geq \tau/2 \} \ ,
	\end{align*}
	it suffices to show
	\begin{align}
		& \bP\Big( \exists \tz \in \bS_{1/2}(\tau,n) : \ \big| \cs_{\hat{\ul{\nu}}_n}(\tz) - \cs_{\ul{\nu}_n}(\tz) \big| \geq \frac{n^{\tilde{\varepsilon}}}{n \Im(\tz)} \Big) \leq \frac{C/2}{n^D} \label{Eq_OuterLaw_Stieltjes_separately}
	\end{align}
	as well as
	\begin{align}
		& \bP\bigg( \exists \tz \in \bS_{1/2}(\tau,n) : \ \Big| \Big\langle v, \Sigma_n^{-1} \big( \tz\Sigma_n^{\frac{1}{2}}\bm{R}_n(\tz)\Sigma_n^{\frac{1}{2}} + \Sigma_n (\Id_d + \cs_{\ul{\nu}_n}(\tz)\Sigma_n)^{-1} \big) \Sigma_n^{-1} w \Big\rangle \Big| \nonumber\\
		& \hspace{4cm} \geq 2n^{\tilde{\varepsilon}} \sqrt{\frac{\Im(\cs_{\ul{\nu}_n}(\tz))}{n \Im(\tz)}} \bigg) \leq \frac{C/2}{n^D} \label{Eq_OuterLaw_Matrix_separately}
	\end{align}
	separately for $\bS_{1}(\tau,n)$ and $\bS_{2}(\tau,n)$. The bounds (\ref{Eq_OuterLaw_Stieltjes_separately}) and (\ref{Eq_OuterLaw_Matrix_separately}) for $\bS_{1}(\tau,n)$ and sufficiently large $C(\tilde{\varepsilon},D,\tau)$ follow directly from Theorem 3.16 (i) and Remark 3.17 of \cite{KnowlesAnisotropicLocalLaws}, so it remains to prove (\ref{Eq_OuterLaw_Stieltjes_separately}) and (\ref{Eq_OuterLaw_Matrix_separately}) for $\bS_{2}(\tau,n)$ with Theorems 3.21 and 3.22 of \cite{KnowlesAnisotropicLocalLaws}, which require two new conditions:
	\begin{itemize}
		\item[i)] There exists a $\tau'>0$ such that $\big|1+\cs_{\ul{\nu}_n}(\tz)\lambda_i(\Sigma_n)\big| \geq \tau'$ for all $n \in \N$, $\tz \in \bS_2(\tau,n)$ and $i \leq d$. This bound is written in (3.20) of \cite{KnowlesAnisotropicLocalLaws}.
		
		\item[ii)] The stability of a re-arranged Marchenko--Pastur equation on $\bS_2(\tau,n)$, as described in Definition 5.4 of \cite{KnowlesAnisotropicLocalLaws}.
	\end{itemize}
	Fortunately, (ii) was already proven to hold for $\tz \in \bS_2(\tau,n)$ with no further assumptions in the (first two paragraphs of the) proof of Lemma A.5 of \cite{KnowlesAnisotropicLocalLaws}, where a stronger property (A.6) in \cite{KnowlesAnisotropicLocalLaws} is shown, which by Definition A.2 of \cite{KnowlesAnisotropicLocalLaws} leads to (ii).
	\\[0.5em]
	The proof of Theorem~\ref{Thm_OuterLaw} thus requires three further steps:
	\begin{itemize}
		\item[a)] proving (\ref{Eq_OuterLaw_Stieltjes_separately}) and (\ref{Eq_OuterLaw_Matrix_separately}) for $\bS_{2}(\tau,n)$, which utilizes Theorems 3.21 and 3.22 of \cite{KnowlesAnisotropicLocalLaws} and thus requires condition (i) to be checked
		
		\item[b)] showing that (\ref{Eq_AnisotropicOuterLaw}) is sufficient for (\ref{Eq_OuterLaw_Matrix})
		
		\item[c)] removing condition (\ref{Eq_Assumption6_TSig}) with the same arguments as used in Section 11 of \cite{KnowlesAnisotropicLocalLaws}.
	\end{itemize}
	
	\subsubsection{Proving (\ref{Eq_OuterLaw_Stieltjes_separately}) and (\ref{Eq_OuterLaw_Matrix_separately}) for $\bS_{2}(\tau,n)$}
	By Lemma 4.10 of \cite{KnowlesAnisotropicLocalLaws} there exists a constant $\mathcal{C}>0$ dependent only on $\tau$ and the asymptotic behavior of $c_n$ such that
	\begin{align}\label{Eq_ImSNu_Comparison}
		& \forall \tz \in \C^+ , \, \tau \leq |\tz| \leq \tau^{-1} : \ \mathcal{C}^{-1} \Im(\tz) \leq \Im(\cs_{\ul{\nu}_n}(\tz)) \leq \mathcal{C}
	\end{align}
	and one may further bound
	\begin{align*}
		& |\Re(\cs_{\ul{\nu}_n}(\tz))| = \bigg|\int_\R \frac{\lambda-\Re(\tz)}{|\lambda-\tz|^2} \, d\ul{\nu}_n(\lambda)\bigg| \leq \int_\R \frac{|\lambda-\Re(\tz)|}{\tau^2} \, d\ul{\nu}_n(\lambda) \leq \frac{\mathcal{C}'}{\tau^2}
	\end{align*}
	for some $\mathcal{C}'>0$, using (b) of Lemma~\ref{Lemma_BasicStieltjesConvergence}.
	Choose $\tau'>0$ small enough such that
	\begin{align}\label{Eq_taubarChoice}
		& \tau'\frac{2\mathcal{C}\,\mathcal{C}'}{\tau^3} \leq 1-\tau' \ ,
	\end{align}
	then for all $i \leq d$ with $\lambda_i(\Sigma_n) \geq \frac{2\tau' \mathcal{C}}{\tau}$ one gets
	\begin{align*}
		& \big|1+\cs_{\ul{\nu}_n}(\tz)\lambda_i(\Sigma_n)\big| \geq \underbrace{|\Im(\cs_{\ul{\nu}_n}(\tz))|}_{\geq \mathcal{C}^{-1} \Im(\tz)} \, \lambda_i(\Sigma_n) \geq \mathcal{C}^{-1} \frac{\tau}{2} \, \frac{2\tau' \mathcal{C}}{\tau} = \tau'
	\end{align*}
	and for all $i \leq d$ with $\lambda_i(\Sigma_n) \leq \frac{2\tau' \mathcal{C}}{\tau}$ further
	\begin{align*}
		& \big|1+\cs_{\ul{\nu}_n}(\tz)\lambda_i(\Sigma_n)\big| \geq 1 - \underbrace{|\Re(\cs_{\ul{\nu}_n}(\tz))|}_{\leq \frac{\mathcal{C}'}{\tau^2}} \, \lambda_i(\Sigma_n) \geq 1 - \frac{\mathcal{C}'}{\tau^2} \, \frac{2\tau' \mathcal{C}}{\tau} \overset{\text{(\ref{Eq_taubarChoice})}}{\geq} 1 - (1-\tau') = \tau' \ ,
	\end{align*}
	which proves condition (i). Theorems 3.21 and 3.22 of \cite{KnowlesAnisotropicLocalLaws} are thus applicable to $\bS_2(\tau,n)$ and yield the averaged and anisotropic local laws (see Definition 3.20 of \cite{KnowlesAnisotropicLocalLaws}), which in the notation of this paper implies the existence of a constant $C''=C''(\tilde{\varepsilon},D,\tau)>0$ such that
	\begin{align*}
		& \bP\Big( \exists \tz \in \bS_{2}(\tau,n) : \ \big| \cs_{\hat{\ul{\nu}}_n}(\tz) - \cs_{\ul{\nu}_n}(\tz) \big| \geq \frac{n^{\tilde{\varepsilon}}}{n \Im(\tz)} \Big) \leq \frac{C''}{n^D}
	\end{align*}
	and
	\begin{align*}
		& \bP\bigg( \exists \tz \in \bS_{2}(\tau,n) : \ \Big| \Big\langle v, \Sigma_n^{-1} \big( z\Sigma_n^{\frac{1}{2}}\bm{R}_n(z)\Sigma_n^{\frac{1}{2}} + \Sigma_n (\Id_d + \cs_{\ul{\nu}_n}(z)\Sigma_n)^{-1} \big) \Sigma_n^{-1} w \Big\rangle \Big|\\
		& \hspace{4cm} \geq n^{\tilde{\varepsilon}} \Big( \sqrt{\frac{\Im(\cs_{\ul{\nu}_n}(\tz))}{n \Im(\tz)}} + \frac{1}{n \Im(\tz)} \Big) \bigg) \leq \frac{C''}{n^D}
	\end{align*}
	hold for all $n \in \N$. Choosing $C\geq 2C''$ immediately yields (\ref{Eq_OuterLaw_Stieltjes_separately}) for $\bS_{2}(\tau,n)$ and the calculation
	\begin{align*}
		& \frac{1}{n \Im(\tz)} \overset{\tz \in \bS_2(\tau,n)}{\leq} \frac{1}{n \tau} \leq \frac{\sqrt{\mathcal{C}}}{\sqrt{n} \tau} \sqrt{\frac{\mathcal{C}^{-1} \Im(\tz)}{n \Im(\tz)}} \overset{\text{(\ref{Eq_ImSNu_Comparison})}}{\leq} \frac{\sqrt{\mathcal{C}}}{\sqrt{n} \tau} \sqrt{\frac{\Im(\cs_{\ul{\nu}_n}(\tz))}{n \Im(\tz)}}
	\end{align*}
	shows that (\ref{Eq_OuterLaw_Matrix_separately}) for $\bS_{2}(\tau,n)$ also follows for sufficiently large $n$. As we have shown (\ref{Eq_OuterLaw_Stieltjes_separately}) and (\ref{Eq_OuterLaw_Matrix_separately}) for $\bS_{1}(\tau,n)$ as well as $\bS_{2}(\tau,n)$, the statements (\ref{Eq_OuterLaw_Stieltjes}) and (\ref{Eq_AnisotropicOuterLaw}) are proved.
	
	\subsubsection{Proving (\ref{Eq_OuterLaw_Matrix}) with (\ref{Eq_AnisotropicOuterLaw})}\label{Subsubsection_EntrywiseLocalLaw_to_GeneralizedStil}
	For $v$ and $w$ in (\ref{Eq_AnisotropicOuterLaw}), insert the vectors
	\begin{align*}
		& v_j = \frac{1}{\sigma} \Sigma_n^{\frac{1}{2}} e_j \ \ \text{ and } \ \ w_j = \frac{1}{\sigma \|g\|_\infty} \Sigma_n^{\frac{1}{2}} g(\Sigma_n) e_j \ ,
	\end{align*}
	where $e_j \in \R^d$ denotes the vector $(\delta_{k,j})_{k \leq d}$ and $g : [0,\infty) \rightarrow [0,\infty)$ is a bounded, measurable function. The term $g(\Sigma_n)$ canonically stands for $U_n \diag(g(\lambda_1(\Sigma_n)),\dots,g(\lambda_d(\Sigma_d))) U_n^*$, where $U_n \diag(\lambda_1(\Sigma_n),\dots,\lambda_d(\Sigma_n))U_n^*$ is the spectral decomposition of $\Sigma_n$. Since $\Sigma_n$ and $g(\Sigma_n)$ commute, one may calculate
	\begin{align*}
		& \Big\langle v_j, \Sigma_n^{-1} \big( \tz\Sigma_n^{\frac{1}{2}}\bm{R}_n(\tz)\Sigma_n^{\frac{1}{2}} + \Sigma_n (\Id_d + \cs_{\ul{\nu}_n}(\tz)\Sigma_n)^{-1} \big) \Sigma_n^{-1} w_j \Big\rangle\\
		& = \frac{1}{\sigma^2 \|g\|_\infty} e_j^\top \Sigma_n^{\frac{1}{2}} \Sigma_n^{-1} \big( \tz\Sigma_n^{\frac{1}{2}}\bm{R}_n(\tz)\Sigma_n^{\frac{1}{2}} + \Sigma_n (\Id_d + \cs_{\ul{\nu}_n}(\tz)\Sigma_n)^{-1} \big) \Sigma_n^{-1} \Sigma_n^{\frac{1}{2}} g(\Sigma_n) e_j\\
		& = \frac{1}{\sigma^2 \|g\|_\infty} \Big( \big( \tz\bm{R}_n(\tz)g(\Sigma_n)\big)_{j,j} + \big( (\Id_d + \cs_{\ul{\nu}_n}(\tz)\Sigma_n)^{-1} g(\Sigma_n) \big)_{j,j} \Big)
	\end{align*}
	and thus
	\begin{align*}
		& \frac{1}{d} \sum\limits_{j=1}^d \Big\langle v_j, \Sigma_n^{-1} \big( \tz\Sigma_n^{\frac{1}{2}}\bm{R}_n(\tz)\Sigma_n^{\frac{1}{2}} + \Sigma_n (\Id_d + \cs_{\ul{\nu}_n}(\tz)\Sigma_n)^{-1} \big) \Sigma_n^{-1} w_j \Big\rangle\\
		& = \frac{1}{\sigma^2 \|g\|_\infty} \Big( \frac{1}{d}\tr\big( \tz\bm{R}_n(\tz)g(\Sigma_n)\big) - \frac{1}{d}\tr\big( -(\Id_d + \cs_{\ul{\nu}_n}(\tz)\Sigma_n)^{-1} g(\Sigma_n) \big) \Big) \ .
	\end{align*}
	Since $|\tz| \geq \tau$ for all $\tz \in \bS(\tau,n)$, the property (\ref{Eq_AnisotropicOuterLaw}) directly yields
	\begin{align}\label{Eq_AnisotropicOuterLaw_inserted}
		& \forall g : [0,\infty) \rightarrow \R \ \text{ bounded and measurable} : \nonumber\\
		& \bP\bigg( \exists \tz \in \bS(\tau,n) : \ \Big| \frac{1}{d} \tr\big( g(\Sigma_n) \bm{R}_n(\tz) \big) - \frac{1}{d} \tr\Big(- \frac{1}{\tz}g(\Sigma_n)(\Id_d + \cs_{\ul{\nu}_n}(\tz)\Sigma_n)^{-1} \Big) \Big| \nonumber\\
		& \hspace{4cm} \geq 2n^{\tilde{\varepsilon}} \frac{\sigma^2\|g\|_\infty}{\tau} \sqrt{\frac{\Im(\cs_{\ul{\nu}_n}(\tz))}{n \Im(\tz)}} \bigg) \leq d\frac{C}{n^D} \ .
	\end{align}
	The calculation
	\begin{align*}
		& \Im(\cs_{\ul{\nu}_n}(\tz)) = \int_\R \Im\Big( \frac{1}{\lambda-\tz} \Big) \, d\ul{\nu}_n(\lambda) = \int_\R \frac{\Im(\tz)}{|\lambda-\tz|^2} \, d\ul{\nu}_n(\lambda) \overset{\tz \in \bS(\tau,n)}{\leq} \frac{\Im(\tz)}{\tau^2}
	\end{align*}
	proves $\sqrt{\frac{\Im(\cs_{\ul{\nu}_n}(\tz))}{\Im(\tz)}} \leq \frac{1}{\tau}$, which turns (\ref{Eq_AnisotropicOuterLaw_inserted}) into (\ref{Eq_OuterLaw_Matrix}).
	
	\subsubsection{Removing condition (\ref{Eq_Assumption6_TSig})}
	So far, Theorem~\ref{Thm_OuterLaw} is proved, when the sample covariance matrix has the form
	\begin{align*}
		& \tilde{\bm{S}}_n = \frac{1}{n} \tilde{\Sigma}_n^{\frac{1}{2}} \bm{X}_n \bm{X}_n^* \tilde{\Sigma}_n^{\frac{1}{2}}
	\end{align*}
	and $\tilde{\Sigma}_n$ is positive definite. It remains to extend the result to sample covariance matrices of the form
	\begin{align*}
		& \bm{S}_n = \frac{1}{n} B_n \bm{X}_n \bm{X}_n^* B_n^*
	\end{align*}
	for general $B_n \in \C^{d \times d}$ such that $\Sigma_n = B_nB_n^*$ may be semi-definite.
	\\[0.5em]
	Let $B_n = U_n D_n V_n^*$ be the singular value decomposition of $B_n$, such that $D_n \in \R^{d \times d}$ is diagonal and $U_n,V_n \in \C^{d \times d}$ are unitary. For any $\varepsilon \in (0,1)$ define $\tilde{\Sigma}_n \coloneq V_n (D_n^2+\varepsilon D_n+ \varepsilon^2 \Id_d) V_n^*$ such that $\tilde{\Sigma}_n^{\frac{1}{2}} = V_n (D_n + \varepsilon \Id_d) V_n^*$. Under Assumption~\ref{Assumption_EigInf_Main} for $\bm{S}_n$, the same assumptions also hold for $\tilde{\bm{S}}_n = \frac{1}{n} \tilde{\Sigma}_n^{\frac{1}{2}} \bm{X}_n \bm{X}_n^* \tilde{\Sigma}_n^{\frac{1}{2}}$, where we must define $\tilde{\sigma}^2 = (\sigma + 1)^2$. The proof thus far will (for any $\tilde{\varepsilon},D,\tau$ as in Theorem~\ref{Thm_OuterLaw}) yield the existence of a constant $C = C(\tilde{\varepsilon},D,\tau)>0$ such that
	\begin{align*}
		& \bP\Big( \exists \tz \in \bS(\tau,n) : \ \big| \cs_{\ul{\tilde{\hat{\nu}}}_n}(\tz) - \cs_{\ul{\tilde{\nu}}_n}(\tz) \big| \geq \frac{n^{\tilde{\varepsilon}}}{n \Im(\tz)} \Big) \leq \frac{C}{n^D} \ ,
	\end{align*}
	where $\ul{\tilde{\hat{\nu}}}_n$ is the spectral distribution $\frac{1}{d} \sum\limits_{j=1}^d \delta_{\lambda_j(\frac{1}{n} \bm{X}_n^* \tilde{\Sigma}_n \bm{X}_n)}$ and $\ul{\tilde{\nu}}_n = (1-c_n) \delta_0 + c_n \tilde{\nu}_n$ for $\tilde{\nu}_n$ the probability distribution on $[0,\infty)$ arising from $c_n$ and $\frac{1}{d} \sum\limits_{j=1}^d \delta_{\lambda_j(\tilde{\Sigma}_n)}$ by Lemma~\ref{Lemma_MP}.
	\\[0.5em]
	The constant $C$ is independent of $\varepsilon \in (0,1)$ and we may thus for each $n \in \N$ let $\varepsilon$ go to zero. The $\omega$-wise convergence $\cs_{\ul{\tilde{\hat{\nu}}}_n}(z) \xrightarrow{\varepsilon \searrow 0} \cs_{\hat{\ul{\nu}}_n}(z)$ is clear, as the spectral norm of the difference of the involved matrices goes to zero for $\varepsilon \searrow 0$. The convergence $\cs_{\ul{\tilde{\nu}}_n}(z) \xrightarrow{\varepsilon \searrow 0} \cs_{\ul{\nu}_n}(z)$ follows from the convergence $\frac{1}{d} \sum\limits_{j=1}^d \delta_{\lambda_j(\tilde{\Sigma}_n)} \xRightarrow{\varepsilon \searrow 0} H_n$ analogously to (c) from Lemma~\ref{Lemma_BasicStieltjesConvergence}. By letting $\varepsilon$ go to zero $n$-wise and $\omega$-wise, the bounds (\ref{Eq_OuterLaw_Stieltjes}) and (\ref{Eq_OuterLaw_Matrix}) follow. \qed
	\\[0.5em]

	By integrating along a curve separating $\bS(\tau,n)$ from the supports of $\ul{\nu}_n$ and $\hat{\ul{\nu}}_n$, one may use Cauchy's integral formula to strengthen the first result of the previous theorem. This will lead to the proof of Theorem~\ref{Corollary_OuterLaw}. First, an auxiliary lemma.
	
	\begin{lemma}[Bound of the largest sample eigenvalue]\label{Lemma_LargestEVBound}\
		\\
		Under Assumption~\ref{Assumption_EigInf_Main} for any $\tau,D > 0$ there exists an $N_0(\tau,D) > 0$ such that
		\begin{align*}
			& \bP\Big( \lambda_{\max}(\bm{S}_n) \leq \sigma^2(1+\sqrt{c_n})^2 + \tau \Big) \geq 1 - n^{-D}
		\end{align*}
		holds for all $n \geq N_0(\tau,D)$.
	\end{lemma}
	\begin{proof}\
		\\
		Let $B_n = U_n \diag(\sigma_{1,n},\dots,\sigma_{d,n}) V_n$ be the singular value decomposition of $B_n$. By assumption (\ref{Eq_Assumption4_sigmaBound}) one has $\sigma_{1,n}^2 \leq \sigma^2$. Since for $\bm{Y}_n = B_n \bm{X}_n$ the difference
		\begin{align*}
			& \frac{1}{n} \bm{X}_n^* V_n^* \diag(\sigma^2,\dots,\sigma^2) V_n \bm{X}_n - \frac{1}{n} \bm{Y}_n^* \bm{Y}_n\\
			& = \frac{1}{n} \bm{X}_n^* V_n^* \diag(\sigma^2-\sigma_{1,n}^2,\dots,\sigma^2-\sigma_{d,n}^2) V_n \bm{X}_n
		\end{align*}
		is positive semi-definite, it must hold that
		\begin{align*}
			& \lambda_{\max}(\bm{S}_n) = \lambda_{\max}\Big( \frac{1}{n} \bm{Y}_n^* \bm{Y}_n \Big)\\
			& \leq \lambda_{\max}\Big( \frac{1}{n} \bm{X}_n^* V_n^* \diag(\sigma^2,\dots,\sigma^2) V_n \bm{X}_n \Big) = \sigma^2 \lambda_{\max}\Big( \frac{1}{n} \bm{X}_n \bm{X}_n^* \Big) \ .
		\end{align*}
		By Theorem 2.10 of \cite{BloemendalIsotropicLocalLaws} (with $\alpha = 1$) one for all $\delta,K'>0$ gets the existence of an $N_0(\delta,K') > 0$ such that
		\begin{align}\label{Eq_BloemendalBound}
			& \bP\Big( \Big| \lambda_{\max}\Big( \frac{1}{n} \bm{X}_n \bm{X}_n^* \Big) - \gamma_{1,n} \Big| \leq n^{\delta-\frac{2}{3}} \Big) \geq 1-n^{-K'} \ ,
		\end{align}
		where $\gamma_{1,n}$ is the (deterministic) classical eigenvalue location for the largest eigenvalue (see 2.17 of \cite{BloemendalIsotropicLocalLaws}). This classical eigenvalue location $\gamma_{1,n}$ must approach the right-hand edge of the deterministic equivalent spectrum ($=(1+\sqrt{c}_n)^2$ by (2.4) and (2.5) of \cite{BloemendalIsotropicLocalLaws}) with rate $\mathcal{O}(\frac{1}{n})$ (though for the sake of this Lemma, any rate $\mathcal{O}(n^{-\varepsilon})$ for $\varepsilon>0$ would suffice).
		For $\delta < \frac{2}{3}$ and sufficiently large $n$, one gets $n^{\delta-\frac{2}{3}} \leq \tau$ and the desired bound follows from (\ref{Eq_BloemendalBound}).
	\end{proof}
	
	\subsection{Proof of Theorem~\ref{Corollary_OuterLaw}}\label{Proof_Corollary_OuterLaw}
	Since $\tilde{\varepsilon}>0$ may be arbitrarily small, it suffices to show
	\begin{align}\label{Eq_OuterLaw_StieltjesInfty_copy2}
		& \bP\Big( \forall \tz \in \hat{\bS}(\tau,n) : \ \big| \cs_{\hat{\ul{\nu}}_n}(\tz) - \cs_{\ul{\nu}_n}(\tz) \big| \leq n^{3\tilde{\varepsilon}-1} \Big) \geq 1 - \frac{C}{n^D}
	\end{align}
	for some $C>0$ dependent on $\tilde{\varepsilon}$, $\tau$ and $D$.
	Without loss of generality assume $\tau < 1$ and that $n$ is large enough to satisfy
	\begin{align}
		& \frac{4\sigma^2(1+\sqrt{c_n})^2)+12\tau}{2\tau^2\pi} \leq n^{\tilde{\varepsilon}} \label{Eq_CorSInf_LargeN_Assumption1}\\
		& \frac{2 n^{\tilde{\varepsilon}} (\sigma^2(1+\sqrt{c_n})^2+3\tau)}{\tau^2\pi n} \leq \frac{n^{2\tilde{\varepsilon}}}{n} \label{Eq_CorSInf_LargeN_Assumption2}\\
		& 0 < \log(\tau) - \log\Big( \frac{n^{\tilde{\varepsilon}}}{2n} \Big) \leq n^{\tilde{\varepsilon}} \label{Eq_CorSInf_LargeN_Assumption3}\\
		& 3 \leq n^{\tilde{\varepsilon}} \label{Eq_CorSInf_LargeN_Assumption4} \ .
	\end{align}
	Since $\ul{\nu}_n$ by definition (see (\ref{Eq_Def_ulNu}) and Lemma~\ref{Lemma_MP}) has compact support and $\supp(\hat{\ul{\nu}}_n) \subset \{0\} \cup \{\lambda_j(\bm{S}_n) \mid j \leq d\}$, the open set
	\begin{align*}
		& U_{n,\tau} \coloneq \bigcup\limits_{x \in \supp(\ul{\nu}_n) \cup \supp(\hat{\ul{\nu}}_n) \cup \{0\}} (x-\tau,x+\tau)
	\end{align*}
	must consist of finitely many disjoint intervals $(a_{1,n},b_{1,n}),\dots,(a_{m,n},b_{m,n})$.
	For each $l \leq m$ let $\gamma_{l,n} : (0,1) \rightarrow \C^+$ be the composite curve $\gamma_3 \circ \gamma_2 \circ \gamma_1$ with:
	\begin{itemize}
		\item $\gamma_1$ going straight up from $b_{l,n}$ to $b_{l,n}+\bm{i}\tau$
		
		\item $\gamma_2$ going straight to the left from $b_{l,n}+\bm{i}\tau$ to $a_{l,n}+\bm{i}\tau$
		
		\item $\gamma_3$ going straight down from $a_{l,n}+\bm{i}\tau$ to $a_{l,n}$
	\end{itemize}
	By Lemma~\ref{Lemma_LargestEVBound} and (b) of Lemma~\ref{Lemma_BasicStieltjesConvergence}, there exists a $C'>0$ such that
	\begin{align}\label{Eq_CorSInf_Event2_Bound}
		& \bP\big( U_{n,\tau} \subset (-\tau,\sigma^2(1+\sqrt{c_n})^2)+2\tau \big) \geq 1 - \frac{C'}{n^D}
	\end{align}
	holds for all $n \in \N$. Note that on this event it must hold that
	\begin{align}\label{Eq_CorSInf_Event2_Consequence}
		& m \leq \frac{\sigma^2(1+\sqrt{c_n})^2+3\tau }{2\tau} \ \ \text{ and } \ \ \sum\limits_{l=1}^m (b_{l,n}-a_{l,n}) \leq \sigma^2(1+\sqrt{c_n})^2+3\tau \ .
	\end{align}
	Cauchy's integral formula yields
	\begin{align*}
		& \cs_{\hat{\ul{\nu}}_n}(\tz) \overset{\text{(\ref{Eq_DefStieltjes})}}{=} \int_{\supp(\hat{\ul{\nu}}_n)} \frac{1}{\lambda - \tz} \, d\hat{\ul{\nu}}_n(\lambda) = \sum\limits_{l=1}^m \int_{(a_{l,n},b_{l,n})} \frac{1}{\lambda - \tz} \, d\hat{\ul{\nu}}_n(\lambda)\\
		& = \sum\limits_{l=1}^m \int_{(a_{l,n},b_{l,n})} \frac{1}{2\pi \bm{i}} \bigg(\int_{\gamma_{l,n}} \frac{1}{v - \tz} \frac{1}{\lambda-v} \, dv - \int_{\gamma_{l,n}} \frac{1}{\ol{v} - \tz} \frac{1}{\lambda-\ol{v}} \, dv \bigg) \, d\hat{\ul{\nu}}_n(\lambda)
	\end{align*}
	and the fact that the integrals
	\begin{align*}
		& \int_{(a_{l,n},b_{l,n})} \frac{1}{2\pi \bm{i}} \bigg(\int_{\gamma_{l',n}} \frac{1}{v - \tz} \frac{1}{\lambda-v} \, dv - \int_{\gamma_{l',n}} \frac{1}{\ol{v} - \tz} \frac{1}{\lambda-\ol{v}} \, dv \bigg) \, d\hat{\ul{\nu}}_n(\lambda)
	\end{align*}
	are zero for $l \neq l'$ together with Fubini's theorem allows for
	\begin{align}\label{Eq_CorSInf_sHatNuCauchy}
		\cs_{\hat{\ul{\nu}}_n}(\tz) & = \sum\limits_{l=1}^m \int_{\R} \frac{1}{2\pi \bm{i}} \bigg(\int_{\gamma_{l,n}} \frac{1}{v - \tz} \frac{1}{\lambda-v} \, dv - \int_{\gamma_{l,n}} \frac{1}{\ol{v} - \tz} \frac{1}{\lambda-\ol{v}} \, dv \bigg) \, d\hat{\ul{\nu}}_n(\lambda) \nonumber\\
		& = \frac{1}{2\pi \bm{i}} \sum\limits_{l=1}^m \int_{\gamma_{l,n}} \frac{\cs_{\hat{\ul{\nu}}_n}(v)}{v - \tz} \, dv - \frac{1}{2\pi \bm{i}} \sum\limits_{l=1}^m \int_{\gamma_{l,n}} \frac{\cs_{\hat{\ul{\nu}}_n}(\ol{v})}{\ol{v} - \tz} \, dv \ .
	\end{align}
	One may in complete analogy also show
	\begin{align}\label{Eq_CorSInf_sNuCauchy}
		& \cs_{\ul{\nu}_n}(\tz) = \frac{1}{2\pi \bm{i}} \sum\limits_{l=1}^m \int_{\gamma_{l,n}} \frac{\cs_{\ul{\nu}_n}(v)}{v - \tz} \, dv - \frac{1}{2\pi \bm{i}} \sum\limits_{l=1}^m \int_{\gamma_{l,n}} \frac{\cs_{\ul{\nu}_n}(\ol{v})}{\ol{v} - \tz} \, dv \ .
	\end{align}
	Further, one from $\dist(\gamma_{l,n},\supp(\hat{\ul{\nu}}_n)) \geq \tau$ and $\dist(\gamma_{l,n},\supp(\ul{\nu}_n)) \geq \tau$ for all $v \in \mathrm{image}(\gamma_{l,n})$ gets
	\begin{align*}
		& |\cs_{\hat{\ul{\nu}}_n}(v)| = \bigg| \int_{\supp(\hat{\ul{\nu}}_n)} \frac{1}{\lambda - v} \, d\hat{\ul{\nu}}_n(\lambda) \bigg| \leq \int_{\supp(\hat{\ul{\nu}}_n)} \frac{1}{|\lambda - v|} \, d\hat{\ul{\nu}}_n(\lambda) \leq \frac{1}{\tau}
	\end{align*}
	and analogously $|\cs_{\ul{\nu}_n}(v)| \leq \frac{1}{\tau}$, which yields
	\begin{align}\label{Eq_CorSInf_sDiffBoundBasic}
		& \forall v \in \mathrm{image}(\gamma_{l,n}) : \ |\cs_{\hat{\ul{\nu}}_n}(v) - \cs_{\ul{\nu}_n}(v)| \leq \frac{2}{\tau} \ .
	\end{align}
	As a final preparation, observe that a curve $\gamma_{l,n}$ can have at most distance $\sqrt{2}\tau$ from $\supp(\ul{\nu}_n) \cup \supp(\hat{\ul{\nu}}_n) \cup \{0\}$, which for all $\tz \in \hat{\bS}(\tau,n)$ and $l \leq m$ implies
	\begin{align}\label{Eq_CorSInf_tzGammaDistance}
		& \dist(\tz,\gamma_{l,n}) \geq 2\tau-\sqrt{2}\tau \geq \frac{\tau}{2} \ .
	\end{align}
	For any $\omega$ from the (high-probability) event
	\begin{align}
		& \bigg\{ \omega \in \Omega \ \bigg| \ \forall v \in \bS(\tau,n) : \ \big| \cs_{\hat{\ul{\nu}}_n}(v) - \cs_{\ul{\nu}_n}(v) \big| \leq \frac{n^{\tilde{\varepsilon}}}{n \Im(v)} \bigg\} \label{Eq_CorSInf_Event1}\\
		& \cap \big\{ \omega \in \Omega \ \big| \ U_{n,\tau} \subset (-\tau,\sigma^2(1+\sqrt{c_n})^2)+2\tau \big\} \label{Eq_CorSInf_Event2}
	\end{align}
	one gets
	\begin{align*}
		& \big| \cs_{\hat{\ul{\nu}}_n}(\tz) - \cs_{\ul{\nu}_n}(\tz) \big|\\
		& \overset{\substack{\text{(\ref{Eq_CorSInf_sHatNuCauchy})} \\ \text{(\ref{Eq_CorSInf_sNuCauchy})}}}{\leq} \frac{1}{2\pi} \sum\limits_{l=1}^m \bigg| \int_{\gamma_{l,n}} \frac{\cs_{\hat{\ul{\nu}}_n}(v)-\cs_{\ul{\nu}_n}(v)}{v - \tz} \, dv \bigg| + \frac{1}{2\pi} \sum\limits_{l=1}^m \bigg| \int_{\gamma_{l,n}} \frac{\cs_{\hat{\ul{\nu}}_n}(v)-\cs_{\ul{\nu}_n}(v)}{v - \ol{\tz}} \, dv \bigg|\\
		& \overset{\text{(\ref{Eq_CorSInf_tzGammaDistance})}}{\leq} \frac{2}{\tau\pi} \sum\limits_{l=1}^m \int_{\gamma_{l,n}} \big| \cs_{\hat{\ul{\nu}}_n}(u) - \cs_{\ul{\nu}_n}(u) \big| \, |du|\\
		& = \frac{2}{\tau\pi} \sum\limits_{l=1}^m \int_{0}^{\tau} \underbrace{\big| \cs_{\hat{\ul{\nu}}_n}(b_{l,n}+\bm{i}t\tau) - \cs_{\ul{\nu}_n}(b_{l,n}+\bm{i}t\tau) \big|}_{ \leq \frac{n^{\tilde{\varepsilon}}}{n t\tau} \land \frac{2}{\tau} \text{ by (\ref{Eq_CorSInf_Event1}) and (\ref{Eq_CorSInf_sDiffBoundBasic})}} \, dt\\
		& \hspace{0.5cm} + \frac{2}{\tau\pi} \sum\limits_{l=1}^m \int_{a_{l,n}}^{b_{l,n}} \underbrace{\big| \cs_{\hat{\ul{\nu}}_n}(t+\bm{i}\tau) - \cs_{\ul{\nu}_n}(t+\bm{i}\tau) \big|}_{\leq \frac{n^{\tilde{\varepsilon}}}{n \tau} \text{ by (\ref{Eq_CorSInf_Event1})}} \, dt\\
		& \hspace{0.5cm} + \frac{2}{\tau\pi} \sum\limits_{l=1}^m \int_{0}^{\tau} \underbrace{\big| \cs_{\hat{\ul{\nu}}_n}(a_{l,n}+\bm{i}t\tau) - \cs_{\ul{\nu}_n}(a_{l,n}+\bm{i}t\tau) \big|}_{\leq \frac{n^{\tilde{\varepsilon}}}{n t\tau} \land \frac{2}{\tau} \text{ by (\ref{Eq_CorSInf_Event1}) and (\ref{Eq_CorSInf_sDiffBoundBasic})}} \, dt\\
		& \leq \frac{4m}{\tau\pi} \int_{0}^{\tau} \frac{n^{\tilde{\varepsilon}}}{n t\tau} \land \frac{2}{\tau} \, dt + \frac{2}{\tau\pi} \frac{n^{\tilde{\varepsilon}}}{n \tau} \sum\limits_{l=1}^m (b_{l,n}-a_{l,n})\\
		& \overset{\text{(\ref{Eq_CorSInf_Event2_Consequence})}}{\leq} \underbrace{\frac{4\sigma^2(1+\sqrt{c_n})^2)+12\tau}{2\tau^2\pi}}_{\leq n^{\tilde{\varepsilon}} \text{ by (\ref{Eq_CorSInf_LargeN_Assumption1})}} \bigg( \int_0^{\frac{n^{\tilde{\varepsilon}}}{2n}} 2 \, dt + \int_{\frac{n^{\tilde{\varepsilon}}}{2n}}^\tau \frac{n^{\tilde{\varepsilon}}}{n t} \, dt \bigg) + \underbrace{\frac{2 n^{\tilde{\varepsilon}} (\sigma^2(1+\sqrt{c_n})^2+3\tau)}{\tau^2\pi n}}_{\leq \frac{n^{2\tilde{\varepsilon}}}{n} \text{ by (\ref{Eq_CorSInf_LargeN_Assumption2})}}\\
		& \leq n^{\tilde{\varepsilon}} \bigg( \frac{n^{\tilde{\varepsilon}}}{n} + \frac{n^{\tilde{\varepsilon}}}{n} \underbrace{\big[ \log(t) \big]_{\frac{n^{\tilde{\varepsilon}}}{2n}}^\tau}_{\leq n^{\tilde{\varepsilon}} \text{ by (\ref{Eq_CorSInf_LargeN_Assumption3})}} \bigg) + \frac{n^{2\tilde{\varepsilon}}}{n} \leq 3\frac{n^{2\tilde{\varepsilon}}}{n} \overset{\text{(\ref{Eq_CorSInf_LargeN_Assumption4})}}{\leq} \frac{n^{3\tilde{\varepsilon}}}{n} \ .
	\end{align*}
	By choosing $C>0$ larger than $C'$ from (\ref{Eq_CorSInf_Event2_Bound}) and $C$ from (\ref{Eq_OuterLaw_Stieltjes}) and also large enough that the bounds (\ref{Eq_CorSInf_LargeN_Assumption1})--(\ref{Eq_CorSInf_LargeN_Assumption4}) hold for all $n \geq C^{\frac{1}{D}}$ the wanted bound (\ref{Eq_OuterLaw_StieltjesInfty_copy2}) follows. \qed

	\section{Consistency of the PLSS estimator}
	
	\subsection{Proof of Corollary~\ref{Cor_PLSS_Consistency}}\label{Proof_Cor_PLSS_Consistency}
	\begin{itemize}
		\item[i)] \textit{Cauchy's integral formula}:\\
		Let $\tilde{\gamma}_n$ be the closed curve going through $\mathrm{image}(\gamma_n) \cup \ol{\mathrm{image}(\gamma_n)}$ with a counter-clockwise orientation.
		Cauchy's integral formula gives
		\begin{align*}
			& L_{n,\gamma_n}(g) = \frac{1}{d} \hspace{-0.5cm} \sum\limits_{\substack{j=1 \\ \lambda_j(\Sigma_n) \in (a_{\gamma_n},b_{\gamma_n})}}^d \hspace{-0.5cm} g(\lambda_j(\Sigma_n)) = \frac{1}{d} \hspace{-0.5cm} \sum\limits_{\substack{j=1 \\ \lambda_j(\Sigma_n) \in (a_{\gamma_n},b_{\gamma_n})}}^d \hspace{-0.5cm} \frac{1}{2\pi \bm{i}} \oint_{\tilde{\gamma}_n} \frac{g(z)}{z-\lambda_j(\Sigma_n)} \, dz\\
			& = \frac{-1}{2\pi \bm{i}} \oint_{\tilde{\gamma}_n} g(z) \frac{1}{d} \hspace{-0.5cm} \sum\limits_{\substack{j=1 \\ \lambda_j(\Sigma_n) \in (a_{\gamma_n},b_{\gamma_n})}}^d \hspace{-0.5cm} \frac{1}{\lambda_j(\Sigma_n)-z} \, dz = \frac{-1}{2\pi \bm{i}} \oint_{\tilde{\gamma}_n} g(z) \cs_{H_n}(z) \, dz \ ,
		\end{align*}
		where the Stieltjes transform $\cs_{H_n}$ may be canonically extended to $\C^- = \{z \in \C \mid \Im(z)<0\}$ with $\cs_{H_n}(\ol{z}) = \ol{\cs_{H_n}(z)}$. This same anti-symmetry further allows for
		\begin{align}\label{Eq_CauchyApplication_new}
			& L_{n,\gamma_n}(g) = \frac{-1}{2\pi \bm{i}} \int_{\gamma_n} \big(g(z) \cs_{H_n}(z) - g(\ol{z}) \ol{\cs_{H_n}(z)}\big) \, dz \ .
		\end{align}
		
		\item[ii)] \textit{Bounding the estimation error}:\\
		By assumption (\ref{Eq_AdmissibleCurve_Cond1}) it on the mentioned event also holds that
		\begin{align}\label{Eq_PLSS_Consisnency_Step2_1}
			& \forall z \in \mathrm{image}(\gamma_n):  \ \big| \hat{s}_n(z) - \cs_{H_n}(z) \big| \leq n^{\varepsilon-1} \ .
		\end{align}
		One may thus bound
		\begin{align*}
			& \big| \hat{L}_{n,\gamma_n}(g) - L_{n,\gamma_n}(g) \big|\\
			& \overset{\text{(\ref{Eq_CauchyApplication_new})}}{\leq} \frac{1}{2\pi} \int_{\gamma_n} (|g(z)|+|g(\ol{z})|) \overbrace{|\hat{s}_n(z) - \cs_{H_n}(z)|}^{\leq n^{4\tilde{\varepsilon}-1} \text{ by (\ref{Eq_PLSS_Consisnency_Step2_1})}} \, |dz|\\
			& \leq \frac{2\|g\|_{\gamma_n}}{2\pi} \int_{\gamma_n} n^{\varepsilon-1} |dz| = n^{\varepsilon-1} \frac{\|g\|_{\gamma_n}}{\pi} \int_{\gamma_n} |dz| \ . \qed
		\end{align*}
	\end{itemize}

	\section{Proof of Theorem~\ref{Thm_GLSS_consistency}}
	Without loss of generality, assume $n$ to be large enough for the following bounds to hold, where $\tilde{\varepsilon}>0$ is a small constant satisfying $\tilde{\varepsilon} < \varepsilon \land (\frac{\varepsilon}{3}+\frac{1}{6})$.\\
	\begin{minipage}{0.35\textwidth}
		\begin{align}
			& \hspace{0.3cm} c_n \leq n  \label{Eq_GLLS_consistency_LargeN_1}\\
			& \hspace{-0cm} n^{\tilde{\varepsilon}-1} \leq \frac{1}{2\kappa} \label{Eq_GLLS_consistency_LargeN_2}\\
			& \hspace{-0.2cm} 2\kappa^2 n^{\tilde{\varepsilon}-1} \leq \frac{\tau}{2} \label{Eq_GLLS_consistency_LargeN_3}\\
			& \hspace{-0.7cm} \frac{\kappa}{c_n\tau^2} + \frac{1+c_n}{c_n\tau} \leq n^{\tilde{\varepsilon}} \label{Eq_GLLS_consistency_LargeN_5}
		\end{align}
	\end{minipage}
	\hspace{0.3cm}\vrule\hspace{0.3cm}
	\begin{minipage}{0.59\textwidth}
		\begin{align}
			& \hspace{-1.9cm} \frac{4\kappa n^{\tilde{\varepsilon}-1}}{\tau^2} + \frac{4\kappa^2 n^{2\tilde{\varepsilon}-1}}{\tau^2} + \frac{2 n^{\tilde{\varepsilon}-1}}{c_n\tau} \leq \frac{n^{3\tilde{\varepsilon}-1}}{2} \label{Eq_GLLS_consistency_LargeN_6}\\
			& \hspace{-2.05cm} \frac{4\kappa^2 n^{\tilde{\varepsilon}-1}}{c_n \tau^3} \Big( 1+c_n+c_n\kappa n^{\tilde{\varepsilon}} + \frac{\kappa}{\tau} \Big) \leq \frac{n^{3\tilde{\varepsilon}-1}}{2} \label{Eq_GLLS_consistency_LargeN_7}\\
			& \hspace{-0.6cm} \frac{2\sigma^2}{\pi\tau^2}n^{\tilde{\varepsilon}-\frac{1}{2}} + n^{3\tilde{\varepsilon}-1} \leq n^{\varepsilon-\frac{1}{2}} \label{Eq_GLLS_consistency_LargeN_4}
		\end{align}
	\end{minipage}
	\begin{itemize}
		\item[i)] \textit{First application of Cauchy's integral formula}:\\
		Let $\tilde{\gamma}_n^{(f)}$ and $\tilde{\gamma}_n^{(g)}$ be the closed curves going through $\mathrm{image}(\gamma_n^{(f)})$ and $\mathrm{image}(\gamma_n^{(g)})$ respectively with a counter-clockwise orientation. Cauchy's integral formula yields
		\begin{align}\label{Eq_GLSS_Consistency_Step1}
			& L_n\big( \overbrace{f \mathbbm{1}_{[a_{\gamma_n^{(f)}},b_{\gamma_n^{(f)}}]}}^{\eqcolon  \tilde{f}_n} \, , \, \overbrace{g \mathbbm{1}_{[a_{\gamma_n^{(g)}},b_{\gamma_n^{(g)}}]}}^{\eqcolon  \tilde{g}_n} \big) = \frac{1}{d} \tr\big( \tilde{f}_n(\bm{S}_n) \tilde{g}_n(\Sigma_n) \big) \nonumber\\
			& = \frac{-1}{2\pi \bm{i}} \oint_{\tilde{\gamma}_n^{(f)}} f(\tz_2) \frac{1}{d} \tr\big( \underbrace{(\bm{S}_n-\tz_2\Id_d)^{-1}}_{\eqcolon \bm{R}_n(\tz_2)} \tilde{g}_n(\Sigma_n) \big) \, d\tz_2 \ .
		\end{align}
		
		\item[ii)] \textit{Approximation of $\frac{1}{d} \tr(\bm{R}_n(\tz_2) \tilde{g}_n(\Sigma_n))$}:\\
		The assumption of $(\gamma_n^{(f)},\gamma_n^{(g)})$ being an admissible pair guarantees the inclusion
		\begin{align}\label{Eq_GLSS_consistency_fCurveIn_bSn}
			& \mathrm{image}(\gamma_n^{(f)}) \subset \bS(\tau,n) \ ,
		\end{align}
		where $\bS(\tau,n)$ is the spectral domain defined in Theorem~\ref{Thm_OuterLaw}. Further, result (\ref{Eq_OuterLaw_Matrix}) of Theorem~\ref{Thm_OuterLaw} yields the existence of a $C'' = C''(\tilde{\varepsilon},D+2K+4,\tau)>0$ such that
		\begin{align}\label{Eq_OuterLaw_Matrix_copy}
			& \forall g : [0,\infty) \rightarrow \R \ \text{ bounded and measurable} : \nonumber\\
			& \bP\bigg( \forall \tz \in \bS(\tau,n) : \ \Big| \frac{1}{d} \tr\big( g(\Sigma_n) \bm{R}_n(\tz) \big) - \frac{1}{d} \tr\Big(- \frac{1}{\tz}g(\Sigma_n)(\Id_d + \cs_{\ul{\nu}_n}(\tz)\Sigma_n)^{-1} \Big) \Big| \nonumber\\
			& \hspace{3cm} \leq  2\frac{\sigma^2 \|g\|_\infty}{\tau^2} n^{\tilde{\varepsilon}-\frac{1}{2}} \bigg) \geq 1 - \frac{d C''}{n^{D+2K+4}} \overset{\text{(\ref{Eq_GLLS_consistency_LargeN_1})}}{\geq} 1 - \frac{C''}{n^{D+2K+2}} \ .
		\end{align}
		An issue with the direct application of (\ref{Eq_OuterLaw_Matrix_copy}) to the approximation of $\frac{1}{d} \tr\big( \tilde{g}_n(\Sigma_n) \bm{R}_n(\tz) \big)$ is that $\tilde{g}_n = g \mathbbm{1}_{[a_{\gamma_n^{(g)}},b_{\gamma_n^{(g)}}]}$ is not deterministic in the sense that the interval $[a_{\gamma_n^{(g)}},b_{\gamma_n^{(g)}}]$ depends on the random curve $\gamma_n^{(g)}$. For this reason, it was assumed in (a) of Definition~\ref{Def_AdmissiblePair} that $a_{\gamma_n^{(g)}},b_{\gamma_n^{(g)}} \in J_{n,K}$. For any two points $a,b$ from the net
		\begin{align*}
			& J_{n,K} \coloneq \Big\{\frac{k}{n^K} \ \Big| \ k \in \Z , \, \frac{k}{n^K} \in [-\kappa,\kappa]\Big\} \ ,
		\end{align*}
		denote $\tilde{g}_{a,b}(x) = g(x) \mathbbm{1}_{x \in [a,b]}$, then the application of (\ref{Eq_OuterLaw_Matrix_copy}) to each such $\tilde{g}_{a,b}$ yields
		\begin{align}\label{Eq_OuterLaw_Matrix_copy2}
			& \bP\bigg( \exists \tz \in \bS(\tau,n) , \, \exists a,b \in J_{n,K} : \nonumber\\
			& \hspace{1cm} \Big| \frac{1}{d} \tr\big( \tilde{g}_{a,b}(\Sigma_n) \bm{R}_n(\tz) \big) - \frac{1}{d} \tr\Big(- \frac{1}{\tz}\tilde{g}_{a,b}(\Sigma_n)(\Id_d + \cs_{\ul{\nu}_n}(\tz)\Sigma_n)^{-1} \Big) \Big| \nonumber\\
			& \hspace{3cm} >  2\frac{\sigma^2 \|\tilde{g}_{a,b}\|_\infty}{\tau^2} n^{\tilde{\varepsilon}-\frac{1}{2}} \bigg) < C''\frac{(\#J_{n,K})^2}{n^{D+(2K+2)}} \leq \frac{C''}{n^{D}} \ .
		\end{align}
		From (\ref{Eq_GLSS_consistency_fCurveIn_bSn}) and (a) of Definition~\ref{Def_AdmissiblePair} it thus follows that
		\begin{align}\label{Eq_OuterLaw_Matrix_copy3}
			& \bP\bigg( \forall \tz \in \mathrm{image}(\gamma_n^{(g)}) : \nonumber\\
			& \hspace{1cm} \Big| \frac{1}{d} \tr\big( \tilde{g}_{n}(\Sigma_n) \bm{R}_n(\tz) \big) - \frac{1}{d} \tr\Big(- \frac{1}{\tz}\tilde{g}_{n}(\Sigma_n)(\Id_d + \cs_{\ul{\nu}_n}(\tz)\Sigma_n)^{-1} \Big) \Big| \nonumber\\
			& \hspace{3cm} \leq  2\frac{\sigma^2 \|\tilde{g}_{n}\|_\infty}{\tau^2} n^{\tilde{\varepsilon}-\frac{1}{2}} \bigg) \geq 1 - \frac{C''}{n^{D}} \ .
		\end{align}

		\item[iii)] \textit{Second application of Cauchy's integral formula}:\\
		A second application of Cauchy's integral formula gives the equality
		\begin{align}\label{Eq_GLSS_trueKernelDef}
			& \frac{-1}{2\pi \bm{i}} \oint_{\tilde{\gamma}_n^{(f)}} f(\tz_2) \frac{1}{d} \tr\Big(- \frac{1}{\tz_2}\tilde{g}(\Sigma_n)(\Id_d + \cs_{\ul{\nu}_n}(\tz_2)\Sigma_n)^{-1} \Big) \, d\tz_2 \nonumber\\
			& = \Big(\frac{-1}{2\pi \bm{i}}\Big)^2 \oint_{\tilde{\gamma}_n^{(f)}} f(\tz_2) \oint_{\tilde{\gamma}_n^{(g)}} g(z_1) \nonumber\\
			& \hspace{2cm} \times \underbrace{\frac{1}{d} \tr\Big(- \frac{1}{\tz_2}(\Sigma_n-z_1\Id_d)^{-1}(\Id_d + \cs_{\ul{\nu}_n}(\tz_2)\Sigma_n)^{-1} \Big)}_{\eqcolon  -k_n(z_1,\tz_2)} \, dz_1 \, d\tz_2 \ .
		\end{align}
		Further examination of $k_n(z_1,\tz_2)$ with the notation $z_2 \coloneq \frac{-1}{\cs_{\ul{\nu}_n}(\tz_2)}$ yields
		\begin{align}\label{Eq_GLSS_theoreticalKernelCalc}
			& k_n(z_1,\tz_2) = \frac{1}{d} \tr\Big(\frac{1}{\tz_2}(\Sigma_n-z_1\Id_d)^{-1}(\Id_d + \cs_{\ul{\nu}_n}(\tz_2)\Sigma_n)^{-1} \Big) \nonumber\\
			& = \frac{1}{\tz_2} \int_\R \frac{1}{(\lambda-z_1) (1+\cs_{\ul{\nu}_n}(\tz_2)\lambda)} \, dH_n(\lambda) \nonumber\\
			& = -\frac{z_2}{\tz_2} \int_\R \frac{1}{(\lambda-z_1) (\lambda-z_2)} \, dH_n(\lambda) \overset{\text{(\ref{Eq_DefStieltjes})}}{=} -\frac{z_2}{\tz_2} \frac{\cs_{H_n}(z_1)-\cs_{H_n}(z_2)}{z_1-z_2} \nonumber\\
			& = \frac{1}{\tz_2} \frac{(1-c_nz_2\cs_{H_n}(z_1)-c_n)--(1-c_nz_2\cs_{H_n}(z_2)-c_n)}{c_n(z_1-z_2)} \nonumber\\
			& = \frac{1}{\tz_2} \frac{(1-c_nz_2\cs_{H_n}(z_1)-c_n)z_2-\overbrace{(1-c_nz_2\cs_{H_n}(z_2)-c_n)z_2}^{= \Phi_{H_n,c_n}(z_2) = \tz_2 \text{ by Lemma~\ref{Lemma_SpaceTransform}}}}{c_n(z_1-z_2)z_2} \nonumber\\
			& = \frac{(1-c_nz_2\cs_{H_n}(z_1)-c_n) - \tz_2/z_2}{c_n(z_1-z_2) \tz_2} = \frac{(1+c_n\frac{\cs_{H_n}(z_1)}{\cs_{\ul{\nu}_n}(\tz_2)}-c_n) + \tz_2\cs_{\ul{\nu}_n}(\tz_2)}{c_n(z_1+\frac{1}{\cs_{\ul{\nu}_n}(\tz_2)}) \tz_2} \ .
		\end{align}
		
		\item[iv)] \textit{Lower bound for $|z+1/\cs_{\ul{\nu}_n}(\tz)|$}:\\
		Assumption (c) of Definition~\ref{Def_AdmissiblePair} yields
		\begin{align}\label{Eq_GLSS_kApprox1}
			& \Big| z+\frac{1}{\cs_{\hat{\ul{\nu}}_n}(\tz)} \Big| \overset{\text{(\ref{Eq_Def_varphi})}}{=} \big| z - \varphi_{\hat{\nu}_n,c_n}(\tz) \big| \geq \tau
		\end{align}
		for all $z \in \mathrm{image}(\gamma_n^{(g)})$ and $\tz \in \mathrm{image}(\gamma_n^{(f)})$. By assumption (\ref{Eq_GLSS_gammafCondition2}), the image $\mathrm{image}(\gamma_n^{(f)})$ lies in the spectral domain $\hat{\bS}(\tau/2,n)$ from Theorem~\ref{Corollary_OuterLaw}, which yields the existence of a constant $C'=C'(\tilde{\varepsilon},\tau/2,D)$ such that
		\begin{align}\label{Eq_OuterLaw_StieltjesInfty_copy}
			& \bP\Big( \forall \tz \in \mathrm{image}(\gamma_n^{(f)}) : \ \big| \cs_{\hat{\ul{\nu}}_n}(\tz) - \cs_{\ul{\nu}_n}(\tz) \big| \leq n^{\tilde{\varepsilon}-1} \Big) \geq 1 - \frac{C'}{n^D} \ .
		\end{align}
		For all $\tz \in \mathrm{image}(\gamma_n^{(f)})$, assumption (\ref{Eq_GLSS_gammafCondition2_2}) gives the bound
		\begin{align}\label{Eq_GLSSCOnsistency_Step4_1}
			& \frac{1}{\kappa} \leq \underbrace{|\cs_{\hat{\ul{\nu}}_n}(\tz)|}_{\overset{\text{(\ref{Eq_Def_varphi})}}{=} 1/|\varphi_{\hat{\nu}_n,c_n}(\tz)|} \leq \frac{1}{\tau} \ ,
		\end{align}
		so on the event $\mathcal{E}_{n,\text{(\ref{Eq_OuterLaw_StieltjesInfty_copy})}}$ from (\ref{Eq_OuterLaw_StieltjesInfty_copy}) it holds that
		\begin{align}\label{Eq_GLSSCOnsistency_Step4_2}
			& |\cs_{\ul{\nu}_n}(\tz)| \geq \overbrace{|\cs_{\hat{\ul{\nu}}_n}(\tz)|}^{\overset{\text{(\ref{Eq_GLSSCOnsistency_Step4_1})}}{\geq} 1/\kappa} - \overbrace{|\cs_{\hat{\ul{\nu}}_n}(\tz) - \cs_{\ul{\nu}_n}(\tz)|}^{\overset{\text{(\ref{Eq_OuterLaw_StieltjesInfty_copy})}}{\leq} n^{\tilde{\varepsilon}-1}} \overset{\text{(\ref{Eq_GLLS_consistency_LargeN_2})}}{\geq} \frac{1}{2\kappa}
		\end{align}
		and also
		\begin{align}\label{Eq_GLSSCOnsistency_Step4_3}
			& |\varphi_{\hat{\nu}_n,c_n}(\tz) - \varphi_{\nu_n,c_n}(\tz)| = \Big| \frac{-1}{\cs_{\hat{\ul{\nu}}_n}(\tz)} - \frac{-1}{\cs_{\ul{\nu}_n}(\tz)} \Big| = \frac{|\cs_{\hat{\ul{\nu}}_n}(\tz)-\cs_{\ul{\nu}_n}(\tz)|}{|\cs_{\hat{\ul{\nu}}_n}(\tz)|\,|\cs_{\ul{\nu}_n}(\tz)|} \nonumber\\
			& \overset{\text{(\ref{Eq_OuterLaw_StieltjesInfty_copy})}}{\leq} \frac{n^{\tilde{\varepsilon}-1}}{|\cs_{\hat{\ul{\nu}}_n}(\tz)|\,|\cs_{\ul{\nu}_n}(\tz)|} \overset{\text{(\ref{Eq_GLSSCOnsistency_Step4_1}), (\ref{Eq_GLSSCOnsistency_Step4_2})}}{\leq} 2\kappa^2 n^{\tilde{\varepsilon}-1} \ \Big[ \overset{\text{(\ref{Eq_GLLS_consistency_LargeN_3})}}{\leq} \frac{\tau}{2} \Big] \ .
		\end{align}
		Finally, one may combine (\ref{Eq_GLSS_kApprox1}) and (\ref{Eq_GLSSCOnsistency_Step4_3}) to on the event $\mathcal{E}_{n,\text{(\ref{Eq_OuterLaw_StieltjesInfty_copy})}}$ see
		\begin{align}\label{Eq_GLSSCOnsistency_Step4_4}
			&\Big| z+\frac{1}{\cs_{\ul{\nu}_n}(\tz)} \Big| = \big| z - \varphi_{\nu_n,c_n}(\tz) \big| \nonumber\\
			& \geq \underbrace{|z - \varphi_{\hat{\nu}_n,c_n}(\tz)|}_{\overset{\text{(\ref{Eq_GLSS_kApprox1})}}{\geq} \tau} - \underbrace{|\varphi_{\hat{\nu}_n,c_n}(\tz) - \varphi_{\nu_n,c_n}(\tz)|}_{\overset{\text{(\ref{Eq_GLSSCOnsistency_Step4_3})}}{\geq} \tau/2} \geq \frac{\tau}{2} \ .
		\end{align}
		for all $\tz \in \mathrm{image}(\gamma_n^{(f)})$.
		
		\item[v)] \textit{Approximation of $k_n(z_1,\tz_2)$ by $\hat{k}_n(z_1,\tz_2)$}:\\
		For ease of notation, write
		\begin{align}\label{Eq_GLSS_kApprox_EaseOfNotation}
			& s_1 \coloneq \cs_{H_n}(z_1) \ \ ; \ \ \hat{s}_1 \coloneq \hat{s}_n(z_1) \ \ ; \ \ s_2 \coloneq \cs_{\ul{\nu}_n}(\tz_2) \ \ ; \ \ \hat{s}_2 \coloneq \cs_{\hat{\ul{\nu}}_n}(\tz_2) \ ,
		\end{align}
		then the difference between $k_n$ and $\hat{k}_n$ as defined in (\ref{Eq_DefGLSS_kernel}) may be bounded by
		\begin{align}\label{Eq_GLSS_kApprox0}
			& \big|k_n(z_1,\tz_2) - \hat{k}_n(z_1,\tz_2)\big| = \Big| \frac{(1+c_ns_1/s_2-c_n) + \tz_2s_2}{c_n(z_1+1/s_2) \tz_2} - \frac{(1+c_n\hat{s}_1/\hat{s}_2-c_n) + \tz_2\hat{s}_2}{c_n(z_1+1/\hat{s}_2) \tz_2} \Big| \nonumber\\
			& = \Big| \frac{(z_1+1/\hat{s}_2) ((1+c_ns_1/s_2-c_n) + \tz_2s_2) - (z_1+1/s_2)((1+c_n\hat{s}_1/\hat{s}_2-c_n) + \tz_2\hat{s}_2)}{c_n(z_1+1/s_2) (z_1+1/\hat{s}_2) \tz_2} \Big| \nonumber\\
			& \leq \Big| \frac{ ((1+c_ns_1/s_2-c_n) + \tz_2s_2) - ((1+c_n\hat{s}_1/\hat{s}_2-c_n) + \tz_2\hat{s}_2)}{c_n(z_1+1/s_2) \tz_2} \Big| \nonumber\\
			& \hspace{1cm} + \Big| \frac{(z_1+1/\hat{s}_2) - (z_1+1/s_2)}{c_n(z_1+1/s_2) (z_1+1/\hat{s}_2) \tz_2} \Big| \, \big| (1+c_n\hat{s}_1/\hat{s}_2-c_n) + \tz_2\hat{s}_2 \big| \nonumber\\
			& = \overbrace{\Big| \frac{ c_n(s_1/s_2 - \hat{s}_1/\hat{s}_2) + \tz_2(s_2-\hat{s}_2)}{c_n(z_1+1/s_2) \tz_2} \Big|}^{\eqcolon  A} \nonumber\\
			& \hspace{1cm} + \underbrace{\Big| \frac{1/\hat{s}_2 - 1/s_2}{c_n(z_1+1/s_2) (z_1+1/\hat{s}_2) \tz_2} \Big| \, \big| (1+c_n\hat{s}_1/\hat{s}_2-c_n) + \tz_2\hat{s}_2 \big|}_{\eqcolon  B} \ .
		\end{align}
		Assumption (\ref{Eq_GLSS_gammafCondition2}) together with $\gamma_n^{(f)} : (0,1) \rightarrow B_{\kappa}^{\C^+}(0)$ in the definition of an admissible pair yields
		\begin{align}\label{Eq_GLSS_kApprox2}
			& \tau < |\tz_2| < \kappa
		\end{align}
		for all $\tz_2 \in \mathrm{image}(\gamma_n^{(f)})$.
		In (a) of Definition~\ref{Def_AdmissiblePair} it is assumed that $\gamma_n^{(g)}$ is admissible, so $\mathrm{image}(\gamma_n^{(g)}) \subset \hat{\bD}(\tau,\kappa,n)$ and Theorem~\ref{Thm_Consistency} yields the existence of a constant $C'''=C'''(\tau,\kappa,\tilde{\varepsilon},D)>0$ such that
		\begin{align}\label{Eq_GLSS_kApprox4}
			& \bP\Big( \forall z \in \mathrm{image}(\gamma_n^{(g)}) : \ \big| \hat{s}_n(z) - \cs_{H_n}(z) \big| \leq n^{\tilde{\varepsilon}-1} \Big) \geq 1 - \frac{C'''}{n^D} \ .
		\end{align}
		The assumption $\mathrm{image}(\gamma_n^{(g)}) \subset \hat{\bD}(\tau,\kappa,n)$ also gives $|\hat{\Phi}_n(z_1)| \overset{\text{(\ref{Eq_Consistency_EmpDomain2})}}{\leq} \kappa$ and $|z_1| \overset{\text{(\ref{Eq_Consistency_EmpDomain1_1})}}{\geq} \tau$, which may be used to bound $|\hat{s}_n(z_1)|$ by the calculation
		\begin{align}\label{Eq_GLSS_kApprox4_2}
			& |\hat{s}_n(z_1)| = \Big| \frac{\hat{\Phi}_n(z_1)/z_1-1+c_n}{c_nz_1} \Big| \overset{\text{(\ref{Eq_Consistency_EmpDomain1_1})}}{\leq} \frac{|\hat{\Phi}_n(z_1)|}{c_n\tau^2} + \frac{1+c_n}{c_n\tau} \nonumber\\
			& \overset{\text{(\ref{Eq_Consistency_EmpDomain2})}}{\leq} \frac{\kappa}{c_n\tau^2} + \frac{1+c_n}{c_n\tau} \overset{\text{(\ref{Eq_GLLS_consistency_LargeN_5})}}{\leq} n^{\tilde{\varepsilon}} \ .
		\end{align}
		With these tools, one may then on the intersection $\mathcal{E}_{n,\text{(\ref{Eq_OuterLaw_StieltjesInfty_copy})}} \cap \mathcal{E}_{n,\text{(\ref{Eq_GLSS_kApprox4})}}$ of the events from (\ref{Eq_OuterLaw_StieltjesInfty_copy}) and (\ref{Eq_GLSS_kApprox4}) bound the summand $A$ from (\ref{Eq_GLSS_kApprox0}) as follows
		\begin{align*}
			& A = \bigg| \frac{ c_n\big(\frac{\cs_{H_n}(z_1)}{\cs_{\ul{\nu}_n}(\tz_2)} - \frac{\hat{s}_n(z_1)}{\cs_{\hat{\ul{\nu}}_n}(\tz_2)}\big) + \tz_2(\cs_{\ul{\nu}_n}(\tz_2)-\cs_{\hat{\ul{\nu}}_n}(\tz_2))}{c_n\big(z_1+\frac{1}{\cs_{\ul{\nu}_n}(\tz_2)}\big) \tz_2} \bigg|\\
			& \overset{\text{(\ref{Eq_GLSSCOnsistency_Step4_4}), (\ref{Eq_GLSS_kApprox2})}}{\leq} \frac{2}{\tau^2}\Big|\frac{\cs_{H_n}(z_1)}{\cs_{\ul{\nu}_n}(\tz_2)} - \frac{\hat{s}_n(z_1)}{\cs_{\hat{\ul{\nu}}_n}(\tz_2)}\Big| + \frac{2}{c_n\tau} \underbrace{\big|\cs_{\ul{\nu}_n}(\tz_2)-\cs_{\hat{\ul{\nu}}_n}(\tz_2)\big|}_{\overset{\text{(\ref{Eq_OuterLaw_StieltjesInfty_copy})}}{\leq} n^{\tilde{\varepsilon}-1}}\\
			& \leq \frac{2}{\tau^2}\frac{|\cs_{\hat{\ul{\nu}}_n}(\tz_2)\cs_{H_n}(z_1) \pm \hat{s}_n(z_1) \cs_{\hat{\ul{\nu}}_n}(\tz_2) - \hat{s}_n(z_1) \cs_{\ul{\nu}_n}(\tz_2)|}{|\cs_{\hat{\ul{\nu}}_n}(\tz_2) \cs_{\ul{\nu}_n}(\tz_2)|} + \frac{2 n^{\tilde{\varepsilon}-1}}{c_n\tau}\\
			& \leq \frac{2}{\tau^2} \frac{\overbrace{|\cs_{H_n}(z_1) - \hat{s}_n(z_1)|}^{\overset{\text{(\ref{Eq_GLSS_kApprox4})}}{\leq} n^{\tilde{\varepsilon}-1}}}{\underbrace{|\cs_{\ul{\nu}_n}(\tz_2)|}_{\overset{\text{(\ref{Eq_GLSSCOnsistency_Step4_2})}}{\geq} 1/2\kappa}} + \frac{2}{\tau^2} \frac{\overbrace{|\hat{s}_n(z_1)|}^{\overset{\text{(\ref{Eq_GLSS_kApprox4_2})}}{\leq} n^{\tilde{\varepsilon}}} \, \overbrace{|\cs_{\hat{\ul{\nu}}_n}(\tz_2) - \cs_{\ul{\nu}_n}(\tz_2)|}^{\overset{\text{(\ref{Eq_OuterLaw_StieltjesInfty_copy})}}{\leq} n^{\tilde{\varepsilon}-1}}}{\underbrace{|\cs_{\hat{\ul{\nu}}_n}(\tz_2) \cs_{\ul{\nu}_n}(\tz_2)|}_{\overset{\text{(\ref{Eq_GLSSCOnsistency_Step4_1}), (\ref{Eq_GLSSCOnsistency_Step4_2})}}{\geq} 1/2\kappa^2}} + \frac{2 n^{\tilde{\varepsilon}-1}}{c_n\tau}\\
			& \leq \frac{4\kappa n^{\tilde{\varepsilon}-1}}{\tau^2} + \frac{4\kappa^2 n^{2\tilde{\varepsilon}-1}}{\tau^2} + \frac{2 n^{\tilde{\varepsilon}-1}}{c_n\tau} \overset{\text{(\ref{Eq_GLLS_consistency_LargeN_6})}}{\leq} \frac{n^{3\tilde{\varepsilon}-1}}{2} \ .
		\end{align*}
		Similarly, the summand $B$ may on the intersection $\mathcal{E}_{n,\text{(\ref{Eq_OuterLaw_StieltjesInfty_copy})}} \cap \mathcal{E}_{n,\text{(\ref{Eq_GLSS_kApprox4})}}$ be bounded by
		\begin{align*}
			& B = \bigg| \frac{\frac{1}{\cs_{\hat{\ul{\nu}}_n}(\tz_2)} - \frac{1}{\cs_{\ul{\nu}_n}(\tz_2)}}{c_n\big(z_1+\frac{1}{\cs_{\ul{\nu}_n}(\tz_2)}\big) \big(z_1+\frac{1}{\cs_{\hat{\ul{\nu}}_n}(\tz_2)}\big) \tz_2} \bigg| \, \Big| \Big(1+c_n\frac{\hat{s}_n(z_1)}{\cs_{\hat{\ul{\nu}}_n}(\tz_2)}-c_n\Big) + \tz_2\cs_{\hat{\ul{\nu}}_n}(\tz_2) \Big|\\
			& \overset{\substack{\text{(\ref{Eq_GLSS_kApprox1}),} \\ \text{(\ref{Eq_GLSSCOnsistency_Step4_4}), (\ref{Eq_GLSS_kApprox2})}}}{\leq} \frac{2}{c_n\tau^3} \underbrace{\Big| \frac{1}{\cs_{\hat{\ul{\nu}}_n}(\tz_2)} - \frac{1}{\cs_{\ul{\nu}_n}(\tz_2)} \Big|}_{\overset{\text{(\ref{Eq_GLSSCOnsistency_Step4_3})}}{\leq} 2\kappa^2 n^{\tilde{\varepsilon}-1}} \, \Big( 1+c_n + c_n\underbrace{\Big|\frac{\hat{s}_n(z_1)}{\cs_{\hat{\ul{\nu}}_n}(\tz_2)}\Big|}_{\overset{\text{(\ref{Eq_GLSS_kApprox4_2}), (\ref{Eq_GLSSCOnsistency_Step4_1})}}{\leq} \kappa n^{\tilde{\varepsilon}}} + \underbrace{|\tz_2\cs_{\hat{\ul{\nu}}_n}(\tz_2)|}_{\overset{\text{(\ref{Eq_GLSSCOnsistency_Step4_1}), (\ref{Eq_GLSS_kApprox2})}}{\leq} \kappa/\tau} \Big)\\
			& \leq \frac{4\kappa^2 n^{\tilde{\varepsilon-1}}}{c_n \tau^3} \Big( 1+c_n+c_n\kappa n^{\tilde{\varepsilon}} + \frac{\kappa}{\tau} \Big) \overset{\text{(\ref{Eq_GLLS_consistency_LargeN_7})}}{\leq} \frac{n^{3\tilde{\varepsilon}-1}}{2} \ .
		\end{align*}
		By combining the last two bounds with (\ref{Eq_GLSS_kApprox0}) it was shown that
		\begin{align*}
			& \forall z_1 \in \mathrm{image}(\gamma_n^{(g)}) \, \forall \tz_2 \in \mathrm{image}(\gamma_n^{(f)}) : \ \big|k_n(z_1,\tz_2) - \hat{k}_n(z_1,\tz_2)\big| \leq A+B \leq n^{3\tilde{\varepsilon}-1}
		\end{align*}
		and the same bounds may be shown analogously for $z_1 \in \ol{\mathrm{image}(\gamma_n^{(g)})}$ or $\tz_2 \in \ol{\mathrm{image}(\gamma_n^{(f)})}$, which leads to the bound
		\begin{align}\label{Eq_GLSSConsistency_Step5_result}
			& \forall z_1 \in \mathrm{image}(\tilde{\gamma}_n^{(g)}) \, \forall \tz_2 \in \mathrm{image}(\tilde{\gamma}_n^{(f)}) : \ \big|k_n(z_1,\tz_2) - \hat{k}_n(z_1,\tz_2)\big| \leq n^{3\tilde{\varepsilon}-1}
		\end{align}
		on the event $\mathcal{E}_{n,\text{(\ref{Eq_OuterLaw_StieltjesInfty_copy})}} \cap \mathcal{E}_{n,\text{(\ref{Eq_GLSS_kApprox4})}}$.
		
		\item[vi)] \textit{Gathering bounds}:\\
		The following calculations take place on the intersection $\mathcal{E}_{n,\text{(\ref{Eq_OuterLaw_Matrix_copy3})}} \cap \mathcal{E}_{n,\text{(\ref{Eq_OuterLaw_StieltjesInfty_copy})}} \cap \mathcal{E}_{n,\text{(\ref{Eq_GLSS_kApprox4})}}$ of the events from (\ref{Eq_OuterLaw_Matrix_copy3}), (\ref{Eq_OuterLaw_StieltjesInfty_copy}) and (\ref{Eq_GLSS_kApprox4}). By step (i), one has
		\begin{align*}
			& L_n\big( f \mathbbm{1}_{[a_{\gamma_n^{(f)}},b_{\gamma_n^{(f)}}]} \, , \, g \mathbbm{1}_{[a_{\gamma_n^{(g)}},b_{\gamma_n^{(g)}}]} \big) \overset{\text{(\ref{Eq_GLSS_Consistency_Step1})}}{=} \frac{-1}{2\pi \bm{i}} \oint_{\tilde{\gamma}_n^{(f)}} f(\tz_2) \frac{1}{d} \tr\big( \bm{R}_n(\tz_2) \tilde{g}_n(\Sigma_n) \big) \, d\tz_2 \ ,
		\end{align*}
		which with step (ii) yields
		\begin{align}\label{Eq_GLSSConsistency_Step6_1}
			& \bigg| L_n\big( f \mathbbm{1}_{[a_{\gamma_n^{(f)}},b_{\gamma_n^{(f)}}]} \, , \, g \mathbbm{1}_{[a_{\gamma_n^{(g)}},b_{\gamma_n^{(g)}}]} \big) \nonumber\\
			& \hspace{2cm} - \underbrace{\frac{-1}{2\pi \bm{i}} \oint_{\tilde{\gamma}_n^{(f)}} f(\tz_2) \frac{1}{d} \tr\Big(- \frac{1}{\tz}\tilde{g}_{n}(\Sigma_n)(\Id_d + \cs_{\ul{\nu}_n}(\tz)\Sigma_n)^{-1} \Big) \, d\tz_2}_{\eqcolon  I_n} \bigg| \nonumber\\
			& \overset{\text{(\ref{Eq_OuterLaw_Matrix_copy3})}}{\leq} \frac{1}{2\pi} \oint_{\tilde{\gamma}_n^{(f)}} |f(\tz_2)| 2\frac{\sigma^2 \|\tilde{g}_n\|_\infty}{\tau^2} n^{\tilde{\varepsilon}-\frac{1}{2}} \, |d\tz_2| \leq \frac{\sigma^2 \|f\|_{\gamma_n^{(f)}} \|\tilde{g}_n\|_\infty}{\pi \tau^2} n^{\tilde{\varepsilon}-\frac{1}{2}} \oint_{\tilde{\gamma}_n^{(f)}} |d\tz_2| \nonumber\\
			& = \frac{2\ell(\gamma_n^{(f)}) \sigma^2 \|f\|_{\gamma_n^{(f)}} \|\tilde{g}_n\|_\infty}{\pi \tau^2} n^{\tilde{\varepsilon}-\frac{1}{2}} \ .
		\end{align}
		Note that $\tilde{g}_n$ was introduced in step (ii) as notation for $g \mathbbm{1}_{[a_{\gamma_n^{(g)}},b_{\gamma_n^{(g)}}]}$. By the maximum modulus principle guarantees $\|\tilde{g}_n\|_\infty \leq \|g\|_{\gamma_n^{(g)}}$.
		In turn, step (iii) gives
		\begin{align*}
			& I_n \overset{\text{(\ref{Eq_GLSS_trueKernelDef})}}{=} \frac{1}{4\pi^2} \oint_{\tilde{\gamma}_n^{(f)}} \oint_{\tilde{\gamma}_n^{(g)}} f(\tz_2) g(z_1) k_n(z_1,\tz_2) \, dz_1 \, d\tz_2
		\end{align*}
		and the approximation of $k_n(z_1,\tz_2)$ by $\hat{k}_n(z_1,\tz_2)$ shown in step (v) is sufficient for
		\begin{align}\label{Eq_GLSSConsistency_Step6_2}
			& \bigg| I_n - \overbrace{\frac{1}{4\pi^2} \oint_{\tilde{\gamma}_n^{(f)}} \oint_{\tilde{\gamma}_n^{(g)}} f(\tz_2) g(z_1) \hat{k}_n(z_1,\tz_2) \, dz_1 \, d\tz_2}^{\overset{\text{(\ref{Eq_DefGLSS_estimator})}}{=} \hat{L}_{n,\gamma_n^{(f)},\gamma_n^{(g)}}(f,g)} \bigg| \nonumber\\
			& \overset{\text{(\ref{Eq_GLSSConsistency_Step5_result})}}{\leq} \oint_{\tilde{\gamma}_n^{(f)}} \oint_{\tilde{\gamma}_n^{(g)}} |f(\tz_2)| \, |g(z_1)| \, n^{3\tilde{\varepsilon}-1} \, |dz_1| \, |d\tz_2| \nonumber\\
			& \leq \|f\|_{\gamma_n^{(f)}} \ell(\gamma_n^{(f)}) \, \|g\|_{\gamma_n^{(g)}} \ell(\gamma_n^{(g)}) \, n^{3\tilde{\varepsilon}-1} \ .
		\end{align}
		Finally, observe that (\ref{Eq_GLSSConsistency_Step6_1}) and (\ref{Eq_GLSSConsistency_Step6_2}) show
		\begin{align*}
			& \Big| \hat{L}_{n,\gamma_n^{(f)},\gamma_n^{(g)}}(f,g) - L_n\big( f \mathbbm{1}_{[a_{\gamma_n^{(f)}},b_{\gamma_n^{(f)}}]} \, , \, g \mathbbm{1}_{[a_{\gamma_n^{(g)}},b_{\gamma_n^{(g)}}]} \big) \Big|\\
			& \leq \frac{2\ell(\gamma_n^{(f)}) \sigma^2 \|f\|_{\gamma_n^{(f)}} \|g\|_{\gamma_n^{(g)}}}{\pi \tau^2} n^{\tilde{\varepsilon}-\frac{1}{2}} + \|f\|_{\gamma_n^{(f)}} \ell(\gamma_n^{(f)}) \, \|g\|_{\gamma_n^{(g)}} \ell(\gamma_n^{(g)}) \, n^{3\tilde{\varepsilon}-1}\\
			& \leq \|f\|_{\gamma_n^{(f)}} \ell(\gamma_n^{(f)}) \, \|g\|_{\gamma_n^{(g)}} \big(\ell(\gamma_n^{(g)})+1\big) \, \Big( \underbrace{\frac{2\sigma^2}{\pi\tau^2}n^{\tilde{\varepsilon}-\frac{1}{2}} + n^{3\tilde{\varepsilon}-1}}_{\overset{\text{(\ref{Eq_GLLS_consistency_LargeN_4})}}{\leq} n^{\varepsilon-\frac{1}{2}}} \Big) \ .
		\end{align*}
		Choosing a $C>0$ greater than $C'+C''+C'''$, and large enough that the bounds (\ref{Eq_GLLS_consistency_LargeN_1})--(\ref{Eq_GLLS_consistency_LargeN_4}) hold for all $n \geq C^{\frac{1}{D}}$ thus yields the wanted result. \qedhere
	\end{itemize}

	\section{Proofs of CLTs}\label{Section_CLT_Proofs}

	The proof of Theorem~\ref{Thm_CLT_Inversion} will require the following three auxiliary lemmas. The second is a straightforward extension of Theorem 2.8 in \cite{FunctionalDeltaReference} and may be proved analogously using Theorem 7.24 of \cite{FunctionalDeltaReference}. The first and third lemmas are proven at the end of this section.
	
	\begin{lemma}[The Marchenko--Pastur law under bounded fourth moments]\label{Lemma_FourthMomentMP_law}\
		\\
		Suppose \ref{EI_ItemAssumption_Asymptotics}--\ref{EI_ItemAssumption_sigmaBound} hold and that there exists a constant $C_4>0$ such that
		\begin{align}\label{Eq_BoundedFourthMoments}
			& \forall n \in \N, \, j\leq d, \, k \leq n : \ \E\big[ |(\bm{X}_n)_{j,k}|^4 \big] \leq C_4 \ .
		\end{align}
		Then the Marchenko--Pastur law (\ref{Eq_MP_law}) holds.
	\end{lemma}
	
	\begin{lemma}[Functional delta method]\label{Lemma_FunctionalDeltaMethod_General}\
		\\
		For normed spaces $D$ and $E$, and a subset $D_\phi \subset D$, let $\phi : D_\phi \rightarrow E$ be Hadamard-differentiable at $x \in D_\phi$ tangentially to a set $D_0 \subset D$, i.e. there exists a continuous linear map $\phi'_x : D \rightarrow E$ such that
		\begin{align*}
			& \frac{\phi(x+t_mh_m) - \phi(x)}{t_m} \xrightarrow{m \to \infty} \phi'_x(h)
		\end{align*}
		holds for any sequences $(t_m)_{m \in\N} \subset (0,\infty)$ and $(h_m)_{m \in \N} \subset D$ with $t_m \to 0$, $h_m \to h \in D_0$ and $x+t_mh_m \in D_\phi$ for all $m \in \N$. Further, suppose there exists a sequence $(X_n)_{n \in \N}$ of $D_\phi$-valued random variables and a sequence of values $(x_n)_{n \in \N} \subset D_\phi$ with
		\begin{align*}
			& x_n \xrightarrow{n \to \infty} x \ \ \text{ in } D_\phi
		\end{align*}
		and
		\begin{align*}
			& n(X_n - x_n) \xrightarrow{n \to \infty}_{\mathcal{D}} G
		\end{align*}
		for some tight process $G$ taking values in $D_0$. Then the weak convergence
		\begin{align*}
			& n(\phi(X_n) - \phi(x_n)) \xrightarrow{n \to \infty}_{\mathcal{D}} \phi'_x(G)
		\end{align*}
		holds.
	\end{lemma}

	\begin{lemma}[Hadamard differentiability of the map $f \mapsto f^{-1}$]\label{Lemma_Hadamard}\
		\\
		Let $W \subset \C$ be an open set and let $K \subset \C$ be a compact set with non-empty interior. Let $D_0=D=H^{\infty}(K)$ and $E=H^{\infty}(W)$, where $H^{\infty}(W)$ denotes the Banach space of all holomorphic and bounded functions $f : W \rightarrow \C$ equipped with the supremum norm. Further, let
		\begin{align}
			& D_\phi = \big\{f \in H^{\infty}(K) \ \big| \ f \text{ is injective with } W \subset f(K) \text{ and } \forall \tilde{v} \in K :\, f'(\tilde{v}) \neq 0 \big\}
		\end{align}
		such that one may define the inversion map
		\begin{align*}
			& \phi : D_\phi \rightarrow E \ \ ; \ \ f \mapsto f^{-1}|_W \ .
		\end{align*}
		For any $f \in D_\phi$, the map $\phi$ is Hadamard-differentiable at $f \in D_\phi$ tangentially to the set $D_0$ with
		\begin{align}\label{Eq_HadamardLemma_phiBarDef}
			& \phi_f'(h)(z) = \frac{h(f^{-1}(z))}{f'(f^{-1}(z))} \ .
		\end{align}
	\end{lemma}

	\subsection{Proof of Theorem~\ref{Thm_CLT_Inversion}}\label{Proof_Thm_CLT_Inversion}
	Under the additional assumption that there exists some $K>\sigma^2(1+\sqrt{c_\infty})^2$ such that
	\begin{align}\label{Eq_StrongerInversion_AdditionalConditionK}
		& \forall n \in \N : \supp(\hat{\nu}_n) \cup \supp(\nu_n) \subset [0,K] \ ,
	\end{align}
	it will be shown here that for any compact subset $J \subset U \cup (\R\setminus I)$,
	where
	\begin{align}\label{Eq_StrongerInversion_Def_I}
		& I \coloneq \big[\varphi_{\nu_\infty,c_\infty}^{-1}(-\delta), \varphi_{\nu_\infty,c_\infty}^{-1}(K+\delta) \big] \subset \R
	\end{align}
	for some $\delta>0$, the uniform convergence in distribution
	\begin{align}\label{Eq_CLTInversion_Result_Stronger}
		& n\big( \hat{s}_n^{(0)}(z) - \cs_{H_n}(z) \big)_{z \in J\setminus\R} \xrightarrow{n \to \infty}_{\mathcal{D}} \Big(-\frac{\tilde{G}(\Phi_{H_\infty,c_\infty}(z))}{c_\infty z^2\cs_{\ul{\nu}_\infty}'(\Phi_{H_\infty,c_\infty}(z))}\Big)_{z \in J\setminus\R} \ \text{ in } \ell^{\infty}(J\setminus\R)
	\end{align}
	holds.
	The interval $I$ is seen to be well-defined by (b) of Lemma~\ref{Lemma_StandardBounds} and the definition
	\begin{align*}
		& \varphi_{\nu_\infty,c_\infty}(\tz) \overset{\text{(\ref{Eq_Def_varphi})}}{=} \frac{-1}{\cs_{\ul{\nu}_\infty}(\tz)} \ .
	\end{align*}
	To prove the original version of Theorem~\ref{Thm_CLT_Inversion}, simply replace both $[0,K]$ and $I$ in the following proof with $\R$.
	\begin{itemize}
		\item[i)] \textit{Preliminaries}:\\
		\begin{itemize}
			\item Together with~\ref{EI_ItemAssumption_Asymptotics}--\ref{EI_ItemAssumption_sigmaBound}, the assumed uniform boundedness of the fourth moments of $(\bm{X}_n)_{j,k}$ is by \cite{MP_RowCorrelation_Bai} and Lemma \ref{Lemma_FourthMomentMP_law} sufficient for the Marchenko--Pastur law to hold, i.e.
			\begin{align}\label{Eq_CLTInversion_MP}
				& 1 = \bP\big( \hat{\nu}_n \xRightarrow{n \to \infty} \nu_\infty \big) \ .
			\end{align}
			
			\item Since $U\subset \bD_{H_\infty,c_\infty}(1)$ holds by construction, $J$ is a compact subset of $\bD_{H_\infty,c_\infty}(1) \cup (\R\setminus I)$ and the entire proof of Lemma~\ref{Lemma_ConsistencyBasic} (see Subsection~\ref{Proof_Lemma_ConsistencyBasic}) is applicable. In part (i) of said proof, it is already shown that $\Phi_{H_\infty,c_\infty}(J)$ is a compact subset of $\C\setminus[0,K]$, leading to
			\begin{align}\label{Eq_CLTInversion_Step1_1}
				& \dist\big( \Phi_{H_\infty,c_\infty}(J), [0,K] \big) > 2\tau
			\end{align}
			for some $\tau>0$.
			By restricting consideration to a smaller subset of $U$, which still contains $J$, one may, without loss of generality, assume that $U$ is bounded and satisfies
			\begin{align}\label{Eq_CLTInversion_UAssumption2}
				& \dist\big( \tilde{U}, [0,K] \big) > \tau \ .
			\end{align}
			
			\item In complete analogy to part (ii) from the proof of Lemma~\ref{Lemma_ConsistencyBasic} (see Subsection~\ref{Proof_Lemma_ConsistencyBasic}), with $\mathrm{closure}(\tilde{U})$ instead of $V_\varepsilon$, one by (\ref{Eq_CLTInversion_MP}) gets
			\begin{align}\label{Eq_CLTInversion_Step1_sApprox}
				& 1 = \bP\big( \sup\limits_{\tz \in \tilde{U}} |\cs_{\hat{{\nu}}_n}(\tz)-\cs_{{\nu}_\infty}(\tz)| \xrightarrow{n \to \infty} 0 \big) \ ,
			\end{align}
			which is by (\ref{Eq_Def_ulNu}) equivalent to
			\begin{align}\label{Eq_CLTInversion_Step1_ulsApprox}
				& 1 = \bP\big( \sup\limits_{\tz \in \tilde{U}} |\cs_{\hat{\ul{\nu}}_n}(\tz)-\cs_{\ul{\nu}_\infty}(\tz)| \xrightarrow{n \to \infty} 0 \big) \ .
			\end{align}
			Similarly, the previously stated assumption $K>\sigma^2(1+\sqrt{c_\infty})^2 = \lim\limits_{n \to \infty} \sigma^2(1+\sqrt{c_n})^2$ with Lemma~\ref{Lemma_BasicStieltjesConvergence} yields
			\begin{align}\label{Eq_CLTInversion_Step1_5}
				& \sup\limits_{\tz \in \tilde{U}} |\cs_{\nu_n}(\tz)-\cs_{\nu_\infty}(\tz)| \xrightarrow{n \to \infty} 0 \ .
			\end{align}
			
			\item The Lipschitz continuity of $\Phi_{H_\infty,c_\infty}$, shown in the proof of (d) of Lemma~\ref{Lemma_StandardBounds}, implies that the boundedness of $U$ is inherited by the set $\tilde{U} = \Phi_{H_\infty,c_\infty}(U)$, i.e.
			\begin{align}\label{Eq_CLTInversion_UAssumption3}
				& \forall \tz \in \tilde{U} : \ |\tz| \leq C
			\end{align}
			holds for some $C>0$.
			Property (\ref{Eq_CLTInversion_UAssumption2}) and result (b) of Lemma~\ref{Lemma_StandardBounds} then yield
			\begin{align}
				& \forall \tz \in \tilde{U} : \ |\cs_{\hat{\ul{\nu}}_n}(\tz)| \geq \int_\R \frac{\tau/2}{|\lambda-\tz|^2} \, d\hat{\ul{\nu}}_n(\lambda) \overset{\text{(\ref{Eq_CLTInversion_UAssumption3})}}{\geq} \int_\R \frac{\tau/2}{(\lambda+C)^2} \, d\hat{\ul{\nu}}_n(\lambda) \label{Eq_CLTInversion_Step1_6_1}\\
				& \forall \tz \in \tilde{U} : \ |\cs_{{\ul{\nu}}_n}(\tz)| \geq \int_\R \frac{\tau/2}{|\lambda-\tz|^2} \, d{\ul{\nu}}_n(\lambda) \overset{\text{(\ref{Eq_CLTInversion_UAssumption3})}}{\geq} \int_\R \frac{\tau/2}{(\lambda+C)^2} \, d{\ul{\nu}}_n(\lambda) \label{Eq_CLTInversion_Step1_6_2}
			\end{align}
			and analogously (d) of Lemma~\ref{Lemma_BasicStieltjesConvergence} and (b) of Lemma~\ref{Lemma_StandardBounds} together show
			\begin{align}\label{Eq_CLTInversion_Step1_7}
				& \forall \tz \in \tilde{U} : \ |\cs_{{\ul{\nu}}_\infty}(\tz)| \geq \underbrace{\int_\R \frac{\tau/2}{(\lambda+C)^2} \, d{\ul{\nu}}_\infty(\lambda)}_{>0} \eqcolon  2\delta' \ .
			\end{align}
			
			\item Since the right-hand side of (\ref{Eq_CLTInversion_Step1_6_1}) by (\ref{Eq_CLTInversion_MP}) almost surely converges to the right-hand side of (\ref{Eq_CLTInversion_Step1_7}), one has
			\begin{align}\label{Eq_CLTInversion_Step1_8}
				& 1 = \bP\big( \exists \hat{N}'>0 \, \forall n \geq \hat{N}' \, \forall \tz \in \tilde{U} : \ |\cs_{\hat{\ul{\nu}}_n}(\tz)| \geq \delta' \big) \ .
			\end{align}
			Analogously, the right-hand side of (\ref{Eq_CLTInversion_Step1_6_2}) by Lemma~\ref{Lemma_BasicStieltjesConvergence} also converges to the right-hand side of (\ref{Eq_CLTInversion_Step1_7}), which yields
			\begin{align}\label{Eq_CLTInversion_Step1_8_2}
				& \exists N'>0 \, \forall n \geq N' \, \forall \tz \in \tilde{U} : \ |\cs_{\ul{\nu}_n}(\tz)| \geq \delta' \ .
			\end{align}
			One may, for the sake of this proof, assume
			\begin{align}\label{Eq_CLTInversion_Step1_8_3}
				& \forall \tz \in \tilde{U} : \ |\cs_{\hat{\ul{\nu}}_n}(\tz)| \geq \delta' \ \ \text{ and } \ \ |\cs_{\ul{\nu}_n}(\tz)| \geq \delta' \ .
			\end{align}
			
			\item By (\ref{Eq_CLTInversion_UAssumption2}) together with (\ref{Eq_CLTInversion_Step1_7}) and (\ref{Eq_CLTInversion_Step1_8_3}), the maps
			\begin{align}\label{Eq_CLTInversion_Step1_varphi_hol}
				& \varphi_{\nu_n,c_n}, \varphi_{\hat{\nu}_n,c_n} \text{ and } \varphi_{\nu_\infty,c_\infty} \text{ are holomorphic on } \tilde{U} \ .
			\end{align}
			It then follows from (d) of Lemma~\ref{Lemma_StandardBounds} that
			\begin{align}\label{Eq_CLTInversion_Step1_varphiInf_DerivBound}
				& \forall \tz \in \tilde{U} : \ |\varphi_{\nu_n,c_n}'(\tz)| \geq \frac{1}{2} \ \ \text{ and } \ \ |\varphi_{\nu_\infty,c_\infty}'(\tz)| \geq \frac{1}{2} \ .
			\end{align}
			
			\item One may also combine (\ref{Eq_CLTInversion_Step1_ulsApprox}) and (\ref{Eq_CLTInversion_Step1_8}) for the result
			\begin{align}\label{Eq_CLTInversion_Step1_9}
				& 1 = \bP\big( \sup\limits_{\tz \in \tilde{U}} |\varphi_{\hat{\nu}_n,c_n}(\tz)-\varphi_{\nu_\infty,c_\infty}(\tz)| \xrightarrow{n \to \infty} 0 \big)
			\end{align}
			and likewise combine (\ref{Eq_CLTInversion_Step1_5}) and (\ref{Eq_CLTInversion_Step1_8_2}) for
			\begin{align}\label{Eq_CLTInversion_Step1_10}
				& \sup\limits_{\tz \in \tilde{U}} |\varphi_{\nu_n,c_n}(\tz)-\varphi_{\nu_\infty,c_\infty}(\tz)| \xrightarrow{n \to \infty} 0 \ .
			\end{align}
			
			\item By restricting consideration to an open subset $\tilde{U}_{\text{new}} \subset \tilde{U}$ which also contains $J$ and which satisfies $\mathrm{closure}(\tilde{U}_{\text{new}}) \subset \tilde{U}$, one may by (\ref{Eq_CLTInversion_Step1_9}), and the fact that the maps $\varphi_{\hat{\nu}_n,c_n}$ and $\varphi_{\nu_\infty,c_\infty}$ are holomorphic on $\C\setminus[0,K]$, without loss of generality assume
			\begin{align}\label{Eq_CLTInversion_Step1_9_Deriv}
				& 1 = \bP\big( \sup\limits_{\tz \in \tilde{U}} |\varphi_{\hat{\nu}_n,c_n}'(\tz)-\varphi_{\nu_\infty,c_\infty}'(\tz)| \xrightarrow{n \to \infty} 0 \big) \ .
			\end{align}
			By (\ref{Eq_CLTInversion_Step1_varphiInf_DerivBound}), one may for the sake of this proof then assume
			\begin{align}\label{Eq_CLTInversion_Step1_varphiN_DerivBound}
				& \forall \tz \in \tilde{U} : \ |\varphi_{\hat{\nu}_n,c_n}'(\tz)| \geq \frac{1}{4} \ .
			\end{align}
			
			\item Since $J$ is a compact subset of the open set $\tilde{U}$, it must hold that
			\begin{align}\label{Eq_CLTInversion_Step1_11}
				& \delta'' \coloneq \dist\big( \Phi_{H_\infty,c_\infty}(J) , \C \setminus \tilde{U} \big) > 0
			\end{align}
			and properties (\ref{Eq_CLTInversion_Step1_9}), (\ref{Eq_CLTInversion_Step1_10}) yield
			\begin{align*}
				& 1 = \bP\big( \exists \hat{N} > 0 \, \forall n \geq \hat{N} : \  \overbrace{\varphi_{\nu_\infty,c_\infty}\big( \Phi_{H_\infty,c_\infty}(J)}^{= J \text{ by Lemma~\ref{Lemma_SpaceTransform}}} \subset \varphi_{\hat{\nu}_n,c_n}\big(\tilde{U}) \big)\\
				& \exists N>0 \, \forall n \geq N : \  \underbrace{\varphi_{\nu_\infty,c_\infty}\big( \Phi_{H_\infty,c_\infty}(J) \big)}_{= J \text{ by Lemma~\ref{Lemma_SpaceTransform}}} \subset \varphi_{\nu_n,c_n}(\tilde{U}) \ .
			\end{align*}
			Once more, for the sake of this proof, assume
			\begin{align}\label{Eq_CLTInversion_Step1_12}
				& J \subset \varphi_{\nu_n,c_n}(\tilde{U}) \ \ \text{ and } \ \ J \subset \varphi_{\hat{\nu}_n,c_n}(\tilde{U}) \ .
			\end{align}
		\end{itemize}

		\item[ii)] \textit{Convergence of $n(\varphi_{\hat{\nu}_n,c_n}-\varphi_{\nu_n,c_n})$ by standard delta method}:\\
		Using Lemma~\ref{Lemma_FunctionalDeltaMethod_General}, it will be shown here that
		\begin{align}\label{Eq_FuncDelta_Step2_1}
			& n \big( \varphi_{\hat{\nu}_n,c_n}(\tz) - \varphi_{\nu_n,c_n}(\tz) \big)_{\tz \in \tilde{U}} \xrightarrow{n \to \infty}_{\mathcal{D}} \Big( \frac{\tilde{G}(\tz)}{\cs_{\ul{\nu}_\infty}(\tz)^2} \Big)_{\tz \in \tilde{U}} \ \ \text{ in } \ell^{\infty}(\tilde{U}) \ .
		\end{align}
		Let $D_0=D=E=\ell^\infty(\tilde{U})$ and $D_\phi = \{f \in D \mid \inf\limits_{\tz \in \tilde{U}} |f(\tz)|>\delta'\}$. The maps
		\begin{align*}
			& x = \cs_{\ul{\nu}_\infty} \ \ , \ \ X_n = \cs_{\hat{\ul{\nu}}_n} \ \ \text{and} \ \ x_n=\cs_{\ul{\nu}_n}
		\end{align*}
		may by (\ref{Eq_CLTInversion_Step1_7}) and (\ref{Eq_CLTInversion_Step1_8_3}) be assumed to lie in $D_\phi \subset D_0$. Define
		\begin{align*}
			& \phi : D_\phi \rightarrow E \ \ ; \ \ f \mapsto \frac{-1}{f}
		\end{align*}
		and observe that it is Hadamard-differentiable at $x$ tangentially to $D_0$ with $\phi'_x(h) = \frac{h}{x^2}$ for any sequences $(t_m)_{m \in\N} \subset (0,\infty)$ and $(h_m)_{m \in \N} \subset D$ with $t_m \to 0$, $h_m \to h \in D_0$ and $x+t_mh_m \in D_\phi$. The fact that $x$ and $x+t_mh_m$ are both in $D_\phi$ is in this case sufficient for
		\begin{align*}
			& \frac{\phi(x+t_mh_m) - \phi(x)}{t_m} = \frac{\frac{1}{x} - \frac{1}{x+t_mh_m}}{t_m} = \frac{h_m}{x(x+t_mh_m)} \xrightarrow{m \to \infty} \frac{h}{x^2} \ .
		\end{align*}
		Thus, Lemma~\ref{Lemma_FunctionalDeltaMethod_General} yields (\ref{Eq_FuncDelta_Step2_1}).
		
		\item[iii)] \textit{Functional delta method}:\\
		Set $D_0=D=H^{\infty}(J)$ and $E=H^{\infty}(\tilde{U})$, where $H^{\infty}(\tilde{U})$ denotes the Banach space of all holomorphic and bounded functions $f : \tilde{U} \rightarrow \C$ equipped with the supremum norm. Further, for
		\begin{align*}
			& D_\phi = \big\{f \in H^{\infty}(\tilde{U}) \ \big| \ f \text{ is injective with } J \subset f(\tilde{U}) \text{ and } \forall \tz \in \tilde{U} :\, f'(\tz) \neq 0 \big\}
		\end{align*}
		one may define the inversion map
		\begin{align*}
			& \phi : D_\phi \rightarrow E \ \ ; \ \ f \mapsto f^{-1}|_{J} \ .
		\end{align*}
		Define the functions
		\begin{align*}
			& x = \varphi_{\nu_\infty,c_\infty} \ \ , \ \ X_n = \varphi_{\hat{\nu}_n,c_n} \ \ \text{and} \ \ x_n = \varphi_{\nu_n,c_n} \ ,
		\end{align*}
		which by (\ref{Eq_CLTInversion_Step1_varphiInf_DerivBound}) and (\ref{Eq_CLTInversion_Step1_varphiN_DerivBound}) have non-zero derivative on $\tilde{U}$. The inclusion $x \in D_\phi$ is guaranteed by Lemma~\ref{Lemma_SpaceTransform} and construction of $\tilde{U}=\Phi_{H_\infty,c_\infty}(U)$. The inclusions $X_n,x_n \in D_\phi$ hold by (\ref{Eq_CLTInversion_Step1_12}).
		Part (ii) of this proof has already proved the weak convergence
		\begin{align*}
			& n(X_n-x_n) \xrightarrow{n \to \infty}_{\mathcal{D}} \Big( \frac{\tilde{G}(\tz)}{\cs_{\ul{\nu}_\infty}(\tz)^2} \Big)_{\tz \in \tilde{U}} \ \ \text{ in } \ell^{\infty}(\tilde{U})
		\end{align*}
		and Lemma~\ref{Lemma_Hadamard} for $W_{\text{Lemma~\ref{Lemma_Hadamard}}}=\tilde{U}$ and $K_{\text{Lemma~\ref{Lemma_Hadamard}}}=J$ proves Hadamard-differentiabi\-lity of $\phi$ at any $f \in D_\phi$ tangentially to $D_0$ with derivative $\phi_f'(h) = \frac{h \circ f^{-1}}{f' \circ f^{-1}}$.
		Lemma~\ref{Lemma_FunctionalDeltaMethod_General} is thus applicable and yields
		\begin{align}\label{Eq_FuncDelta_Step3_1}
			& n\big( \varphi_{\hat{\nu}_n,c_n}^{-1}(z) - \Phi_{H_n,c_n}(z) \big)_{z \in J} = n\big( \phi(X_n) - \phi(x_n) \big) \nonumber\\
			& \xrightarrow{\substack{\text{Lemma~\ref{Lemma_FunctionalDeltaMethod_General}} \\ n \to \infty}}_{\mathcal{D}} \Bigg( \frac{\frac{\tilde{G}(\Phi_{H_\infty,c_\infty}(z))}{\cs_{\ul{\nu}_\infty}(\Phi_{H_\infty,c_\infty}(z))^2}}{\varphi_{\nu_\infty,c_\infty}'(\Phi_{H_\infty,c_\infty}(z))} \Bigg)_{z \in J} = \Big( \frac{\tilde{G}(\Phi_{H_\infty,c_\infty}(z))}{\cs_{\ul{\nu}_\infty}'(\Phi_{H_\infty,c_\infty}(z))} \Big)_{z \in J} \ ,
		\end{align}
		where the convergence in distribution is as $\ell^{\infty}(J)$-valued random variables and the identities $\varphi_{\nu_n,c_n}^{-1} = \Phi_{H_n,c_n}$ and $\varphi_{\nu_\infty,c_\infty}^{-1} = \Phi_{H_\infty,c_\infty}$, known by Lemma~\ref{Lemma_SpaceTransform}, were used.
		\\[0.5em]
		By Lemma~\ref{Lemma_ConsistencyBasic}, there almost surely exists an $\hat{N}_J>0$ such that
		\begin{align*}
			& \forall n \geq \hat{N}_J \, \forall \tz \in J : \ \hat{s}_n \text{ exists} \ .
		\end{align*}
		Consequently, for all $n \geq \hat{N}_J$ the map
		\begin{align*}
			& \hat{\Phi}_{n,c_n} : J \rightarrow \C \ \ ; \ \ z \mapsto (1-c_nz\hat{s}_n(z)-c_n)z
		\end{align*}
		exists and by (\ref{Eq_RevMPE_Equivalence1}) satisfies
		\begin{align}\label{Eq_FuncDelta_Step3_2}
			& \hat{\Phi}_{n,c_n} = \varphi_{\hat{\nu}_n,c_n}^{-1}
		\end{align}
		on $J$. For the map $\hat{\Phi}_{n,c_n}^{(0)}(z) = (1-c_nz\hat{s}_n^{(0)}(z)-c_n)z$ it then follows that
		\begin{align}\label{Eq_FuncDelta_Step3_3}
			& n\big( \hat{\Phi}_{n,c_n}^{(0)}(z) - \Phi_{H_n,c_n}(z) \big)_{z \in J} \xrightarrow{n \to \infty}_\mathcal{D} \Big( \frac{\tilde{G}(\Phi_{H_\infty,c_\infty}(z))}{\cs_{\ul{\nu}_\infty}'(\Phi_{H_\infty,c_\infty}(z))} \Big)_{z \in J} \ \ \text{in } \ell^{\infty}(J) \ .
		\end{align}
		
		\item[iv)] \textit{Standard delta method (second application)}:\\
		Another application of Lemma~\ref{Lemma_FunctionalDeltaMethod_General} will convert the result of (iii) into the wanted (\ref{Eq_CLTInversion_Result_Stronger}). Let $D_0=D=D_\phi=E=\ell^\infty(J)$ and define
		\begin{align*}
			& \phi : D_\phi \rightarrow E \ \ ; \ \ f \mapsto \frac{-f}{\Id^2} \ .
		\end{align*}
		such that $\phi(\hat{\Phi}_n^{(0)}) = c_n\hat{s}_n^{(0)} + \frac{c_n-1}{z}$ and $\phi(\Phi_{H_n,c_n})=c_n\cs_{H_n} + \frac{c_n-1}{z}$.
		One verifies that $\phi$ is Hadamard-differentiable at every $f \in D_\phi$ tangentially to $D_0$ with $\phi'_f(h)(z) = -\frac{h}{z^2}$ by the calculation
		\begin{align*}
			& \frac{\phi(f+t_mh_m) - \phi(f)}{t_m} = -\frac{h_m}{\Id^2} \xrightarrow{m \to \infty} -\frac{h}{\Id^2} = \phi'_f(h)
		\end{align*}
		for any sequences $(t_m)_{m \in\N} \subset (0,\infty)$ and $(h_m)_{m \in \N} \subset D$ with $t_m \to 0$, $h_m \to h \in D_0$ and $f+t_mh_m \in D_\phi$. The convergence holds uniformly in $z$ by the fact that $J$ is bounded away from zero and $\|h_m-h\|_\infty \to 0$. Lemma~\ref{Lemma_FunctionalDeltaMethod_General} is thus applicable for $x_n = \Phi_{H_n,c_n}$, $x = \Phi_{H_\infty,c_\infty}$ and $X_n = \hat{\Phi}_{n}^{(0)}$. It follows that
		\begin{align*}
			& n\big( c_n\hat{s}_{n}^{(0)}(z) - c_n\cs_{H_n}(z) \big)_{z \in J} \xrightarrow{n \to \infty}_{\mathcal{D}} \Bigg(-\frac{\frac{\tilde{G}(\Phi_{H_\infty,c_\infty}(z))}{\cs_{\ul{\nu}_\infty}'(\Phi_{H_\infty,c_\infty}(z))}}{z^2}\Bigg)_{z \in J}\\
			& \hspace{6cm} = \Big(-\frac{\tilde{G}(\Phi_{H_\infty,c_\infty}(z))}{z^2\cs_{\ul{\nu}_\infty}'(\Phi_{H_\infty,c_\infty}(z))}\Big)_{z \in J}
		\end{align*}
		as $\ell^{\infty}(J)$-valued random variables. A final application of Slutsky's lemma with $c_n \to c_\infty$ turns the above convergence into (\ref{Eq_CLTInversion_Result}). \qed
	\end{itemize}
	
	\begin{corollary}[CLT in the Gaussian case]\label{Cor_GaussianCLT}\
		\\
		Suppose Assumption~\ref{Assumption_EigInf_Main} holds, and that the entries of $\bm{X}_n$ are additionally i.i.d. real ($\beta=1$) or complex ($\beta=2$) standard normal.
		For any compact set $J \subset \bD_{H_\infty,c_\infty}(1)$, the convergence in distribution
		\begin{align}\label{Eq_GaussianCLT_Result}
			& n\big( \hat{s}_{n}^{(0)}(z) - \cs_{H_n}(z) \big)_{z \in J} \xrightarrow{n \to \infty}_{\mathcal{D}} G(z)_{z \in J} \ \ \text{in } \ell^\infty(J)
		\end{align}
		holds for a Gaussian process $G(z)_{z \in J}$ with mean
		\begin{align}\label{Eq_GaussianCLT_Result_Mean}
			& \E\big[ G(z) \big] = - \mathbbm{1}_{\beta=1} \, \frac{\Phi_\infty''(z)}{2c_\infty \Phi_\infty'(z)}
		\end{align}
		and covariance
		\begin{align}\label{Eq_GaussianCLT_Result_Cov}
			& \Cov\big[ G(z_1), G(z_2) \big] = \frac{1}{\beta c_\infty^2} \Big(\frac{1}{(z_1-z_2)^2} - \frac{\Phi_{H_\infty,c_\infty}'(z_1) \Phi_{H_\infty,c_\infty}'(z_2)}{(\Phi_{H_\infty,c_\infty}(z_1)-\Phi_{H_\infty,c_\infty}(z_2))^2}\Big)
		\end{align}
		as well as
		\begin{align*}
			& \Cov\big[ G(z_1), \ol{G(z_2)} \big] = \Cov\big[ G(z_1), G(\ol{z_2}) \big] \ .
		\end{align*}
		As in Theorem~\ref{Thm_CLT_Inversion}, $\hat{s}_n^{(0)}(z)$ is defined as $\hat{s}_n(z)$, whenever the latter exists, and zero otherwise.
	\end{corollary}
	\begin{proof}\
		\\
		For ease of notation write $\Phi_{\infty} = \Phi_{H_\infty,c_\infty}$. By Lemma~\ref{Lemma_LargestEVBound} and (b) of Lemma~\ref{Lemma_BasicStieltjesConvergence}, it will for sufficiently large $n$ hold that
		\begin{align}\label{Eq_StrongerGaussCLT_AdditionalConditionK}
			& \supp(\hat{\nu}_n) \cup \supp(\nu_n) \subset [0,\underbrace{\sigma^2(1+\sqrt{c_\infty})^2+\tau}_{\eqcolon K}] \ .
		\end{align}
		It will be shown here that for any compact subset $J \subset U \cup (\R\setminus I)$,
		where
		\begin{align}\label{Eq_StrongerGaussCLT_Def_I}
			& I \coloneq \big[\varphi_{\nu_\infty,c_\infty}^{-1}(-\delta), \varphi_{\nu_\infty,c_\infty}^{-1}(K+\delta) \big] \subset \R
		\end{align}
		for some $\delta>0$, the uniform convergence in distribution
		\begin{align}\label{Eq_GaussCLT_Result_Stronger}
			& n\big( \hat{s}_{n}^{(0)}(z) - \cs_{H_n}(z) \big)_{z \in J\setminus\R} \xrightarrow{n \to \infty}_{\mathcal{D}} G(z)_{z \in J\setminus\R} \ \ \text{in } \ell^\infty(J\setminus\R)
		\end{align}
		holds.
		The interval $I$ is seen to be well-defined by (b) of Lemma~\ref{Lemma_StandardBounds} as well as the definition $\varphi_{\nu_\infty,c_\infty}(\tz) = \frac{-1}{\cs_{\ul{\nu}_\infty}(\tz)}$. The above result is stronger than the statement of Corollary~\ref{Cor_GaussianCLT}.
		\begin{itemize}
			\item[i)] \textit{Preliminaries}:\\
			Together with~\ref{EI_ItemAssumption_Asymptotics}--\ref{EI_ItemAssumption_sigmaBound}, the assumed uniform boundedness of the fourth moments of $(\bm{X}_n)_{j,k}$ is by \cite{MP_RowCorrelation_Bai} and Lemma \ref{Lemma_FourthMomentMP_law} sufficient for the Marchenko--Pastur law to hold, i.e.
			\begin{align}\label{Eq_GaussCLT_MP}
				& 1 = \bP\big( \hat{\nu}_n \xRightarrow{n \to \infty} \nu_\infty \big) \ .
			\end{align}
			Since $U\subset \bD_{H_\infty,c_\infty}(1)$ holds by construction, $J$ is a compact subset of $\bD_{H_\infty,c_\infty}(1) \cup (\R\setminus I)$ and the entire proof of Lemma~\ref{Lemma_ConsistencyBasic} (see Subsection~\ref{Proof_Lemma_ConsistencyBasic}) is applicable. In part (i) of said proof, it is already shown that $\Phi_{H_\infty,c_\infty}(J)$ is a compact subset of $\C\setminus[0,K]$, leading to
			\begin{align*}
				& \dist\big( \Phi_\infty(J), [0,K] \big) > 4\tau
			\end{align*}
			for some small $\tau>0$.
			Let $\tilde{U}$ be a bounded open set containing $\Phi_{\infty}(J)$ such that
			\begin{align}\label{Eq_GaussianCLT_Step2_01}
				& \dist\big( \tilde{U}, [0,K] \big) > 3\tau \ .
			\end{align}
			By (\ref{Eq_StrongerGaussCLT_AdditionalConditionK}), (b) of Lemma~\ref{Lemma_StandardBounds} and (d) of Lemma~\ref{Lemma_BasicStieltjesConvergence} it also holds that
			\begin{align}\label{Eq_GaussianCLT_Step2_03}
				& \supp(\ul{\nu}_\infty) \subset [0,K] \ .
			\end{align}
			Define the contour $\mathcal{C}$ as on page 561 of \cite{BaiCLT} with parameters
			\begin{align}\label{Eq_GaussianCLT_Step2_04}
				& (x_l,v_0,x_r)_{\text{\cite{BaiCLT}}} = \big(-\tau,\tau,\underbrace{\sigma^2(1+\sqrt{c_\infty})^2+2\tau}_{=K+\tau}\big) \ .
			\end{align}
			Not only will $\mathcal{C}$ by (\ref{Eq_GaussianCLT_Step2_01}) and (\ref{Eq_GaussianCLT_Step2_03}) separate the support $\supp(\nu_\infty)$ from $\tilde{U}$, it will by construction also hold that
			\begin{align}\label{Eq_GaussianCLT_Step2_05}
				& \dist(\mathcal{C} \cup \ol{\mathcal{C}}, [0,K]) \geq \tau
			\end{align}
			and
			\begin{align}\label{Eq_GaussianCLT_Step2_06}
				& \dist(\mathcal{C} \cup \ol{\mathcal{C}}, \tilde{U}) \geq \tau \ .
			\end{align}
			The last property also yields
			\begin{align}\label{Eq_GaussianCLT_Step2_07}
				& \Big|\frac{1}{\tz-\tilde{u}_1} - \frac{1}{\tz-\tilde{u}_2}\Big| = \frac{|\tilde{u}_1-\tilde{u}_2|}{|\tz-\tilde{u}_1| \, |\tz-\tilde{u}_2|} \leq \frac{|\tilde{u}_1-\tilde{u}_2|}{\tau^2}
			\end{align}
			for all $\tz \in \mathcal{C}$ and $\tilde{u}_1,\tilde{u}_2 \in \tilde{U}$.
			
			\item[ii)] \textit{Application of Lemma 1.1 from \cite{BaiCLT}}:\\
			Let $\mathcal{C}$ be the contour as constructed in (\ref{Eq_GaussianCLT_Step2_04}).
			On page 561 of \cite{BaiCLT}, Bai and Silverstein further define the $n$-dependent contour
			\begin{align}\label{Eq_GaussCLT_Step2_01}
				& \mathcal{C}_n = \{\tz \in \mathcal{C} \mid \Im(\tz) \geq \varepsilon_n/n\}
			\end{align}
			for some sequence $(\varepsilon_n)_{n \in \N} \in (0,\infty)$ with $\varepsilon_n \searrow 0$, and also define the process $\hat{M}_n(\tz)_{\tz \in \mathcal{C}}$ in such a way that it satisfies
			\begin{align}\label{Eq_GaussCLT_Step2_02}
				& \tz \in \mathcal{C}_n : \ \hat{M}_n(\tz) = n(\cs_{\hat{\ul{\nu}}_n}(\tz) - \cs_{\ul{\nu}_n}(\tz))
			\end{align}
			and
			\begin{align}\label{Eq_GaussCLT_Step2_03}
				& \sup\limits_{\tz \in \mathcal{C}} |\hat{M}_n(\tz)| = \sup\limits_{\tz \in \mathcal{C}_n} |\hat{M}_n(\tz)| \overset{\text{(\ref{Eq_GaussCLT_Step2_02})}}{=} \sup\limits_{\tz \in \mathcal{C}_n} |n(\cs_{\hat{\ul{\nu}}_n}(\tz) - \cs_{\ul{\nu}_n}(\tz))| \ .
			\end{align}
			Lemma 1.1 of \cite{BaiCLT} proves the convergence
			\begin{align}\label{Eq_GaussCLT_Step2_Lemma11}
				& \hat{M}_n(\tz)_{\tz \in \mathcal{C}} \xrightarrow{n \to \infty}_{\mathcal{D}} \tilde{G}(\tz)_{\tz \in \mathcal{C}} \ \ \text{in } \ell^\infty(\mathcal{C})
			\end{align}
			for a Gaussian process $\tilde{G}(\tz)_{\tz \in \mathcal{C}}$ whose mean and covariance structure of $\tilde{G}$ is given by
			\begin{align}\label{Eq_GaussianCLT_BaiMean}
				& \E\big[ \tilde{G}(\tz) \big] = \mathbbm{1}_{\beta=1} \frac{c_\infty\int_\R \frac{\cs_{\ul{\nu}_\infty}(\tz)^3 \lambda^2}{(1+\lambda\cs_{\ul{\nu}_\infty}(\tz))^3} \, dH_\infty(\lambda)}{\big(1-c_\infty\int_\R \frac{\cs_{\ul{\nu}_\infty}(\tz)^2 \lambda^2}{(1+\lambda\cs_{\ul{\nu}_\infty}(\tz))^2} \, dH_\infty(\lambda)\big)^2}
			\end{align}
			and
			\begin{align}\label{Eq_GaussianCLT_BaiCov}
				& \Cov\big[ \tilde{G}(\tz), \tilde{G}(\tz') \big] = \E\big[ \big(\tilde{G}(\tz)-\E[\tilde{G}(\tz)]\big) \big(\tilde{G}(\tz')-\E[\tilde{G}(\tz')]\big) \big] \nonumber\\
				& = \frac{1}{\beta} \Big(\frac{\cs_{\ul{\nu}_\infty}'(\tz) \cs_{\ul{\nu}_\infty}'(\tz')}{(\cs_{\ul{\nu}_\infty}(\tz) - \cs_{\ul{\nu}_\infty}(\tz'))^2} - \frac{1}{(\tz-\tz')^2}\Big)
			\end{align}
			as well as
			\begin{align}\label{Eq_GaussianCLT_BaiCov2}
				& \Cov\big[ \tilde{G}(\tz), \ol{\tilde{G}(\tz')} \big] = \Cov\big[ \tilde{G}(\tz), \tilde{G}(\ol{\tz'}) \big] \ .
			\end{align}
			
			\item[iii)] \textit{Cauchy's integral formula}:\\
			Since the contour $\mathcal{C}$ separates $[0,K]$ from $\tilde{U}$, it for all $\tilde{u} \in \tilde{U}$ by Cauchy's integral formula holds that
			\begin{align}\label{Eq_CLTInversion_Step3_01}
				n\big(\cs_{\hat{\ul{\nu}}_n}(\tilde{u})-\cs_{\ul{\nu}_n}(\tilde{u})\big) & = \frac{-1}{2\pi \bm{i}} \int_{\mathcal{C}} \frac{1}{\tz-\tilde{u}} \, n\big(\cs_{\hat{\ul{\nu}}_n}(\tz)-\cs_{\ul{\nu}_n}(\tz)\big) \, d\tz \nonumber\\
				& \hspace{0.5cm} + \frac{1}{2\pi \bm{i}} \int_{\mathcal{C}} \frac{1}{\ol{\tz}-\tilde{u}} \, n\big(\ol{\cs_{\hat{\ul{\nu}}_n}(\tz)-\cs_{\ul{\nu}_n}(\tz)}\big) \, d\tz
			\end{align}
			where the integral goes through $\mathcal{C}$ with a counter-clockwise orientation. Consequently, one may bound
			\begin{align}\label{Eq_CLTInversion_Step3_02}
				& \Big| \frac{-1}{2\pi \bm{i}} \int_{\mathcal{C}} \frac{1}{\tz-\tilde{u}} \, n\big(\cs_{\hat{\ul{\nu}}_n}(\tz)-\cs_{\ul{\nu}_n}(\tz)\big) \, d\tz - \frac{-1}{2\pi \bm{i}} \int_{\mathcal{C}} \frac{1}{\tz-\tilde{u}} \, \hat{M}_n(\tz) \, d\tz \Big| \nonumber\\
				& \overset{\text{(\ref{Eq_GaussCLT_Step2_02})}}{\leq} \frac{1}{2\pi} \int_{\mathcal{C}\setminus\mathcal{C}_n} \overbrace{\frac{1}{|\tz-\tilde{u}|}}^{\leq \frac{1}{\tau} \text{ by (\ref{Eq_GaussianCLT_Step2_06})}} \, \Big| n\big(\cs_{\hat{\ul{\nu}}_n}(\tz)-\cs_{\ul{\nu}_n}(\tz)\big) - \hat{M}_n(\tz)\Big| \, d\tz \nonumber\\
				& \overset{\text{(\ref{Eq_GaussCLT_Step2_01})}}{\leq} \frac{\varepsilon_n/n}{\pi \tau} \sup_{\tz \in \mathcal{C}\setminus\mathcal{C}_n} \Big| n\big(\cs_{\hat{\ul{\nu}}_n}(\tz)-\cs_{\ul{\nu}_n}(\tz)\big) - \hat{M}_n(\tz)\Big| \nonumber\\
				& \overset{\text{(\ref{Eq_GaussCLT_Step2_03})}}{\leq} \frac{\varepsilon_n/n}{\pi \tau} 2\sup_{\tz \in \mathcal{C}} \Big| n\big(\cs_{\hat{\ul{\nu}}_n}(\tz)-\cs_{\ul{\nu}_n}(\tz)\big)\Big| \nonumber\\
				& \leq \frac{2\varepsilon_n}{\pi \tau} \sup\limits_{\tz \in \mathcal{C}} \big( \underbrace{|\cs_{\hat{\ul{\nu}}_n}(\tz)| + |\cs_{\ul{\nu}_n}(\tz)|}_{\leq \frac{2}{\tau} \text{ by (\ref{Eq_GaussianCLT_Step2_05}) and (\ref{Eq_StandardBounds_c2})}} \big) \leq \frac{4\varepsilon_n}{\pi\tau^3} \xrightarrow{n \to \infty} 0
			\end{align}
			and in complete analogy
			\begin{align}\label{Eq_CLTInversion_Step3_03}
				& \Big| \frac{1}{2\pi \bm{i}} \int_{\mathcal{C}} \frac{1}{\ol{\tz}-\tilde{u}} \, n\big(\ol{\cs_{\hat{\ul{\nu}}_n}(\tz)-\cs_{\ul{\nu}_n}(\tz)}\big) \, d\tz - \frac{1}{2\pi \bm{i}} \int_{\mathcal{C}} \frac{1}{\ol{\tz}-\tilde{u}} \, \ol{\hat{M}_n(\tz)} \, d\tz \Big| \leq \frac{4\varepsilon_n}{\pi\tau^3} \xrightarrow{n \to \infty} 0 \ .
			\end{align}

			\item[iv)] \textit{Application of the continuous mapping theorem and Slutsky's lemma}:\\
			The map
			\begin{align*}
				F : \ell^\infty(\mathcal{C}) & \rightarrow \ell^\infty(\tilde{U})\\
				M & \mapsto \Big( \tilde{u} \mapsto \frac{-1}{2\pi \bm{i}} \int_{\mathcal{C}} \frac{1}{\tz-\tilde{u}} \, M(\tz) \, d\tz + \frac{1}{2\pi \bm{i}} \int_{\mathcal{C}} \frac{1}{\ol{\tz}-\tilde{u}} \, \ol{M(\tz)} \, d\tz\Big)
			\end{align*}
			is continuous by the bound
			\begin{align*}
				& \|F(M_1)-F(M_2)\|_{\ell^\infty(\tilde{U})} = \sup\limits_{\tilde{u} \in \tilde{U}} |F(M_1)(\tilde{u})-F(M_2)(\tilde{u})|\\
				& \leq \frac{1}{2\pi} \sup\limits_{\tilde{u} \in \tilde{U}} \int_{\mathcal{C}} \underbrace{\Big(\frac{1}{|\tz-\tilde{u}|} + \frac{1}{|\ol{\tz}-\tilde{u}|}\Big)}_{\leq \frac{2}{\tau} \text{ by (\ref{Eq_GaussianCLT_Step2_06})}} \, \big|M_1(\tz) - M_2(\tz)\big| \, |d\tz|\\
				& \leq \frac{1}{\pi \tau} \Big(\int_{\mathcal{C}} |d\tz|\Big) \, \sup\limits_{\tz \in \mathcal{C}} \big|M_1(\tz) - M_2(\tz)\big| = \frac{1}{\pi \tau} \Big(\int_{\mathcal{C}} |d\tz|\Big) \, \|M_1-M_2\|_{\ell^\infty(\mathcal{C})} \ .
			\end{align*}
			The continuous mapping theorem then (cf. Theorem 7.24 in \cite{FunctionalDeltaReference}) turns the convergence (\ref{Eq_GaussCLT_Step2_Lemma11}) into
			\begin{align}\label{Eq_CLTInversion_Step4_01}
				& F(\hat{M}_n) \xrightarrow{n \to \infty}_{\mathcal{D}} F(\tilde{G}) \ \ \text{ in } \ell^\infty(\tilde{U}) \ .
			\end{align}
			The bounds (\ref{Eq_CLTInversion_Step3_02}) and (\ref{Eq_CLTInversion_Step3_03}) then show
			\begin{align*}
				& \|F(\hat{M}_n) - F\big( n(\cs_{\hat{\ul{\nu}}_n}-\cs_{\ul{\nu}_n}) \big)\|_{\ell^\infty(\tilde{U})} \xrightarrow{n \to \infty} 0
			\end{align*}
			and Slutsky's lemma yields
			\begin{align}\label{Eq_CLTInversion_Step4_02}
				& \underbrace{F\big( n(\cs_{\hat{\ul{\nu}}_n}-\cs_{\ul{\nu}_n}) \big)}_{= n(\cs_{\hat{\ul{\nu}}_n}-\cs_{\ul{\nu}_n}) \text{ by (\ref{Eq_CLTInversion_Step3_01})}} \xrightarrow{n \to \infty}_{\mathcal{D}} F(\tilde{G}) \ \ \text{ in } \ell^\infty(\tilde{U}) \ .
			\end{align}
			Finally, $F(\tilde{G}|_{\mathcal{C}})=\tilde{G}|_{\tilde{U}}$ may either be calculated from the given mean and covariance structure (\ref{Eq_GaussianCLT_BaiMean})--(\ref{Eq_GaussianCLT_BaiCov2}), or may be seen to necessarily hold from the fact that Lemma 1.1 of \cite{BaiCLT} also holds for contours $\mathcal{C}'$ going through $\tilde{U}$. This proves the convergence
			\begin{align}\label{Eq_CLTInversion_Step4_03}
				& n\big(\cs_{\hat{\ul{\nu}}_n}(\tz)-\cs_{\ul{\nu}_n}(\tz)\big)_{\tz \in \tilde{U}} \xrightarrow{n \to \infty}_{\mathcal{D}} \tilde{G}(\tz)_{\tz \in \tilde{U}} \ \ \text{ in } \ell^\infty(\tilde{U}) \ .
			\end{align}
			
			\item[v)] \textit{Application of Theorem~\ref{Thm_CLT_Inversion}}:\\
			The conditions of Theorem~\ref{Thm_CLT_Inversion} follow from (\ref{Eq_CLTInversion_Step4_03}) and the conditions of this corollary, where the fourth moments of the entries of $\bm{X}_n$ are either all $3$ in the real case ($\beta=1$), or all $2$ in the complex case $\beta=2$. Since $U$ was constructed to satisfy $J \subset U\cup(\R\setminus I)$ and $\tilde{U}=\Phi_{\infty}(U)$, the slightly stronger version of Theorem~\ref{Thm_CLT_Inversion} as shown in Subsection~\ref{Proof_Thm_CLT_Inversion} directly yields
			\begin{align}\label{Eq_CLTInversion_Result_copy}
				& n\big( \hat{s}_n^{(0)}(z) - \cs_{H_n}(z) \big)_{z \in J\setminus\R} \xrightarrow{n \to \infty}_{\mathcal{D}} \Big(\underbrace{-\frac{\tilde{G}(\Phi_\infty(z))}{c_\infty z^2\cs_{\ul{\nu}_\infty}'(\Phi_\infty(z))}}_{\eqcolon  G(z)}\Big)_{z \in J\setminus\R} \hspace{-0.1cm} \text{ in } \ell^{\infty}(J\setminus\R) \ .
			\end{align}
			The mean and covariance structure of $G$ is by construction
			\begin{align}\label{Eq_GaussianCLT_MeanFirst}
				& \E\big[ G(z) \big] = -\frac{\E[\tilde{G}(\Phi_\infty(z))]}{c_\infty z^2\cs_{\ul{\nu}_\infty}'(\Phi_\infty(z))}
			\end{align}
			and
			\begin{align}\label{Eq_GaussianCLT_CovFirst}
				& \Cov\big[ G(z_1), G(z_2) \big] = \frac{\Cov[\tilde{G}(\Phi_\infty(z_1)), \tilde{G}(\Phi_\infty(z_2))]}{c_\infty^2z_1^2z_2^2\cs_{\ul{\nu}_\infty}'(\Phi_\infty(z_1)) \cs_{\ul{\nu}_\infty}'(\Phi_\infty(z_2))} \ .
			\end{align}
			It remains to show that this mean and covariance structure coincides with (\ref{Eq_GaussianCLT_Result_Mean}) and (\ref{Eq_GaussianCLT_Result_Cov}).
			
			\item[vi)] \textit{Calculating the mean and covariance structure}:\\
			By Lemma~\ref{Lemma_SpaceTransform}, equation (\ref{Eq_GaussianCLT_BaiMean}) yields
			\begin{align*}
				\E\big[ \tilde{G}(\Phi_\infty(z)) \big] & \overset{(\ref{Eq_GaussianCLT_BaiMean})}{=} \mathbbm{1}_{\beta=1} \frac{c_\infty\int_\R \frac{\cs_{\ul{\nu}_\infty}(\Phi_\infty(z))^3 \lambda^2}{(1+\lambda\cs_{\ul{\nu}_\infty}(\Phi_\infty(z)))^3} \, dH_\infty(\lambda)}{\big(1-c_\infty\int_\R \frac{\cs_{\ul{\nu}_\infty}(\Phi_\infty(z))^2 \lambda^2}{(1+\lambda\cs_{\ul{\nu}_\infty}(\Phi_\infty(z)))^2} \, dH_\infty(\lambda)\big)^2}\\
				& = \mathbbm{1}_{\beta=1} \frac{c_\infty\int_\R \frac{\lambda^2}{(\lambda-\frac{-1}{\cs_{\ul{\nu}_\infty}(\Phi_\infty(z))})^3} \, dH_\infty(\lambda)}{\big(1-c_\infty\int_\R \frac{\lambda^2}{(\lambda-\frac{-1}{\cs_{\ul{\nu}_\infty}(\Phi_\infty(z))})^2} \, dH_\infty(\lambda)\big)^2}\\
				& \overset{\text{(\ref{Eq_Def_varphi})}}{=} \mathbbm{1}_{\beta=1} \frac{c_\infty\int_\R \frac{\lambda^2}{(\lambda-z)^3} \, dH_\infty(\lambda)}{\big(1-c_\infty\int_\R \frac{\lambda^2}{(\lambda-z)^2} \, dH_\infty(\lambda)\big)^2} \ ,
			\end{align*}
			which with
			\begin{align}\label{Eq_PhiDeriv_new}
				& \Phi_{H_\infty,c_\infty}'(z) = \frac{\partial}{\partial z} \big[ (1 - c_\infty z\cs_{H_\infty}(z) - c_\infty) z \big] \overset{\text{(\ref{Eq_DefStieltjes})}}{=} \frac{\partial}{\partial z} \bigg[ \bigg(1 - c_\infty\int_\R \frac{\lambda}{\lambda-z} dH_\infty(\lambda)\bigg) z \bigg] \nonumber\\
				& = \bigg(1 - c_\infty\int_\R \frac{\lambda}{\lambda-z} dH_\infty(\lambda)\bigg) + \bigg(- c_\infty\int_\R \frac{\lambda}{(\lambda-z)^2} dH_\infty(\lambda)\bigg) z \nonumber\\
				& = 1 - c_\infty\int_\R \frac{\lambda (\lambda-z)}{(\lambda-z)^2} dH_\infty(\lambda) - c_\infty\int_\R \frac{\lambda z}{(\lambda-z)^2} dH_\infty(\lambda) \nonumber\\
				& = 1 - c_\infty\int_\R \frac{\lambda^2}{(\lambda-z)^2} dH_\infty(\lambda)
			\end{align}
			and the additional calculation
			\begin{align}\label{Eq_GaussianCLT_varphiDerivCalc}
				&1 = \frac{d}{dz} \varphi_{\nu_\infty,c_\infty}\big( \Phi_\infty(z) \big) = \varphi_{\nu_\infty,c_\infty}'\big( \Phi_\infty(z) \big) \, \Phi_\infty'(z) \nonumber\\
				& \overset{\text{(\ref{Eq_Def_varphi})}}{=} \frac{\cs_{\ul{\nu}_\infty}'(\Phi_\infty(z))}{\cs_{\ul{\nu}_\infty}(\Phi_\infty(z))^2} \, \Phi_\infty'(z) \overset{\text{(\ref{Eq_Def_varphi})}}{=} z^2 \cs_{\ul{\nu}_\infty}'(\Phi_\infty(z)) \, \Phi_\infty'(z) \ ,
			\end{align}
			allows for
			\begin{align*}
				& \E\big[ G(z) \big] \overset{\text{(\ref{Eq_GaussianCLT_MeanFirst})}}{=} -\frac{\E[\tilde{G}(\Phi_\infty(z))]}{c_\infty z^2\cs_{\ul{\nu}_\infty}'(\Phi_\infty(z))} \overset{\text{(\ref{Eq_PhiDeriv_new})}}{=} - \mathbbm{1}_{\beta=1} \frac{\frac{1}{2}\Phi_\infty''(z)}{c_\infty z^2\cs_{\ul{\nu}_\infty}'(\Phi_\infty(z))\Phi_\infty'(z)^2}\\
				& \overset{\text{(\ref{Eq_GaussianCLT_varphiDerivCalc})}}{=} - \mathbbm{1}_{\beta=1} \, \frac{\Phi_\infty''(z)}{2c_\infty \Phi_\infty'(z)} \ ,
			\end{align*}
			thus proving (\ref{Eq_GaussianCLT_Result_Mean}). Similarly, equation (\ref{Eq_GaussianCLT_BaiCov}) becomes
			\begin{align*}
				& \Cov\big[ \tilde{G}(\Phi_\infty(z_1)), \tilde{G}(\Phi_\infty(z_2)) \big]\\
				& = \frac{1}{\beta} \Big(\frac{\cs_{\ul{\nu}_\infty}'(\Phi_\infty(z_1)) \cs_{\ul{\nu}_\infty}'(\Phi_\infty(z_2))}{(\cs_{\ul{\nu}_\infty}(\Phi_\infty(z_1)) - \cs_{\ul{\nu}_\infty}(\Phi_\infty(z_2)))^2} - \frac{1}{(\Phi_\infty(z_1)-\Phi_\infty(z_2))^2}\Big)\\
				& \overset{\text{(\ref{Eq_Def_varphi})}}{=} \frac{1}{\beta} \Big(\frac{\cs_{\ul{\nu}_\infty}'(\Phi_\infty(z_1)) \cs_{\ul{\nu}_\infty}'(\Phi_\infty(z_2))}{(\frac{1}{z_2}-\frac{1}{z_1})^2} - \frac{1}{(\Phi_\infty(z_1)-\Phi_\infty(z_2))^2}\Big)\\
				& = \frac{1}{\beta} \Big(\frac{z_1^2\cs_{\ul{\nu}_\infty}'(\Phi_\infty(z_1)) z_2^2\cs_{\ul{\nu}_\infty}'(\Phi_\infty(z_2))}{(z_1-z_2)^2} - \frac{1}{(\Phi_\infty(z_1)-\Phi_\infty(z_2))^2}\Big) \ .
			\end{align*}
			Combining this with (\ref{Eq_GaussianCLT_CovFirst}) yields
			\begin{align*}
				& \Cov\big[ G(z_1), G(z_2) \big]\\
				& = \frac{1}{\beta c_\infty^2} \Big(\frac{1}{(z_1-z_2)^2} - \frac{1}{(\Phi_\infty(z_1)-\Phi_\infty(z_2))^2 z_1^2z_2^2\cs_{\ul{\nu}_\infty}'(\Phi_\infty(z_1)) \cs_{\ul{\nu}_\infty}'(\Phi_\infty(z_2))}\Big)\\
				& \overset{\text{(\ref{Eq_GaussianCLT_varphiDerivCalc})}}{=} \frac{1}{\beta c_\infty^2} \Big(\frac{1}{(z_1-z_2)^2} - \frac{\Phi_\infty'(z_1) \Phi_\infty'(z_2)}{(\Phi_\infty(z_1)-\Phi_\infty(z_2))^2}\Big) \ ,
			\end{align*}
			which proves (\ref{Eq_GaussianCLT_Result_Cov}). \qedhere
		\end{itemize}
	\end{proof}

	\subsection{Proof of Corollary~\ref{Cor_GaussPLSSCLT}}\label{Proof_Cor_GaussPLSSCLT}
	\begin{itemize}
		\item[i)] \textit{Cauchy's integral formula}:\\
		Let $\tilde{\gamma}$ be the closed curve going counter-clockwise through $\mathrm{image}(\gamma) \cup \ol{\mathrm{image}(\gamma)}$, then Cauchy's integral formula for any $g \in \operatorname{Hol}(\gamma_n)$ gives
		\begin{align}\label{Eq_GaussPLSSCLT_Step1_1}
			& \int_\R g(\lambda) \, dH_n(\lambda) = \int_\R \frac{1}{2\pi \bm{i}} \oint_{\tilde{\gamma}} g(z) \frac{1}{z-\lambda} \, dz \, dH_n(\lambda) \nonumber\\
			& = \frac{-1}{2\pi \bm{i}} \oint_{\tilde{\gamma}} g(z) \int_\R \frac{1}{\lambda-z} \, dH_n(\lambda) \, dz = \frac{-1}{2\pi \bm{i}} \oint_{\tilde{\gamma}} g(z) \cs_{H_n}(z) \, dz \nonumber\\
			& = \frac{-1}{2\pi \bm{i}} \int_{\gamma} g(z)\cs_{H_n}(z) - g(\ol{z})\cs_{H_n}(\ol{z}) \, dz \ .
		\end{align}
		The left hand side of (\ref{Eq_GaussPLSSCLT_Result2}) is thus equal to
		\begin{align}\label{Eq_GaussPLSSCLT_Step1_2}
			& n\Big( \hat{L}_{n,\gamma}^{(0)}(g) - \int_\R g \, dH_n \Big) \nonumber\\
			& = \frac{-1}{2\pi \bm{i}} \int_{\gamma} g(z) n\big(\hat{s}_n^{(0)}(z) - \cs_{H_n}(z)\big) - g(\ol{z}) n\big(\hat{s}_n^{(0)}(\ol{z}) - \cs_{H_n}(\ol{z})\big) \, dz \ .
		\end{align}
		Note that $\hat{s}_n$ and thus also $\hat{s}_n$ was (below Definition~\ref{Def_StieltjesEstimator}) canonically extended to $\C^-$ such that $\hat{s}_n^{(0)}(\ol{z}) = \ol{\hat{s}_n^{(0)}(z)}$.
		
		\item[ii)] \textit{Extension of Corollary~\ref{Cor_GaussianCLT} to $J=\mathrm{closure}(\mathrm{image}(\gamma))$}:\\
		The proof of Corollary~\ref{Cor_GaussianCLT} actually proved the stronger convergence (\ref{Eq_GaussCLT_Result_Stronger}) under the auxiliary support bound (\ref{Eq_StrongerGaussCLT_AdditionalConditionK}) and with the interval $I$ from (\ref{Eq_StrongerGaussCLT_Def_I}).
		\\[0.5em]
		By construction, the image of $\gamma : (0,1) \rightarrow \bD_{H_\infty,c_\infty}(1)$ lies in $\bD_{H_\infty,c_\infty}(1)$.
		Assumption (\ref{Eq_GaussPLSSCLT_Cond1}) guarantees that the compact set
		\begin{align}\label{Eq_GaussPLSSCLT_Step2_1}
			& J\coloneq\mathrm{closure}(\mathrm{image}(\gamma))
		\end{align}
		for small enough $\delta>0$ does not overlap the set $I$ as defined in (\ref{Eq_StrongerGaussCLT_Def_I}), so
		\begin{align*}
			& J \subset \bD_{H_\infty,c_\infty}(1) \cup (\R\setminus I)
		\end{align*}
		and result (\ref{Eq_GaussCLT_Result_Stronger}) applies. It thus from (\ref{Eq_GaussCLT_Result_Stronger}) and $\mathrm{image}(\gamma) \overset{\text{(\ref{Eq_GaussPLSSCLT_Step2_1})}}{=} J\setminus\R$ follows that
		\begin{align}\label{Eq_GaussCLT_Result_Stronger_copy}
			& n\big( \hat{s}_{n}^{(0)}(z) - \cs_{H_n}(z) \big)_{z \in \mathrm{image}(\gamma)} \xrightarrow{n \to \infty}_{\mathcal{D}} G(z)_{z \in \mathrm{image}(\gamma)} \ \ \text{in } \ell^\infty(\mathrm{image}(\gamma)) \ .
		\end{align}
		
		\item[iii)] \textit{Continuous mapping theorem}:\\
		As the map
		\begin{align*}
			& F : \ell^{\infty}\big(\mathrm{image}(\gamma)\big) \rightarrow \R \ \ ; \ \ M \mapsto \frac{-1}{2\pi \bm{i}} \int_{\gamma} g(z) M(z) - g(\ol{z}) \ol{M(z)} \, dz
		\end{align*}
		is continuous, the continuous mapping theorem turns (\ref{Eq_GaussCLT_Result_Stronger_copy}) into
		\begin{align}
			& n\Big( \hat{L}_{n,\gamma}^{(0)}(g) - \int_\R g \, dH_n \Big) \overset{\text{(\ref{Eq_GaussPLSSCLT_Step1_2})}}{=} F\big( n(\hat{s}_n^{(0)}-\cs_{H_n}) \big) \xrightarrow{n \to \infty}_{\mathcal{D}} F(G) \ .
		\end{align}
		Setting $Z(g) = F(G)$ together with the observations $\bm{e}(z) = \E[G(z)]$ and $\bm{c}(z_1,z_2) = \Cov[G(z_1),G(z_2)]$ completes the proof by calculations
		\begin{align*}
			& \E[F(G)] = \frac{-1}{2\pi \bm{i}} \int_{\gamma} g(z) \underbrace{\E[G(z)]}_{=\bm{e}(z)} - g(\ol{z}) \underbrace{\E[G(\ol{z})]}_{=\bm{e}(\ol{z})} \, dz
		\end{align*}
		and, using $G(\ol{z}) = \ol{G(z)}$, also
		\begin{align*}
			& \Cov[F(G),F(G)] = \frac{-1}{4\pi^2} \int_{\gamma} \int_{\gamma} \Cov\big[g(z_1) G(z_1) - g(\ol{z}_1) G(\ol{z}_1) ,\\
			& \hspace{7cm} g(z_2) G(z_2) - g(\ol{z}_2) G(\ol{z}_2)\big] \, dz_1 \, dz_2\\
			& = - \frac{1}{4\pi^2} \int_{\gamma} \int_{\gamma} g(z_1) g(z_2) \Cov[G(z_1),G(z_2)] \, dz_1 \, dz_2\\
			& \hspace{0.5cm} + \frac{1}{4\pi^2} \int_{\gamma} \int_{\gamma} g(z_1) g(\ol{z}_2) \Cov[G(z_1),G(\ol{z}_2)] \, dz_1 \, dz_2\\
			& \hspace{0.5cm} + \frac{1}{4\pi^2} \int_{\gamma} \int_{\gamma} g(\ol{z}_1) g(z_2) \Cov[G(\ol{z}_1),G(z_2)] \, dz_1 \, dz_2\\
			& \hspace{0.5cm} - \frac{1}{4\pi^2} \int_{\gamma} \int_{\gamma} g(\ol{z}_1) g(\ol{z}_2) \Cov[G(\ol{z}_1),G(\ol{z}_2)] \, dz_1 \, dz_2
		\end{align*}
		as well as
		\begin{align*}
			& \Cov[F(G),\ol{F(G)}] = \frac{1}{4\pi^2} \int_{\gamma} \int_{\gamma} \Cov\big[g(z_1) G(z_1) - g(\ol{z}_1) G(\ol{z}_1) ,\\
			& \hspace{7cm} \ol{g(z_2) G(z_2) - g(\ol{z}_2) G(\ol{z}_2)}\big] \, dz_1 \, dz_2\\
			& = \frac{1}{4\pi^2} \int_{\gamma} \int_{\gamma} g(z_1) \ol{g(z_2)} \Cov[G(z_1),G(\ol{z_2})] \, dz_1 \, dz_2\\
			& \hspace{0.5cm} - \frac{1}{4\pi^2} \int_{\gamma} \int_{\gamma} g(z_1) \ol{g(\ol{z}_2)} \Cov[G(z_1),G(z_2)] \, dz_1 \, dz_2\\
			& \hspace{0.5cm} - \frac{1}{4\pi^2} \int_{\gamma} \int_{\gamma} g(\ol{z}_1) \ol{g(z_2)} \Cov[G(\ol{z}_1),G(\ol{z_2})] \, dz_1 \, dz_2\\
			& \hspace{0.5cm} + \frac{1}{4\pi^2} \int_{\gamma} \int_{\gamma} g(\ol{z}_1) \ol{g(\ol{z}_2)} \Cov[G(\ol{z}_1),G(z_2)] \, dz_1 \, dz_2 \ . \qed
		\end{align*}
	\end{itemize}
	
	\subsection{Proof of Lemma \ref{Lemma_FourthMomentMP_law}}
	This lemma is a basic application of Theorem 1.1 of \cite{MP_RowCorrelation_Bai}, though the prerequisites of said theorem must first be shown to follow from the assumptions of this lemma. The prerequisites to Theorem 1.1 of \cite{MP_RowCorrelation_Bai} in this papers notation translate to:
	\begin{itemize}
		\item[1)] Suppose $\E[\ol{(\bm{Y}_n)_{j,k}} (\bm{Y}_n)_{l,k}] = \Sigma_{l,j}$ holds for all all $k \leq n$ and $j,l \leq d$, and with the notation $\bm{y}^{(n)}_k = (\bm{Y}_n)_{\bullet,k} \in \C^d$ suppose
		\begin{align}\label{Eq_BaiMP_condition}
			& \frac{1}{n^2} \sup\limits_{\substack{A \in \C^{d \times d} \\ ||A|| \leq 1}} \E\big[ \big| (\bm{y}^{(n)}_k)^* A \bm{y}^{(n)}_k - \tr(A \Sigma_n) \big|^2 \big] \xrightarrow{n \to \infty} 0 \ .
		\end{align}
		
		\item[2)] Suppose the convergence $\frac{d}{n} \to c$ holds for some $c \in (0,\infty)$.
		
		\item[3)] Suppose there exists a constant $\sigma^2>0$ with $||\Sigma_n|| \leq \sigma^2$ for all $n \in \N$, and also suppose $H_n \xRightarrow{n \to \infty} H$ for some probability measure $H$.
	\end{itemize}
	Note that (2) holds immediately by \ref{EI_ItemAssumption_Asymptotics} and (3) holds by \ref{EI_ItemAssumption_sigmaBound} and \ref{EI_ItemAssumption_PopConv}. Since \ref{EI_ItemAssumption_X_Structure} together with $\bm{Y}_n = B_n \bm{X}_n$ implies $\E[\ol{(\bm{Y}_n)_{j,k}} (\bm{Y}_n)_{l,k}] = \Sigma_{l,j}$, it only remains to show (\ref{Eq_BaiMP_condition}).
	\\[0.5em]
	With the notation $\bm{x}^{(n)}_k = (\bm{X}_n)_{\bullet,k} \in \C^d$, one by $\bm{Y}_n = B_n \bm{X}_n$ and $\Sigma_n = B_n B_n^*$ has
	\begin{align*}
		& (\bm{y}^{(n)}_k)^* A \bm{y}^{(n)}_k - \tr(A \Sigma_n) = (\bm{x}^{(n)}_k)^* B_n^* A B_n \bm{x}^{(n)}_k - \tr(B_n^* A B_n)\\
		& = \sum\limits_{i,j=1}^d \ol{(\bm{X}_n)_{i,k}} (B_n^* A B_n)_{i,j} (\bm{X}_n)_{j,k} - \sum\limits_{i=1}^d (B_n^* A B_n)_{i,i} \\
		& = \underbrace{\sum\limits_{i=1}^d \big( \ol{(\bm{X}_n)_{i,k}} (B_n^* A B_n)_{i,i} (\bm{X}_n)_{i,k} - (B_n^* A B_n)_{i,i} \big)}_{\eqcolon S_1} + \underbrace{\sum\limits_{\substack{i,j=1 \\ i \neq j}}^d \ol{(\bm{X}_n)_{i,k}} (B_n^* A B_n)_{i,j} (\bm{X}_n)_{j,k}}_{\eqcolon S_2} \ .
	\end{align*} 
	The second moment of the "diagonal sum" $S_1$ may be bounded by
	\begin{align*}
		& \E\big[ |S_1|^2 \big] \overset{\text{\ref{EI_ItemAssumption_X_Structure}}}{=} \bV\Big[ \sum\limits_{i=1}^d \ol{(\bm{X}_n)_{i,k}} (B_n^* A B_n)_{i,i} (\bm{X}_n)_{i,k} \Big] \overset{\text{indep.}}{=} \sum\limits_{i=1}^d \bV\big[ \ol{(\bm{X}_n)_{i,k}} (B_n^* A B_n)_{i,i} (\bm{X}_n)_{i,k} \big]\\
		& = \sum\limits_{i=1}^d \big| (B_n^* A B_n)_{i,i} \big|^2 \big( \E[|(\bm{X}_n)_{i,k}|^4] - \E[|(\bm{X}_n)_{i,k}|^2]^2 \big) \overset{\text{(\ref{Eq_BoundedFourthMoments})}}{\leq} C_4 \sum\limits_{i=1}^d \big| (B_n^* A B_n)_{i,i} \big|^2
	\end{align*}
	second moment of the "off-diagonal sum" $S_2$ may in turn be bounded by
	\begin{align*}
		& \E\big[ |S_2|^2 \big] = \sum\limits_{\substack{i,j=1 \\ i \neq j}}^d \sum\limits_{\substack{i',j'=1 \\ i' \neq j'}}^d (B_n^* A B_n)_{i,j} \ol{(B_n^* A B_n)_{i',j'}} \overbrace{\E\big[ \ol{(\bm{X}_n)_{i,k}} (\bm{X}_n)_{j,k} (\bm{X}_n)_{i',k} \ol{(\bm{X}_n)_{j',k}} \big]}^{=0 \text{ unless } (i,j) = (i',j') \text{ or } (i,j) = (j',i')}\\
		& = \sum\limits_{\substack{i,j=1 \\ i \neq j}}^d (B_n^* A B_n)_{i,j} \ol{(B_n^* A B_n)_{i,j}} \E[|(\bm{X}_n)_{i,k}|^2] \E[|(\bm{X}_n)_{j,k}|^2]\\
		& \hspace{0.5cm} + \sum\limits_{\substack{i,j=1 \\ i \neq j}}^d (B_n^* A B_n)_{i,j} \ol{(B_n^* A B_n)_{j,i}} \ol{\E[(\bm{X}_n)_{i,k}^2]} \E[(\bm{X}_n)_{j,k}^2]\\
		& \leq \sum\limits_{\substack{i,j=1 \\ i \neq j}}^d |(B_n^* A B_n)_{i,j}|^2 + \sum\limits_{\substack{i,j=1 \\ i \neq j}}^d |(B_n^* A B_n)_{i,j}| \, |(B_n^* A B_n)_{j,i}| \overset{\text{C.S.}}{\leq} 2 \sum\limits_{\substack{i,j=1 \\ i \neq j}}^d |(B_n^* A B_n)_{i,j}|^2 \ .
	\end{align*}
	Combining the two previous bounds, one has
	\begin{align*}
		& \E\big[ |S_1|^2 + |S_2|^2 \big] \leq (2 \lor C_4) \sum\limits_{i,j=1}^d \big| (B_n^* A B_n)_{i,j} \big|^2 = (2 \lor C_4) \tr\big( B_n^* A B_n B_n^* A^* B_n \big)\\
		& = (2 \lor C_4) \tr\big( A \Sigma_n A^* \Sigma_n \big) \leq (2 \lor C_4) d ||A \Sigma_n A^* \Sigma_n|| \overset{\text{\ref{EI_ItemAssumption_sigmaBound}}}{\leq} (2 \lor C_4) \sigma^4 d ||A||^2 \ ,
	\end{align*}
	which proves (\ref{Eq_BaiMP_condition}) by the calculation
	\begin{align*}
		& \frac{1}{n^2} \sup\limits_{\substack{A \in \C^{d \times d} \\ ||A|| \leq 1}} \E\big[ \big| (\bm{y}^{(n)}_k)^* A \bm{y}^{(n)}_k - \tr(A \Sigma_n) \big|^2 \big] = \frac{1}{n^2} \sup\limits_{\substack{A \in \C^{d \times d} \\ ||A|| \leq 1}} \E\big[ \big| S_1 + S_2 \big|^2 \big]\\
		& \leq \frac{1}{n^2} \sup\limits_{\substack{A \in \C^{d \times d} \\ ||A|| \leq 1}} 2\E\big[ |S_1|^2 + |S_2|^2 \big] \leq \frac{1}{n^2} \sup\limits_{\substack{A \in \C^{d \times d} \\ ||A|| \leq 1}} 2(2 \lor C_4) \sigma^4 d ||A||^2 = 2(2 \lor C_4) \sigma^4 \frac{d}{n^2} \xrightarrow{n \to \infty} 0 \ .
	\end{align*}
	\qed
	
	\subsection{Proof of Lemma~\ref{Lemma_Hadamard}}\label{Proof_Lemma_Hadamard}
	For any sequences $(t_m)_{m \in\N} \subset (0,\infty)$ and $(h_m)_{m \in \N} \subset D$ with $t_m \to 0$, $h_m \to h \in D_0$ and $f+t_mh_m \in D_\phi$, the convergence
	\begin{align}\label{Eq_CLTInversion_FunkDelt_HadamardConv_ToShow}
		&  \frac{\phi(f+t_mh_m) - \phi(f)}{t_m} \xrightarrow{m \to \infty} \phi_f'(h)
	\end{align}
	in the Banach space $H^{\infty}(W)$ is to be shown.
	As $f' \neq 0$ on the compact set $K$, there must exist a $C > 0$ such that
	\begin{align}\label{Eq_CLTInversion_FunkDelt_alpha}
		& \forall \tz \in K : \ \frac{1}{C} \leq |f'(\tz)| \leq C \ .
	\end{align}
	The continuous map
	\begin{align*}
		& h_2 : K \times K \rightarrow [0,\infty) \ \ ; \ \ \tilde{v},\tilde{w} \mapsto \begin{cases}
			\frac{|f(\tilde{v})-f(\tilde{w}) - f'(\tilde{w})(\tilde{v}-\tilde{w})|}{|\tilde{v}-\tilde{w}|^2} & \text{, if } \tilde{v} \neq \tilde{w}\\
			\frac{1}{2} |f''(\tilde{w})| & \text{, if } \tilde{v} = \tilde{w}
		\end{cases}
	\end{align*}
	must take its maximum on the compact domain $K \times K$, which yields the uniform bound
	\begin{align}
		& |f(\tilde{v})-f(\tilde{w}) - f'(\tilde{w})(\tilde{v}-\tilde{w})| \leq C'|\tilde{v}-\tilde{w}|^2 \label{Eq_CLTInversion_FunkDelt_UniformTaylor2}
	\end{align}
	for some $C'>0$.
	For each $z \in W$ define
	\begin{align}\label{Eq_CLTInversion_FunkDelt_tzmDef}
		& \tz_m = \tz_m(z) \coloneq (f+t_mh_m)^{-1}(z) \in K \ \ \text{ and } \ \ \tz = \tz(z) \coloneq f^{-1}(z) \in K \ .
	\end{align}
	Suppose $\max\limits_{z \in W}|(f+t_mh_m)^{-1}(z)-f^{-1}(z)|$ does not converge to zero, then after transition to a sub-sequence there must exist an $\varepsilon > 0$ and $z_m \in W$ such that $|(f+t_mh_m)^{-1}(z_m)-f^{-1}(z_m)| \geq \varepsilon$ for all $m \in \N$. Compactness of $K$ yields that again after transition to a sub-sequence, one may assume $(f+t_mh_m)^{-1}(z_m) \xrightarrow{m \to \infty} \tz_*$ and $f^{-1}(z_m) \xrightarrow{m \to \infty} \tz_{**}$ for some $\tz_*,\tz_{**} \in K$, which by the previous construction may not be the same. It however must hold that
	\begin{align*}
		& f(\tz_*) - f(\tz_{**}) = \lim\limits_{m \to \infty} f((f+t_mh_m)^{-1}(z_m)) - f(f^{-1}(z_m))\\
		& = \lim\limits_{m \to \infty} z_m -z_m - t_mh_m(z_m) = 0 \ .
	\end{align*}
	This is a contradiction to the injectivity of $f \in D_\phi$, which proves
	\begin{align}\label{Eq_CLTInversion_FunkDelt_maxDiff}
		& \max\limits_{z \in W}|\tz_m(z)-\tz(z)| = \max\limits_{z \in W}|(f+t_mh_m)^{-1}(z)-f^{-1}(z)| \xrightarrow{m \to \infty} 0 \ .
	\end{align}
	As a first consequence, one may without loss of generality assume $m$ to be large enough that
	\begin{align}
		& |\tz_m-\tz| \leq \frac{1}{2CC'} \overset{\text{(\ref{Eq_CLTInversion_FunkDelt_alpha})}}{\leq} \frac{|f'(\tz_m)|}{2C'} \label{Eq_CLTInversion_FunkDelt_mLarge1} \ .
	\end{align}
	As a second consequence, the uniform convergence $\|h_m-h\|_\infty \to 0$ with (\ref{Eq_CLTInversion_FunkDelt_alpha}) ensures the convergence
	\begin{align}\label{Eq_CLTInversion_FunkDelt_hmConvergence}
		& \sup\limits_{z \in W} \Big| \underbrace{\frac{h(\tz)}{f'(\tz)}}_{=\phi_f(h)(z)} - \frac{h_m(\tz_m)}{f'(\tz_m)} \Big| \xrightarrow{m \to \infty} 0 \ .
	\end{align}
	The bound
	\begin{align}\label{Eq_CLTInversion_FunkDelt_Bound1}
		& \big| t_mh_m(\tz_m) - f'(\tz_m)(\tz-\tz_m) \big| \nonumber\\
		& = \big| \overbrace{f(\tz)}^{=z} - \overbrace{(f+t_mh_m)(\tz_m)}^{=z} + t_mh_m(\tz_m) - f'(\tz_m)(\tz-\tz_m) \big| \nonumber\\
		& = \big| f(\tz) - f(\tz_m) - f'(\tz_m)(\tz-\tz_m) \big| \overset{\text{(\ref{Eq_CLTInversion_FunkDelt_UniformTaylor2})}}{\leq} C' |\tz-\tz_m|^2
	\end{align}
	may be combined with (\ref{Eq_CLTInversion_FunkDelt_mLarge1}) to by triangle inequality see
	\begin{align*}
		& \frac{1}{2} \big| f'(\tz_m)(\tz-\tz_m) \big| \leq t_m|h_m(\tz_m)| \leq t_m \underbrace{\sup\limits_{M \in \N} \|h_M\|_\infty}_{\eqcolon  C'' < \infty} \ ,
	\end{align*}
	thus proving
	\begin{align}\label{Eq_CLTInversion_FunkDelt_DffBound_tm}
		& |\tz_m-\tz| \leq t_m \frac{2C''}{|f'(\tz_m)|} \overset{\text{(\ref{Eq_CLTInversion_FunkDelt_alpha})}}{\leq} t_m 2CC'' \ .
	\end{align}
	Finally, Hadamard-differentiability (\ref{Eq_CLTInversion_FunkDelt_HadamardConv_ToShow}) follows from the calculation
	\begin{align*}
		& \Big| \frac{\phi(f+t_mh_m) - \phi(f)}{t_m}(z) - \frac{h_m(\tz_m)}{f'(\tz_m)} \Big| \overset{\text{(\ref{Eq_CLTInversion_FunkDelt_tzmDef})}}{=} \Big| \frac{\tz_m-\tz}{t_m} - \frac{h_m(\tz_m)}{f'(\tz_m)} \Big|\\
		& \overset{\text{(\ref{Eq_CLTInversion_FunkDelt_alpha})}}{\leq} \frac{C}{t_m} \big| f'(\tz_m)(\tz_m-\tz) - t_mh_m(\tz_m) \big| \overset{\text{(\ref{Eq_CLTInversion_FunkDelt_Bound1})}}{\leq} \frac{CC'}{t_m} |\tz_m-\tz|^2 \overset{\text{(\ref{Eq_CLTInversion_FunkDelt_DffBound_tm})}}{\leq} t_m 4C^3C'(C'')^2 \to 0
	\end{align*}
	together with (\ref{Eq_CLTInversion_FunkDelt_hmConvergence}). \qed

	\newpage
	\cleardoublepage
	\section{List of symbols}\label{ListOfSymbols}
	\begin{itemize}
		
		\item[] $a_{\gamma}$ \tabto{1.7cm} end point of the curve $\gamma$ (see (\ref{Eq_AdmissibleCurve_Cond2}) or (\ref{Eq_GLSS_gammafCondition1}))
		
		\item[] $B_n$ \tabto{1.7cm} a deterministic $(d \times d)$ matrix (see~\ref{EI_ItemAssumption_X_Structure} of Assumption~\ref{Assumption_EigInf_Main})
		
		\item[] $B_{\varepsilon}^{\C}(z)$ \tabto{1.7cm} the open $\varepsilon$-neighborhood around $z$ in $\C$ (see Subsection~\ref{Subsection_ModelAndNotation})
		
		\item[] $B_{\varepsilon}^{\C^+}(z)$ \tabto{1.7cm} the open $\varepsilon$-neighborhood around $z$ in $\C^+$ (see Subsection~\ref{Subsection_ModelAndNotation})
		
		\item[] $b_{\gamma}$ \tabto{1.7cm} starting point of the curve $\gamma$ (see (\ref{Eq_AdmissibleCurve_Cond2}) or (\ref{Eq_GLSS_gammafCondition1}))
		
		\item[] $\C^+$ \tabto{1.7cm} the set of complex numbers with positive imaginary part (see Subsection~\ref{Subsection_ModelAndNotation})
		
		\item[] $\C^-$ \tabto{1.7cm} the set of complex numbers with negative imaginary part (see Subsection~\ref{Subsection_ModelAndNotation})
		
		\item[] $c$ \tabto{1.7cm} a constant $c>0$ (see Lemma~\ref{Lemma_MP})
		
		\item[] $c_{n/\infty}$ \tabto{1.7cm} the ratio $c_n=\frac{d}{n}$ or its limit $c_\infty = \lim\limits_{n \to \infty} \frac{d}{n}$ (see~\ref{EI_ItemAssumption_Asymptotics} of Assumption~\ref{Assumption_EigInf_Main})
		
		\item[] $\mathrm{closure}$ \tabto{1.7cm} the closure of a set (see Subsection~\ref{Subsection_ModelAndNotation})
		
		\item[] $\bm{D}(\tau)$ \tabto{1.7cm} a spectral domain (see Theorem~\ref{Thm_OuterLaw})
		
		\item[] $\bD_{H,c}(\theta)$ \tabto{1.7cm} a certain spectral domain (see (\ref{Eq_DefD}))
		
		\item[] $\hat{\bD}(\tau,\kappa,n)$ \tabto{1.7cm} an empirical spectral domain (see Theorem~\ref{Thm_Consistency})
		
		\item[] $d$ \tabto{1.7cm} the data-dimension (see the introduction)
		
		\item[] $\diag$ \tabto{1.7cm} an operator for turning a sequence of numbers into the diagonal matrix\\
		\tabto{1.7cm} with said numbers on its diagonal
		
		\item[] $\dist$ \tabto{1.7cm} the distance between subsets of $\C$ (see Subsection~\ref{Subsection_ModelAndNotation})
		
		\item[] $\delta$ \tabto{1.7cm} symbol used to denote small constants in various settings
		
		\item[] $\delta_x$ \tabto{1.7cm} point-probability measure with all its mass at $x \in \R$
		
		\item[] $\varepsilon$ \tabto{1.7cm} symbol used to denote small constants in various settings
		
		\item[] $\tilde{\varepsilon}$ \tabto{1.7cm} symbol used to denote small constants in various settings
		
		\item[] $F_{n,\bullet}$ \tabto{1.7cm} a function used in the application of Rouch\'{e}'s theorem (see (\ref{Eq_BasicConsistency_Step4_6}) or (\ref{Eq_Consistency_Step4_1}))
		
		\item[] $f$ \tabto{1.7cm} placeholder for a function
		
		\item[] $G$ \tabto{1.7cm} a limiting Gaussian process (see Corollary~\ref{Cor_GaussianCLT})
		
		\item[] $\tilde{G}$ \tabto{1.7cm} a limiting Gaussian process (see Theorem~\ref{Thm_CLT_Inversion})
		
		\item[] $G_n$ \tabto{1.7cm} a function used in the application of Rouch\'{e}'s theorem (see (\ref{Eq_BasicConsistency_Step4_5}) or (\ref{Eq_Consistency_Step4_2}))
		
		\item[] $g$ \tabto{1.7cm} placeholder for a function
		
		\item[] $\gamma$ \tabto{1.7cm} placeholder for a curve through $\C$
		
		\item[] $\gamma_n$ \tabto{1.7cm} an admissible curve (see Definition~\ref{Def_AdmissibleCurve})
		
		\item[] $\gamma_n^{(f)}$ \tabto{1.7cm} part of an admissible pair $(\gamma_n^{(f)},\gamma_n^{(g)})$ (see Definition~\ref{Def_AdmissiblePair})
		
		\item[] $\gamma_n^{(g)}$ \tabto{1.7cm} part of an admissible pair $(\gamma_n^{(f)},\gamma_n^{(g)})$ (see Definition~\ref{Def_AdmissiblePair})
		
		\item[] $H$ \tabto{1.7cm} a probability measure (see Lemma~\ref{Lemma_MP})
		
		\item[] $H_n$ \tabto{1.7cm} the population spectral distribution (see (\ref{Eq_Def_Hn}))
		
		\item[] $H_\infty$ \tabto{1.7cm} the limiting population spectral distribution (see~\ref{EI_ItemAssumption_PopConv} of Assumption~\ref{Assumption_EigInf_Main})
		
		\item[] $\operatorname{Hol}(\gamma)$ \tabto{1.7cm} a certain space of holomorphic functions (see Subsection~\ref{Subsection_ModelAndNotation})
		
		\item[] $\Im(z)$ \tabto{1.7cm} the imaginary part of a complex number $z \in \C$
		
		\item[] $\mathrm{image}(\gamma)$ \tabto{1.7cm} the image of a curve $\gamma$ in $\C$ (see Subsection~\ref{Subsection_ModelAndNotation})
		
		\item[] $\bm{i}$ \tabto{1.7cm} the imaginary unit
		
		\item[] $J$ \tabto{1.7cm} placeholder for a compact subset of $\C$
		
		\item[] $K$ \tabto{1.7cm} a constant such that $\supp(\hat{\nu}_n)$ and/or $\supp(\nu_n)$ is contained in $[0,K]$\\
		\tabto{1.7cm} (see Lemma~\ref{Lemma_StandardBounds})
		
		\item[] $\kappa$ \tabto{1.7cm} symbol used to denote large constants in various settings
		
		\item[] $L_n(g)$ \tabto{1.7cm} a population linear spectral statistic (see the introduction)
		
		\item[] $\hat{L}_{n,\gamma_n}(g)$ \tabto{1.7cm} the PLSS estimator (see Definition~\ref{Def_PLSS_estimator})
		
		\item[] $L_n(f,g)$ \tabto{1.7cm} a generalized linear spectral statistic (see (\ref{Eq_DefGLSS_theoretical}))
		
		\item[] $\hat{L}_{n,\gamma_n^{(f)},\gamma_n^{(g)}}(f,g)$ \tabto{3.0cm} the GLSS estimator (see Definition~\ref{Def_GLSS_estimator})
		
		\item[] $\ell^\infty(S)$ \tabto{1.7cm} the function space of essentially bounded Borel-measurable $f : S \rightarrow \C$
		
		\item[] $\ell(\gamma)$ \tabto{1.7cm} the length of a curve $\gamma$
		
		\item[] $\lambda$ \tabto{1.7cm} an integration variable (see e.g. (\ref{Eq_DefStieltjes}))
		
		\item[] $\lambda_j$ \tabto{1.7cm} an eigenvalue of $\Sigma_n$ (notation only used on page 18)
		
		\item[] $\lambda_j(A)$ \tabto{1.7cm} an eigenvalue of a Hermitian matrix $A$ (see above (\ref{Eq_Def_Hn}))
		
		\item[] $\mu$ \tabto{1.7cm} placeholder for a measure
		
		\item[] $n$ \tabto{1.7cm} the sample size (see the introduction)
		
		\item[] $\nu$ \tabto{1.7cm} a probability measure on $[0,\infty)$ (see Lemma~\ref{Lemma_MP})
		
		\item[] $\hat{\nu}_n$ \tabto{1.7cm} the sample spectral distribution (see Subsection~\ref{Subsection_ModelAndNotation})
		
		\item[] $\nu_n$ \tabto{1.7cm} the deterministic equivalent of $\hat{\nu}_n$ (see Subsection~\ref{Subsection_ModelAndNotation})
		
		\item[] $\nu_\infty$ \tabto{1.7cm} the limiting spectral distribution (see below (\ref{Eq_MP_law}))
		
		\item[] $\hat{\ul{\nu}}_n$ \tabto{1.7cm} a transform of $\hat{\nu}_n$ (see (\ref{Eq_Def_ulNu}))
		
		\item[] $\ul{\nu}_n$ \tabto{1.7cm} a transform of $\nu_n$ (see (\ref{Eq_Def_ulNu}))
		
		\item[] $\ul{\nu}_\infty$ \tabto{1.7cm} a transform of $\nu_\infty$ (see below (\ref{Eq_Def_ulNu}))
		
		
		\item[] $\textrm{PD}_d$ \tabto{1.7cm} the set of all positive semi-definite $(d \times d)$-matrices (see lemma~\ref{Lemma_BorelMeasurability})
		
		\item[] $\Phi_{H,c}$ \tabto{1.7cm} a certain holomorphic map (see Lemma~\ref{Lemma_SpaceTransform})
		
		\item[] $\hat{\Phi}_n$ \tabto{1.7cm} the map $z \mapsto (1-c_nz\hat{s}_n(z)-c_n)z$ (see e.g. Theorem~\ref{Thm_Consistency})
		
		\item[] $\varphi_{\nu,c}$ \tabto{1.7cm} a certain holomorphic map (see Lemma~\ref{Lemma_SpaceTransform})
		
		\item[] $\Re(z)$ \tabto{1.7cm} the real part of a complex number $z \in \C$
		
		\item[] $\bm{S}_n$ \tabto{1.7cm} the $(d \times d)$ sample covariance matrix (see~\ref{EI_ItemAssumption_X_Structure} of Assumption~\ref{Assumption_EigInf_Main})
		
		\item[] $\bS(\tau,n)$ \tabto{1.7cm} a spectral domain (see Theorem~\ref{Thm_OuterLaw})
		
		\item[] $\hat{\bS}(\tau,n)$ \tabto{1.7cm} a spectral domain (see Theorem~\ref{Corollary_OuterLaw})
		
		\item[] $\cs_\bullet(z)$ \tabto{1.7cm} the Stieltjes transform of a measure at position $z \in \C$ (see (\ref{Eq_DefStieltjes}))
		
		\item[] $\hat{s}_n(z)$ \tabto{1.7cm} the population Stieltjes transform estimator (see Definition~\ref{Def_StieltjesEstimator})
		
		\item[] $\hat{s}_n^{(0)}(z)$ \tabto{1.7cm} a version of the population Stieltjes transform estimator,\\
		\tabto{1.7cm} which takes the value $0$ when $\hat{s}_n(z)$ does not exist (see Theorem~\ref{Thm_CLT_Inversion})
		
		\item[] $\supp$ \tabto{1.7cm} the support of a measure (see Subsection~\ref{Subsection_ModelAndNotation})
		
		\item[] $\Sigma_n$ \tabto{1.7cm} the $(d \times d)$ population covariance matrix (see~\ref{EI_ItemAssumption_X_Structure} of Assumption~\ref{Assumption_EigInf_Main})
		
		\item[] $\sigma^2$ \tabto{1.7cm} upper bound for $\|\Sigma_n\|$ (see~\ref{EI_ItemAssumption_sigmaBound} of Assumption~\ref{Assumption_EigInf_Main})
		
		\item[] $\tr$ \tabto{1.7cm} the trace of a square matrix
		
		\item[] $\tau$ \tabto{1.7cm} symbol used to denote small constants in various settings
		
		\item[] $\theta$ \tabto{1.7cm} a parameter (usually $1$ or $\infty$) used in the definition (\ref{Eq_DefD})
		
		\item[] $U$ \tabto{1.7cm} placeholder for a neighborhood or an open set
		
		\item[] $U_n$ \tabto{1.7cm} placeholder for a unitary matrix
		
		\item[] $V_n$ \tabto{1.7cm} placeholder for a unitary matrix
		
		\item[] $\bm{X}_n$ \tabto{1.7cm} a $(d \times n)$ matrix with independent entries (see~\ref{EI_ItemAssumption_X_Structure} of Assumption~\ref{Assumption_EigInf_Main})
		
		\item[] $\xi_{c,z}$ \tabto{1.7cm} a function introduced for proofs of measurability (see Lemma~\ref{Lemma_BorelMeasurability})\\
		\tabto{1.7cm} it by construction holds that $\hat{s}_n^{(0)}(z) = \xi_{c_n,z}(\bm{S}_n)$
		
		\item[] $\bm{Y}_n$ \tabto{1.7cm} the $(d \times n)$ data matrix (see~\ref{EI_ItemAssumption_X_Structure} of Assumption~\ref{Assumption_EigInf_Main})
		
		\item[] $z$ \tabto{1.7cm} a complex number (from $H$-space, see Lemma~\ref{Lemma_SpaceTransform})
		
		\item[] $\tz$ \tabto{1.7cm} a complex number (from $\nu$-space, see Lemma~\ref{Lemma_SpaceTransform})
		
		\item[] $\ol{z}$ \tabto{1.7cm} complex conjugate of a $z \in \C$ (see Subsection~\ref{Subsection_ModelAndNotation})
		
		\item[] $\|\cdot\|$ \tabto{1.7cm} spectral norm of a matrix
		
		\item[] $\|\cdot\|_\gamma$ \tabto{1.7cm} supremum norm on the image of the curve $\gamma$ (see Subsection~\ref{Subsection_ModelAndNotation})
		
		\item[] $\xrightarrow{n \to \infty}_\mathcal{D}$ \tabto{1.7cm} convergence in distribution
		
		\item[] $\xRightarrow{n \to \infty}$ \tabto{1.7cm} weak convergence of measures

	\end{itemize}

	\section*{Acknowledgements}
	The presented results are part of the author’s PhD thesis. Support of the Research Unit 5381, DFG-grant RO3766/8-1, is gratefully acknowledged.
	
	\bibliographystyle{apalike}
	\bibliography{literature}
	
\end{document}